\documentclass[11pt]{article}
\usepackage{amsfonts}
\usepackage{amscd}
\usepackage[all]{xy}

\def\proof{\vspace{2ex}\noindent{\bf Proof.} }

\def\endproof{\relax\ifmmode\expandafter\endproofmath\else
  \unskip\nobreak\hfil\penalty50\hskip.75em\hbox{}\nobreak\hfil\bull
  {\parfillskip=0pt \finalhyphendemerits=0 \bigbreak}\fi}
\def\endproofmath$${\eqno\bull$$\bigbreak}
\def\bull{\vbox{\hrule\hbox{\vrule\kern3pt\vbox{\kern6pt}\kern3pt\vrule}\hrule}}
\def\msigma{\ifmmode{\cal M}^\sigma\else {${\cal M}^\sigma$}\fi}

\newtheorem{theorem}{Theorem}[subsection]
\newtheorem{proposition}[theorem]{Proposition}
\newtheorem{lemma}[theorem]{Lemma}
\newtheorem{claim}[theorem]{Claim}

\newtheorem{corollary}[theorem]{Corollary}
\newtheorem{assumption}[theorem]{Assumption}
\newtheorem{D}[theorem]{Definition}
\newenvironment{defn}{\begin{D} \rm }{\end{D}}
\newtheorem{addendum}[theorem]{Addendum}
\newtheorem{R}[theorem]{Remark}
\newenvironment{remark}{\begin{R}\rm }{\end{R}}

\def\ov{\overline}
\def\spcheck{^{\vee}}
\def\frak{\mathfrak}
\newsymbol\rtimes 226F

\title{Almost commuting elements in  compact Lie
  groups}  
\author{Armand Borel\\
Institute for Advanced Study, Princeton NJ 08540 
\and Robert
Friedman\thanks{Research partially supported by NSF Grant
DMS-96-22681} \\ Department of
Mathematics, Columbia University, New York, NY 10027
\and John W. Morgan\thanks{Research partially 
supported by NSF Grant DMS-97-04507} \\ Department of
Mathematics, Columbia University, New York, NY 10027}

\begin{document}
\maketitle
\tableofcontents

\section{Introduction}

Let $K$ be a compact, connected and semisimple Lie group. This paper
describes the moduli space of  isomorphism classes of flat
connections on principal
$K$-bundles over the two-torus and the three-torus. There are two
motivations for this study. The first is the relation between flat
$K$-bundles over the two-torus and holomorphic  principal bundles
over an elliptic curve with structure group  the complexification of
$K$. The second  is to give a proof of a conjecture of
Witten concerning the moduli space of  flat
$K$-bundles over the three-torus, used in \cite{Witten}. 

 Of
course, a flat bundle is completely determined by its holonomy
representation, so that the problem of classifying flat bundles over
two- and three-tori becomes the question of classifying ordered
pairs and triples of commuting elements in $K$, up to simultaneous
conjugation. If two elements commute in $K$, then any lifts of
them to the universal cover of $K$ commute up to an element of the
center. Thus we shall work in the simply connected cover $G$ of $K$ and
study pairs and triples of elements in $G$ which commute up to the
center, hence the name ``almost commuting." This  is the  form in which
we attack the question. Our point of view is that the extended Dynkin
diagram of
$G$, the action of
$\pi_1(K)$ on this diagram, and the coroot integers associated to this
diagram completely determine the answer in a manner which we shall
describe below.

\medskip

{\sl Notation used throughout this paper}

\noindent
$\bullet$ $G$ is a compact, connected
and simply connected Lie group, and in particular $G$ is semi-simple.

\noindent
$\bullet$ ${\cal C}G$  denotes the center of $G$.

\noindent
$\bullet$ For any subset $X\subseteq G$, $Z_G(X)$
denotes the centralizer of $X$ in $G$. We will denote $Z_G(X)$ by
$Z(X)$ when $G$ is clear from the context.

\noindent
$\bullet$ If $S$ is a torus in $G$, not necessarily maximal, we let
$W(S,G)$ be the finite group $N_G(S)/Z_G(S)$ and we call it the {\sl
Weyl group of $S$ in $G$}.  If $\frak s ={\rm Lie}(S)$, then we set
$W(\frak s,G)=W(S,G)$.

\noindent
$\bullet$ Given $x,y\in G$, we define the commutator $[x,y] =
xyx^{-1}y^{-1}$ and denote  conjugation by $x$ as
$i_x(y) = {}^xy = xyx^{-1}$. 

\medskip

One convention concerning subtori used throughout the paper is
the following.  Let ${\frak t}$ be a vector space with a positive
definite inner product $\langle\cdot,\cdot \rangle$, let
$\Lambda\subseteq {\frak t}$ be a  lattice, such that $\langle v,
w\rangle\in {\bf Z}$ for all
$v,w\in \Lambda$, and let
$T={\frak t}/\Lambda$ be the associated torus.  
For any subtorus $S\subseteq T$  let
${\frak s}\subseteq {\frak t}$ be its tangent space.
Let ${\frak s}^\perp$ be the
perpendicular subspace to ${\frak s}$.
Then ${\frak s}^\perp$ is the tangent space of another 
subtorus $S'\subseteq T$ and $F_S=S\cap S'$ is a finite group. We
denote by
$\ov S$ the quotient torus $S/F_S= T/S'$. Clearly $S = \frak s
/(\frak s \cap \Lambda)$ and $\ov S = \frak s/\pi (\Lambda)$,
where
$\pi\colon
\frak t\to \frak s$ is orthogonal projection.

\subsection{Preliminaries}
Let $\Phi$ be a root system on a vector space ${\frak t}$.
Fix  a set of
simple roots $\Delta$ for $\Phi$. Denote by $Q= Q(\Phi)$ the
root lattice and by
$P$ the lattice of weights with basis
$\{\varpi_\alpha\}_{\alpha\in \Delta}$.  Further, we denote the
inverse root system by 
$\Phi\spcheck$; it is a subset of ${\frak t}$. We have
$Q\spcheck=Q(\Phi\spcheck)$ the coroot lattice dual to $P$ and
$P\spcheck$ the  
coweight lattice dual to $Q(\Phi)$.
We denote the Weyl group by $W(\Phi)$ or simply by $W$ if $\Phi$
is clear from the context. We fix an inner product
$\langle\,\cdot,\cdot\,\rangle$ on ${\frak t}$ invariant under
the Weyl  group  and use it to identify
$\frak t$ with
$\frak t^*$. We choose this inner product so that in each
irreducible factor the shortest length of a coroot is 
$\sqrt2$. Coroots of this length are called {\sl short coroots\/}
whereas coroots of longer lengths are called {\sl long coroots}. 
Roots inverse to short coroots are thus {\sl long roots\/} and
roots inverse to long coroots are {\sl  short roots}.

Assume now that $\frak t ={\rm Lie}(T)$, where  
$T$ is a maximal torus of $G$, and let $\Phi=\Phi(T,G)$  be
the corresponding root system. 
In this case,  the exponential map identifies
${\frak t}/Q\spcheck$ with $T$ and ${\frak t}/P\spcheck$ with the
maximal torus in the  adjoint form of $G$, and $P\spcheck/Q\spcheck$
with ${\cal C}G$. More generally, for an abstract root system
$\Phi$, we define   ${\cal C}\Phi$ to be the
finite group $P\spcheck(\Phi)/Q\spcheck(\Phi)$.

Associated to  $a,b\in\Phi$ is the Cartan
integer defined by
\begin{equation}\label{Cartan}
n(a,b)=\frac{2\langle a,b\rangle }{\langle b,b\rangle}.
\end{equation}
The Cartan integers for pairs of elements in $\Delta$ determine and
are determined by the Dynkin diagram $D(G)$ whoses nodes 
are indexed by the $a\in \Delta$ and whose bonds together with their
multiplicities and arrows \cite{Bour}.  The set
$\Delta\spcheck\subseteq \Phi\spcheck$ consists of the coroots 
$a\spcheck$ inverse to each root $a\in\Delta$.
The Cartan
integers $n(a\spcheck,b\spcheck)$ for $a\spcheck,b\spcheck\in
\Delta\spcheck$ are described by a
Dynkin diagram $D\spcheck(G)$, {\sl the coroot diagram} for $G$. Its nodes
are identified in the 
obvious way with the nodes of $D(G)$. Its bonds, including the
multiplicities, are exactly the 
same as the bonds in $D(G)$, but the direction of every arrow
is reversed.

Suppose that $\Phi$ is irreducible.
Let $d$ be the highest root of $\Phi$ with respect to $\Delta$.
Set $\tilde a=-d$ and let  $\widetilde \Delta=\Delta\cup\{\tilde
a\}$ be the {\sl extended set of
simple roots}. Let $C_0\subseteq {\frak t}$ be the positive Weyl
chamber associated to
$\Delta$ and let $A\subseteq C_0$ be the unique alcove in $C_0$
containing the origin. The walls of $A$ are given by
$\{a=0\}_{a\in\Delta}$ and $\{\tilde a=-1\}$. Hence there is a natural
bijection between $\widetilde \Delta$ and the walls of $A$.
The set $\widetilde \Delta$ is the set of nodes for the extended
Dynkin diagram $\widetilde D(G)$.
The Cartan integers $n(a,b)$ for $a,b\in \widetilde \Delta$ are 
recorded in the multiplicities of the bonds  and the directions
of the arrows of
 $\widetilde D(G)$, by exactly the same rules as given in the case 
of $D(G)$. In the case of $\widetilde A_1$, we shall always make
the convention that the two nodes are connected by two single
bonds, so that the diagram is a cycle. Dually, there is the
extended coroot diagram
$\widetilde D\spcheck(G)$ whose nodes are the coroots
$\widetilde\Delta\spcheck$ inverse to the roots in $\widetilde
\Delta$. As in the case of $D(G)$, the diagram $\widetilde
D\spcheck(G)$ is obtained from
$\widetilde D(G)$ by reversing the directions of all the arrows
on the multiple bonds. In case $\Phi =
\coprod_i\Phi_i$ is reducible and the $\Phi_i$ are the irreducible
factors, we define the set of extended roots
$\widetilde \Delta$ as $\coprod _i\widetilde \Delta _i$, and define the
extended root and coroot diagrams as the disjoint union of the
corresponding diagrams of the factors.

Assuming again that $\Phi$ is irreducible, there is a single linear
relation among the roots of
$\widetilde
\Delta$, namely 
\begin{equation}\label{highest}
\sum_{a\in \widetilde \Delta}h_aa=0
\end{equation}
for positive integers $h_a$, with $h_{\tilde a}=1$.
The $h_a$ will be  called the {\sl root integers}.
The sum $\sum_{a\in \widetilde \Delta}h_a=h$ is the {\sl Coxeter number}
of the group. 
Dually,  there is a single linear relation between the coroots
inverse to the roots in $\widetilde
\Delta$. It takes the form
\begin{equation}
\sum_{a\in \widetilde \Delta}g_aa\spcheck=0
\end{equation}
where the $g_a$ are all positive integers and the coefficient of
the coroot $\tilde a\spcheck$ is one.
(N.B. The coroot inverse to the highest root is  the
highest short coroot, and is equal to the highest coroot if and only if
$\Phi$ is simply laced.) The integers $g_a$ are called the {\sl coroot
integers} and the sum
$g=\sum_{a\in \widetilde \Delta}g_a$ is called the {\sl dual
Coxeter  number}.
Since $\tilde a$ is a long root, it follows that $g_a|h_a$ for 
every
$a\in \widetilde \Delta$, with equality exactly for the long
roots  in
$\widetilde \Delta$. It will be convenient to view the coroot
integers as defining a function 
${\bf g}\colon \widetilde \Delta
\to {\bf N}$ by the formula ${\bf g}(a) = g_a$.

There is an action 
 of ${\cal C}G$ as a group of affine isometries of $\frak t$
normalizing the alcove $A$, which will be described in more detail 
in  Section~\ref{action1}.
For each
$c\in {\cal C}G$,   the differential  $w_c$ 
of its action on $A$ is a linear map which is an 
element of $ W$   normalizing $\widetilde \Delta\subset {\frak t}^*$.
Since the action of $W$ preserves the Cartan integers, it
follows that the resulting action of
$w_c$ on the nodes of $\widetilde D(G)$ is a  diagram
automorphism. The action preserves the root integers $h_a$, in the
sense that $h_{w_c\cdot a} = h_a$. Of course, $w_c$
also induces a diagram automorphism of $\widetilde D \spcheck(G)$
preserving the coroot
integers $g_a$. This defines a faithful representation of ${\cal
C}G$ as  a group of diagram automorphisms of $\widetilde D(G)$
and of
$\widetilde D\spcheck(G)$.

Let ${\cal C}$ be a subgroup of ${\cal C}G$. By restriction, 
${\cal C}$ acts on
${\frak t}$, and the linear part $w_{\cal C}$ of this action 
normalizes the subsets
$\widetilde \Delta$ and $\widetilde
\Delta\spcheck$ and defines an action of ${\cal C}$ by
diagram automorphisms on $\widetilde D(G)$ and $\widetilde
D\spcheck(G)$.  This action  preserves the coroot integers
$g_a$. For each orbit $\ov a\in \widetilde
\Delta/{\cal C}$ we define $g_{\ov a}$ to be $n_{\ov a}g_a$
where $n_{\ov a}$ is the cardinality of the orbit $\ov a$ and $g_a$ is
the coroot integer associated to any $a$ in this orbit. Let
$\widetilde\Delta_{\cal C}$, resp.\ $\widetilde\Delta_{\cal
C}\spcheck$, be the quotient of $\widetilde\Delta$, resp.\ 
$\widetilde\Delta\spcheck$,  under the action of ${\cal C}$. Let
$S^{w_{\cal C}}$ be the subtorus of $T$ whose Lie algebra is
${\frak t}^{w_{\cal C}}$, the  subspace of ${\frak t}$ pointwise
fixed under the action of $w_{\cal C}$. Let $\pi\colon \frak t
\to {\frak t}^{w_{\cal C}}$ be orthogonal projection.

\subsection{The case of commuting pairs in a simply connected group}

In \cite{Tohoku}, the first author proved
that, if $x$ and 
$y$ are commuting 
elements in a compact, simply connected group $G$, then there is a
maximal torus $T\subseteq G$ containing both $x$ and $y$. Furthermore,
two pairs of elements in $T$ are conjugate in $G$ if and only if
they are conjugate by an element of $W$. Thus, the
moduli space of conjugacy classes of commuting pairs of elements
in $G$ is identified with 
$(T\times T)/W$.
The torus  $T$ is the quotient of ${\frak t}$ by the lattice
$Q\spcheck$ 
generated by $\widetilde \Delta\spcheck$ and $W$ is the group of
isometries of ${\frak t}$ generated by reflections in the
hyperplanes through the origin defined by the elements of
$\widetilde \Delta\spcheck$. The group $W$ acts on ${\frak t}$
preserving the lattice $Q\spcheck$ and hence there is an induced
$W$-action on $T$.
This is the model result that we carry over into all the
other cases we study.

\subsection{$c$-pairs}

Next, let us consider a compact, connected, simple,
but not necessarily simply
connected group $K$ with $G$ as simply connected covering. The
first invariant of a $K$-bundle 
$\xi$ over $T^2$ is the characteristic class $w(\xi)\in
H^2(T^2;\pi_1(K))=\pi_1(K)$.
We identify $\pi_1(K)$ with a
subgroup of the center ${\cal C}G$ of $G$, so that $w(\xi)\in
{\cal C}G$. If $\xi$ has a flat 
connection and if $x,y\in K$ are the
holonomy images of the standard generators of the fundamental group of
the two-torus, then for any lifts $\tilde x\in G$, resp.\  $\tilde
y\in G$, of $x$, resp.\  $y$, we have $[\tilde x,\tilde
y]=w(\xi)$. Our classification results in this case are
simplified by three assumptions:
\begin{enumerate}
\item We fix the topological type of the bundle $\xi$, or
equivalently the class $w(\xi) = c\in {\cal C}G$.
\item We assume that $\xi$ does not lift to a bundle over any
  non-trivial covering group of $K$.
\item We classify flat connections on $\xi$ up to restricted gauge
equivalence, i.e., up to $G$-gauge equivalence. In this case, it turns
out that restricted gauge
equivalence is the same as $K$-equivalence.
\end{enumerate}

Translating these conditions gives the following equivalent group
theoretic problem.

{\sl Let $G$ be  simple, and let $c\in {\cal C}G$.
A pair of elements $(x,y)$ in $G$ is said to be a $c$-pair if 
$[x,y]=c$. We classify $c$-pairs up to simultaneous conjugation
by elements of $G$.}

Note that $\frak t^{w_{\langle c \rangle}} = \frak t^{w_c}$.
Thus we define $S^{w_c} = S^{w_{\langle c \rangle}}$.

\begin{theorem}\label{c-pairsthm}
Let $G$  be  simple, and let $c\in {\cal C}G$. Then the  moduli
space of conjugacy classes of $c$-pairs of
  elements in 
$G$ is homeomorphic to 
$(\ov S^{w_c}\times \ov S^{w_c})/W(S^{w_c},G)$.
\end{theorem}

In a very closely related form, this theorem was first proved by
Schweigert \cite{Schweig}. In Theorem~\ref{diagram1} below, we
shall describe $\ov S^{w_c}$ and $W(S^{w_c},G)$ in terms of the
extended coroot diagram of $G$ and the action of $c$ on this
diagram.

\subsection{Commuting triples}

Next, we let $G$ be simple and we turn to flat $G$-bundles over the $3$-torus.
The holonomy of such a connection around the standard basis of the fundamental
group of the torus is a commuting
triple $(x,y,z)$ in $G$ well-defined up to simultaneous conjugation.
 Let ${\cal T}_G$ denote the moduli space of
conjugacy classes of commuting triples in $G$.
In general, ${\cal T}_G$ has several
components even though there is only one topological type for
a $G$-bundle over $T^3$.

\begin{theorem}\label{commuttrip}
Let $G$ be simple.
 For any $k\ge 1$ dividing
at least one of the coroot integers $g_a$ we set 
$\widetilde{I}(k) =\{a\in \widetilde{\Delta}: k\not|g_a\}$, and we let
$S(k)\subseteq T$ be the subtorus with
$${\rm Lie}(S(k))={\frak
t}(k)= \bigcap _{a\in \widetilde{I}(k)}{\rm Ker}\,a.$$
Note that $\dim S(k)$ is one less than the number of $a$ such
that $k|g_a$. Then: 
\begin{enumerate}
\item For each commuting triple $(x,y,z)$ in $G$, there is a
unique integer $k\ge 1$ dividing
at least one of the coroot integers $g_a$ such that $S(k)$ is 
conjugate  to a maximal torus for $Z(x,y,z)$.
The integer $k$ is 
called the {\rm order} of $(x,y,z)$.
\item The order is a conjugacy class invariant and defines a
locally constant function on ${\cal T}_G$.  We define the order
of a component $X$ of ${\cal T}_G$ to be the value of this
function on $X$. 
\item  If $k\geq 1$ divides at least one of the $g_a$,
there  are  exactly $\varphi(k)$ components of ${\cal
T}_G$ of order $k$, where
$\varphi$  is the Euler $\varphi$-function. Given a component $X$
of ${\cal T}_G$, let $d_X = \frac13\dim X + 1$. Then
$$\sum _X d_X = g.$$
\item  Each component of ${\cal T}_G$ of order $k$ is
homeomorphic to 
$$\left(\ov S(k)\times \ov S(k)\times \ov S(k)\right)/W(S(k), G).$$
\item Let $\pi_k$ be orthogonal projection from $\frak t$ to
$\frak t(k)$.  For 
$a\in\widetilde \Delta$,  $\pi_k( a\spcheck)$ is non-zero if and
only if $a\notin
\widetilde I(k)$. Thus $\pi_k$ determines an embedding of
$\widetilde
\Delta\spcheck-\widetilde I\spcheck(k)$ into ${\frak t}(k)$.   
The torus $\ov S(k)$ is the quotient of ${\frak t}(k)$ by the
lattice $\pi_k(Q\spcheck(k))$, which is spanned by the image of
 $\widetilde
\Delta\spcheck-
\widetilde I\spcheck(k)$. The group  $W(S(k),G)$ is the group of
isometries of 
${\frak t}(k)$ generated by reflections in the 
hyperplanes $\pi_k(a\spcheck)^\perp$ for $a\in \widetilde
\Delta -\widetilde I (k)$.  
\end{enumerate}
\end{theorem}

Results along these lines have been obtained independently by 
Kac-Smilga \cite{Kac-Smilga}.

In Section~\ref{quot} we give a result which shows that the image
under $\pi_k$ of 
$\Delta\spcheck- \widetilde I(k)\spcheck$ in ${\frak t}(k)$ is the
extended set of simple coroots  of a root system $\Phi(\frak
t(k))$ of 
${\frak
  t}(k)$ and derives its extended coroot diagram from the 
extended coroot diagram $\widetilde D\spcheck (G)$ and the coroot
integers $g_a$.

\subsection{$C$-triples}

Let $K$ be a  connected, simple group 
with simply connected covering $G$.
Now we consider principal $K$-bundles with flat connection over the
three-torus. 
Given a basis for the fundamental group of the $3$-torus, a flat
connection is determined by 
the holonomy image $x_1,x_2,x_3\in K$ for the given basis of the 
fundamental group of the torus. These elements are defined up to
simultaneous  conjugation in $K$. Lifting these to elements
$\tilde x_1,\tilde x_2,\tilde x_3$ in 
$G$, we have elements $c_{ij}=[\tilde
x_i,\tilde x_j]$ in $\pi_1(K)\subseteq {\cal C}G$.
These  elements determine the topological type of the $K$-bundle. 
We record these elements by constructing 
$3\times 3$ matrix $C=(c_{ij})$ with values in ${\cal C}G$.
The entries of $C$ satisfy: $c_{ii} =1, c_{ij}=c_{ji}^{-1}$ for
$i\neq j$. Such a matrix will be called {\sl antisymmetric}.Let
$\langle C\rangle$ be the subgroup of ${\cal C}G$ generated by
the entries of $C$.  A triple of elements $(x_1, x_2,
x_3)$ in
$G$ such that $c_{ij}=[ 
x_i,  x_j]$ is a {\sl
$C$-triple}. If
$c_{12} = c$ and $c_{13} =c_{23} =1$, we  refer to $(x_1, x_2,
x_3)$ as a  {\sl $c$-triple}. As long as 
$\langle C\rangle$ is cyclic and generated by $c$, the moduli space of
$C$-triples can be identified with the moduli space of $c$-triples.  As
in the case of the   two-torus, we assume that the bundle does not
lift to a proper covering group of
$K$. This is equivalent to saying that the group $\langle C\rangle $
generated by the entries of $C$
is all of $\pi_1(K)$. Also as before, we
classify these bundles up to restricted gauge equivalence, i.e., 
up to conjugation by elements of $G$. In this case, this
classification differs slightly from the classification up to
conjugation by $K$. The classification of commuting triples in
$K$ follows easily from the results of this paper, but we shall
not give the details here.

Given $C$, we denote by $w_C$ the subgroup $w_{\langle C\rangle}$
of $W$. Following this convention,  $\frak t^{w_C}$,
$S^{w_C}$,
$\widetilde \Delta_C$, and $\widetilde \Delta\spcheck_C$ are
defined similarly.

\begin{theorem}\label{ctrip}
Let $G$ be simple and let $C$ be an antisymmetric $3\times 3$-matrix
with values in ${\cal C}G$. Let ${\cal T}_G(C)$ be the moduli
space of conjugacy classes of $C$-triples in $G$.  Let $\ov {\bf
g}\colon \widetilde
\Delta _C \to {\bf N}$ be the function defined by $\ov {\bf g}
(\ov a) = g_{\ov a}$. For any
$k\ge 1$ dividing at least one of  the $g_{\ov  a}$,   let
$\widetilde I_C(k)$ be the set of $\ov{a}\in \widetilde
\Delta_C$ such that $k \not| g_{\ov{a}}$ and let
$S^{w_C}(\ov {\bf g},k)\subseteq S^{w_C}$ be the 
  subtorus whose  Lie   algebra is
$${\frak t}^{w_C}(\ov {\bf g}, k)=\bigcap_{\ov{a} \in
\widetilde I_C(k)}
{\rm Ker}\,\ov a \subseteq {\frak t}^{w_C}.$$
Let  $\pi^C_k$ be orthogonal projection from $\frak t$ to ${\frak
t}^{w_C}(\ov {\bf g}, k)$. Then:
\begin{enumerate}
\item For every $C$-triple $(x,y,z)$, there is an integer
$k\ge 1$ dividing at least one of the $g_{\ov a}$, called the
{\rm order} of $(x,y,z)$,  such that $S^{w_C}(\ov {\bf g}, k)$ is
conjugate to a maximal torus for $Z(x,y,z)$.   The dimension of
$S^{w_C}(\ov {\bf g}, k)$ is one less than the number of $\ov a$
such that
$k|g_{\ov a}$.
\item The order is a conjugacy class invariant and defines a
locally constant function on ${\cal T}_G(C)$. We define the order
of a component $X$ of ${\cal T}_G(C)$ to be the value of this
function on  $X$. 
\item For any
$k\ge 1$ dividing at least one of  the $g_{\ov  a}$, there are
exactly 
$\varphi(k)$
  components of ${\cal T}_G(C)$ of order $k$. Given a component
$X$ of ${\cal T}_G(C)$, let $d_X = \frac13\dim X + 1$. Then
$$\sum _X d_X = g.$$
\item   In the case where $\langle C\rangle$  is cyclic, i.e.
the case of
$c$-triples, each  component of
${\cal T}_G(C)$ of order $k$ is
  homeomorphic to either
$$\left(\ov S\times \ov S\times
    \ov S\right)/W(S,G)
\ \ {\textit{or}}\ \  
\left(\ov S\times \ov S\times
     S\right)/W(S,G),$$
where $S=S^{w_c}(\ov {\bf g}, k)$.
We will describe the possibilities precisely in
Section~\ref{ctripsect}. In case $\langle C\rangle$ is not cyclic,
each component of the moduli space is homeomorphic to
$$\left((S\times S\times S)/F)/W(S,G)\right),$$ 
where $S = S^{w_C}(\ov {\bf g}, k)$ and $F$ is a finite subgroup
of
$S\times S\times S$ which will be described explicitly in Section
12.  
\item  For $a\in \widetilde \Delta$, $\pi_k^C(a\spcheck) \neq
0$ depends only on $\ov a\in \widetilde \Delta_C$ and is 
nonzero if and only if
$\ov a \notin
\widetilde I_C(k)$.  Thus
$\pi^C_k$ determines an embedding of
$\widetilde \Delta\spcheck_C-\widetilde I\spcheck_C(k)$
into ${\frak t}^{w_C}(\ov {\bf g}, k)$.
The torus $\ov 
S^{w_C}(\ov {\bf g}, k)$ is the quotient of ${\frak t}^{w_C}(\ov
{\bf g},k)$ by the  lattice $Q\spcheck_C(k)$ spanned by image of 
$\widetilde 
\Delta\spcheck_C-
\widetilde I\spcheck_C(k)$. The group  $W(S^{w_C}(\ov {\bf
g},k),G)$ is the  group of isometries of
${\frak t}^{w_C}(\ov {\bf g},k)$ generated by reflections in the
hyperplanes $\pi^C_k( a\spcheck)^\perp$ for
$\ov a\in \widetilde
\Delta_C-\widetilde I_C(k)$. 
\end{enumerate}
\end{theorem}

In Section~\ref{quot} we state a result which shows that the
image of 
$\widetilde
\Delta_C\spcheck-
\widetilde I_C\spcheck(k)$ in ${\frak t}^{w_C}(\ov {\bf g},k)$ is
the extended set of simple coroots    of a root system, and which
describes its extended coroot  diagram in terms of the
extended coroot diagram $\widetilde D\spcheck (G)$,  the
action of
$\langle C\rangle$ on this diagram, and the coroot integers $g_a$.

\subsection{Quotients of diagram automorphisms}

Let $\widetilde D$ be an extended coroot diagram.
A diagram automorphism 
$\sigma\colon \widetilde D\to \widetilde D$ automatically 
preserves the  coroot integers.
Suppose that ${\cal C}$ is a group of diagram
automorphisms of
$\widetilde D$.

\begin{defn}\label{defofdiagram}
We form the quotient diagram $\widetilde
D/{\cal C}$ as follows. The nodes of $\widetilde
D/{\cal C}$ are the ${\cal C}$-orbits of nodes of
$\widetilde D$.  The orbit of the node $v$ is
denoted
$\ov v$.

There is one case that we handle directly.
If $\widetilde D$ is $\widetilde A_n$ and ${\cal C}$ acts
transitively on the set of nodes of 
$\widetilde D$, then the quotient has one element. 

Ruling out this degenerate case, 
there are two types of ${\cal C}$-orbits on an 
extended diagram $\widetilde D$: an orbit consisting of nodes, no 
two of which are connected by a bond, and an orbit consisting of
subdiagrams of type $A_2$, no two of which are connected by a
bond. The first type of orbit is called {\sl ordinary} and the
second {\sl exceptional}. Let
$\epsilon(\ov v)=1$ if $\ov v$ is an ordinary orbit and let
$\epsilon(\ov v)=2$ if $\ov v$ is exceptional. We describe the
bonds of $\widetilde D/{\cal C}$, their multiplicities
and their arrows by giving the Cartan  integers 
$n(\ov u,\ov v)$. If
$\ov u$ and
$\ov v$ are distinct orbits  such that $n(u_i,v_j)=0$
for all $u_i\in
\ov u$ and
$v_j\in \ov v$, then $n(\ov u,\ov v)=0$ and
there is no bond connecting $\ov u$ and
$\ov v$ in the quotient diagram. 
If there are nodes $u\in \ov u$ and $v\in \ov v$ which are
connected by a bond, then either ${\rm Stab}(u)\subseteq {\rm
  Stab}(v)$ or ${\rm Stab}(v)\subseteq {\rm Stab}(u)$.
In the first case we define 
$$n(\ov u,\ov v)=\epsilon(\ov
v)n(u,v),$$ and in the second case 
we define
$$n(\ov u,\ov v)=\epsilon(\ov v)\frac{n_{\ov
v}}{n_{\ov u}}n(u,v),$$ 
where as before $n_{\ov u}$ is the cardinality of the orbit
$\ov u$.   It is easy to see that these numbers are well-defined
integers and  satisfy:  $n(\ov u,\ov v)=0$
implies
$n(\ov v,\ov u)=0$; otherwise
$n(\ov u,\ov v)<0$ for $\ov u\not=\ov v$.
Thus, these numbers determine the bonds of the quotient diagram
$\widetilde D/{\cal C}$,  their multiplicities and the directions
of their arrows.

We have the set $\{g_v\}$ of coroot integers on the
original diagram. 
We define integers $g_{\ov v}$ on the nodes of $\widetilde D/{\cal
C}$ as follows. In all cases, including the degenerate ones, the
induced coroot integers $g_{\ov v}$ on  the quotient diagram are
defined by
$g_{\ov v}=\sum_vg_v$ where $v$ ranges over the nodes in the
orbit $\ov v$. Since diagram automorphisms preserve the integers
$g_v$, we have
$g_{\ov v}=n_{\ov v}g_v$ for any $v\in \ov v$.
\end{defn}

Here is the theorem which describes the fixed tori and their Weyl
groups in terms of diagram automorphisms.

\begin{theorem}\label{diagram1} Let $G$ be simple and let ${\cal
C}$ be a subgroup of ${\cal C}G$.  Let $\pi\colon
  {\frak t}\to   {\frak 
    t}^{w_{\cal C}}$  be orthogonal projection. 
Then:
\begin{enumerate}
\item  Restriction of  $\pi$  to $\widetilde \Delta\spcheck$
factors to induce   an embedding of
  $\widetilde\Delta_{\cal C}\spcheck$ in ${\frak t}^{w_{\cal
C}}$.  This embedding identifies the set of nodes $\widetilde
\Delta\spcheck$ of $\widetilde D\spcheck/{\cal C}$  with an
extended set of simple coroots for a root system $\Phi(w_{\cal
C})$. 
\item   $\widetilde D\spcheck/{\cal C}$ is the extended coroot
diagram of 
$\Phi(w_{\cal C})$.
\item The coroot lattice of $\Phi(w_{\cal C})$ is equal to
$\pi(Q\spcheck)$ , which by definition is the fundamental group 
of the torus
$\ov S^{w_{\cal C}}$, and the group $W(S^{w_{\cal C}},G)$ is the
Weyl group of $\Phi(w_{\cal C})$. 
\end{enumerate}
\end{theorem}

This result together with Theorem~\ref{c-pairsthm} yields an 
explicit description of the 
moduli space of $c$-pairs in terms of the extended coroot diagram of
$G$ and the action of $w_c$ on this diagram.

\subsection{Description of $\ov S(k)$ and $\ov
S^{w_C}(\ov {\bf g}, k)$}\label{quot}

Let $\widetilde D$ be a connected extended coroot diagram, whose
set of nodes is
$\widetilde \Delta\spcheck$,  and let ${\bf n}\colon \widetilde
\Delta \to {\bf N}$ be a function of the form $n_0{\bf g}$ for
some positive integer $n_0$.  Given an integer
$k\ge 1$  dividing at least one of the ${\bf n}(a)$, let
$\widetilde I ({\bf n}, k) =\{a\in
\widetilde \Delta : k\not|{\bf n}(a)\}$. 
We let $\widetilde D'({\bf n},k)$ be the largest subdiagram of
$\widetilde D$ having $\widetilde I\spcheck({\bf n},k)$ as its
set of nodes. Fix a length function $\ell\colon A\to {\bf R}^+$
such that if 
$v,v'\in 
\widetilde \Delta\spcheck$
are connected by a bond of $\widetilde D$ and $\ell(v)\ge
\ell(v')$ then the multiplicity of the bond between $v$ and $v'$
is 
$\ell(v)^2/\ell(v')^2$ and the arrow (if the multiplicity is not
one) points toward $v'$.

\begin{proposition}\label{typesofconnections}
Let $v \in \widetilde \Delta\spcheck- \widetilde I\spcheck({\bf
n}, k)$. Then exactly one of the following holds.
\begin{enumerate}
\item[\rm (Type $\infty$)] $v$ is the only element of $\widetilde
\Delta\spcheck- \widetilde I\spcheck({\bf n},k)$.
\item[\rm (Type 1)] The node $v$ is not connected by a bond of
$\widetilde D$ to a node of $\widetilde D'({\bf n},k)$, but is
not of Type $\infty$.
\item[\rm (Type 2)(i)] The node $v$ is connected by bonds of
$\widetilde D$ to exactly two components
of $\widetilde D'({\bf n},k)$, each of
  which is of  type $A_1$, with nodes  $v_1$ and $v_2$, say,
 and $\ell(v)\le {\rm min}(\ell(v_1),\ell(v_2))$.
\item[\rm (Type 2)(ii)] The node $v$ is connected by bonds of
$\widetilde D$ to exactly one node $v_1$ of  
$\widetilde D'({\bf n},k)$, and $\ell(v) < \ell (v_1)$.
\item[\rm (Type 3)] The node $v$ is connected by bonds of $\widetilde D$
to exactly two  components
of $\widetilde D'({\bf n},k)$, each of which is of type  $A_2$.
\item[\rm (Type 4)] The node $v$ is connected by bonds of $\widetilde D$
to exactly two  components, and exactly one of the following holds:
\begin{enumerate} 
\item[\rm (i)] Both of these components are of
type $A_1$, with nodes $v_1$ and $v_2$, say, and 
$\ell(v)>{\rm min}(\ell(v_1),\ell(v_2))$;
\item[\rm (ii)] Each component is of type $A_3$;
\item[\rm (iii)] One component is of type $A_3$ and
the  other is of type $A_1$.
\end{enumerate}
\end{enumerate}
\end{proposition}

\begin{defn}\label{newlengths}
Define a new length function $\ell_k\colon \widetilde
\Delta\spcheck- \widetilde I\spcheck({\bf n},k)\to {\bf R}^+$ by setting
$\ell_k(v)=\ell(v)/\sqrt{r}$ if $v$ is of Type $r$
according to  the above proposition. We will see later that, if
$\ell_k(v)
\geq \ell_k(w)$, then in fact $\ell_k(v)^2/\ell_k(w)^2 \in {\bf Z}$. 

Let   $\widetilde D({\bf n},k)$ be the unique diagram with nodes $\widetilde
\Delta\spcheck- \widetilde I\spcheck({\bf n},k)$
satisfying the following three conditions:
\begin{enumerate}
\item Nodes $v,v'\in \widetilde
\Delta\spcheck- \widetilde I\spcheck({\bf n},k)$ are connected by a bond in 
$\widetilde D({\bf n},k)$ if and only if either $v$ and $v'$ are connected by a
bond in $\widetilde D$ or $v$ and $v'$ are connected by bonds of
$\widetilde D$ to the same component of $\widetilde D'({\bf n},k)$.
\item Suppose that $v$ and $v'$ are distinct nodes of $\widetilde
\Delta\spcheck- \widetilde I\spcheck({\bf n},k)$
which are connected by a bond in
$\widetilde D({\bf n},k)$ and suppose that
neither of these nodes is connected by a bond of 
$\widetilde D({\bf n},k)$ to any other 
node of $\widetilde
\Delta\spcheck- \widetilde I\spcheck({\bf n},k)$. 
If $\ell_k(v)=\ell_k(v')$, then the subdiagram of $\widetilde D({\bf n},k)$
spanned by $v$ and $v'$ is of type $\widetilde A_2$.
\item Suppose that $v$ and $v'$ are distinct nodes of $\widetilde
\Delta\spcheck- \widetilde I\spcheck({\bf n},k)$
connected by a bond of $\widetilde D({\bf n},k)$ and suppose that $v$ and
$v'$ are not as in the previous case. Lastly, supppose that $\ell_k(v)
\ge \ell_k(v')$. Then the multiplicity of the bond in $\widetilde D({\bf n},k)$
connecting $v$ to $v'$ is $\ell_k(v)^2/\ell_k(v')^2$ and the arrow (if 
this multiplicity is greater than one) points toward $v'$. 
\end{enumerate}
\end{defn}

\begin{proposition}\label{newcorootdiag} 
The diagram $\widetilde D({\bf n},k)$ defined above is the coroot diagram of a
{\rm(}possibly non-reduced{\rm)} root system.
\end{proposition}

In particular, taking the quotient
diagram
$\widetilde D\spcheck/\langle C\rangle$ as our extended coroot
diagram  and the function $\ov {\bf g}$   defined in
Theorem~\ref{ctrip} as 
${\bf n}$ defines an extended coroot diagram $(\widetilde
D\spcheck/\langle C\rangle)(\ov {\bf g},k)$. The next theorem
shows that this diagram is in fact the extended coroot of a root
system on ${\frak t}^{w_C}(\ov {\bf
g},k)$.

\begin{theorem}\label{quotdiag}
Let $G$ be simple, let $C$ be an antisymmetric $3\times 3$
matrix with values in ${\cal C}G$, and let
$k\ge 1$ be an integer dividing at least one of the $g_{\ov a}$.
Orthogonal projection from $\frak t$ to ${\frak t}^{w_C}(\ov {\bf
g},k)$ factors to induce an embedding of $\widetilde
\Delta_C\spcheck-\widetilde I_C\spcheck(k)$ in ${\frak
t}^{w_C}(\ov {\bf g},k)$. This embedding identifies the nodes $\widetilde
\Delta_C\spcheck-\widetilde I_C\spcheck(k)$ of $(\widetilde
D\spcheck/\langle C\rangle)(\ov {\bf g},k)$ with an
extended set of simple coroots of a root system
$\Phi(w_C,k)$.
The diagram
$(\widetilde D\spcheck/\langle C\rangle)(\ov {\bf g},k)$
 is the extended coroot diagram 
of $\Phi(w_C,k)$. 
\end{theorem}

This theorem, together with Theorem~\ref{commuttrip} and
Theorem~\ref{ctrip} lead to an explicit description of the 
components of the moduli spaces of commutative triples and
$C$-triples in terms of the diagrams $\widetilde D\spcheck(G)$
and 
$\widetilde D\spcheck(G)/\langle C\rangle$ and the integers 
$g_{\ov a}$.

\subsection{Chern-Simons invariants and Witten's ``Clockwise Symmetry
 Conjecture''}

Assume that  $G$ is simple. 
 A  connection on a principal
$G$-bundle $\xi$ over a three-manifold $M$ has a Chern-Simons
invariant, see \cite{ChS},
which measures its difference from a trivial connection. This
invariant is well-defined modulo ${\bf Z}$ on isomorphism classes
of connections.  The Chern-Simons function is constant on
continuous paths of  flat
connections, so that we can view it as a function from the
components of the moduli space of gauge equivalence classes of
flat connections on
$\xi$ to
${\bf R}/{\bf Z}$. In case $M$ is the three-torus, this moduli
space is identified with the  space of
conjugacy classes of commuting triples in $G$, so that the
Chern-Simons invariant defines a function from the set of
components of  ${\cal T}_G$  to ${\bf R}/{\bf Z}$.

Let us consider now principal $K$-bundles, where $K$ is compact,
connected and simple. 
Let $G$ be the universal covering group of $K$. 
The topological type of a $K$-bundle $\xi$ over $T^3$ is determined by
its characteristic class $w(\xi)\in H^2(T^3; \pi_1(K))$. We construct
an antisymmetric matrix $C(\xi)=(c_{ij}(\xi))$ where $c_{ij}(\xi)$
is the value of $w(\xi)$ on the coordinate two-torus $T_{ij}$ in
the $(ij)^{\rm th}$-coordinate directions.

In this case, there is only a relative Chern-Simons invariant
which is  well-defined modulo $(1/n){\bf Z}$ for
some integer $n\ge 1$ depending on $K$.
To get a Chern-Simons invariant well-defined modulo ${\bf Z}$ we 
need to consider enhanced $K$-bundles over a three-manifold $M$,
where by an enhanced $K$-bundle $\Xi$ we mean  an underlying  
$K$-bundle
$\xi$ together with a lifting to $G$ of the structure group of
$\xi$ over the one-skeleton of $M$. We define $C(\Xi) = C(\xi)$.
Given an antisymmetric matrix $C$ with coefficients in
$\pi_1(K)$,  there is, up to isomorphism, a unique enhanced
$K$-bundle
$\Xi$ with
$C(\Xi) = C$. Under automorphisms  of
$\Xi$, the Chern-Simons invariant is well-defined modulo ${\bf Z}$. 

Let  $A$ be a flat connection on a $K$-bundle $\xi$ over
$T^3$. The holonomy of $A$ is identified with a conjugacy class
of commuting triples in
$K$.  Given an enhanced structure $\Xi$ on $\xi$, by a flat connection on
$\Xi$ we simply mean a flat connection $A$ on $\xi$. The connection
$A$ then lifts over  the $1$-skelton to a flat $G$-connection, and this
 defines a lifting of the holonomy of $A$ to a   conjugacy class of
$C(\xi)$-triples in $G$, called the {\sl $G$-holonomy\/} of $A$ on $\Xi$.
 The $G$-holonomy determines a bijection between the set of
isomorphism classes of flat connections on enhanced $K$-bundles $\Xi$
over $T^3$ with $C(\Xi)=C$
and the moduli space ${\cal T}_G(C)$ of conjugacy classes of
$C$-triples in
$G$.

\begin{theorem}\label{CS}
Let $G$ be   simple,  let $C$ be an
antisymmetric matrix with entries in ${\cal C}G$, and let $K=G/\langle
C\rangle$. Let $X_1$ be the unique component of order $1$ of the moduli
space of $C$-triples. Let
$\Xi$ be an enhanced
$K$-bundle with
$C(\Xi)=C$.
Let $\Gamma_1$ be a flat connection  on $\Xi$
whose  $G$-holonomy is a conjugacy class in $X_1$.
 Let $A$ be a flat connection on $\Xi$. Then ${\rm
CS}_{\Gamma_1}(A)\bmod {\bf Z}$ is independent of the choice of
$\Gamma_1$. Let
${\rm CS}(A)$ be the class of ${\rm CS}_{\Gamma_1}(A)\bmod {\bf Z}$. 
\begin{enumerate}
\item The function ${\rm CS}(A)$
is a well-defined function on the set of isomorphism classes 
of flat connections on $\Xi$, or equivalently 
on ${\cal T}_G(C)$, 
to ${\bf R}/{\bf Z}$. 
\item Viewing ${\rm CS}$ as a function on ${\cal T}_G(C)$, it is
constant on  components. 
\item The function ${\rm CS}$ induces a bijection from the set of
components of ${\cal T}_G(C)$ of order $k$ to the set of points
of order
$k$ in ${\bf R}/{\bf Z}$. 
\end{enumerate}
\end{theorem}

There is a refinement of this theorem, leading to a
surprising symmetry involving the Chern-Simons invariants and the
dimensions of the components of ${\cal T}_G(C)$. It was first
discovered by Witten from considerations in quantum field theory.

\begin{theorem}\label{clock}
Let $G$ be simple and let $C$ be an
antisymmetric matrix with entries in ${\cal C}G$. 
Let $g$ be the dual Coxeter number of $G$.
For each component
$X$ of ${\cal T}_G(C)$, let $d_X$ be $\frac{1}{3}{\rm
dim}(X)+1$, and let ${\rm CS}(X)$ be the value of
the Chern-Simons functional on this  component.
Let $J(X)\subseteq {\bf R}/{\bf Z}$ be a subset of $d_X$ equally spaced
points centered at ${\rm CS}(X)$ with spacing $1/g$. Then the $J(X)$
are disjoint subsets of ${\bf R}/{\bf Z}$. 
Then
$$J=\bigcup_XJ(X)\subset \frac{1}{2g}{\bf Z}/{\bf Z},$$
and $J$ consists of $g$ points invariant under translation by $1/g$.
Thus, identifying ${\bf R}/{\bf Z}$ with the unit circle, the set
$J$ is invariant under by a rotation through angle
$2\pi/g$.  
\end{theorem}

\subsection{Outline of the paper}

We begin in Section 2 with a discussion of almost commuting 
$n$-tuples  in a compact group $K$. We define the rank of such an
$n$-tuple to be the dimension of a maximal torus of the
centralizer of the $n$-tuple. We show a general finiteness
statement: for a compact group $K$, there are only finitely many
conjugacy classes of rank zero almost commuting triples in   $K$.
More generally, given a torus $S$ in
$K$, we give an explicit description of the subspace of the moduli
space of conjugacy classes of  almost commuting  $n$-tuples
${\bf x}$ such that $S$ is conjugate to a maximal torus of the 
centralizer of ${\bf x}$.
Each such component  is a quotient of  the product of $n$ copies
of $S$ by a finite group. The problem now is to determine the
possible tori and finite groups which arise in the cases we study in more
detail. 

Section 3 contains various characterizations of groups of
$A_n$-type.
These results are applied to deduce the following property of the
coroot integers $g_a$ associated to the simple roots $a$: The
subdiagram $D_k$ of the extended Dynkin diagram of $G$ consisting
of all nodes corresponding to roots  $a$ for which $g_a$ is not
divisible by a  fixed positive integer $k$ is a disjoint union of
diagrams of
$A_{n_i}$-type, for integers $n_i$ with the property that 
$(n_i+1)|k$. Furthermore, if $I_i$ denotes the set of simple
roots corresponding to the nodes of one of the components of
$D_k$, then $\sum_{a\in
  I_i}(g_a/k)a\spcheck$ lies in the center of the corresponding 
simple subgroup of $G$ and generates the center. This result
implies that the set of integers which occur as $g_a$ form an
interval $[1,N]$ and 
that, for the Euler $\varphi$-function, one has
$\varphi(N)=2$, so that $N\in
\{1,2,3,4,6\}$.  These properties, as well as the
characterization of groups of $A_n$-type, are important in
the study of $c$-pairs and
$c$-triples.

Section 4 takes up the case of $c$-pairs and contains a proof of
the first item in Theorem~\ref{c-pairsthm}.
Following the pattern laid down in Section 2, we first consider 
rank zero $c$-pairs and then pass to the general case of higher
rank.

Section 5 considers the case of commuting triples in $G$, once
again taking up the case of commuting triples of rank zero first.
This section contains a proof of the first four parts of
Theorem~\ref{commuttrip}. While this section will be for the most
part subsumed in Section 9, it seemed worthwhile to give the
reasonably straightforward arguments needed to handle this case. 

Sections 6, 7, and 8 are preparatory for the study of
$C$-triples. They examine questions which, to us, are interesting 
in their own right. Section 6 considers a group $\tau$ of affine
automorphisms of a vector space which normalizes an alcove of a
root system on that vector space. Such automorphisms are
equivalent to diagram automorphisms of the extended Dynkin
diagram of the root system. We show that the walls of the
original root system divide the fixed point set of $\tau$ into
alcoves and that these are the alcoves of a Coxeter group acting
on the fixed point set. Thus, there is a reduced root system on
the underlying vector space of the fixed point set whose alcove
structure is identified with the given alcove structure on the
affine subspace fixed by
$\tau$. We use this to study two closely related root systems on
the linear subspace fixed by the differential of $\tau$. One root
system, the restricted root system, consists of the non-trivial
restrictions of roots of the original system to the fixed
subspace. The other, the projection system, is the inverse root
system to  the root system obtained by taking the nonzero
orthogonal projections of the coroots of the original system.
These root systems are not in general equal, nor are they always
inverse systems, nor is either equal in general to the root
system on the fixed point set of the affine automorphism. It is
also not true in general that these roots systems are reduced.
Nevertheless, all three of them have the same Weyl groups, and
hence the same set of roots up to positive multiples.

Section 7 applies the work of Section 6 to the torus $S^{w_c}$
fixed by the $w_c$-action on $T$. The main result is the 
completion of the proof of Parts 2 and 3 of
Theorem~\ref{c-pairsthm} describing the Weyl group of this torus
and the fundamental group of $\ov S^{w_c}$. The results of
Section 6 are also used here to prove Theorem~\ref{diagram1}.

Section 8 is concerned with the centralizer $Z(x,y)$ of a
$c$-pair
$(x,y)$. We describe the root system corresponding to the   Lie
algebra of $Z(x,y)$, the fundamental group, and
the component group of $Z(x,y)$. We study $Z(x,y)$ by viewing
conjugation by $y$ as an automorphism of the compact group
$Z(x)$. It is natural to consider the most general problem along
these lines: let $H$ be an arbitrary compact connected group and
let $\sigma$ be an automorphism of $H$. We study the Lie algebra
and the component group of the fixed subgroup $H^\sigma$ of
$\sigma$. Our results on the Lie algebra of $H^\sigma$ generalize
those of Kac \cite{Kac1} in case $\sigma$ has finite order. The
description of the component group of $H^\sigma$ contains, in
particular, a generalization of the result of the first author in
\cite{Tohoku}, that $H^\sigma$ is connected if $H$ is simply
connected.

With the preliminary work about centralizers of $c$-pairs
established, we turn in Section 9  to $c$-triples, for $c$
non-trivial. Following the 
general pattern, we first consider the case of $c$-triples of rank
zero. It turns out there there is a finite and short list of
simple groups that have rank zero
$c$-triples. We then go on to establish Parts 1 through 4 of
Theorem~\ref{ctrip} in the case when $\langle C\rangle$ is
cyclic. 

In Section 10 we turn to the tori $\ov S(k)$ and $\ov
S^{w_c}(\ov {\bf g}, k)$ and compute their Weyl groups and
fundamental groups. The result is a proof of
Theorem~\ref{quotdiag}, of Part 3 of Theorem~\ref{commuttrip} and
of Part 3 of Theorem~\ref{ctrip} in the case when $C$ is cyclic.
Once again, we find a related Coxeter group by arguments which
are formally similar to those used to study the fixed subspace of
an affine automorphism.

Section 11 considers the Chern-Simons invariant of a flat 
connection with holonomy a given $c$-triple.
This invariant is identified with an invariant of the $c$-triple
defined using the Weyl invariant inner product on the   Cartan
subalgebra of $G$. We then prove Theorem~\ref{CS}.
Lastly, we establish Witten's clockwise symmetry statement,
Theorem~\ref{clock}, in the case of $c$-triples.

In Section 12 we consider the case  when
the group $\langle C\rangle$ is non-cyclic, establishing both
Theorems~\ref{ctrip} and~\ref{clock} in these cases by
explicit computation.

At the end of the paper, we give a list of the possible coroot
diagrams and quotient coroot diagrams, and give tables
summarizing the tori and  root systems we have defined, as well as
other relevant information.  

A guiding principle of this paper has been to avoid
classification and case-by-case checking wherever possible. We
have preferred to give more general, conceptual arguments, even
at the cost of increasing the length.

\subsection{History}

Questions related to the ones   considered here have a long
history. One motivation was to understand the experimental
connection between the torsion primes of $G$, i.e. the primes $p$
for which there is $p$-torsion in the integral homology of $G$,
and those primes $p$ for which $G$ has an elementary abelian
$p$-group of rank three contained in a torus. The first author
and J-P. Serre proved that if $p$ is not a torsion prime, then
every elementary $p$-group in $G$ is contained in a torus. The
converse was checked, case-by-case,  in \cite{Tohoku} and in the
earlier references cited in that paper. From the point  of view of
this paper, the subgroups constructed in \cite{Tohoku} are simply
the  commuting triples of rank zero. It was also checked
that a torsion prime divides one of the root integers
$h_a$, but the converse does not quite hold. R. Steinberg pointed
out later  that the torsion primes listed in \cite{Tohoku}
are exactly the prime divisors of the coroot integers
$g_a$. In \cite{Stein}, he shows that the prime $p$ divides one of the
$g_a$ if and only if $G$ contains an elementary abelian
$p$-subgroup not contained in a torus.  The methods of the
present paper are, in part, quite similar to those of \cite{Stein}.  
More recently, Griess \cite{Griess} classified the possibilities
for elementary abelian
$p$-groups in $G$, apparently by completely different methods
from those of this paper. The last two authors of the present
paper were led to this circle of problems from another direction,
the study of holomorphic principal $G$-bundles over an elliptic
curve. In this case, if the bundle is flat, then its
holonomy defines a $c$-pair up to conjugation. The
singularities of the moduli space of such bundles are
then closely related to the component groups of centralizers of
$c$-pairs, and thus ultimately to $c$-triples.  We should
also add that   Kac-Smilga in \cite{Kac-Smilga} as well as
Keurentjes \cite{Keur1,Keur2} have independently established
results overlapping significantly with ours in the case of
commuting triples. Their approach seems very similar to ours.

It is natural to ask about the moduli space of almost commuting
$N$-tuples for $N > 3$. We believe that the methods developed here
can attack this question as well. However, and perhaps not
surprisingly, there seem to be very few essentially new cases.

The impetus for the study we carry out here
of commuting triples and $c$-triples was questions, conjectures, 
and statements of Witten about the nature of the moduli space, 
and especially the number and structure of its components. He was
led to these questions by studying the quantum field theory of
gauge theories over the three torus, in particular the
$R$-symmetries of these theories. With this heuristic guide, he 
conjectured the clockwise-symmetry statement and checked it in
many cases. It is our pleasure to thank Witten for pointing out
these questions to us, and for many stimulating discussions on
these and other related matters.

\section{Almost commuting $N$-tuples}

An ordered $N$-tuple ${\bf x}=(x_1,\ldots,x_N)$  of elements in $G$ is
{\sl almost commuting} if $[x_i,x_j]\in {\cal C}G$ for every $1\le
i,j\le N$.
Notice that ${\bf x}$ is almost commuting if and only if the image
ordered $N$-tuple $\ov {\bf x}=(\ov x_1,\ldots,\ov x_N)$ in the 
adjoint form of $G$ is a commuting ordered $N$-tuple.
 Given an almost commuting $N$-tuple ${\bf x}$ let $c_{ij}=[x_i,x_j]$ and
 let $C$ be the $N\times N$ matrix $C=(c_{ij})$
of elements in ${\cal C}G$. The matrix $C$ is antisymmetric in the  sense
that $c_{ii}=1$ and $c_{ij}=c_{ji}^{-1}$ for all $i,j$. We say that
${\bf x}$ is an {\sl ordered $N$-tuple of $C$-type}. 

Clearly, the space $\widetilde{\cal M}_G(C)$ of all ordered
$N$-tuples in  $G$ of $C$-type is
identified with a closed subspace of $\prod_{i=1}^NG$ and thus
is a compact Hausdorff space.
The compact group $G$ acts on $\prod_{i=1}^NG$ by simultaneous conjugation
normalizing the subspace $\widetilde{\cal M}_G(C{})$. Hence, the
quotient  ${\cal M}_G(C{})=\widetilde 
{\cal M}_G(C{})/G$ is a compact Hausdorff space. It is the space of
conjugacy classes of ordered $N$-tuples in $G$ of $C{}$-type. When $G$
is clear from 
context, we shall denote this space by ${\cal M}(C{})$. Our goal in this
section is to prove a very general qualitative result concerning ${\cal
M}(C{})$. We shall show that its connected components are homeomorphic
to quotients of products of subtori of $T$ by finite groups.

\subsection{An invariant for almost commuting $N$-tuples}

We define the {\sl rank} of an ordered $N$-tuple  ${\bf x}$, denoted
${\rm rk}({\bf x})$, to be the rank of 
$Z({\bf x})$, the centralizer of ${\bf x}$ in $G$.
Notice that ${\bf x}$ is of rank zero if and only if $Z({\bf x})$ is a
finite group.
There is a related but finer invariant of ${\bf x}$ derived from any
maximal torus of $Z({\bf x})$ which we shall now describe.

For any subset $I\subseteq \Delta$ let ${\frak t}_I=\bigcap_{a\in
 I}{\rm Ker}\,a\subseteq {\frak t}$, and let  
 $S_I$ be the subtorus of $T$ with Lie algebra ${\frak t}_I$. We
denote by $L_I$ the derived group of  $Z(S_I)$:
$$L_I=DZ(S_I).$$
Since $Z(S_I)$ is the centralizer of a torus, it is connected, and thus
$L_I$ is also connected.

\begin{lemma}\label{subset}
\begin{enumerate}
\item $L_I$ is non-trivial if and only if $I\neq \Delta$.
\item If $L_I$ is non-trivial, then $I$ forms a set of simple roots for
$L_I$.
\item $L_I$ is simply connected.
\item $S_I$ is the component of the identity of the center of its
  centralizer. 
\item Conversely, if $S$ is a torus in $G$ and is equal to the
  component of the identity of the center of its centralizer, then
  there is a subset   $I\subseteq \Delta$ such that $S$ is conjugate to
  $S_I$. 
\end{enumerate}
\end{lemma}

\begin{proof}
All of these elementary facts are proved in \cite{Bour}.
\end{proof}

\begin{lemma}
Let ${\bf x}$ be an ordered $N$-tuple and let $S$ be a maximal torus of
$Z({\bf x})$. 
There is a unique subset $I\subseteq
\Delta$ such that $S$ is conjugate to $S_I$.
\end{lemma}

\begin{proof}
Let $S'$ be the component of the identity  of the center  of $Z(S)$.
Clearly, since $S$ is contained in the center of $Z(S)$, $S\subseteq
S'$. Since ${\bf x}\subseteq Z(S)$, we have 
$S'\subseteq Z({\bf x})$. Since $S$ is a maximal torus of $Z({\bf x})$,
this implies that $S'=S$. Since $S$ is the identity component of its
centralizer,  $S$ is conjugate to $S_I$ for some $I\subseteq \Delta$.
Uniqueness of $I$ is clear.
\end{proof}

The subset $I\subseteq \Delta$ given in the last lemma is an invariant
of the ordered $N$-tuple  ${\bf x}$ and is denoted by $I({\bf x})$. The
cardinality of $\Delta-I({\bf x})$ is the rank of ${\bf x}$. 

\begin{lemma} The torus $S$ is a maximal torus of $G$ if and only if the
components $x_i$ are contained in $S$, and in this case $C ={\rm Id}$.
\end{lemma}
\begin{proof} If the $x_i$ are all contained in the torus $S$, then they
are mutually commuting. Moreover, if $T$ is a maximal torus containing
$S$, then $T\subseteq Z({\bf x})$. Thus, since $S$ is a maximal torus of 
$Z({\bf x})$, $S=T$.
Conversely, suppose that $S$ is a maximal torus of $G$. Since each $x_i$
commutes with $S$, it lies in $Z_G(S) =S$.
\end{proof}

\subsection{The case of rank zero}

There is a general finiteness  result in this case:

\begin{proposition}\label{modfinite}
Let $K$ be a compact group with finite center.
{\rm (}We do not assume that $K$ is connected nor that the component of
the identity is simply connected.{\rm )}
Up to conjugation, there are only finitely many  almost commuting
ordered $N$-tuples of   rank zero in $K$.
Thus, for each anti-symmetric $N\times N$-matrix $C=(c_{ij})$
with coefficients in ${\cal C}K$, the moduli space ${\cal M}_K(C)$ is a
finite set. 
\end{proposition}

\begin{proof}
The proof is by induction on $N$. The case $N=1$ is deduced from the
following lemma.

\begin{lemma}\label{f}
Suppose that ${\cal C}K$ is
finite. There are only finitely 
many conjugacy classes of elements $g\in K$ for which the center of
$Z(g)$  is finite.
In case $K$ is connected and simply connected, and thus semi-simple, the
center of
$Z_K(x)$ is finite if and only if $x$ is conjugate to the exponential of a
vertex of an alcove $A$.
\end{lemma}

\begin{proof}
Fix a maximal torus $T\subseteq K^0$. Let $W(T,K)$ be the normalizer of
$T$ in $K$ modulo its centralizer. It is a finite group, whose action on
$T$ is covered by a   linear action on ${\mathfrak t}$. 

Fix $g\in K$. 
According to \cite{deS}, II \S 3 Proposition $2$ (see also
\cite{BM})  there is a regular element in 
${\frak g}$ fixed under ${\rm Ad}\,g$. Thus, after conjugation we can
assume that $g$ normalizes $T$ and
a positive Weyl chamber  $C_0\subseteq {\mathfrak t}$. This
implies  that 
$T_0=(T^g)^0$ contains a   regular point of $T$.  The torus $T_0$ is a
maximal torus 
of $K^g$ and the normalizer of $T_0$ in $K$ is contained in
 the normalizer of $T$. Let us suppose that the center of $Z(g)$ is finite.
This implies that $W_0=W(T_0, Z(g))$ acts on $T_0$ and that there are
only finitely many points fixed by the action of $W_0$. Equivalently,
$W_0$ acts on the quotient torus $\ov T_g=T/({\rm Id}-{\rm Ad}\,g)T$ with
only finitely many points fixed by  $W_0$.

It suffices to show that there are only finitely many conjugacy classes
of elements $g'\in K$ which (i) are congruent to $g$ modulo $K^0$, (ii)
such that ${\rm Ad}\,g'|T={\rm Ad}\,g|T$ and
(iii) have the property that $W_0'=W(T_0, K^{g'})$ is equal to $W_0$.

Any element $g'\in K$ satisfying (i) and (ii)  is  an element
of the form $tg$ for some $t\in T$. The conjugacy class of $tg$ in $K$
depends only on $[t]\in \ov T_g$.

For each $w\in W_0$ there is an element $h_w\in N_K(T_0)\subseteq
N_K(T)$ such that 
$[h_w,g]=1$ and $[h_w]=w\in W_0$.
Suppose that $g' = tg$ satisfies (iii) for some $t\in T$. Then, for all
$w\in W_0$ there exists $h_w'=t_wh_w$, $t_w\in T$, commuting with
$tg$. Then 
$$t_wh_wtgh_w^{-1}t_w^{-1}=tg$$
or
$$t_w{}^w(t){}^g(t_w^{-1})=t.$$
Thus  $[t]\in \ov T_g$ is fixed by $w$ for every $w\in W_0$.
This means that there are only
finitely many possibilities for $[t]\in T_g$ and hence only finitely
many possibilities for the conjugacy class of $tg$ in $K$.

In the case when $K$ is connected and simply connected, $Z(x)$ is
connected.
Thus, it has a finite center if and only if it is semi-simple.
This occurs exactly when $x$ is conjugate to the exponential of a vertex
of an alcove $A$.
\end{proof}

There is the following which we shall need later:

\begin{lemma}\label{vertex}
Let $K$ be  connected and simply connected,  let
$(x_1,\ldots,x_N)$ be a rank zero subset of almost commuting elements of
$K$, and suppose that $[x_1, x_i] =1$ for all $i$. Then $Z(x_1)$ is
semi-simple and $x_1$ is conjugate in
$K$ to the exponential of a vertex of an alcove. 
\end{lemma}

\begin{proof}
The $(N-1)$-tuple $(x_2,\ldots,x_N)$ is a rank zero
almost commuting $(N-1)$-tuple in 
$Z(x_1)$. Hence the center of $Z(x_1)$ is finite, so that $Z(x_1)$ is
semi-simple. It follows from Lemma~\ref{f} that $x_1$ is conjugate to the
exponential of a vertex of an alcove.
\end{proof}

\noindent{\bf Proof of Proposition~\ref{modfinite}.}

\noindent
{\bf Case $N=1$.}
An almost commuting $1$-tuple is  a single element  $g\in K$.
Rank zero means that 
$Z(g)$ is finite and {\sl a fortiori\/} that its center is finite.
Thus, the result in this case is immediate from Lemma~\ref{f}.

\noindent
{\bf General Case.}
We prove the general case by induction on $N$.
Suppose that the result is known for all groups with finite center
and for all almost commuting $k$-tuples for $k<N$. Consider an ordered
almost commuting $N$-tuple 
${\bf x}=(x_1,\ldots,x_N)\subseteq K$ of rank zero.
Let $\hat Z(x_N)$ be the subgroup of elements of $K$ whose
commutator with $x_N$ lies in ${\cal C}K$.
Since ${\bf x}\subseteq \hat Z(x_N)$ and ${\bf x}$ is of rank
zero,  we see that the center of $\hat Z(x_N)$ must be finite.
Let $\ov K$ be the quotient of $K$ by its center and let $\ov x_N$ be
the image of $x_N$ in $\ov K$. Then we have an exact sequence
$$\{1\}\to {\cal C}K\to \hat Z(x_N)\to Z_{\ov K}(\ov x_N)\to \{1\}$$
and it follows that the center of $Z_{\ov K}(\ov x_N)$ is finite.
Applying the previous lemma to $\ov K$, we see that there are only
finitely many possibilities for $\ov x_N$ up to conjugation in $\ov K$.
Hence, there are only finitely many 
possibilities  for $x_N\in K$ up to conjugation.
Let ${\bf x}'=(x_1,\ldots,x_{N-1})$. 
This is an ordered almost commuting $(N-1)$-tuple in 
$\hat Z(x_N)$.   
Consider the center of the centralizer $Z$ of ${\bf x}'$ in $\hat
Z(x_N)$. Clearly, $Z({\bf x})=Z\cap Z(x_N)$ is a subgroup of finite
index of $Z$. 
Since the center of $Z({\bf x})$ is finite, it follows that the center
of $Z$ is finite. This means that ${\bf x}'$ is of rank zero.
Thus, by the inductive hypothesis, there are only
finitely many possibilities for ${\bf x}'$ up to conjugation in $\hat
Z(x_N)$, and hence only finitely many possibilities for ${\bf x}'$ in
$\hat Z(x_N)$ up to conjugation by $Z(x_N)$. This completes the
inductive step. 
\end{proof}

\subsection{The case of arbitrary rank}

In this section we return to the group $G$, which   is  connected
and  simply connected and thus semi-simple. Let $C=(c_{ij})$ be an
antisymmetric
$N\times N$ matrix with coefficients in ${\cal C}G$. 

Fix a subset $I\subseteq \Delta$.
Suppose that $L_I$ contains
all the elements $c_{ij},\, 1\le i,j\le N$.
Let ${\cal M}^0_{L_I}(C{})\subseteq {\cal M}_{L_I}(C{})$ be the subset of
conjugacy classes 
of rank zero ordered $N$-tuples of $C$-type in $L_I$.
Let $\widetilde{\cal M}^0 _{L_I}(C{})$ be a set of representatives for
the finite set ${\cal M}^0_{L_I}(C{})$.

Let $F_I=S_I\cap L_I$. This is a finite subgroup of the center of $L_I$.
Consider the action of $F_I^N$ on $L_I^N$ given by
$$(f_1,\ldots,f_N)\cdot (y_1,\ldots,y_N)=(f_1y_1,\ldots,f_Ny_N).$$
Clearly, this operation does not change the pairwise commutators, nor
does it change the centralizer in $L_I$ of the subset. Hence, it
defines an action of $F_I^N$ on ${\cal M}^0_{L_I}(C{})$.

Then we have a map
$$S_I^N\times {\cal M}^0_{L_I}(C{})) \to {\cal M}_G(C{})$$
which associates to $(s_1,\ldots,s_N)\times r$ the
conjugacy class of the ordered $N$-tuple 
$${\bf x}=(s_1\tilde r_1,\ldots,s_N\tilde r_N),$$
where $(\tilde r_1,\ldots,\tilde r_N)$ is the chosen representative in
$\widetilde{\cal M}_G(C)^0$
for the element $r\in {\cal M}^0_{L_I}(C{})$.
It is clear from the definitions that this map is independent of the
choice of representatives $\widetilde {\cal M}^0_{L_I}(C{})$ for the
conjugacy 
classes ${\cal M}^0_{L_I}(C{})$ and that it factors to define a
continuous map
$$p\colon S_I^N\times_{F_I^N} {\cal M}^0_{L_I}(C) \to
{\cal M}_G(C{}).$$

The group $N_G(S_I)$ acts by conjugation normalizing $S_I$ and hence $L_I$.
Thus, it acts on ${\cal M}^0_{L_I}(C{})$. The group $Z_G(S_I)$ acts on
$L_I$ by inner automorphisms of $L_I$ and hence acts trivially on ${\cal
M}^0_G(C{})$. Thus, we have an induced action of
$W(S_I,G)=N_G(S_I)/Z_G(S_I)$ on 
$$S_I^N\times_{F_I^N} {\cal M}^0_{L_I}(C{}).$$ 
Clearly, the map $p$ factors through this action.

\begin{theorem}\label{main}
  For  $I\subseteq \Delta$, let
${\cal M}^I_G(C)\subseteq {\cal M}_G(C{})$ be the subspace of
conjugacy classes of ordered $N$-tuples  ${\bf x}$ of
$C$-type in $G$ whose centralizer has a maximal torus which is
conjugate to $S_I$.
The map $p$ induces a homeomorphism
$$\ov p\colon \left(S_I^N\times_{F_I^N} {\cal
    M}^0_{L_I}(C{}))\right)/W(S_I,G)\to {\cal M}_G^I(C).$$
In particular, ${\cal M}^I_G(C)$ is a
compact Hausdorff space and hence a closed subset of ${\cal M}_G(C)$.
\end{theorem} 

\begin{corollary}\label{cor}
  For each $I\subseteq \Delta$ the subset ${\cal M}^I_G(C{})$  is a union
of components of ${\cal M}_G(C{})$. In particular, 
the rank of ${\bf x}$ and the subset $I({\bf x})$ are locally constant
functions on ${\cal M}_G(C{})$. 
\end{corollary}

\noindent
{\bf Proof of Corollary~\ref{cor}.}
According to the theorem the subset ${\cal M}^I_G(C{})$ 
is compact and hence is a closed subset of ${\cal M}_G(C{})$.
Since ${\cal M}_G(C{})$ is a disjoint union of the ${\cal M}^I_G(C{})$
for the various $I\subseteq \Delta$ and since these subsets are
closed and finite in number, each is  a union of components.
\relax\ifmmode\expandafter\endproofmath\else
  \unskip\nobreak\hfil\penalty50\hskip.75em\hbox{}\nobreak\hfil\bull
  {\parfillskip=0pt \finalhyphendemerits=0 \bigbreak}\fi

\noindent{\bf Proof of Theorem~\ref{main}.}
  Fix a subset $I\subseteq \Delta$.
Suppose that ${\bf x}$ is an ordered $N$-tuple of  $C{}$-type in $G$ and
that a maximal torus $S$ 
for $Z({\bf x})$ is conjugate to $S_I$.
Conjugating ${\bf x}$ we can assume that $S=S_I$.
Then ${\bf x}\subseteq Z(S_I)=S_I\cdot L_I$.
Thus, we can write
${\bf x}=(s_1y_1,\ldots,s_Ny_N)$ with $s_i\in S_I$ and $y_i\in L_I$.
It follows that $c_{ij}\in L_I$ for all $1\le i<j\le N$ and that
${\bf x}'=(y_1,\ldots,y_N)$ is an ordered $N$-tuple of $C$-type in $L_I$. 
Since $S\cdot Z_{L_I}({\bf x}')\subseteq Z({\bf x})$, the $N$-tuple ${\bf
x}'$  is of rank zero in $L_I$.

This shows that the map $\ov p$ is onto the subset ${\cal M}^I_G(C{})$
of ${\cal M}_G(C{})$. 
To prove that $\ov p$ is a homeomorphism onto ${\cal M}^I_G(C)$ we
need only show that $\ov p$ is one-to-one. 
Suppose that $\ov p((s_1,\ldots,s_N),{\bf x})
=\ov p((s_1',\ldots,s_N'),{\bf x}')$ where ${\bf x}=(y_1,\ldots, y_N)$
and ${\bf x}'=(y_1',\ldots, y_N')$ are elements of $\widetilde
{\cal M}^0_{L_I}(C{})$.
Then there is  $g\in G$  such that
$$g(s_1y_1,\ldots, s_Ny_N) g^{-1}=(s_1'y_1',\ldots,s_N'y_N').$$
Thus $gS_Ig^{-1}$ and $S_I$ are both maximal tori for
$Z=Z(s_1'y_1',\ldots,s_N'y_N')$ and hence there is an element in
$h\in Z$ such that $hS_Ih^{-1}=gS_Ig^{-1}$. Replacing $g$ by $h^{-1}g$
allows us to assume that $g\in N_G(S_I)$. 
Then $g{\bf x} g^{-1}$ is a rank-zero ordered $N$-tuple of  $C$-type in
$L_I$. Also,  
$gy_ig^{-1}=(gs_i^{-1}g^{-1}s_i')y_i'$, implying that
$f_i=gs_i^{-1}g^{-1}s_i'\in F_I$. 
Let $f\in F_I^N$ be the given by $(f_1,\ldots, f_N)$.
Then there is an element of $L_I$ which
conjugates $g{\bf x} g^{-1}$ to the given representative for the conjugacy
class of $f\cdot {\bf x}'$ in ${\cal M}^0_{L_I}(C{})$.
This proves that $((s_1,\ldots,s_N),{\bf x})$ and
$((s_1',\ldots,s_N'),{\bf x}')$ have the same image in 
$(S_I^N\times_{F_I^N}{\cal M}^0_{L_I}(C{}))/W(S_I,G)$.
\relax\ifmmode\expandafter\endproofmath\else
  \unskip\nobreak\hfil\penalty50\hskip.75em\hbox{}\nobreak\hfil\bull
  {\parfillskip=0pt \finalhyphendemerits=0 \bigbreak}\fi

\begin{corollary}\label{determ}
There is a component of ${\cal M}_G(C{})$ consisting of conjugacy
classes of  ordered $N$-tuples ${\bf x}$ of $C{}$-type
with 
$I({\bf x})=I$ if and only if {\rm (i)} $c_{ij}\in L_I$ for all $1\le
i,j\le N$ and {\rm (ii)} $L_I$ contains an ordered $N$-tuple of 
$C{}$-type of rank zero. The set of these components is identified
with $\left({\cal M}^0_{L_I}(C{})/F_I^N\right)/W(S_I,G)$. The
dimension of each such 
component is $N\cdot\#(\Delta-I)$.
  For each component $X$ of ${\cal M}^I_G(C{})$
let  $F_I^N(X)\subseteq F_I^N$  be the stabilizer of the
corresponding point of ${\cal M}^0_{L_I}(C{})$. 
Then  $X$ is homeomorphic to
$$\left(S_I^N/F_I^N(X)\right)/W(S_I,G).$$
\end{corollary}

There is a similar but somewhat more involved statement of the corollary
in case $G$ is not connected.

\section{A characterization of groups of type $A_n$}

\subsection{Action of ${\cal C}G$ on an alcove}\label{action1}
Let $G$ be  simple, and let $A\subseteq{\frak
t}$ be the alcove containing the origin determined by the set of simple
roots 
$\Delta$. Recall that the walls of $A$ correspond to the extended simple
roots
$\widetilde
\Delta$. We then have the following well-known lemma, whose proof is left
to the reader:

\begin{lemma}\label{alcove}
Let $\tilde x \in A$ and let $x =\exp \tilde x \in T$. Let $\Phi(x)$
be the set of
$a\in \Phi$ such that $a(\tilde x) \in {\bf Z}$. Finally let $\widetilde
I(x)$ be the set of extended simple roots $a$ such that $\tilde x$
lies in the wall $W_a$ of $A$ corresponding to $a$. Then
$\Phi(x)$ is a closed sub-root system of
$\Phi$, and $\widetilde I(x)$ is a set of simple roots for
$\Phi(x)$. Thus $\widetilde I( x)$ corresponds to a proper subdiagram of
$\widetilde D(G)$, and moreover every proper subdiagram of $\widetilde
D(G)$ corresponds to a subset $\widetilde I(x)$ for some $x=\exp
\tilde x$ with $\tilde x
\in A$.
\end{lemma}

There is the following lemma on the relationships of the coroot lattices.

\begin{lemma}\label{coroottors}
Let $\widetilde I$ be a proper subset of $\widetilde\Delta$, and let
$Q_{\widetilde I}\spcheck$ be the sublattice of $Q\spcheck$ spanned by
the coroots $a\spcheck$, $a\in \widetilde I$. Let $k=\{\, \gcd g_a:
a\notin
\widetilde I\,\}$. Then the torsion subgroup of
$Q\spcheck/Q_{\widetilde I}\spcheck$ is cyclic of order $k$, and a
generator is
$$\zeta = -\frac{1}{k}\sum_{a\in \widetilde I}g_aa\spcheck.$$
\end{lemma}

\begin{corollary}\label{pi1}
Let $x={\rm exp}(\tilde x)$ for some $\tilde x\in A$. Let
$\widetilde  I(x)$ be the subset of $\widetilde \Delta$
consisting of all $a\in \widetilde \Delta$ such that the
corresponding wall
$W_a$ contains $\tilde x$.  Then $\widetilde  I(x)$ is a set
of simple roots for $DZ(x)$. Moreover $\pi_1(DZ(x))$ is a cyclic group of
order
$k={\rm gcd}_{a\in
\widetilde \Delta- \widetilde I(x)}(g_a)$.
A generator for the fundamental group, viewed as a central element in the
simply connected covering group $\widetilde {DZ}(x)$, is ${\rm
exp}(\zeta)$ where  
$$\zeta=-\frac{1}{k}\sum_{a\in \widetilde
I(x)}g_aa\spcheck.$$  
\end{corollary}

The affine Weyl group $W_{\rm aff}(G)$ of $G$ with
respect to $T$ acts as a group of affine linear transformations of
${\frak t}$ with the alcove $A$ as fundamental domain.
  For each element
$c\in {\cal C}G$ there is a unique point $\zeta_c\in A$ with ${\rm
  exp}(\zeta_c)=c\in T$. 
In fact, $\zeta_c$ is a vertex of $A$ and if $\{a=0\}$
defines the wall opposite this vertex then, in Equation~\ref{highest},
the root integer $h_a$ is equal to $1$. 
This identifies the center of $G$ with the subset of $a\in
\Delta$ for which $h_a=1$.

Let $A'=A-\zeta_{c^{-1}}$. Then $A'$ is another alcove containing
$0$. Hence there is a unique element $w_c\in W(T,G)$ with the
property that
$w_c\cdot A=A'$.
The map 
$$\varphi_c(t)= w_c\cdot( t-\zeta_{c^{-1}}), $$
is an affine linear transformation of ${\frak t}$ carrying $A$ to
itself. 
We denote its fixed point set by ${\frak t}^c$ and we set
$A^c=A\cap{{\frak t}}^c$.
The map $c\mapsto \varphi_c$ 
defines a homomorphism  from ${\cal C}G$ to the group of
affine   linear automorphisms of $A$. This map is called {\sl the
  action of   ${\cal C}G$ on the alcove $A$}.
The element $w_c\in W$ is called the {\sl Weyl part of the action of
  $c$}.
We define a homomorphism $\nu = \nu_A\colon {\cal C}G\to W(T,G)$ by
associating to $c$ the element $w_c$.
Notice that $\nu$ is an injective
homomorphism ${\cal C}G\to W(T,G)$. Its image is the
stabilizer in $W$  of $\widetilde \Delta$ and acts simply
transitively on the set of $a\in \widetilde \Delta$ such that
$h_a=1$. If $\gamma_c\in N_G(T)$ projects to $w_c\in W$ and
$x=\exp\tilde x$ for some $\tilde x\in A^c$, then
\begin{equation}\label{conjeqn}
\gamma_c\cdot x \cdot \gamma_c^{-1} = c^{-1}x.
\end{equation}
If $A'$ is another alcove containing the origin, then $\nu_{A'}$
is conjugate to $\nu_A$. More generally, for any alcove
$A'$, not necessarily containing the origin, there is an action of ${\cal
C}G$ on $A'$ by elements of the affine Weyl group.

Notice that $A^c$ contains the barycenter of the alcove $A$ and hence
is non-empty. Thus,
the dimension of $A^c$ is equal to the dimension of ${\frak t}^c$
which in turn is equal to the dimension of the 
fixed point subspace ${\frak t}^{w_c}$ for the linear action of $w_c$
on ${\frak t}$.
The affine action $\varphi_c$ permutes the codimension-one faces of
$A$, and hence $w_c\in W$ permutes the roots in $\widetilde \Delta$ and
induces a diagram automorphism of $\widetilde D(G)$. In fact, every
diagram automorphism of $\widetilde D(G)$ is induced by a permutation of 
$\widetilde \Delta$ of the form $w_c\circ \sigma$, where $c\in {\cal C}G$
and $\sigma$ is an outer automorphism of $G$, which can be identified
with a diagram automorphism of $D(G) \subseteq \widetilde D(G)$.   Since
there is only one integral linear relation among the $a\in
\widetilde \Delta$ with positive  coefficients with no common factor, 
namely $\sum_{a\in\widetilde \Delta}h_aa=0$,
the action of $w_c$ preserves this
relation.
Thus, $h_a=h_{w_c\cdot a}$ for all $a\in \widetilde \Delta$.
  For the same reason, the action of $w_c$ preserves the linear relation
$\sum_{a\in\widetilde \Delta}g_aa\spcheck$, i.e. $g_a=g_{w_c\cdot
a}$ for all $a\in \widetilde \Delta$.

The above constructs an action of the center of $G$, or equivalently the
fundamental group of the adjoint form of $G$, on the alcove $A$. More
generally, let
$H$ be a compact connected group with maximal torus
$T$. Let $\Lambda = \pi_1(T)$. Then the coroots inverse to the roots
$\Phi_H$ of
$H$ with respect to $T$ span a subspace ${\frak d}\subseteq {\frak
t}$, and $\Phi_H$ is a root system on $\frak d$. Let $Q\spcheck_H$ be
the sublattice of $\Lambda$ spanned by the coroots and let $P\spcheck
_H$ be the corresponding coweight lattice. Then $\pi_1(H) \cong
\Lambda/Q\spcheck_H$, and this group can be identified with a subgroup
of the center of
$\widetilde H$, the universal covering group of $H$. If
$\pi$ is orthogonal projection from $\frak t$ to $\frak d$ under a Weyl
invariant inner product, then $\pi(\Lambda) \subseteq P\spcheck_H$. Thus
there is an induced homomorphism
$$\Lambda/Q\spcheck_H  \to P\spcheck_H/Q\spcheck_H.$$

We
consider the  decomposition
${\cal A}_H$ of
${\frak t}$ into alcoves under the walls of the affine Weyl group $W_{\rm
aff}(H)$.  If ${\frak z}$ is the Lie
algebra of the center of $H$, then the alcoves in ${\cal A}_H$ are of
the form $A'\times {\frak z}$ where $A'$ is an alcove for the
affine Weyl group of the derived group $DH$ with respect to its
maximal torus $T\cap DH$.  
Thus, the alcoves in ${\frak t}$ are compact if and only if $H$ is
semi-simple. 
Since $W_{\rm aff}(H)=W_{\rm aff}(DH)$, it follows that this group
acts simply transitively on the set of alcoves in ${\frak t}$.

  Fix a set of simple roots $\Delta_H$ for the root system of $H$ with
respect to
$T$ and let $\widetilde\Delta_H$ be the corresponding set of extended
roots. Let $A\subseteq {\frak t}$ be the alcove determined by $\Delta_H$.
This alcove containes the origin.
Conversely, given any alcove $A'$ containing the origin, it
corresponds to a set of simple roots $\Delta'$ for $\Phi_H$. 
Exactly as in the semi-simple case, the center of $\widetilde H$
acts as a group of affine linear isometries of $A$.
The linear part of this
automorphism defines a homomorphism $\nu \colon {\cal C}\widetilde
H\to W(H)$. 
The only difference with the semi-simple case is that $\nu$
is not injective -- its kernel is the identity component of   ${\cal
C}\widetilde H$.
  For each $c\in {\cal C}\widetilde H$, we let $w_c=\nu(c)$.
The element $w_c$ acts on the roots of $H$, normalizing the set
$\widetilde\Delta_H$. The action on $\widetilde\Delta_H$ is the one
induced by the action of
$c$ on the walls of $A$. Using the homomorphism $\Lambda \to {\cal
C}\widetilde H$ described above, we can also view $\Lambda$ as acting on
$A$.

\begin{lemma}\label{char}
Let $H$ be a
compact group with maximal torus $T$ and let 
$\varphi\colon {\frak t}\to {\frak t}$ be an affine linear map whose
translation part is given by an element  $v\in {\frak t}$ which 
exponentiates to $c\in {\cal C}\widetilde H$ and whose linear
part is an  element of the Weyl group
$W(T,H)$. 
If there is an alcove $A\subseteq {\frak t}$ for the affine Weyl 
group of
$H$  such that $\varphi(A) =A$, then $\varphi$ is the action of $c$ on
this alcove. In particular, its linear part is $w_c$.
\end{lemma}

\begin{proof}
Consider the composition of $\varphi^{-1}$ and the affine linear map
which is the action of $c$ on $A$. The translation part of this map is
given by an element of ${\frak t}$ which exponentiates to the identity
in $\widetilde H$
and hence is contained in the coroot lattice for $\widetilde H$. The 
linear part is a composition of elements of the Weyl group and hence is
an element of the Weyl group. That is to say this composition is an
element of the affine Weyl group of $H$. Of course, it sends $A$ to $A$.
This means that it is the trivial element of the affine Weyl group.
\end{proof}

  For future reference, we shall also need the following lemma on the
stabilizer of a point of $T$ under the action of the Weyl group.

\begin{lemma}\label{stabilizer}
Fix an alcove $A\subset \frak t$. Let $t\in T$ and let $\tilde t\in A$
be a lift of $t$. Let
$W= W(T,H)$ and let
$W(t)$ be the Weyl group generated by the roots $a$ such that $a(\tilde
t) \in {\bf Z}$, in other words by the root system $\Phi_H(t)$. Then the
map
$\nu$ induces an isomorphism
$${\rm Stab}_{\Lambda/Q_H\spcheck}(\tilde t) \cong {\rm Stab}_W(t)/W(t).$$
Moreover, if $I$ is the image of ${\rm Stab}_{\Lambda/Q_H\spcheck}(\tilde
t)$ under
$\nu$, then $I = {\rm Stab}_W(t) \cap {\rm Im} \,\nu$ and ${\rm Stab}_W(t)
\cong W(t)
\rtimes I$.
\end{lemma}
\begin{proof} Suppose $c \in \Lambda/Q_H\spcheck$. If
$\varphi_c$ stabilizes $\tilde t$, then $\tilde t =
w_c(\tilde t-\zeta_{c^{-1}})$ for some $\zeta_{c^{-1}} \in
\Lambda$ projecting to
$c^{-1}$, and so 
$w_c(t) = t$. Thus $\nu$ defines a homomorphism from ${\rm
Stab}_{\Lambda/Q_H\spcheck}(\tilde t)$ to ${\rm Stab}_W(t)$. If in the
above notation
$w_c\in W(t)$, then by Lemma~\ref{alcove}, $w_c$ is in the group
generated by reflections in the roots which are integral on 
$\tilde t$. Thus $w_c(\tilde t) = \tilde t + \lambda$, where
$\lambda \in Q_H\spcheck$. From
$$w_c(\tilde t) = \tilde t +w_c(\zeta_{c^{-1}})  = \tilde t +
\lambda,$$ it follows that $w_c(\zeta_{c^{-1}}) \in Q_H\spcheck$
and hence that $\zeta_{c^{-1}}\in Q_H\spcheck$. Thus, $c^{-1}\in
\Lambda/Q_H\spcheck$ is trivial, so that $c$ is trivial as well.
This shows that
$\nu|{\rm Stab}_{\Lambda/Q_H\spcheck}(\tilde t)$ is injective, and, if
$I$ is its image, then $I\cap W(t) = \{1\}$.

Now suppose that $w\in W$ fixes $t$. Thus $w(\tilde t) =\tilde t+ \tilde
c$ for some $\tilde c \in \Lambda$. Hence $A$ and $A'=w(A)-\tilde c$ are
two alcoves, both of which contain $\tilde t$. After transforming $A'$
by an element in the group generated by reflections in the walls of
$A$ containing $\tilde t$, we can then assume that $A'=A$. By
Lemma~\ref{char}, it follows that $w=w_c$. This says that every element
of ${\rm Stab}_W(t)$ can be written as a product of an element in ${\rm
Im}\,\nu$ times an element of the group generated by reflections in the
walls of
$A$ containing $\tilde t$. But by Lemma~\ref{alcove}, this second group
is exactly $W(t)$. This proves that $I = {\rm Stab}_W(t) \cap {\rm Im}\,
\nu$ and that ${\rm Stab}_W(t) = W(t)\cdot I$. Since $W(t)$ is a normal
subgroup of
${\rm Stab}_W(t)$ and the product decomposition is unique, by the first
paragraph of the proof, we see that 
${\rm Stab}_W(t) \cong W(t) \rtimes I$.
\end{proof}

Clearly, if $\lambda \in \Lambda/Q_H\spcheck$ has infinite order, then it
does not stabilize any point. Thus ${\rm
Stab}_{\Lambda/Q_H\spcheck}(\tilde t)$ is identified with the stabilizer
of $\tilde t$ in the torsion subgroup of $\Lambda/Q_H\spcheck$.

Finally, we note that all of the results of this subsection go over, with
essentially identical proofs, to the case where $\frak t=\frak d \oplus 
\frak z$ is a real vector space, $\Phi$ is a 
root system on $\frak d$, not necessarily reduced, and $\Lambda \subseteq
\frak t$ is a lattice such that
$$Q\spcheck(\Phi) \subseteq \Lambda \cap \frak d \subseteq
\pi_{\frak d}(\Lambda)\subseteq P\spcheck(\Phi),$$
where $\pi_{\frak d}$ is projection to the factor $\frak d$.

\subsection{A first characterization of groups of type $A_n$}

\begin{lemma}\label{root}
Let $\Phi$ be irreducible, but not necessarily reduced, with $d$ as
highest root and let
$e=\sum_{a\in\Delta} e_aa$ be the highest
short root. {\rm (}Of course, $d=e$ if and only if $\Phi$ is simply
laced.{\rm )}
\begin{enumerate}
\item $d$ and $e$ are orthogonal to all but either one or two simple
  roots.
\item The following are equivalent:
\begin{enumerate}
\item[(i)] $R$ is of type $A_n$ for some $n\ge 2$;
\item[(ii)] $d$ is not orthogonal to two simple roots;
\item[(iii)] $e$ is not orthogonal to two simple roots.
\end{enumerate}
\item If  $R$ is not of type $A_n$ for any $n\ge 1$, then
$e=\varpi_a$ for some $a\in \Delta$.
\end{enumerate}
\end{lemma}

\begin{proof} 
We have
$$2=n(d,d)=\sum_{a\in \Delta} h_an(a,d)=n(e,e)=\sum_{a\in
\Delta}e_an(a,e).$$ The coefficients $h_a$ and $e_a$, are
positive integers. Since $d$ is the  
highest root, resp.\ $e$ is the highest short root, the $n(a,d)$,
resp.\
$n(a,e)$ are also nonnegative.
Statement 1 follows.

Clearly, (i) implies (ii) and (iii). Moreover,
since the passage to an 
inverse system is conformal and permutes the roles of $d$ and $e$, (ii) is
equivalent to (iii). It remains to see that (ii) implies (i). Assume
(ii) holds. Let $a,b$ be the simple roots not orthogonal to $d$. There
exists a root $f$ which is a sum of distinct simple roots including $a,b$.
(Namely, $f$ is the sum of the simple roots corresponding to the nodes
of the interval connecting $a$ and $b$ in the Dynkin diagram for $\Phi$.) 
Then $n(f,d)=2$. But since $d$ is a root at least as long as all other
roots, the Cartan integer $n(f,d)=2$ implies that $f=d$, and
(i) follows. This proves Part (2).

By Statement 2, there is a unique simple root $a$ for
which $n(e,a)\not= 0$ meaning that $e$ is a multiple of $\varpi_a$.
Of course, $n(e,a)\ge 0$.
Since $e$ is short, $n(e,a)\in\{1,2\}$, with $n(e,a)=2$ if
and only if $a=e$. Since $R$ is not of type $A_1$, this can never
happen, proving Part (3).
\end{proof}

\begin{proposition}\label{first}
\begin{enumerate}
\item Let $G$ be simple.
Then there exists a fundamental  weight
$\varpi_a$
such that $\varpi_a({\cal C}G)=1$ if and only if $G$ is not
of type $A_n$ for any $n$. 
\item Let
$G=\prod_{i=1}^s G_i$ be the decomposition of
$G$ into simple
  factors, and let $c\in {\cal C}G$ be the product $c=c_1\cdots c_s$
  with $c_i\in {\cal C}G_i$. Write ${\rm log}(c)=\sum_{a\in \Delta}
  \lambda_aa\spcheck$.  If no $\lambda_a$ is integral, then for each
  $i,1\le i\le s$, we have
$G_i=SU(n_i)$ for some $n_i\ge 2$ and
  the element  $c_i$ generates ${\cal C}G_i$. 
\item Suppose that $G$ is of type $A_n$ and that $c$ generates ${\cal
C}G$. Write ${\rm log}(c)=\sum_{a\in \Delta}
  \lambda_aa\spcheck$. Then 
$$\{\lambda_a \bmod {\bf Z}, a\in \Delta\} = \left\{\frac{k}{n+1}, 1\leq k
\leq n\right\}.$$
In particular, no $\lambda _a$ is integral.
\end{enumerate}
\end{proposition}

\begin{proof}
(1) If $G$ is simple and not of type $A_n$, then by Lemma~\ref{root}
there is a
fundamental  weight $\varpi_a$ which is a root and hence
kills ${\cal C}G$.  Conversely, if $G=SU(n+1)$, the
fundamental representations are the exterior powers $\bigwedge ^i{\bf
C}^{n+1}$, $1\leq i\leq n$, of the standard representation, and their
highest weights are nontrivial on ${\cal C}SU(n+1)$.

If ${\rm log}(c)=\sum_a\lambda_aa\spcheck$, then $\varpi_a(c)={\rm
  exp}(2\pi\sqrt{-1}\lambda_a)$. Thus, 
(2) follows from (1) and the fact, again easily verified by direct
inspection,  that, for a proper subgroup $C$ 
of the center of $SU(n)$,  there is a fundamental  weight
which  kills $C$.

Finally, (3) follows by examining a generator for the center of $SU(n+1)$.
\end{proof}

The proof actually establishes the following:

\begin{addendum}\label{corootno} Let $\Phi$ be a not necessarily reduced
root system with irreducible factors $\Phi_i$. Let $\Delta$ be a set of
simple roots for
$\Phi$. Suppose that $\zeta = \sum _{a\in \Delta} \lambda_aa\spcheck \in
P\spcheck (\Phi)$, where $\lambda_a \in {\bf Q}$. Then no 
$\lambda_a$ is integral if and only if every $\Phi_i$ is of type 
$A_{n_i}$ for some integer $n_i$, and, for every $i$, the
projection of
$\zeta
\in P\spcheck (\Phi)/Q\spcheck(\Phi) \cong \bigoplus _iP\spcheck
(\Phi_i)/Q\spcheck(\Phi_i)$ to the factor  $P\spcheck
(\Phi_i)/Q\spcheck(\Phi_i)$ generates this factor.
\end{addendum}

\subsection{Subgroups associated with elements of the
  center}\label{center} 

 Let $c\in {\cal C}G$ and let
      $\lambda\in {\frak t}$ be such that ${\rm exp}(\lambda)=c$. Write
\begin{equation}
\label{lambdaeqn}
\lambda=\sum_{a\in \Delta}\lambda_aa\spcheck.
\end{equation}
Then the $\lambda_\alpha$ are rational numbers.
Let $\Delta(c)=\{a\in \Delta|\lambda_a\not\in {\bf Z}\}$. The set
$\Delta(c)$ depends only on $c$ and not on the choice of a lift $\lambda$.
In the notation of Lemma~\ref{subset}, we set ${\frak t}_c= {\frak
t}_{\Delta(c)}$,
$S_c\subset T$ equal to the subtorus with ${\rm Lie}(S_c)={\frak
  t}_c$, and $L_c=L_{\Delta(c)}=DZ(S_c).$

\begin{lemma}\label{c}
With the previous notation we have
\begin{enumerate}
\item $\varpi_a(c)=1$ if and
  only if $a\not\in \Delta(c)$. 
\item $c\in L_c$.
\item If $I\subseteq \Delta$ has the property that $c\in L_I$, then
  $\Delta(c)\subseteq I$, and hence $L_c\subseteq L_I$.
\end{enumerate}
\end{lemma}

\begin{proof}
Recall that
$\varpi_a(c)={\rm exp}(2\pi\sqrt{-1}\lambda_a)$ and that $a\spcheck\in
{\frak t}\cap {\rm Lie}\, L_I$ if and only if $a\in \Delta(c)$. From these
facts (1) and (2) are clear.   

If $c\in L_I$, then $c={\rm exp}(\lambda')$ for some element
$\lambda'$ in the real linear span of the coroots $a\spcheck$ for 
$a\in \Delta_I$.
Since the element $\lambda'$ differs by an element of 
the coroot lattice from $\lambda$, we see that $\lambda_a\in {\bf Z}$
for all $a\notin I$. That is to say $\Delta(c)\subseteq I$.
\end{proof}

\subsection{A further characterization of products of groups of
type $A_n$}

\begin{proposition}\label{finite}
Let $G=\prod_{i=1}^sG_i$, with the $G_i$ simple. 
Let $c\in {\cal C}G$. The following conditions are equivalent.
\begin{enumerate}
\item  The fixed point set $T^{w_c}$ of the $w_c$-action on $T$ is finite.
\item  For each $i,\, 1\le i\le s$, the group $G_i=SU(n_i)$ and the 
projection
  $c_i$ of $c$ into $G_i$ generates ${\cal C}G_i$.
\item  $\Delta(c)=\Delta$, i.e. in Equation~\ref{lambdaeqn}, no
coefficient $\lambda_a$ is integral.
\end{enumerate}
\end{proposition}

\begin{proof}
We have $c=c_1 \cdots c_s$ with $c_i\in {\cal C}G_i$, 
$w_c=w_{c_1}\cdots w_{c_s}$, and $A=A_1\times\cdots
\times A_s$ where $A_i$ is an alcove for $G_i$.
The condition that $T^{w_c}$ be finite implies that for each
$i$ the action of $c_i$ on the alcove $A_i$
permutes transitively the vertices of $A_i$.
Hence, every vertex
of the alcove $A_i$ is an element of the center 
of $G_i$.
This means that the highest root is
the sum of the simple roots.
The only groups with this property are groups of type $A_n$.
This shows that (1) implies (2).

It is clear that (2) implies (1) and (3). The fact that (3)
implies (2) follows from Proposition~\ref{first} and
Lemma~\ref{c}. 
\end{proof}

\begin{corollary}\label{LisLc} Let $c\in {\cal C}G$, and let $L$ be a
subgroup of $G$ of the form
$L_I$. Then $L=L_c$ if and only if 
\begin{enumerate}
\item $c\in L$;
\item  $L$ is a product of simple factors $L_i\cong SU(n_i)$ for some
$n_i$;
\item The projection of $c$ to $L_i$ generates the center of $L_i$.
\end{enumerate}
\end{corollary}
\begin{proof} 
First suppose that $L=L_c$. By Subsection~\ref{center}, $\Delta(c)$ is a
set of simple roots for $L_c$  and  $c\in L_c$. Of course, in the
expression 
$c=\sum_{a\in \Delta(c)}\lambda_aa\spcheck$
all the coefficients $\lambda_a$ are non-integral. Hence by
Proposition~\ref{finite} $L_c$ is a product of groups $\prod_{i=1}^sL_i$
where $L_i$ is isomorphic to $SU(n_i)$ for some integer $n_i\ge
2$ and $c=c_1\cdots c_s$ where $c_i$ generates the center of $L_i$.

Conversely, suppose that $L=L_I$ for some $I$, that $c\in L$ and that $L$
is a product of groups
$\prod_{i=1}^sL_i$ where $L_i$ is isomorphic to $SU(n_i)$ for some
integer $n_i\ge 2$ and $c=c_1\cdots c_s$ where $c_i$ generates the center
of $L_i$. By Lemma~\ref{c}, $\Delta(c) \subseteq I$. On the other hand, by
Part (3) of  Proposition~\ref{first},
no coefficient of $c$, expressed as a linear combination 
of the $a\spcheck, a\in I$, is integral, and hence $I\subseteq
\Delta(c)$. Thus $I=
\Delta(c)$ and so $L=L_c$.
\end{proof}

\begin{defn} Fix an element $c\in {\cal C}G$.
 Let $w_c$ be the Weyl part of the action of
$c$ on the alcove $A$. Let $T^{w_c}$ be the fixed points of the action
of $w_c$ on $T$. Let $S^{w_c}$ be the identity component of   $T^{w_c}$,
and let $\frak t^{w_c}$ be the fixed subspace for the action of $w_c$ on
$\frak t$. Clearly, ${\rm Lie}(S^{w_c}) = \frak t^{w_c}$.
\end{defn}

\begin{proposition}\label{torus}
The torus $S^{w_c}$ is conjugate to 
$S_c$.
\end{proposition}

\begin{proof}
Let $\frak t_{L_c}=\frak t\cap {\rm Lie}(L_c)$. It is the Lie algebra of a
maximal torus for
$L_c$. Let $A'$ be the alcove in ${\frak t}_{L_c}$ associated
with the set of simple roots $I_c\subseteq \Delta$  for $L_c$ with
respect to
$\frak t_{L_c}$. Let $\varphi_c'$ be the element of $W_{\rm aff}(L_c)$ whose
restriction to $A'$ induces the action of $c$ on $A'$.
By Proposition~\ref{finite} $\varphi_c'$ fixes a unique
point, say $\hat p$, of $A'$.
Let $\widetilde\varphi_c$ be the extension of $\varphi_c'$ to $\frak t =
{\frak t}_{L_c}\oplus{\frak t}_c$ by the identity on ${\frak t}_c$. Thus
$\widetilde\varphi_c$ is the image of
$\varphi_c'$ in
$W_{\rm aff}(G)$.

A root of
$G$ which is integral on the affine space
${\frak t}_c+\hat p$ must vanish on ${\frak t}_c$ and hence be a root of
$L_c$. But since $\hat p$ is a regular element for $L_c$, it follows that
there are no roots of $G$ taking integral values on ${\frak t}_c+\hat p$.
Thus there is an open dense subset of ${\frak t}_c+\hat p$
consisting of regular elements for $G$. In particular, there is  $v\in
{\frak t}$ with the following three properties:
\begin{enumerate}
\item ${\rm exp}\,v$  is a regular element for $G$;
\item $v$ is fixed by $\widetilde\varphi_c$;
\item  the unique alcove $A$ for the affine Weyl group of $G$ containing
$v$ also contains $\hat p$. 
\end{enumerate}

The point $\hat p$ lies in the alcove $A'\subset \frak t_{L_c}$, and $A'$
contains the origin. Condition 3 above implies that
$A$ contains the origin. It follows from  Conditions 1 and 2 that
$\widetilde\varphi_c$ sends
$A$ to itself. By Lemma~\ref{char} we see that $\widetilde\varphi_c$ is
the action of $c$ on
$A$. 
Thus, the fixed point set of the Weyl part $w'_c$ of $\varphi_c$ is
exactly  ${\frak t}_c$ and exponentiates onto $S_c$. The proposition now
follows since the Weyl part $w_c$ of the action of $c$ on $A$ is
conjugate to $w'_c$, and hence $S^{w_c}$ is conjugate
to $S_c$.
\end{proof}

\subsection{A consequence of Proposition~\protect{\ref{finite}}}

\begin{theorem}\label{Antype}
With notation as above, 
fix an integer $k>1$ dividing at least 
one of the $g_a$ for $a\in \widetilde\Delta$.
Let $\widetilde{I}(k)=\{\, a\in\widetilde \Delta: k \not|g_a\,\}$.
Let $H(k)$ be the closed, connected subgroup of  $G$ whose
complexified
Lie algebra is
generated by  $\left\{(\mathfrak g^a\oplus\mathfrak
  g^{-a})\right\}_{a\in \widetilde{I}(k)}$.     
Then $H(k)$ is isomorphic to
$$\left(\prod_{i=1}^r H_i\right)/\langle c\rangle$$
where:
\begin{enumerate}
\item $c$ has order $k$;
\item for each $i,\ 1\le i\le r$, the group $H_i$ is isomorphic to
  $SU(n_i)$ for some   $n_i|k$;
\item for at  least one $i$ between $1$ and $r$ we have $n_i=k$;  
\item $c=\prod_{i=1}^rc_i$ and $c_i\in H_i$ generates the
  center of $H_i$.
\end{enumerate}
\end{theorem}

\begin{proof}
The group $H(k)$ is a semi-simple subgroup of $G$ for which
$\widetilde{I}(k)$ is a set of simple roots. 
Let
${\frak t}_{H(k)}\subseteq \frak t$ be the subspace spanned by the
coroots $a\spcheck$ for $a\in \widetilde{I}(k)$. Then ${\frak t}_{H(k)}$ is 
the Lie algebra of a maximal torus of $H(k)$. The element 
\begin{equation}\label{lambda1}
\lambda=\frac{1}{k}\sum_{a\in \widetilde{I}(k)} g_aa\spcheck
\end{equation} is contained in ${\frak t}_{H(k)}$ and
is also in $Q\spcheck$ since it is equal to
$$-\frac{1}{k}\sum_{a\in \widetilde \Delta-\widetilde{I}(k)}g_aa\spcheck$$
and by definition $k|g_a$ for every $a\in \widetilde
\Delta-\widetilde{I}(k)$.  This means that every root of $H(k)$ takes
integral values on
$\lambda$, and hence that $\lambda$ exponentiates to an element
$c$ contained in the center of  
the simply connected form  $\widetilde H(k)$ of $H(k)$. By
Equation~\ref{lambda1}, 
$c^k=1$.

The definition of $\lambda$ implies that, when $\lambda$ is
expressed as a linear combination of the basis $a\spcheck, a\in
\widetilde{I}(k)$,  all the
coefficients of the simple coroots for $H(k)$  are non-integral.
Hence, by Proposition ~\ref{first},
$\widetilde H(k)$ is a product  $\prod_{i=1}^rH_i$, where for each $i$ 
the group  $H_i$ is isomorphic to 
$SU(n_i)$ for some $n_i\ge 2$,  and $c$ is of the form
$c_1\cdots c_r$ where for each $i$ the element $c_i$ generates the
center of $H_i$. 
Since $c^k=1$, each of the $c_i$ has order dividing $k$.
Since $c_i$ generates the center of $H_i$, its order is $n_i$. We
conclude that $n_i|k$ for each $i$.

Consider now the component of $\widetilde D(k)$ that contains 
$\tilde a=-\tilde d$.
 (Since $g_{\tilde a}=1$, we have $\tilde a\in \widetilde D(k)$.)
We index the $H_i$  so that this component corresponds to $H_1$.
Since the expression for $\lambda$ as a linear combination of
coroots has $1/k$ as the coefficient of $\tilde a$, it
follows that the order of $c_1$ in $H_1$ is divisible by $k$,
i.e., that $k|n_1$. Since we have already shown the opposite
divisibility, it  must  be the case that $n_1=k$, showing that $H_1$
is isomorphic to $SU(k)$. Moreover, the order of $c$ is divisible by $k$
and hence is equal to $k$, and $k$ is the least common multiple of the
$n_i$. 

The fundamental group of $H(k)$ is cyclic, by
Lemma~\ref{coroottors}, and contains the element
$c$, which is of order $k$. Thus $k$ divides the order of
$\pi_1(H(k))$.  On the other hand, $\pi_1(H(k))$ is identified
with  a cyclic subgroup of
$\prod_i{\cal C}SU(n_i) =
\prod_i{\bf Z}/n_i{\bf Z}$, and hence its order divides the least
common multiple of the $n_i$, namely
$k$. Since $c$ has order $k$, it generates the fundamental group of
$H(k)$. 
\end{proof}

\begin{corollary}\label{values}
If $N$ is the maximal value for $g_a$ for
$a\in \Delta$, then for each $k,\ \  1\le k\le N$, there is at least one
$a\in \widetilde\Delta$ for which $g_a=k$. In fact there is a simply laced
chain of length $N$ in $\widetilde D(G)$ containing $\tilde a$ as one end 
so that the $g_a$, in order, along this chain are  $1,2,\ldots, N$.
\end{corollary}

\begin{corollary}\label{gcd}
If $\Phi$ is an irreducible root system and $k$ is a positive integer
dividing at least one of the coroot integers $g_a$, then
$\gcd\{\, g_a: k|g_a\,\} = k$.
\end{corollary}

\subsection{Application to generalized Cartan matrices and affine
diagrams}\label{affine}

In this subsection, we apply the above results on groups of type 
$A_n$ to establish numerology concerning coroot integers, root
integers, and more generally integers which are the coefficients
of linear relations for nodes of diagrams of affine type. This
numerology will be crucial for the proof of the Clockwise
Symmetry result.

Let $\Psi$ be an irreducible, but not necessarily reduced, root system
and suppose that $\Upsilon$ is a set of simple roots for $\Psi$. As in the
case of reduced root systems, there is the extended set of simple roots
$\widetilde
\Upsilon$, obtained by adding minus the highest root to $\Upsilon$.
There is an extended Dynkin diagram $\widetilde D(\Psi)$ whose nodes are the
extended set of simple roots $\widetilde \Upsilon$, and whose
bonds, with their multiplicities   and arrows are
determined by the Cartan integers $n(a,b)$ for $a,b\in \widetilde \Upsilon$
exactly as in the case of the ordinary diagram.
Dually, there is the extended diagram $\widetilde D\spcheck(\Psi)$
associated with the coroots inverse to the extended roots; it is obtained
from $\widetilde D(\Psi)$ by reversing the direction of every arrow.

More generally, suppose that we are given 
a finite set $\widetilde \Upsilon$ and an integral  matrix $N=(n(a,b))$,
where
$a,b\in
\widetilde \Upsilon$.  The matrix $N$ is called a {\sl generalized Cartan
matrix\/} if:
\begin{enumerate}
\item $n(a,a)=2$ for all $a\in S$.
\item $n(a,b)\le 0$ for all $a\not= b\in S$.
\item $n(a,b)=0$ implies $n(b,a)=0$.
\end{enumerate} 
If in addition the set $\widetilde \Upsilon$ cannot be divided into two
disjoint non-empty 
subsets $S_1$ and $S_2$ such that $n(a,b)=0$ for all $(a,b)\in
S_1\times S_2$, then we call $N$ {\sl indecomposable}. We will
assume throughout that $N$ is indecomposable.

We can form a diagram associated with a generalized
Cartan matrix whose nodes are indexed by $\widetilde \Upsilon$. It has
bonds with multiplicities and arrows determined by the same rules as in
the case of Dynkin diagrams. Because of indecomposability, the diagram
associated with a generalized Cartan matrix is connected.
Notice that one can reconstruct the generalized Cartan matrix from its
diagram. 

Given a
real vector space $\widetilde V$ of dimension $d+1$ and a positive
semidefinite bilinear form
$\langle
\cdot, \cdot
\rangle$, suppose that $\widetilde \Upsilon$ is a basis of
$\widetilde V$ such that $\langle v, v\rangle \neq 0$ for all $v \in
\widetilde
\Upsilon$. Then we can define $n(v,w)$ by the usual formula
\begin{equation}
\label{Cartan2} n(v,w)=\frac{2\langle v,w\rangle }{\langle w,w\rangle}. 
\end{equation}
If
$N=(n(v,w))$ is a generalized Cartan matrix and  there is a vector
$u=\sum_{v\in \widetilde \Upsilon}n_vv$ in $\widetilde V$ such that 
$n_v>0$ for all $v$ and such that $\langle u, x\rangle =0$ for all
$x\in V$, then the generalized Cartan matrix is
said to be {\sl of 
affine type}. In this case the associated diagram is called {\sl an affine
diagram}, and, according to a theorem of Kac \cite{Kac2}, the
diagram is either the extended root or the extended coroot diagram of
a possibly non-reduced root system, and the coefficients of the vector
$u$ are a fixed positive integral multiple of the (root or coroot) integers on
the extended diagram. Moreover, every proper subdiagram is the
Dynkin diagram of a root system.  Let $\frak t(\widetilde \Upsilon)$ be
the quotient of
$\widetilde V$ by the one-dimensional radical of the semidefinite form and let
$Q(\widetilde \Upsilon)$ be the lattice in $\frak t(\widetilde \Upsilon)$ 
spanned by the image of $\widetilde \Upsilon$. Note that $\frak 
t(\widetilde
\Upsilon)$ is a $d$-dimensional vector space with a positive definite inner
product. Projection induces a bijection from $\widetilde \Upsilon$ to a
spanning set of
$\frak t(\widetilde \Upsilon)$ of cardinality $d$, and the Cartan 
integers  are determined by  inner products of their images in $\frak
t(\widetilde
\Upsilon)$ by Equation~\ref{Cartan2}. It follows from 
the theorem of Kac that the reflections in the $v\in \widetilde
\Upsilon$ generate a Weyl group acting by isometries on $\frak t(\widetilde
\Upsilon)$ and that the lattice
$Q(\widetilde
\Upsilon)$ spanned by $\widetilde\Upsilon$ is invariant
under this group. It also follows that there is one linear relation
between the vectors of $\widetilde \Upsilon$ and that the coefficients of
this relation can be chosen to be  positive integers.

This construction can be reversed: suppose that $V$ is a vector space of
dimension
$d$ with a positive definite inner product and let $\widetilde
\Upsilon$ be a subset of cardinality $d+1$ spanning $V$ such that the
Cartan numbers defined by Equation~\ref{Cartan} are integers and determine
an indecomposable generalized Cartan matrix. Then this matrix  is  of
affine type.

For example, if $G$ is   simple  of rank $r$, then
$\widetilde \Delta\spcheck$ is a subset of ${\frak t}$ of
cardinality 
$r+1$. This embedding induces an identification of ${\frak t}(\widetilde
\Delta\spcheck)$ with ${\frak t}$, and further identifies the
coroot  lattice
$Q\spcheck$ with the lattice $Q(\widetilde \Delta\spcheck)$
and the Weyl group of $G$ with the group generated by the
reflections in $\widetilde
\Delta\spcheck$. 

In this general context of affine diagrams, we have the following
generalization of Corollary~\ref{values}. 

\begin{proposition}\label{numerology} Let $D(\widetilde \Upsilon)$ 
be a connected affine diagram  whose nodes are indexed by $v\in
\widetilde
\Upsilon$. Let ${\bf n}\colon \widetilde \Upsilon \to {\bf N}$ be
a function such that
$\sum {\bf n}(v)v =0$ in $\frak t(\widetilde \Upsilon)$.We denote
${\bf n}(v)$ by $n_v$.  Let
$k$ be a positive integer which divides at least one of the $n_v$.
Let 
$\widetilde I({\bf n},k) = \{\, v: k \not|n_v\,\} \subseteq
\widetilde
\Upsilon$. Then there exist cyclic subgroups $C_i\subseteq {\bf Z}/k{\bf
Z}$, not necessarily distinct, and a bijection
$$\phi_k \colon \widetilde I({\bf n},k) \to \coprod_i(C_i
-\{0\})$$ such that $n_v \equiv \phi_k(v) \bmod k$ for all $v\in
\widetilde I({\bf n},k)$.
\end{proposition}
\begin{proof} Let 
$$\zeta = -\frac{1}{k}\sum _{v\in \widetilde I({\bf n},k)}n_vv.$$
Then $\zeta$ lies in the ${\bf R}$-span $\frak t(\widetilde
I({\bf n},k))$ of
$\widetilde I({\bf n},k)$ in $\frak t(\widetilde \Upsilon)$ as well as in
$Q(\widetilde \Upsilon)$. Since
$\widetilde I({\bf n},k)$ is a proper subset of $\widetilde\Upsilon$, it is a set
of simple roots for a root system on $\frak t(\widetilde I({\bf n},k))^*$ with
the property that the given inner product on $\frak t(\widetilde
\Upsilon)$ restricts to a Weyl invariant inner product on $\frak
t(\widetilde I({\bf n},k))$. Let $\Psi$ be the inverse root system on $\frak
t(\widetilde I({\bf n},k))$, so that $\widetilde I({\bf n},k)$ is
a set of simple coroots for $\Psi$. Since $\zeta \in Q(\widetilde
\Upsilon)$, $\zeta$ has integral inner product with the lattice
spanned by $\widetilde I({\bf n},k)$ and hence $\zeta \in
P\spcheck(\Psi)$. Moreover, all of the coefficients of
$\zeta$ are non-integral with respect to the set of simple coroots given
by $\widetilde I({\bf n},k)$. Thus, by Addendum~\ref{corootno}, $\Psi$ is a
product of irreducible root systems $\Psi_i$ of type $A_{N_i}$ for some
integers
$N_i$ and
$\zeta$ projects into each factor $P\spcheck(\Psi_i)/Q\spcheck(\Psi_i)$
as a generator. Since $k\zeta \in  Q\spcheck(\Psi)$, it
follows that $(N_i+1)|k$. Let
$C_i$ be the cyclic subgroup of ${\bf Z}/k{\bf
Z}$ of order $N_i +1$. The result now follows immediately from Part 3 of
Proposition~\ref{first}.
\end{proof}

\subsection{Numerical consequences of
Proposition~\protect{\ref{numerology}}}

Let $\widetilde \Upsilon$ be a finite set and let ${\bf n}\colon
\widetilde
\Upsilon \to {\bf N}$ be a function. We denote ${\bf n}(v)$ by
$n_v$. For every  positive integer
$k$ which divides at least one of the $n_v$, let 
$\widetilde I({\bf n},k) = \{\, v: k \not|n_v\,\} \subseteq \widetilde
\Upsilon$.  We suppose
throughout that the pair $(\widetilde \Upsilon,{\bf n})$
satisfies the conclusions of Proposition~\ref{numerology}:

\begin{assumption}\label{assump} For every  positive integer
$k$ which divides at least one of the $n_v$,  there exist cyclic
subgroups $C_i\subseteq {\bf Z}/k{\bf Z}$, not necessarily distinct, and
a bijection
$$\phi_k \colon \widetilde I({\bf n},k) \to \coprod_i(C_i -\{0\})$$
such that $n_v \equiv \phi_k(v) \bmod k$. 
\end{assumption}

Let $n_0 = \gcd \{n_v: v \in \widetilde
\Upsilon\}$. We call the pair $(\widetilde
\Upsilon,{\bf n})$ {\sl reduced\/} if $n_0=1$. In general, 
if we define $n'_v = n_v/n_0$, the pair $(\widetilde \Upsilon,n')$ also
satisfies the conclusions of Assumption~\ref{assump} and is reduced.

Let us introduce the following notation:
\begin{eqnarray*}
i(x) &=& \#\{v\in \widetilde \Upsilon: n_v =x\}; \\
i(x,k) &=& \#\{v\in \widetilde \Upsilon: n_v \equiv x\bmod k\}= \sum _{\ell
\in {\bf Z}}i(x+\ell k);\\ 
N &=& \max\{ n_v: v\in \widetilde \Upsilon\}:\\
g &=& \sum _{v\in \widetilde \Upsilon}n_v = \sum _{x\geq 1}xi(x).
\end{eqnarray*}

Note for example that $i(x,N) = i(x)$ for all $x$ such that $1\leq x\leq N$.

The following is a consequence of
Assumption~\ref{assump}:

\begin{lemma}\label{compare}
Suppose $r,s\in {\bf Z}$ are not divisible by $k$. If $\langle s
\rangle \subseteq \langle r \rangle$ as subgroups of  ${\bf
Z}/k{\bf Z}$,  $i(r,k) \leq i(s,k)$. 
Hence, if $\langle s \rangle = \langle r
\rangle$, then $i(r,k) = i(s,k)$.
\end{lemma}

\begin{lemma}\label{string}
The $x$ such that $i(x)\neq 0$, in other words the
integers of the form $n_v$, are exactly the positive multiples of $n_0$
less than or equal to $N$.
\end{lemma}
\begin{proof} It suffices to consider the case where  $(\widetilde
\Upsilon,{\bf n})$ is reduced and to show that the $x$ such that
$i(x)\neq 0$ are exactly the integers $x$ such that $1\leq x \leq
N$. Let $\ell$ be the smallest positive integer such that
$i(\ell) \neq 0$. Thus $i(t)=0$ for
$t<\ell$, and hence, by Lemma~\ref{compare}, $i(t) =0$ for $N-\ell<t<N$. If
$\ell \neq 1$, since 
$(\widetilde
\Upsilon,{\bf n})$ is reduced, there exists an $x$ with $i(x)\neq
0$ and with $x$ not divisible by $\ell$. Choose $x$ to be the
smallest such positive integer. In particular, $x>\ell$, and, by
Lemma~\ref{compare},
$i(N-x) = i(x)\neq 0$. Since
$$-(x-\ell) = \ell -x \equiv N-x \bmod {(N-\ell)},$$
it follows that $i(x-\ell, N-\ell) \neq 0$. For $d\geq 1$,
$$d(N-\ell) +x-\ell \geq N-\ell  +x-\ell > N-\ell,$$
and so if $t\equiv x-\ell\bmod {(N-\ell)}$ and $i(t)\neq 0$, then $t=x-\ell$.
This says that $i(x-\ell) = i(x-\ell, N-\ell) \neq 0$. But
$x-\ell < x$ and
$\ell{\not|}x-\ell$, contradicting the choice of $x$. Hence $\ell =1$.
Applying Lemma~\ref{compare} with $k=N$ and $r=1$, we see that 
$i(s) \geq i(1)\geq 1$ for all $s$ with $1\leq s\leq N$.
\end{proof}

\begin{lemma}\label{atmost6}
Suppose that $(\widetilde
\Upsilon,{\bf n})$ is reduced. Then $\varphi(N) \leq 2$, where
$\varphi$ is the Euler $\varphi$-function, and hence $N\in
\{1,2,3,4,6\}$.
\end{lemma}
\begin{proof} Suppose that $1\leq x \leq N-1$ and that $x$ is relatively
prime to $N$. Then by Lemma~\ref{compare} and Lemma~\ref{string}, $i(x) = i(1)
\geq 1$. Assume that $x\neq 1, N-1$. Then
$$i(1)+1\leq i(1)+i(N) \leq i(1, N-1) \leq i(x, N-1) = i(x) = i(1),$$
a contradiction. Thus $x=1$ or $x=N-1$, and so $\varphi(N) \leq 2$.
\end{proof}

\begin{lemma}\label{atleast2}
Let $\ell>1$ be a positive integer such that $i(\ell) \neq 0$. Then either
$i(t\ell) = 0$ for $t>1$ or $i(\ell) \geq 2$.
\end{lemma}
\begin{proof} We may assume that $(\widetilde
\Upsilon,{\bf n})$ is reduced and that $2\ell \leq N$. Since
$\ell +2 \leq 2\ell
\leq N$, $i(1, \ell +1) \geq 2$. First suppose that $\ell \geq 3$.  By
Lemma~\ref{compare}, $i(1, \ell+1) =i(\ell, \ell +1)$. Since $2\ell +1 >
6\geq N$, by Lemma~\ref{atmost6}, it follows that $i(\ell) = i(\ell, \ell
+1) = i(1, \ell+1)$, which as we have just seen is at least $2$. Now
suppose that
$\ell =2$. Then $N=4$ or $6$. If $N=4$, then as before $2\ell +1 > N$ and
so $i(2)= i(2, 3)
\geq 2$. If $N=6$, then $i(1,5) \geq 2$. But $i(1,5) = i(2,5) =i(2)$ and
so again
$i(2)\geq 2$.
\end{proof}

In case $(\widetilde
\Upsilon, {\bf n})$ is   reduced, it is easy to check that
necessary and sufficient conditions on the integers $i(x)$ for
the pair 
$(\widetilde
\Upsilon, {\bf n})$ to satisfy Assumption~\ref{assump} are as
follows: 
\begin{itemize}
\item For $N=1$ or $2$ there is no condition on the $i(x)$. In this case
$g=i(1)$, resp.\ $g=i(1)+2i(2)$. 
\item For $N=3$, a necessary and sufficient condition is $i(1)= i(2)$. In
this case
$g= 3(i(1)+i(3))$.
\item For $N=4$, necessary and sufficient conditions are:  $i(1) = i(3)$;
$i(1)+i(4) = i(2)$. In this case $g=6(i(1)+i(4))$.
\item For $N=6$, necessary and sufficient conditions are:  $i(1)
=i(5)=i(6)$, $i(2)=i(3)=i(4)=2i(1)$. In this case $g=30i(1)$.
\end{itemize}

Next we define
$$d_x  =  d_x(\widetilde \Upsilon, {\bf n})= \#\{v\in \widetilde
\Upsilon: x|n_v\} = 
\sum_{\ell \geq 1}i(\ell x).$$
Note that $d_x\neq 0$ if and only if $x|n_v$ for some $v$.

\begin{lemma}\label{standardfact}
$$\sum _{x\leq N}\varphi(x)d_x = g.$$
\end{lemma} 
\begin{proof} Using the identity $\sum
_{d|n}\varphi(d) = n$, we have
$$\sum _{x\leq N}\varphi(x)d_x = \sum _{v\in
\widetilde\Upsilon}\sum _{d|n_v}\varphi(d) = \sum _{v\in
\widetilde\Upsilon}n_v = g.$$
\end{proof}

We come now to one form of the statement of clockwise symmetry:

\begin{theorem}\label{clocked} Suppose that the pair $(\widetilde
\Upsilon, {\bf n})$
 satisfies Assumption~\ref{assump}. For each $x\leq N$ such that
$d_x\neq 0$,  $x|2g$. For
each such $x\leq N$, and for each $r \leq x$ and   relatively prime to
$x$,   consider the  subset of the integers mod
$2g$ given by 
$$J(x,r)= \left\{\frac{2gr}{x}-d_x+1, \frac{2gr}{x}-d_x+3, \dots,
\frac{2gr}{x}+d_x -3, \frac{2gr}{x}+d_x-1\right\}.$$
Thus $J(x,r)$ consists of $d_x$ integers, centered at $2gr/ x$ and with
spacing $2$. Then for distinct pairs $(x,r)\neq (y,s)$, the sets
$J(x,r)$ and $J(y,s)$ are disjoint, and
$\bigcup_{x,r} J(x,r)\subset {\bf Z}/2g{\bf Z}$  is either $\{0,2,
\dots, 2g-2\}$ or $\{1, 3, \dots, 2g-1\}$.
\end{theorem}
\begin{proof}  Let us first show that it suffices to consider
the  the case where $(\widetilde
\Upsilon, {\bf n})$ is reduced. For a general pair $(\widetilde
\Upsilon, {\bf n})$, let $(\widetilde
\Upsilon, {\bf n}')$ be the associated reduced pair. Thus $g =
n_0g'$. For each
$x$, write $x=\ell m$, where $\ell =\gcd(x,n_0)$. Then $d_x= d'_m$, where
$d'_m= d_m(\widetilde
\Upsilon, n')$. Moreover, it is easy to see that $d_x= 0$ for all other
$x$. An elementary argument shows that the set of rational numbers of the
form
$2g'r/m + a$, with $1\leq r\leq m$ and $r$ relatively prime to $m$, and
$0\leq a
< n_0$, is exactly the set of rational numbers of the form $2n_0g's/x$,
with $1\leq s\leq x$ and $s$ relatively prime to $x$ and with $x=\ell m$,
where $\ell =\gcd (x, n_0)$. Thus, $\bigcup_{x,r} J(x,r)\subset {\bf
Z}/2g{\bf Z}$ for $(\widetilde
\Upsilon, {\bf n})$ is invariant under translation by $2g'$ and
the image of
$\bigcup_{x,r} J(x,r)$ in ${\bf Z}/2g'{\bf Z}$ is the corresponding
subset for the pair $(\widetilde
\Upsilon, {\bf n}')$. Hence it suffices to consider the reduced
case.

Let ${\cal F}_N =\{0/1, 1/N, 1/(N-1), \dots\}$
be the Farey sequence of rational numbers between $0$ and $1$ whose
denominator is at most $N$, written in increasing order. We call integers
$x$ and $y$ {\sl adjacent\/} with respect to $N$ if there exist
$r,s$ with $(r,x) = (s,y) =1$ such that $r/x$ and $s/y$ are
consecutive terms in ${\cal F}_N$. If
$r/x$ and $s/y$ are consecutive terms in ${\cal F}_N$, then it is a
standard fact that $sx-ry=1$. The conclusions of Theorem~\ref{clocked}
are easily seen to be equivalent to the following statement:

For all  consecutive terms $r/x$ and $s/y$ in ${\cal
F}_N$,
$$\frac{2gs}{y} = \frac{2gr}{x} + d_x+d_y.$$
Using the fact that $sx-ry=1$, this condition is equivalent to:

For all integers $x$ and $y$ which are adjacent  with respect to $N$,
$$g =\frac{xy}{2}(d_x+d_y).$$

Another way to write the conclusions of Theorem~\ref{clocked} is as
follows: for all $x\leq N$ and $r\leq x$ with $(r,x)=1$, 
$$\frac{2gr}{x} = d_1 + \sum_{\scriptstyle(y,t)=1\atop\scriptstyle 2gt/y
< 2gr/x}2d_y + d_x.$$
By the symmetry $r\mapsto -r$ it is sufficient to check these conditions
for $r/x\leq 1/2$. The case $r/x =1/2$ follows from $\sum _{x\leq
N}\varphi(x)d_x = g$. Thus, for $ N \leq 4$, it is  enough to
check the first two conditions:
\begin{eqnarray*}
\frac{2g}{N}&=& d_1+d_N;\\
\frac{2g}{N-1}&=& d_1+2d_N+d_{N-1}.
\end{eqnarray*}
Let us consider  the first condition. By Assumption~\ref{assump},
$$d_1+ d_N = \sum _{t=1}^{N-1}i(t) + 2i(N) =
\sum_{\scriptstyle t|N\atop\scriptstyle t<N}\varphi(N/t)i(t) + 2i(N).$$
Since $g = \sum _{t=1}^Nti(t) = \sum _{t=1}^{N-1}ti(t) + Ni(N)$,
we see that it suffices to show that
$$\sum _{t=1}^{N-1}ti(t) = \sum_{\scriptstyle t|N\atop\scriptstyle
t<N}\frac{N}{2}\varphi(N/t)i(t).$$
On the other hand, by Assumption~\ref{assump},
$$\sum _{t=1}^{N-1}ti(t) = \sum_{\scriptstyle t|N\atop\scriptstyle
t<N}\left(\sum_{\scriptstyle s<N\atop\scriptstyle
\langle s \rangle = \langle t \rangle}s\right)i(t).$$
The condition then follows from the elementary lemma:

\begin{lemma}   For every positive integer $N$, and every positive divisor
$t$ of $N$,
$$\sum_{\scriptstyle s<N\atop\scriptstyle
\langle s \rangle = \langle t \rangle}s = \frac{N}{2}\varphi(N/t).$$
\end{lemma}
\begin{proof} Fix $t|N$. Then
that 
$$\sum_{\scriptstyle s<N\atop \scriptstyle\langle s \rangle = \langle t
\rangle}s = t\sum _{\scriptstyle u<N/t \atop \scriptstyle(u, N/t)=1}u=
t\frac{N}{2t}\varphi(N/t)=\frac{N}{2}\varphi(N/t),$$ where the second
equality follows from
$$M\varphi(M) = \sum_{\scriptstyle u<M \atop \scriptstyle(u, M)=1}M =
\sum_{\scriptstyle u<M \atop \scriptstyle(u, M)=1}(u+(M-u))
=2\sum_{\scriptstyle u<M \atop \scriptstyle(u, M)=1}u.$$ 
\end{proof}

This proves that the first condition holds under under
Assumption~\ref{assump}.  A very similar argument handles the second
condition. The result follows for $N\leq 4$. The case
$N=6$ can be checked directly.
\end{proof}

\begin{remark} One can ask if, given a positive integer $N$ there are
collections of (not necessarily positive) integers $i(x)$, not
satisfying Assumption~\ref{assump}, but such that the corresponding
integers
$d_x$  satisfy the conclusions of
Theorem~\ref{clocked}. It is easy to see by the proof of
Theorem~\ref{clocked} that, for $N\in
\{1,2,3,4,6\}$, Theorem~\ref{clocked} is equivalent to
Assumption~\ref{assump}. For $N=5$, fixing a positive integer $d$ and
setting $i(1)= 2d, i(2)=3d, i(3)=3d, i(4)=2d, i(5)=d$, the
corresponding integers $d_x$ satisfy Theorem~\ref{clocked} but of
course the 
$i(x)$ cannot satisfy Assumption~\ref{assump}, by
Lemma~\ref{atmost6}. These are in fact the only nonzero examples.
\end{remark}

\section{$c$-pairs}

Let
$C$ be an antisymmetric $2\times 2$-matrix with entries in ${\cal C}G$.
Then $C$ is completely specified by  $c=c_{12}=c_{21}^{-1}$. We consider 
ordered pairs of elements $(x,y)$ 
in $G$ satisfying $[x,y]=c$ and
call such pairs {\sl $c$-pairs}.

\subsection{The rank zero case}

Following our discussion of the structure of the moduli space of
almost commuting $N$-tuples in $G$, our first task is to
determine the set of rank zero $c$-pairs in $G$.

\begin{proposition}\label{rank0pair}
Suppose that  $c\in {\cal C}G$ is an
element of order $k>1$. Let
$(x,y)$ be a $c$-pair of rank zero in $G$. Then:
\begin{enumerate}
\item Both $x$ and $y$ are regular elements of $G$ which are conjugate;
\item The group $G$ is a product of $r$ simple factors $G_i$,
where each $G_i$ is isomorphic to $SU(n_i)$ for some $i\ge 2$;
\item $c=c_1\cdots c_r$  where $c_i$ generates the center of $G_i$;
\item the subgroup of $G/\langle c\rangle $ generated by $x,y$ is
isomorphic to $({\bf Z}/k{\bf Z})^2$ where $k$ is the order of $c$;
\item All $c$-pairs in $G$ are conjugate;
\item  $Z(x,y)={\cal C}G$.
\end{enumerate}

Conversely, if $G$ is as in {\rm (2)} and $c\in {\cal C}G$ is as in {\rm
(3)}, then there is a rank zero $c$-pair in $G$.
\end{proposition}

\begin{proof}
Let $(x,y)$ be a $c$-pair in $G$.
Conjugation by the element $y$ normalizes the connected group $Z(x)$.
Thus, by \cite{deS}  II \S 2, since $Z(x,y)$ is finite,  $Z(x)$ is
a torus. Hence $x$ is a regular element of $G$. By symmetry, $y$
is also regular.

Recall that $A\subseteq {\frak t}$ is the alcove containing the origin
associated to the set of simple roots $\Delta$ for $\Phi$. 
By conjugation we can assume that $x\in
T$, so that $T=Z(x)$. Conjugation by $y$
normalizes $Z(x)=T$, and hence
$y\in N_G(T)$. Finally, conjugation by an element of $N_G(T)$ makes
$x$ the image under the exponential mapping of  a point $\tilde x\in
A$.  Since $x$ is regular, $\tilde x$ is an interior point of $A$.
Let $w\in   W(T,G)$ be the  Weyl element 
defined by conjugation by $y$.
The relation
$yxy^{-1}=xc^{-1}$ yields
$w\cdot\tilde x-\xi=\tilde x$ for some $\xi\in {\frak t}$  such that
${\rm exp}\,\xi=c^{-1}$. We denote by $\varphi$ the affine linear
map
$v\mapsto w\cdot v-\xi$. The map $\varphi$ normalizes the alcove
structure for
$W_{\rm aff}(\Phi)$ and $\varphi(\tilde x) =\tilde x$, where $\tilde x$
is an interior point of $A$. Thus,
$\varphi(A)=A$.  It follows that $\xi\in A$ and, by
Lemma~\ref{char},  that
$\varphi$ is the action of $c$ on $A$. This means that $w=w_c$.
Since $T^{w_c}\subseteq Z(x,y)$, the fact that $(x,y)$ is of rank 
zero implies that ${\frak t}^{w_c}$ is a single point.
Thus, Proposition~\ref{finite} implies that $G$ is a product of groups
$\prod_{i=1}^r G_i$, where, for each 
$i$, the group  $G_i$ is isomorphic to $SU(n_i)$ for some $n_i>1$, and
$c$ projects to a generator of ${\cal C}G_i$.
The alcove $A$ is a product of alcoves $A_i$ for the simple factors
$G_i$ of $G$. The unique fixed point of the $c$-action on $A$ is
the product of the barycenters of the $A_i$.
Thus, $x$ is the image under the exponential mapping of the product of
the barycenters of the $A_i$.

Let $\langle x,y\rangle $ be the subgroup of $G$ generated by $x,y$.
In fact, it is a subgroup of $N_G(T)$ and $\langle x,y\rangle/(\langle
x,y\rangle\cap T)$ is the cyclic group generated by $[y]=w_c\in W$. This
element is of order $k$, the order of $c\in {\cal C}G$.  Also, $x\in T$
is a product of barycenters in the $A_i$, and so by inspection $x$ has
order
$k$ modulo
$\langle c\rangle$. From this it is clear that $\langle x,y\rangle\subseteq
G/\langle c\rangle $ is isomorphic to $({\bf Z}/k{\bf Z})^2$.

Lastly, reversing the roles of $x$ and $y$ and replacing $c$ by
$c^{-1}$ we see that $y$ is also conjugate to the product of the
barycenters of the $A_i$. Hence, $x$ and $y$ are conjugate in $G$.

Let $(x',y')$ be any $c$-pair in $G$ and let $S$ be a maximal torus of
$Z(x',y')$. By conjugation we can assume that $S=S_I$ for some subset
$I\subseteq \Delta$. Then $c\in L_I$, and hence, by Lemma~\ref{c},
$\Delta(c)\subseteq I$. Since $G$ is a product of simple factors
of type $A_n$ and $c$ projects to a generator of every factor, it
follows from   Proposition~\ref{finite}
that $\Delta(c)=\Delta$. Thus, $I=\Delta$ and consequently, $S$ is
trivial, which means that $(x',y')$ is of rank zero.
By what we proved above, we see that $(x',y')$ is conjugate to a pair
$(x,y'')$ where $x$  is the image under the exponential mapping of
the product of the barycenters of the $A_i$ and $y''\in w_cT$.
But $T$
operating by inner automorphism on $w_cT$ is transitive since
the component of the identity of $T^{w_c}$ is trivial. [Proof:
Consider the map $\mu$ defined by $t\mapsto 
tw_ct^{-1}$. The isotropy group of $w_c$ is $\{t|w_c^{-1}tw_c=t\}$,
which is finite.]
Thus, all such pairs are conjugate by elements of $T$.

Suppose that $z\in Z(x,y)$. Since $x$ is a regular element of $T$, $z\in
T$ and $z=yzy^{-1}=w_cz$. Since $G$ is a product of groups of type $A_n$
and the image of $c$ in each factor generates the center of that factor,
it follows by inspection that $z\in {\cal C}G$.

Conversely, suppose that $G$ is a product of simple groups isomorphic to
$SU(n_i)$ and $c$ projects to a generator of the center of each factor.
Set
$x$ equal to the image under the exponential mapping of the product of
the barycenters of the alcoves  for the simple factors of $G$ and
take $y\in w_cT$. Then $(x,y)$ is a $c$-pair of rank zero.
\end{proof}

\begin{corollary}\label{orders}
If $(x,y)$ is a rank zero $c$-pair in $SU(n)$, then 
$(x,y)$ is conjugate to a pair of the form $(x_0,y_0)$ where $x_0$ is
the image under the exponential mapping of the barycenter of $A$ and
$y_0$ normalizes $T$ and projects to $w_c$ in $W(T,G)$.
In particular, both $x$ and $y$ are regular elements and are conjugate
in $SU(n)$.
If $n$ is odd, each of $x$ and $y$ has  
order $n$, which is the order of $c$. If $n$ is even, then each of
$x$ and $y$ has order $2n$, whereas the order of $c$ is $n$. In fact, in
this case
$x^n = y^n = c^{n/2}$.
\end{corollary}
\begin{proof} 
All of these statements were established in the course of the proof of 
the previous proposition, except the statements about the orders of $x$
and
$y$. This statement follows by inspection of the order of the image under
the exponential mapping of the barycenter in a group of
type $A_n$.
\end{proof}

\subsection{The general case}

The next step is to determine the maximal torus for the centralizer of
a $c$-pair $(x,y)$.

\begin{proposition}\label{S_c}
Let $G$ be simple, let $c\in {\cal C}G$, and
let ${\bf x}=(x,y)$ be a $c$-pair. 
We write
$$c={\rm exp}\left(\sum_{a\in \Delta}r_aa\spcheck\right).$$
Then $I({\bf x})\subseteq \Delta$ is equal to
$$I_c=\{a\in \Delta|r_a\notin{\bf Z}\}.$$
Thus any maximal torus of $Z(x,y)$ is conjugate in $G$ to 
$S_c$.
Finally, there is a unique rank zero $c$-pair up to conjugation in
$L_c=DZ(S_c)$.
\end{proposition}

\begin{proof}
Set $I=I({\bf x})$.
According to Corollary~\ref{determ}, $c\in L_I$ and there is a rank
zero $c$-pair in $L_I$.
Since $L_I$ is simply connected, it follows from
Proposition~\ref{rank0pair}  that there are
integers $n_1,\ldots,n_r$ such that $L_I$ is
isomorphic to $\prod_{i=1}^rSU(n_i)$ and under this isomorphism
$c=\prod_{i=1}^rc_i$ where $c_i$ generates the center of $SU(n_i)$.
According to Proposition~\ref{finite}, writing
$$c={\rm exp}(\sum_{a\in I}r_aa\spcheck)$$
no $r_a$ is an integer. It then follows from Lemma~\ref{c} that
$I=I_c$, and hence that any maximal torus of $Z(x,y)$ is conjugate to 
$S_c$.  It also follows from Proposition~\ref{rank0pair} that the
$c$-pair in
$L_c$ is unique up to conjugation. 
\end{proof}

\begin{corollary}\label{maxtorus}
Let $c\in {\cal C}G$ be given. 
Fix a $c$-pair $(x_0,y_0)$ in $L_c$, and define a map
$S_c\times S_c$ to the space of $c$-pairs
in $G$ by sending $(s_1,s_2)$ to $(s_1x_0,s_2y_0)$. This map factors
to induce a homeomorphism from 
$(\ov S_c\times \ov S_c)/W(S_c,G)$ to the  moduli space of conjugacy
classes of $c$-pairs in $G$.
\end{corollary}

\begin{proof}
According to  Proposition~\ref{rank0pair}, the moduli space  ${\cal
M}^0_{L_I}(c)$ of $c$-pairs of rank zero in $L_I$ is empty unless
$I=I_c$ and the moduli space
${\cal M}^0_{L_{I_c}}(c)$ is a single point. Thus, the actions
of $(F(S_c))^2$ and of $W(S_c, G)$ on ${\cal M}^0_{L_{I_c}}(c)$ are
trivial. The result is now immediate from Corollary~\ref{determ}.
\end{proof}

\noindent
{\bf Proof  of Theorem~\ref{c-pairsthm}.}
Theorem~\ref{c-pairsthm} is an immediate consequence of
Proposition~\ref{torus}  and
Corollary~\ref{maxtorus}. 

\medskip

Prior to Lemma~\ref{c}, given $c$, we have defined a  torus
$S_c$ and the group
$DZ(S_c) = L_c$. There is also a
description of the abstract group $L_c$ in terms of the action of
$w_c$:

\begin{proposition} The group $L_c$ is isomorphic to 
$$\prod
_{\ov a \in \widetilde \Delta/\langle w_c\rangle}SU(n_{\ov a}),$$
where the
$\ov a$ are the orbits of $w_c$ acting on $\widetilde \Delta$ and $n_{\ov
a}$ is the number of elements in $\ov a$.
\end{proposition}
\begin{proof} We know by Proposition~\ref{rank0pair} that $L_c$ is
isomorphic to a product of groups of the form $SU(n_i)$ and that, if
$c_i$ is the component of $c$ in the $i^{\rm th}$ factor, then $c_i$
generates the center of $SU(n_i)$. The vector space
$\frak t$ is a quotient of the vector space $\bigoplus _{a\in \widetilde
\Delta}{\bf R}\cdot a\spcheck$  by the one-dimensional space spanned by
$\sum _ag_aa\spcheck$. The element
$w_c$ acts on $\bigoplus _{a\in \widetilde \Delta}{\bf R}\cdot a\spcheck$
by permuting the $a\spcheck$ and is the identity on $\sum _ag_aa\spcheck$.
For each orbit
$\ov a$, the subspace $\bigoplus _{b\in \ov a}{\bf R}\cdot a\spcheck$ is
$w_c$-invariant, and the eigenvalues of $w_c$ on this subspace are
$\zeta _{n_{\ov a}}^i, i=0,
\dots, n_{\ov a} - 1$, where $\zeta _{n_{\ov a}}$ is a primitive $n_{\ov
a}^{\rm th}$ root of unity.  On the
other hand,
$w_c$ is conjugate to the product of the $w_{c_j}$, and the eigenvalues
of $w_{c_j}$ on the subspace of $\frak t$ corresponding to the
simple factor $SU(n_j)$ are $\zeta _{n_j}^i, i=0, \dots, n_j-1$.
The proposition follows by comparing the two forms for
the set of   eigenvalues.
\end{proof}

We conclude with a preliminary normal form for a $c$-pair; we shall give a
more precise form in Section 6.

\begin{defn}
A $c$-pair $(x,y)$ is said to be {\sl in weak normal form\/} (with respect
to  the maximal torus  $T$ and the alcove $A$) if $x\in T$, $x$ is
  the image under the exponential mapping of a point $\tilde x\in \frak
t^c$
  and $y\in N_G(T)$ projects to $w_c$ in $W(T,G)$.
\end{defn}

\begin{corollary}\label{weaknormalform}
Every $c$-pair in $G$ is conjugate to one in weak normal form. 
\end{corollary}

\begin{proof}
By Corollary~\ref{maxtorus}, after conjugation, we can assume that the
maximal torus of 
$Z(x,y)$ is $S_c$. There is a rank zero $c$-pair $(x_0,y_0)$ in
$L_c=DZ(S_c)$ and elements $(s_1,s_2)\in S_c\times S_c$ such that
$(x,y)=(s_1x_0,s_2y_0)$ . The intersection $T\cap
L_c$ is the 
maximal torus of 
$L_c$. According to Corollary~\ref{orders} we can assume that 
$x_0$ is the barycenter of an alcove $A'$ of $T\cap L_c$ and that $y_0$
normalizes this torus and has image in $W(L_c)$ equal to the Weyl
part $w$ of the action of $c$ on the alcove $A'$ for $L_c$.
Thus, $x\in T$ and $y$ normalizes $T$ and projects to 
the image of $w\in W(L_c)$ in $W(G)$. Since this image is the Weyl
part of the action of $c$ with respect to any alcove $A$ for $G$
containing
$A'$, the result now follows by conjugating $x$ and $y$ by an element of
$W(T,G)$ which sends $A'$ to $A$.
\end{proof}

\section{Commuting triples}\label{commtripsect}

In this section $G$ is  simple   and
$c_{ij}=1$ for  
$1\le i<j\le 3$. In other words, we consider conjugacy classes
of commuting triples $(x,y,z)$ in $G$.
We denote the moduli space of conjugacy classes of such
triples by ${\cal T}_G$.

\subsection{Commuting triples of rank zero}

As usual, let $A$ be the  alcove  $T$ containing the origin
coresponding to the set of simple roots  $\Delta$.

\begin{lemma}\label{rank0trip}
Let $(x,y,z)$ be a commuting triple in $G$ of rank zero.
Then $x$ is conjugate in $G$ to the image under the exponential
mapping of a vertex $v$ of the alcove $A$.
Let $a\in \Delta$ be the root with the property that $\{a=0\}$ defines
the wall of $A$ opposite $v$.
Then $g_a=h_a$ and the order of $x$ in $G$ is
$g_a$. For every $b\in 
\widetilde \Delta-\{a\}$ we have that $g_a\not|g_b$.
Each of $y$ and $z$ is conjugate in $G$ to $x$ and 
each has order equal to $g_a$. Conversely, let $k$ be a positive
integer such that $k|g_a$ for exactly one $a$. Then $k=g_a$ and there
exists a commuting triple $(x,y,z)$ in $G$ of rank zero such that the
order of $x$ is $k$.
\end{lemma}

\begin{proof}
Let $(x,y,z)$ have rank zero.
By Lemma~\ref{vertex},  $Z(x)$ is semi-simple and 
$x$ is conjugate to the image under the exponential
mapping of a vertex  $v\in A$.
Letting $a$ be the simple root such that $\{a=0\}$ defines the wall
of
$A$ opposite
$v$, we see that the order of $x$ modulo ${\cal C}G$ is $h_a$. In
particular,
$h_a$ divides the order of $x$.

The pair $(y,z)$ is a commuting pair  in
$Z(x)$ of rank zero.
Lift $y,z$ to elements $\tilde y,\tilde z$ in the universal covering
$\widetilde Z(x)$ of $Z(x)$. Then $(\tilde y,\tilde z)$ is a rank zero
$c$-pair for some 
$c\in \pi_1(Z(x))\subseteq {\cal C}\widetilde Z(x)$. 
According to Lemma~\ref{coroottors}, the group $\pi_1(Z(x))$ 
is cyclic of order $g_a$ and
is generated by $c_0={\rm exp}(\zeta)$, where
$$\zeta=-\frac{1}{g_a}\sum_{b\in\widetilde\Delta-\{a\}}g_bb\spcheck.$$ 
By  Proposition~\ref{rank0pair}, the existence of a rank zero
$c$-pair implies that $c$ generates $\pi_1(Z(x))$ and  that, if ${\rm
exp}(\zeta')= c$, then,  
in the expression of $\zeta'$ as a linear combination of the
simple coroots, all coefficients are non-integral. Since $c$ is a
power of $c_0$, the same is true for $\zeta$. This implies that, for each
$b\in \widetilde \Delta-\{a\}$, the integer $g_b$ is not divisible by $g_a$.

By Proposition~\ref{rank0pair}, the elements $\tilde y$ and $\tilde
z$ are conjugate in $\widetilde Z(x)$, and hence $y$ and $z$ are 
conjugate in  $Z(x)$. It  follows from the same proposition that
the subgroup of $Z(x)$ 
  generated by $y,z$ is isomorphic to
$({\bf Z}/g_a{\bf Z})^2$. In particular, $y$ and $z$ have the same
order $g_a$  in $G$.

Interchanging the roles of $x$ and $y$ in this construction, we see
that $x$ and $z$ are conjugate.
Thus, $x,y,z$ are all conjugate in $G$ and hence all have the same
order, $g_a$.
Since we have already shown that the order of $x$ is divisible by
$h_a$ and  since $g_a|h_a$, it follows  that $g_a=h_a$.

To see the converse, suppose that $k|g_a$ for exactly one $a$.
Then,  by Lemma~\ref{atleast2}, $k=g_a$.  Let $x$ be the image
under the exponential map of the vertex of the alcove opposite
the face
$\{a=0\}$. By Theorem~\ref{Antype}, the
universal cover $\widetilde Z(x)$ is a product of groups of
type $A_n$, and   $Z(x) =\widetilde Z(x)/\langle\zeta\rangle$,
where $\zeta$ has order $k$ and projects to a generator of every
factor. By Proposition~\ref{rank0pair}, there is a rank zero
$\zeta$-pair $(\tilde y, \tilde z)$ in $\widetilde Z(x)$. It  
suffices to take
$(y,z)$ to be the image in
$Z(x)$ of $(\tilde y, \tilde z)$.
\end{proof}

We define the {\sl order of a commuting triple of rank zero} to be the
common order of each of its elements.

\begin{corollary}\label{pi1order}
Suppose that $(x,y,z)$ is a rank zero commuting triple in $G$ of order
$k$. Then 
$\pi_1(Z(x))$ is a cyclic group of order $k$. If the image of 
$\sum_{a\in \Delta}r_aa\spcheck$
under the exponential mapping lies in $\pi_1(Z(x))\subseteq{\cal
C}\widetilde Z(x)$  and generates this group, then $r_a\notin {\bf
Z}$ for every $a\in \Delta$. 
\end{corollary}

\begin{remark}\label{notAn}
If $G=SU(n+1)$, then there does not exist  a rank zero
commuting triple in $G$. This follows by Lemma~\ref{rank0trip},
since all of the $g_a$ are one in this case. Of course, it is 
elementary that every commuting $N$-tuple in
$SU(n+1)$ is contained in a maximal torus, and hence can never be
of rank zero.
\end{remark}

\begin{remark}
Suppose that $(x,y,z)$ has rank zero 
and order  $k>1$.
The lifts $\tilde y,\tilde z$ in the universal covering $\widetilde
Z(x)$ form a $c$-pair for some $c\in \pi_1(Z(x))\subseteq {\cal
C}G$. It is easy to see directly that they generate 
a subgroup of $\widetilde Z(x)$ which modulo the center of $\widetilde
Z(x)$ is isomorphic to $({\bf Z}/k{\bf Z})^2$ and which meets the
center in a subgroup of $\pi_1(Z(x))$. Thus, 
the elements $y,z$ generate a subgroup of $Z(x)$ which is
isomorphic to $({\bf Z}/k{\bf Z})^2$ and has trivial intersection with
the center of $Z(x)$. Since $x$ lies in
the center of $Z(x)$ and has order $k$, it follows that the group
generated by $x,y,z$ is
isomorphic to $({\bf Z}/k{\bf Z})^3$.
\end{remark}

\begin{proposition}\label{varphi}
There is a rank zero commuting triple of order $k$ in $G$
if and only if
$k=g_a$ for exactly one $a$. In this case, there are exactly
$\varphi(k)$ conjugacy classes of rank zero commuting triples of order
$k$ in $G$. If $(x,y,z)$ has rank zero and order
$k$, then the other conjugacy classes of such triples are 
represented by $(x,y,z^\ell)$  for $1\le \ell<k$ and $\ell$
relatively prime to 
$k$.
\end{proposition}

\begin{proof}
The first statement follows from Lemma~\ref{rank0trip}.
Suppose that $(x,y,z)$ is rank zero commuting
triple of order $k$. 
According to Lemma~\ref{rank0trip}, $x$ is conjugate to 
the image under the exponential mapping of the vertex of the alcove
opposite the face of $A$ defined by 
$\{a=0\}$ where $a$ is the unique element of $\widetilde \Delta$ such that
$k=g_a$. 
Then $(y,z)$ is a  rank zero commuting pair in
$Z(x)$. Let
$(\tilde y,\tilde z)$ be a lift 
of $(y,z)$ 
to the universal covering $\widetilde Z(x)$ of $Z(x)$. This is  a
rank zero $c$-pair for some $c$ generating $\pi_1(Z(x))\subseteq {
\cal C}\widetilde Z(x)$. 
It follows that $\widetilde Z(x)$ is $\prod_{i=1}^rG_i$ with
$G_i$ isomorphic to $ SU(n_i)$ for an integer $n_i|k$ and 
$c=c_1\cdots c_r$ where $c_i$ generates the center of $G_i$.
The element $c\in \pi_1(Z(x))$ dependes only on the conjugacy
class of $(y,z)$ in $Z(x)$.
By Proposition~\ref{rank0pair},
$c$ determines the conjugacy class of $(\tilde y,\tilde z)$ in 
$\widetilde Z(x)$ and hence 
the conjugacy class of  $(y,z)$
in $Z(x)$. On the other hand, again by
Proposition~\ref{rank0pair}, for each
$c'\in \pi_1(Z(x))\subseteq{\cal C}\widetilde Z(x)$ there is
a $c'$-pair $(\tilde y',\tilde z')$ in $\widetilde Z(x)$.
Moreover, the 
$c'$-pair $(\tilde y',\tilde z')$ is of rank zero if and only
if $c'$ generates $\pi_1(Z(x))$, and in this case $(\tilde
y',\tilde z')$ is unique up to conjugation in $\widetilde Z(x)$.
The image
$(y',z')$ is a rank zero commuting pair in $Z(x)$. The group
$\pi_1(Z(x))$  is cyclic of order $k$ and hence has
$\varphi(k)$ generators.
This shows that there are exactly $\varphi(k)$
conjugacy classes of commuting pairs of rank zero in $Z(x)$, and hence
$\varphi(k)$ conjugacy classes of rank zero commuting triples in $G$. 
 
Clearly, if $[\tilde y,\tilde z]=c$ then $[\tilde y,\tilde
z^\ell]=c^\ell$. This proves the last statement.
\end{proof}

\subsection{A list of all simple groups with rank zero commuting 
triples}

Suppose that $G$ is simple and contains a rank zero commuting
triple of order
$k$. Then by Lemma~\ref{rank0trip}, there is exactly one coroot
integer $g_a$ which is divisible by
$k$, and in fact $g_a=k$. Conversely, if 
there is exactly one coroot integer $g_a$ equal to
$k$, then by Lemma~\ref{atleast2}, none of the other coroot integers
is divisible by
$k$ and $G$ contains a rank zero commuting
triple of order
$k$.  Examining the coroot integers on the Dynkin diagrams of the
simple groups, one sees that the following are the only
possibilities:
\begin{enumerate}
\item $k=1$: $G$ is the trivial group.
\item $k=2$: $G$ is of type $D_4$, of type $B_3$, or of type
$G_2$.
\item $k=3$: $G$ is of type $E_6$ or $F_4$.
\item $k=4$: $G$ is of type $E_7$.
\item $k=5$ or $k=6$: $G$ is of type $E_8$.
\end{enumerate}

\subsection{Action of the outer automorphism group of $G$}

\begin{proposition}\label{1234}
Let ${\cal T}^0_G$ be the space of conjugacy classes of commuting
triples of rank zero in $G$. 
Then the action of every automorphism $\sigma$ of $G$ on ${\cal
T}^0_G$ is trivial. 
\end{proposition}

\begin{proof}
Let $(x,y,z)$ be a rank zero commuting triple  in
$G$.  Let $k$ be the order of $(x, y, z)$. There is a
unique $a\in \widetilde \Delta$ such that $k|g_a$. 
After conjugation we can assume that $x$ is the 
exponential   of the vertex $\tilde x$ of the alcove
$A$ opposite the face $\{a=0\}$.
After composing $\sigma$ with a suitable inner automorphism, we can
assume that
$\sigma$ normalizes 
$\Delta$. The action of $\sigma$
on the set of simple roots 
preserves the integers $\{g_b\}_{b\in \widetilde \Delta}$. Hence
$\sigma(a) = a$ and therefore $\sigma (\tilde x) = \tilde
x$ and $\sigma(x) = x$. 
Thus
$\sigma$ acts  on $\pi_1(Z(x))$, and by
Corollary~\ref{pi1}, since $\sigma$ preserves the coroot integers
$g_b$,  the action is trivial. Hence
$\sigma$ acts  trivially on the conjugacy class of $(y,z)$ in
$Z(x)$, by Proposition~\ref{rank0pair}, and therefore on the
conjugacy class of
$(x,y,z)$ in $G$.
\end{proof}

\subsection{Action of the center of $G$}

There is an action of $({\cal C}G)^3$ on the space of conjugacy
classes of commuting triples defined by 
$(\gamma_1,\gamma_2,\gamma_3)\cdot(x,y,z)=(\gamma_1x,\gamma_2y,
\gamma_3z)$.  Clearly,
$Z(x,y,z)=Z(\gamma_1x,\gamma_2y,\gamma_3z)$ so that this action
preserves the subspace of conjugacy classes of commuting triples
of rank zero.

\begin{proposition}\label{action}
The induced action of $({\cal C}G)^3$ on ${\cal T}^0_G$ is trivial.
\end{proposition}

\begin{proof}
Let $(x,y,z)$ be a rank zero commuting triple  of order $k>1$.
First consider the action of ${\cal C}G$ on $x$.
We can assume that $x$ is the image under the exponential mapping of the
vertex of the alcove $A$ 
opposite the face $\{a=0\}$ where $a\in\widetilde \Delta$ is the unique
element with $k=g_a$.
For any $\gamma\in{\cal C}G$ let $w_\gamma$ be the Weyl element
which is the linear part of the action of $c$ on $A$.
The Weyl element $w_\gamma$ normalizes $\widetilde \Delta$
and, according to Subsection~\ref{action1}, preserves the $g_b$ in the
sense that $g_{w_\gamma\cdot b}=g_b$.
This implies that $w_\gamma\cdot a=a$.
Thus, the affine automorphism $\varphi_\gamma$ of $A$ fixes $x$.
By Equation~\ref{conjeqn}, this means that 
if $h\in N_G(T)$ projects to  $w_\gamma\in W$, then
$hxh^{-1}=x\gamma^{-1}$. Thus the triples $(x\gamma, y,z)$ and
$(x, {}^hy, {}^hz)$ are conjugate. 

\begin{claim}
Conjugation by $h$ normalizes $Z(x)$ and induces the identity
automorphism of its fundamental group.
\end{claim}

\begin{proof}
Since $Z(x)=Z(x\gamma)$ the first statement is clear.
By Lemma~\ref{pi1}, the exponential map identifies  the
subgroup of 
$(Q\spcheck\otimes {\bf Q})/Q\spcheck$ generated by $ \zeta
=(-1/g_a)\sum_{b\in\widetilde \Delta-\{a\}}g_bb\spcheck$ with
the fundamental group of $Z(x)$.
The element $w_\gamma$   normalizes the set $\widetilde
\Delta-\{a\}\subset {\frak t}^*$ and the subset $\widetilde
\Delta\spcheck-\{a\spcheck\}\subset {\frak t}$.
Since 
$g_{w_\gamma b}=g_b$, for $b\in \widetilde \Delta$,
it is clear that $w_\gamma$ fixes  $\zeta\in{\frak t}$, and 
hence acts trivially on $\pi_1(Z(x))$.
\end{proof}

It follows immediately from the claim that $({}^hy, {}^hz)$ lifts
to a $c$-pair in $\widetilde Z(x)$, where $c=[\tilde y, \tilde
z]$  for any two lifts of $y,z$ to $\widetilde Z(x)$.  By
Proposition~\ref{rank0pair}, $({}^hy, {}^hz)$ is conjugate in
$Z(x)$ to $(y,z)$. Hence $(x\gamma, y,z)$ is conjugate to
$(x,y,z)$, and so the action of
${\cal C}G$ on the first factor of commuting triples induces the
trivial action on the space of conjugacy classes of commuting
triples of rank zero. The situation is completely symmetric in
$x,y,z$ and it then follows that the action of ${\cal C}G$ on the
space of conjugacy classes of commuting triples of rank zero is
trivial.
\end{proof}

\subsection{The general case}

If  $k\ge 1$ divides at least one of the $g_a$, we set
$\widetilde I(k)\subseteq\widetilde \Delta =\{\,a\in\widetilde
\Delta:k\not|g_a\}$. Let $f(k)$ be the face of $A$ which is the
intersection of all the walls of $A$ defined by  the $a\in
\widetilde I(k)$ and let
$${\frak t}(k)=\bigcap_{a\in\widetilde I(k)}{\rm Ker}\, a.$$
Of course ${\frak t}(k)$ is the linear space parallel to $f(k)$.
By convention, if $k=1$ then we let $\frak t(k) = \frak t$. In all
other cases $\frak t(k)$ is a proper subspace of $\frak t$. Let 
$S(k)$ be the subtorus of $T$ whose Lie algebra is ${\frak t}(k)$ and
let 
$L(k)$ be the derived group of $Z(S(k))$.

\begin{proposition}\label{component}
Let $(x,y,z)$ be a commuting triple. Then there is a unique
positive integer
$k\ge 1$ dividing at least one of the $g_a$ with the following 
properties:
\begin{enumerate}
\item $S(k)$ is conjugate to a maximal
torus of $Z(x,y,z)$.
\item After conjugation, we can find a decomposition
$(x,y,z)=(s_1x_0,s_2y_0,s_3z_0)$ where 
$s_i\in S(k)$ and $(x_0,y_0,z_0)$ is a commuting triple of rank zero
and of order $k$ in $L(k)$.
\item The element $x$ is conjugate to the exponential of an element
of $f(k)$.
\end{enumerate}
We call $k$ the {\rm order} of the commuting triple $(x,y,z)$.
\end{proposition}

\begin{proof}
Let $(x,y,z)$ be a commuting triple,
let $S$ be a maximal torus of $Z(x,y,z)$ and let ${\frak s}={\rm
Lie}(S)$.
By Theorem~\ref{main}, there is a commuting triple $(x_0,y_0,z_0)$ of
rank zero 
in $L=DZ(S)$ such that $x\in S\cdot x_0$. Let $k$ be the order of 
$(x_0,y_0,z_0)$. If $k=1$, then $L$
is   trivial.  In this case,  $x,y,z\in S$ and hence they are  
all contained in a maximal torus. Thus $S$ itself is a
maximal torus for $G$. Conversely, if
$x,y,z$ are contained in a maximal torus $T$, then $T$ is a
maximal torus for $Z(x,y,z)$ and hence $L$ is trivial. For the rest
of the proof, we assume that
$L$ is not trivial, or equivalently that
$S$ is not a maximal torus.

\begin{claim}\label{whoknows}  $L$ is a simply connected,
simple group, not of type $A_n$ for any $n$. If $G$ is not simply
laced, then
$L$ is not simply laced.
\end{claim}

\begin{proof} Let $L=\prod_i L_i$ be the decomposition of $L$ into its simple
factors. Since $L=DZ(S)$, $L$ is simply connected and the Dynkin
diagram of $L$ is identified with a subdiagram of the Dynkin
diagram of $G$. Thus, at most
one of the components $L_i$ is not of type $A_n$.
Since $L$ has a commuting triple of rank zero, it follows that each
of the $L_i$ has such a triple. By Remark~\ref{notAn} this implies
that no $L_i$ is of type $A_n$. Thus, $L$ is simple. The last
statement is clear since, if $G$ is not simply laced, then $D(G)$
is a chain, and hence every simply laced subdiagram is of type
$A_n$ for some $n$. 
\end{proof}

\begin{claim}\label{generic}
For a generic $x'\in S\cdot x_0$, $DZ(x')=Z_{L}(x_0)$. In particular,
$Z_L(x_0)$ is semi-simple.
\end{claim}

\begin{proof}
Suppose that $x'$ is generic in $S\cdot x_0$.
The roots in $DZ(x')$ are the roots that annihilate $x'$.
Since $x'$ is generic, a root annihilates $x'$ if and only if it 
annihilates $S$ and annihilates $x_0$. This shows that the roots of
$G$ annihilating $x'$ are exactly the roots of $L$ annihilating
$x_0$ and hence that $DZ(x')=DZ_{L}(x_0)$.
Since $(x_0,y_0,z_0)$ is a rank zero commuting triple in $L$,
the center of $Z_{L}(x_0)$ is finite. Since $Z_L(x_0)$ is connected,
it is semi-simple and equal to its own derived group. Thus
$DZ(x')=DZ_{L}(x_0) = Z_L(x_0)$.
\end{proof}

Let $x'\in S\cdot x_0$ be generic.
After conjugation, we can assume that $x'$ is the
the image under the exponential mapping of a point in $A$. 
We let $\widetilde I(x')\subseteq \widetilde \Delta$
be the subset consisting of all the roots vanishing on $x'$. 
The subset $\widetilde I(x')$ forms a set of simple roots for the
root system of
$DZ(x')$ with respect to the maximal torus $T\cap DZ(x')$.

It now follows from  Corollary~\ref{pi1order} that 
$\pi_1(Z_L(x_0))$ is a cyclic group of order $k$, and hence
$\pi_1(DZ(x'))$ is cyclic of order $k$. By Corollary~\ref{pi1},
this means that $k$ divides $g_a$ for every $a\in  \widetilde
\Delta-\widetilde I(x')$.  That is to say  $\widetilde I(k)\subseteq
\widetilde I(x')$. Hence $x'$ is the image under the exponential
of a point of $f(k)$. Since $x'$ is generic, the same conclusion
holds up to conjugation for   $x$.

Since $(x_0,y_0,z_0)$ is a rank zero triple in $L$, it follows from
Corollary~\ref{pi1order}  that  if $I$ is any set of simple roots for
the root system $Z_L(x_0)$ and if 
$c=\sum_{a\in I}r_aa\spcheck$ has the property that its image 
under the exponential mapping is a central element in
$\widetilde Z_L(x_0)$ generating
$\pi_1(Z_{L}(x_0))$, then no $r_a$ is an integer.
But according to Corollary~\ref{pi1}, $\widetilde
I(x')$ is a set of simple roots for the root system $DZ(x')=Z_L(x_0)$
and the exponential of 
$$-\frac{1}{k}\sum_{a\in\widetilde I(x')}g_aa\spcheck$$
generates the fundamental group of $DZ(x')=Z_L(x_0)$. It follows
that 
$k\not|g_a$ for every $a\in \widetilde I(x')$. This
implies that
$\widetilde I(x')\subseteq \widetilde I(k)$, and hence that
$\widetilde I(x')=\widetilde I(k)$.

Hence $x'$ is contained in the
interior of the face $f(k)$ of $A$ and 
$DZ(x')$ has $\widetilde I(k)$
as a set of simple roots. Consequently, $Z_L(x_0)$ has $\widetilde
I(k)$ as a set of simple roots. Since $Z_L(x_0)$  is semi-simple,
${\rm Lie}(Z_L(x_0))\cap {\frak t}={\rm Lie}(L)\cap {\frak t}$. 
Thus, ${\frak s}={\rm Lie}(S)$ is the perpendicular space to
$\{a\spcheck\}_{a\in\widetilde I(k)}$. It follows  
that
${\frak s}={\frak t}(k)$ and $L=L(k)$.
\end{proof}

We now establish a converse to Proposition~\ref{component}:

\begin{proposition}\label{converse}
Suppose that $k\geq 1$ is a positive integer dividing at least one of
the $g_a$. Then there exists a commuting triple of order $k$ in $G$.
\end{proposition}
\begin{proof} We begin with the following result about root systems:

\begin{proposition}\label{kroots}
Let $\Phi$ be a reduced and irreducible root system on a vector
space $\frak t$, and suppose that
$k>1$ is an integer such that $k|g_a$ for some $a$. Define
$\widetilde I(k)$ and $\frak t(k)$ as before, and let $\Phi(k)$ be
the set of all roots which annihilate $\frak t(k)$. Then $\Phi(k)$
is an irreducible root system. Moreover, if the coroot integers for
$\Phi(k)$ are of the form $m_b, b\in \widetilde \Delta (\Phi(k))$,
then $k|m_b$ for exactly one $b$, and in this case $k=m_b$.
\end{proposition}
\begin{proof} We may assume that $\Phi$ is the root system of a
simple and simply connected group $G$. Thus there is a torus $S(k)$ 
corresponding to $\frak t(k)$. Let $L(k) = Z(S(k))$. Then
$\Phi(k)$ is the set of roots for $L(k)$, and in particular it is a
root system. Let
$Q_{L(k)}\spcheck$ be the sublattice of $Q\spcheck$ generated by the
coroots of  $L(k)$. Then $Q_{L(k)}\spcheck$ is a primitive
sublattice of $Q\spcheck$, by Lemma~\ref{subset}. Let $Q_{\widetilde
I(k)}\spcheck$ be the lattice spanned by $\widetilde I(k)$. Then
$Q_{\widetilde I(k)}\spcheck$ is a sublattice of $Q_{L(k)}\spcheck$.
By Lemma~\ref{coroottors},
${\rm Tor}(Q\spcheck/Q_{\widetilde I(k)}\spcheck)$ is a cyclic group
of order equal to the gcd of the $g_a$ such that $k|g_a$, and by
Corollary~\ref{gcd}, this gcd is $k$. Since $Q_{L(k)}\spcheck$ is a
primitive sublattice of $Q\spcheck$, it follows that
$Q_{L(k)}\spcheck/Q_{\widetilde I(k)}\spcheck$ is also cyclic of
order $k$.

Next, we have the following description of the root system $\Phi(k)$:

\begin{claim}\label{kroots2}
Let $\Phi^+$ be the set of positive roots for $\Phi$ corresponding
to $\Delta$ and let $I(k) = \widetilde I(k)\cap \Delta$. Then
$\Phi^+(k)=\Phi^+\cap  
\Phi(k)$ is a set of positive roots for $\Phi(k)$. Let $\Delta(k)$ be
the corresponding set of simple roots. Then $\Delta(k)=
I(k)\cup \{b\}$ for some root $b\in \Phi(k)$. The root system
$\Phi(k)$ is irreducible and the highest root $d$ for $\Phi$  is
also a highest root for $\Phi(k)$.
\end{claim}

\begin{proof} Choose any $\tilde p$ contained in the interior
of
$A$. Then the roots of $\Phi(k)$ which are positive on $\tilde p$ are
exactly those in $\Phi^+(k)$. Thus, $\Phi^+(k)$ is
a set of positive 
roots with respect to some set of simple roots of $\Phi(k)$.
Clearly, the elements of $I(k)\subseteq \Phi(k)$ are positive roots.
Since none of these can be written as a non-trivial linear
combination of positive roots of
$\Phi$, {\it a fortiori} none of these can be written as a
non-trivial linear combination of positive roots of $\Phi(k)$. Thus
$I(k)$ is a subset of the set of simple roots $\Delta(k)$ determined
by $\Phi^+(k)$. Since   $g_{\tilde a}=1$,  $\tilde
a\in\widetilde I(k)$ implying that  the cardinality of $I(k)$ is one 
less than the dimension of the span of $\Phi(k)$. Thus, there is a
root $b\in
\Phi^+(k)$ with the property that
$\Delta(k)=I(k)\cup \{b\}$.

Let $d$ be the highest root of $\Phi$ with respect to the positive
roots
$\Phi^+$. Since $\tilde a \in \widetilde I(k)$, $d=-\tilde a \in
\Phi^+(k)$. Write $d\spcheck=\sum_{a\in
I(k)}m_aa\spcheck+m_bb\spcheck$. Since
$\{a\spcheck, a\in I(k)\} \cup
\{b\spcheck\}$ is a basis for
$Q_{L(k)}\spcheck$ and $\{ a\spcheck, a\in I(k)\} \cup
\{d\spcheck\}$ is a basis for $Q_{\widetilde
I(k)}\spcheck$, it follows that $Q_{L(k)}\spcheck/Q_{\widetilde
I(k)}\spcheck$ is  cyclic of order $m_b$. Thus $m_b=k$. Since 
$k\not|g_a$ for all $a\in I(k)$, it follows that $m_a\neq 0$ for
all $a\in I(k)$.  This proves that all 
the coefficients is this expression are non-trivial, and hence
$\Phi(k)$ is irreducible.
Since the sum of $d$ 
and any positive root in $\Phi$ is not a root of $\Phi$, it follows
that $d$ is the highest root of $\Phi(k)$ with respect to the set of
simple roots
$\Delta(k)$.
\end{proof} 

Returning to the proof of Proposition~\ref{kroots}, we see that we
have proved that $\Phi(k)$ is irreducible and 
that $k=m_b$ and $k\not|m_a$ for $a\neq b$. This completes the
proof of Proposition~\ref{kroots}.
\end{proof}

Finally, let us finish the proof of Proposition~\ref{converse}.
Since $L(k)$ is a simple group and $k$ is equal to exactly one of
the coroot integers of $L(k)$, it follows by 
Proposition~\ref{varphi} that $L(k)$ contains a commuting rank zero
triple of order $k$. Of course, such a triple will also be a
commuting triple of order $k$ in $G$.
\end{proof}

\begin{theorem}\label{commtrip}
Let $G$ be simple.
Let $k\ge 1$ be an integer.
\begin{enumerate}
\item If $(x,y,z)$ is a commuting triple of order $k$ in $G$, then
$k$ divides at least one of the coroot integers $g_a$ and 
$S(k)$ is conjugate  to a maximal torus for $Z(x,y,z)$. 
\item The order is a conjugacy class invariant and defines a
locally constant function on ${\cal T}_G$.
\item  If $k$ divides at least one of the $g_a$, there  are 
exactly $\varphi(k)$ components of ${\cal T}_G$ consisting of
conjugacy classes of commuting triples of order $k$, where
$\varphi$  is the Euler $\varphi$-function. Given a component $X$
of ${\cal T}_G$, let $d_X = \frac13\dim X + 1$. Then
$$\sum _X d_X = g.$$
\item  Each component consisting of commuting triples of
order $k$ in $G$ is homeomorphic to 
$$\left(\ov S(k)\times \ov S(k)\times \ov S(k)\right)/W(S(k), G).$$
\end{enumerate}
\end{theorem}

\begin{proof}
The first statement follows from Proposition~\ref{component}.
Clearly, the order is a conjugacy class invariant, and it is locally
constant on ${\cal T}_G$ by Corollary~\ref{cor}. Now suppose that
$k$ divides at least one of the $g_a$. By
Proposition~\ref{converse}, there is a commuting triple ${\bf x}$ of
order
$k$ in $G$. By Part (1), we can assume after conjugation that $S(k)$
is a maximal torus of $Z({\bf x})$. By Part (2) of
Proposition~\ref{component}, there is a rank zero commuting triple
of order $k$ in $L(k)$.    It then follows from Lemma~\ref{varphi}
there are exactly
$\varphi(k)$ conjugacy classes of commuting triples of rank zero in
$L(k)$. By Proposition~\ref{action},   the 
center of $L(k)$ acts trivially on the space of conjugacy classes of
commuting triples in $L(k)$ or, in 
the notation of Corollary~\ref{determ}, the group  $F^3$ acts
trivially on the space of conjugacy classes of commuting triples in
$L(k)$. By Proposition~\ref{1234}, the Weyl group of $W(S(k), G)$
acts trivially on the set of conjugacy classes of rank zero
commuting triples in $L(k)$.
Corollary~\ref{determ} now implies that there are exactly $\varphi(k)$
components of ${\cal T}_G$ of triples of order $k$, and each of
these components is homeomorphic to
$$(\ov S(k)\times \ov S(k)\times \ov S(k))/W(S(k),G),$$ proving
the first sentence in Part  (3) and Part (4) of the theorem. Let
$X$ be a component of ${\cal T}_G$ of order $k$. By Part (4),
$$\frac13\dim X + 1 = \dim S(k) + 1.$$
It follows directly from the definition of $S(k)$ that $\dim S(k)
+ 1$ is equal to the number of $a$ such that $k|g_a$. The second
statement of Part (3) then follows from
the first statement of Part (3) and Lemma~\ref{standardfact}.
\end{proof}

The first four parts of Theorem~\ref{commuttrip} are contained in
the statement of Theorem~\ref{commtrip}.  We shall prove the last
item of Theorem~\ref{commuttrip}  in Section~\ref{sd}.

\begin{remark} Assume that $G\neq L(k)$, in other words that the
$c$-triple $(x,y,z)$ has positive rank. We have defined
$\ov S(k)$ to be
$S(k)/(S(k)\cap L(k))$. Here $S(k)\cap L(k)\subseteq {\cal
C}L(k)$, and is easily checked to be
$\pi_k(Q\spcheck)/Q_{L(k)}\spcheck$, where
$Q_{L(k)}\spcheck$ is the coroot lattice of $L(k)$ and $\pi_k$ is
orthogonal projection onto the real vector space spanned by
$Q_{L(k)}\spcheck$. Using this remark, it is not difficult to
check that
$S(k)\cap L(k)= {\cal C}L(k)$ except for the case where $G$ is of
type $D_n$ for
$n>4$ and $k=2$, so that $L(2)$ is of type $D_4$. In this case,
$S(2)\cap  L(2)$ has order $2$.
\end{remark}

\section{Some results on diagram automorphisms and associated root
systems}

Let $\Phi$ be a reduced but not necessarily irreducible root system on the
vector space
$\frak t$ with a basis $\Delta$ a set of simple roots.  Let
${\cal  A}$ be the decomposition of $\frak t$ into alcoves determined by
the set of affine walls ${\cal W}$ associated to 
$\Phi$.  Suppose that
$\tau$ is a group of  affine isometries of $\frak
t$ normalizing the alcove decomposition ${\cal A}$. Suppose that
$A_0$ is an an alcove such that $\tau(A_0) = A_0$. After
conjugating by an element of the affine Weyl group group,  we can
assume that
$ A_0 =A$ is the alcove associated to  $\Delta$. If $\Phi$ is
irreducible, then $A$ is a simplex and every group $\tau$ of
affine isometries of $A$ fixes the barycenter of $A$, which is an
interior point of $A$. In general, the group $\tau$ fixes the
product of the barycenters of the  factors of $A$. Let
$\ell$ be the associated group of linear isomorphisms of $\frak t$ and
let $\frak t^\ell$ be the fixed subspace of $\ell$.  Clearly, the group
$\ell$ normalizes $\widetilde \Delta$ and  defines a group of diagram
automorphisms of $\widetilde D(\Phi)$. Conversely, every such group of    
diagram automorphisms leads to a group of affine automorphisms $\tau$
as above. We denote 
by $\widetilde \Delta/\ell$ the quotient set. Note that $\ell$
acts on the set $\Phi$ of all roots as well.

 The purpose of this
section is to study the set of nonzero restrictions of elements of $\Phi$
to
$\frak t^\ell$. We show that this set forms a root system $\Phi^{\rm
res}(\ell)$ whose Weyl group is the Weyl group of $\frak t^\ell$ in $W$,
and explicitly identify the inverse coroots. There are two other
closely related root systems with  the same Weyl group which we
also study. Related results, in a more general context, have been
given in \cite{FSS}  and 
\cite{FRS}.

\subsection{A chamber structure and a Coxeter group on the
fixed subspace}\label{Linear}

\begin{lemma}
\begin{enumerate}
\item
No wall $W$ of ${\cal W}$ contains the fixed subspace ${\frak
t}^\tau$ of $\tau$.
\item If $W$ is a wall of ${\cal W}$ meeting ${\frak
t}^\tau$, then the intersection $W^\tau=W\cap {\frak t}^\tau$ is a
codimension-one affine subspace of ${\frak t}^\tau$.
\item  The walls $W^\tau$ divide ${\frak t}^\tau$ into compact convex
subsets with nonempty interior. We denote this collection of subsets by
${\cal A}^\tau$.
\item The elements of ${\cal A}^\tau$ are exactly the subsets ${\frak
t}^\tau\cap A'$, where $A'$ is an alcove of ${\cal A}$ such that $\tau(A')
=A'$.
\end{enumerate}
\end{lemma}
 
\begin{proof}
Since ${\frak t}^\tau$ contains an interior point of an alcove, no wall of
${\cal W}$ can contain ${\frak t}^\tau$. The second and third statements are
now clear. As to the last, if $A'$ is an alcove of ${\cal A}$
normalized by $\tau$, then $A'\cap {\frak t}^\tau$ contains an interior
point of $A'$, and hence $A'\cap {\frak t}^\tau$ is the closure of its
interior in ${\frak t}^\tau$. Clearly, $A'\cap {\frak t}^\tau \in {\cal
A}^\tau$. Conversely, let
$B\in{\cal A}^\tau$. Since $B$ contains a non-empty open subset of
${\frak t}^\tau$, it contains a   element of $\frak t-\bigcup_{W\in \cal
W}W$. This shows that $B$ is contained in a unique alcove $A'$ of ${\cal
A}$. Clearly,
$B={\frak t}^\tau\cap A'$ and $\tau (A')=A'$.
\end{proof}

\begin{lemma}\label{restr}
Two elements of $\widetilde \Delta$ have the same restriction to ${\frak
t}^{\ell}$ if and only if they lie in the same $\ell$-orbit. 
If $a\in
\widetilde \Delta$, then $a|{\frak
t}^{\ell} \neq 0$.
\end{lemma}

\begin{proof}
Clearly, it suffices to establish this result in the case that $\Phi$ is
irreducible, so that $A$ is a simplex. 
Suppose that $a,a'\in\widetilde \Delta$ are in the same $\ell$-orbit. Then
their restrictions to ${\frak t}^\ell$ are equal. Furthermore,
the walls $W_a$ and $W_{a'}$ of $A$ determined by $a$ and $a'$ are in
the same $\tau$-orbit. This means that $W_a\cap {\frak
  t}^\tau=W_{a'}\cap {\frak t}^\tau$.

Since the restrictions of the walls of $A$ to ${\frak t}^\tau$ cut out
a compact convex body, there must be at least ${\rm dim}({\frak
  t}^\tau)+1$ distinct and non-parallel walls. But this is exactly the
cardinality of 
$\widetilde \Delta/\ell$. Hence, it follows that
distinct $\ell$-orbits in $\widetilde \Delta$ cut out distinct 
and non-parallel walls in ${\frak t}^\tau$, and hence have distinct, nonempty
restrictions to ${\frak t}^\ell$.
\end{proof}

Note: it is not in general true that a general $a\in \Phi$ has
nonzero restriction to ${\frak t}^\ell$, or that, if two elements
of
$\Phi$ have the same (nonzero) restriction to ${\frak t}^\ell$,
then  they lie in the same $\ell$-orbit. See for example
Lemma~\ref{assumptionOK} below.

\begin{lemma}\label{cent=norm}
\begin{enumerate}
\item
The group $\tau$ normalizes the affine Weyl group
$W_{\rm aff}(\Phi)$.
\item The centralizer $Z_{W_{\rm aff}(\Phi)}(\tau)$ is  equal to
the normalizer $N_{W_{\rm aff}(\Phi)}(\frak t^\tau)$  of 
${\frak t}^\tau$ in $W_{\rm aff}(\Phi)$ and acts simply
transitively on  ${\cal A}^\tau$.
\end{enumerate}
\end{lemma}

\begin{proof}
Since $\tau$ normalizes ${\cal A}$, it normalizes the group generated
by reflections in the walls of ${\cal W}$, i.e., the affine Weyl group
$W_{\rm aff}(\Phi)$.
To establish (2), clearly $Z_{W_{\rm aff}(\Phi)}(\tau)\subseteq N_{W_{\rm
aff}(\Phi)}(\frak t^\tau)$. Conversely, suppose that $g\in N_{W_{\rm
aff}(\Phi)}(\frak t^\tau)$. Then $B =gA$ is an element of ${\cal A}$
meeting $\frak t^\tau$ in an interior point. Thus for all $f\in \tau$,
$f^{-1}gfA = B = gA$. By Part 1, $f^{-1}gf\in W_{\rm aff}(\Phi)$
and hence $f^{-1}gf = g$. Since this is true for all $f\in \tau$, $g\in
Z_{W_{\rm aff}(\Phi)}(\tau)$.  The final statement is now clear since
$W_{\rm aff}(\Phi)$ acts simply transitively on the set of all alcoves.
\end{proof}

\begin{lemma}\label{wall}
For each wall $W\cap {\frak t}^\tau$ there is an element of $Z_{W_{\rm
aff}(\Phi)}(\tau)$ which is a geometric reflection in this wall.
The group $Z_{W_{\rm aff}(\Phi)}(\tau)$ is generated by the
reflections in the walls of any given element of ${\cal A}^\tau$. 
Thus, $Z_{W_{\rm aff}(\Phi)}(\tau)$ is an affine Coxeter group with
fundamental domain $B=A\cap {\frak t}^\tau$.
\end{lemma}

\begin{proof}
The wall $W\cap {\frak t}^\tau$ is a common wall between two alcoves
$B_1,B_2$ of ${\cal A}^\tau$. Let
$A_1,A_2$ be 
the $\tau$-invariant alcoves of ${\cal A}$ containing $B_1,B_2$. Let
$g\in Z_{W_{\rm aff}(\Phi)}(\tau)$ be the unique element carrying $A_1$
to $A_2$. Then $g$ is a product of reflections about walls separating
$A_1$ and $A_2$, and hence it is the identity on $A_1\cap A_2$ and {\it a
fortiori\/} on $B_1\cap B_2$. Since $B_1\cap B_2$ contains a nonempty
open subset of $W\cap {\frak t}^\tau$, 
$g|{\frak t}^\tau$ is an isometry fixing $W\cap {\frak t}^\tau$ and
sending
$B_1$ to
$B_2$. It is then the reflection in $W\cap {\frak t}^\tau$.
\end{proof}

\begin{proposition}\label{multcoroots} Let $\Phi$ be a reduced
root system on $\frak t$ and let
$\tau$ be a group of affine isometries of $\frak t$ normalizing the
alcove decomposition associated to $\Phi$. Suppose that $\tau$ normalizes
the alcove $A$. Then there is a point
$v\in B=A\cap {\frak t}^\tau$ which is a vertex of the alcove
decomposition
${\cal A}^\tau$ so that, using $v$ to  identify ${\frak  t}^\tau$ with
${\frak t}^\ell$,  there is a uniquely determined reduced root system
$\Phi^\tau$  on ${\frak t}^\ell$ whose
 alcove structure is ${\cal A}^\tau$. The affine Weyl group of
$\Phi^\tau$ is
$Z_{W_{\rm aff}(\Phi)}(\tau)$. The root system $\Phi^\tau$ is irreducible
if 
$\Phi$ is irreducible. There is a set of extended simple roots $\widetilde
\Delta^\tau$ for
$\Phi^\tau$  and a bijection  $\iota\colon 
\widetilde \Delta/\ell\to \widetilde \Delta^\tau$ such
that  the restriction mapping
$\widetilde\Delta/\ell\to \left({\frak t}^\ell\right)^*$ sends
$\ov a\in \widetilde\Delta/\ell$ to a positive multiple of $\iota(\ov
a)\in\widetilde \Delta^\tau$. 
\end{proposition}

\begin{proof}
We may assume that $\Phi$ is irreducible. By Lemma~\ref{wall},  there
is a Coxeter group with
$A\cap {\frak t}^\tau$ as fundamental domain. By the general
classification result  for such Coxeter groups \cite{Bour}, it
follows that this Coxeter group is  isomorphic to the affine Weyl
group of a (reduced) root system
$\Phi^\tau$. That is to say  there is a linear structure on ${\frak
t}^\tau$ compatible with its  given  affine structure, such that under
this identification the  Coxeter group becomes the affine Weyl group of a
root system. Of course this linear structure is determined by choosing a
point
$v\in{\frak t}^\tau$ to identify ${\frak t}^\tau$ with ${\frak
t}^\ell$. The point $v$ must be a vertex of an alcove. In fact, we can 
choose it to  be a vertex of the alcove $A\cap {\frak t}^\tau$. When we
do this, the restricted roots  defining the walls of  $A\cap{\frak
t}^\tau$  become the set of extended simple roots. By Lemma~\ref{restr}
these  roots are exactly the walls of the restrictions of the orbits 
$\widetilde\Delta/\ell$. The proof
of Lemma~\ref{restr} shows that $A\cap {\frak t}^\tau$ has $\dim \frak
t^\tau + 1$ walls. Hence
$A\cap {\frak t}^\tau$ is a simplex, so that $\Phi^\tau$ is irreducible.
The last statement follows  since a root is
determined up to a  multiple by the wall it defines, and it is easy to
see in this case that the multiple must be positive.
\end{proof}

\begin{corollary}\label{Weylgroupgen}
The Weyl group $W(\Phi^\tau)$ is the group of isometries of ${\frak
t}^{\ell}$ generated by the reflections in $\{\ov a\}_{\ov a\in
\widetilde\Delta/\ell}.$ 
\end{corollary}

As a first application to the study of $c$-pairs we have the following:

\begin{lemma} Let $(x,y)$ be a $c$-pair in $G$. Then up to conjugation we
can assume that $x=\exp \tilde x$, where $\tilde x$ lies in the 
fixed set $A^c$ of the alcove
$A$, and $y\in N_G(T)$ is an element projecting to $w_c\in
W(T,G) =W$.
\end{lemma}
\begin{proof} By Corollary~\ref{weaknormalform}, we may assume that $x$
is  the image under the exponential mapping of a point
$\tilde x\in \frak t^c$ and that $y\in N_G(T)$ projects to $w_c$ in
$W$. By Proposition~\ref{multcoroots} applied to the affine
automorphism
$\varphi_c$ of Section~\ref{action1}, there is a $\gamma\in W_{\rm
aff}(G)$, commuting with
$\varphi_c$, such that $\gamma\cdot \tilde x \in A^c$. Let $h\in
N_G(T)$ project to the element of $W$ which is the linear part of
$\gamma$. Then
$hxh^{-1}$,  $hyh^{-1}$ satisfy the conclusions of the lemma.
\end{proof}

\begin{defn}\label{normalform}
A $c$-pair $(x,y)$ is said to be in {\sl normal form\/} (with
respect to  the maximal torus  $T$ and the alcove $A$) if $x\in
T$, $x$ is the image under the exponential mapping of a point
$\tilde x\in A^c$  and $y\in N_G(T)$ projects to $w_c$ in
$W(T,G)$. By the above lemma, every $c$-pair is conjugate to one
in normal form.
\end{defn}

\subsection{The restricted root system and the projection
root system}\label{Digression}

In this subsection we keep the notation of   Subsection~\ref{Linear}. Let
$\widetilde D(\Phi)$ be the extended diagram of $\Phi$. The action of
$\tau$ permutes the connected components of $\widetilde D(\Phi)$. If $D$
is  a connected component of $\widetilde D(\Phi)$, then the stabilizer
$\tau_D$ acts as a group of diagram automorphisms of $D$. Moreover, given
$a\in D$, the
$\tau$-orbit of $a$ is a disjoint union of mutually perpendicular copies
of the
$\tau_D$-orbit of $a$. Likewise, $\frak t$ is an
orthogonal direct sum of subspaces $\frak t_D$ indexed by the components
of
$\widetilde D(\Phi)$, and $\frak t^\ell$ is an orthogonal direct sum  over
the
$\tau$-orbits of the set of components of $\widetilde D(\Phi)$, of the
fixed spaces $\frak t_D^{\ell_D}$, where $\ell_D$ is the linear group
corresponding to $\tau_D$.
Thus, in the proofs that follow, we will be able to reduce to the case
where $\Phi$ is irreducible.

\begin{lemma} 
Let ${\cal O}$ be an orbit of
$\ell$ acting on $\widetilde\Delta$.  Then exactly one of the following 
holds:
\begin{enumerate}
\item ${\cal O}$ is a union of components $D$ of $\widetilde D(\Phi)$ of
type $A_n$ and the stabilizer $\tau_D$ of   a component $D$ of ${\cal O}$
acts transitively on $D$. In this case, if $\frak t_D$ is the
subspace of $\frak t$ corresponding to a component $D$ of ${\cal O}$ and
$\ell_D$ is the linear group corresponding to the stabilizer $\tau_D$ of
$D$ in $\tau$, then 
${\frak t_D}^{\ell_D}=\{0\}$. 
\item For all $a,b\in {\cal O}$ with $a\neq b$, $a$ and $b$ are
orthogonal. In this case we say that ${\cal O}$ is a {\rm ordinary} orbit.
\item ${\cal O}=\coprod_i{\cal O}_i$, where each ${\cal O}_i$ is of
cardinality $2$, say ${\cal O}_i=\{a_{i,1},a_{i,2}\}$. Furthermore, 
$a_{i,1}$ and $a_{i,2}$ have the same lengths and 
$n(a_{i,1},a_{i,2})=-1$. Lastly, elements from distinct ${\cal O}_i$
are orthogonal.
We say that ${\cal O}$ is an {\rm exceptional} orbit and that each of
the ${\cal O}_i$ is an {\rm exceptional pair}. If moreover $\Phi$ is
irreducible, not of type $A_n$ and $\tau$ is cyclic, then there is at
most one exceptional orbit and, if it exists, it has exactly two elements.
\end{enumerate}
\end{lemma} 

\begin{proof} 
By the remarks before the statement of the lemma, it suffices to
consider the case where $\Phi$ is irreducible. If
$\Phi$ is of type
$A_n$, then the corresponding group of diagram automorphisms is a
cyclic group or a dihedral group and the lemma follows by
inspection. Otherwise,  
$\widetilde D(\Phi)$ is a connected, contractible diagram which is not a
single orbit. Let
$D$ be the proper subdiagram of $\widetilde D(\Phi)$ defined by ${\cal O}$
and suppose that
$D$ is a union of $k$ connected subdiagrams $D_i$. Then the $D_i$ are
all isomorphic and the stabilizer of $D_i$ in $\ell$ acts transitively
on the nodes of $D_i$.
It follows that each $D_i$ is of type
$A_1$ or to $A_2$.
The first case corresponds to an ordinary orbit, the second to an
exceptional orbit. The final statement of Case 3 follows easily from the
Lefschetz fixed point formula.
\end{proof}

We remark that Case 1 above occurs if and only if $\tau_D$
contains a rotation of order
$n+1$,  or $n+1=2k$ is even, $\tau_D$ contains a rotation of order
$k$, and an involution with no fixed points. 

\begin{lemma}\label{samewall} 
Let $B=A\cap {\frak
t}^\tau$. Let $W$ be the wall of $B$ corresponding to the orbit
${\cal O}\subseteq \widetilde\Delta$. Suppose that  $a\in \Phi$
and that the wall
$W_a\cap \frak t^\tau$ is equal to $W$. Then either $\pm a\in 
{\cal O}$ or
${\cal O} $ is exceptional and $\pm a=a_1+a_2$ for an exceptional 
pair
$\{a_1,a_2\}$ in ${\cal O}$.
\end{lemma}

\begin{proof} Choose a point $\tilde x\in B\cap W$ which is
contained in no other wall. By Lemma~\ref{restr}, the set of
roots in $\widetilde\Delta$ which take integral values on $\tilde
x$ is exactly ${\cal O}$. By Lemma~\ref{alcove}, since
$\tilde x\in A$, the set of roots
$a$ in $\Phi$ such that $a$ is integral on $\tilde  x$ is a root system with
simple roots equal to ${\cal O}$. If ${\cal O}$ is ordinary, then this root
system exactly is a product of root systems of type $A_1$ and  one of
$a,-a$ is contained in ${\cal O}$. Otherwise ${\cal O}$ is exceptional
and this root system is a product of root systems of type $A_2$. Every
root of this system, up to sign, is either in ${\cal O}$ or is the sum of
an exceptional pair in ${\cal O}$.
\end{proof}

Using the previous lemma and Weyl invariance, we can extend the description of
the possible orbit types from the set of orbits of extended roots to all
orbits.

\begin{corollary}\label{orbittype}
Let $a\in \Phi$ and suppose that $a|\frak t^\ell \neq 0$. Let ${\cal O}$ be the
orbit of $a$. Then either the elements of ${\cal O}$ are mutually
orthogonal or ${\cal O}$ is a disjoint union of mutually orthogonal
subsets ${\cal O}_i$ where each ${\cal O}_i =\{a_{i,1}, a_{i,2}\}$
and where  $a_{i,1}$ and $a_{i,2}$ have the same length and $n(a_1, a_2)
= -1$.
\end{corollary}

\begin{proof} Since $a|\frak t^\ell \neq 0$, the wall $W_a$ corresponding to
$a$ meets $\frak t^\tau$ in a hyperplane. Thus there is $g\in Z_{W_{\rm
aff}(\Phi)}(\tau)$ such that $g\cdot W_a$ defines a wall of $B$.
Suppose that
$g\cdot W_a$ corresponds to the orbit ${\cal O}'$ of simple roots. If
$w$ is the linear part of $g$, then $w$ commutes with $\ell$ and
hence sends $\ell$-orbits to $\ell$-orbits. Then result now follows
from Lemma~\ref{samewall}.
\end{proof}

In the first case of the corollary, we call ${\cal O}$  {\sl ordinary} and
in the second case we call ${\cal O}$ {\sl exceptional} and the
subsets ${\cal O}_i$ {\sl exceptional pairs}.

For $a\in \Phi$, let $\ov a$ be the
$\ell$-orbit of $a$. For $a\in \Phi$, $a|\frak t^\ell$ is a linear form,
depending only on the orbit $\ov a$.
Given an orbit
$\ov a$, let $n_{\ov a}$ be the number of elements of $\ov a$.
Define the {\sl restricted roots\/} $\Phi^{\rm res}(\ell)\subseteq ({\frak t}^\ell)^*$ to be the set of nonzero linear maps of the form
$a|{\frak t}^\ell$ for $a\in \Phi$.
Note that distinct orbits may define the same restricted root, although
this does not happen for orbits contained in $\widetilde \Delta$.

For an orbit $\ov a$ in $\Phi/\ell$ we define  $\epsilon(\ov a) = 1$ if
$\ov a$ is ordinary and $\epsilon(\ov a) = 2$ if  $\ov a$ is exceptional.
Now we define the coroots inverse to the elements of $\Phi^{\rm
  res}(\ell)$ as follows:
For each $u\in \Phi^{\rm res}(\ell)$ we choose $a\in \Phi$ such that
$a|{\frak t}^\ell =u$ and define the inverse coroot 
\begin{equation}\label{rooteqn}u\spcheck=\epsilon(\ov
a)\sum_{a'\in\ov a}(a')\spcheck.
\end{equation}

\begin{claim}
$u\spcheck$ is independent of the choice of $a\in \Phi$ restricting to
give $u$ and  $u\spcheck = \epsilon(\ov
a)n_{\ov a}\pi(a\spcheck)$,
where $\pi$ is orthogonal projection ${\frak t}\to {\frak t}^\ell$.
\end{claim}

\begin{proof}
Elements $a,b\in \Phi$ restrict to give the same root in $\Phi^{\rm
res}(\ell)$ if and only if $\pi(a\spcheck)=\pi(b\spcheck)$. Clearly,
$\epsilon(\ov a)\sum_{a'\in \ov a}(a')\spcheck$ is a positive real
multiple of
$\pi(a\spcheck)$ and by Corollary~\ref{orbittype} $\langle
\pi(a\spcheck),\epsilon(\ov a)\sum_{a'\in \ov a}(a')\spcheck\rangle =2$.
The proves the first statement. The second follows from the fact that
$\pi(a\spcheck) = (1/n_{\ov a})\sum_{a'\in\ov a}(a')\spcheck$.
\end{proof}

\begin{proposition}\label{Wres=Wtau}
The set $\Phi^{\rm res}(\ell)$ is a  possibly nonreduced root system in
$\frak t^\ell$. It is irreducible if $\Phi$ is irreducible. For $u\in
\Phi^{\rm res}(\ell)$, the coroot inverse to $u$ is $u\spcheck$
as given in Equation~\ref{rooteqn}.  The Weyl group of
$\Phi^{\rm res}(\ell)$ is equal to $W(\Phi^\tau)$. 
\end{proposition}

 We call
$\Phi^{\rm res}(\ell)$ the {\sl restricted root system}. 

\begin{proof} Clearly the $\{\ov a: a\in \Phi\}$ span the dual space to $\frak
t^\ell$ and hence the roots of $\Phi^{\rm res}(\ell)$ span this space.
As we saw in the above claim, for any $u\in \Phi^{\rm res}(\ell)$ the
inner product $\langle u, u\spcheck\rangle$ is $2$.
It is clear from the definitions that for $u,v\in \Phi^{\rm
res}(\ell)$ we have
$\langle u, v\spcheck\rangle \in {\bf Z}$.
Thus it suffices to show that, for all $u, v\in \Phi^{\rm
res}(\ell)$, we have
$$r_{u}(v) = v - \langle v, u\spcheck\rangle u
\in \Phi^{\rm res}(\ell).$$
We fix $a\in \Phi$, resp.\ $b\in \Phi$, such that restriction of
$a$, resp.\ $b$ to ${\frak t}^\ell$ is $u$, resp.\ $v$.
Then
\begin{eqnarray*}
r_u(v) & = & v-\langle v,u\spcheck\rangle u\, = \, v-\langle
v,\epsilon(\ov a)\sum_{a'\in \ov a}(a')\spcheck\rangle u \\
& = & v-\langle b,\epsilon(a)\sum_{a'\in\ov a}(a')\spcheck\rangle u \\
& = & \left(b-\langle b,\epsilon(a)\sum_{a'\in\ov
a}(a')\spcheck\rangle a\right)|{\frak t}^\ell
\end{eqnarray*}

First assume that $\ov a$ is an
ordinary orbit. Since $a'|{\frak t}^\ell=a|{\frak t}^\ell$ for all
$a'\in \ov a$ and since $\epsilon(\ov a)=1$, we have
$$r_u(v)=\left(b-\sum_{a'\in\ov a}\langle b,(a')\spcheck\rangle
a'\right)|{\frak t}^\ell.$$ 
Suppose that $\ov a =\{a_1,\ldots, a_n\}$ where the $a_i$ are pairwise
distinct. Since the
$a_i\in\ov a$ are mutually orthogonal, this last equation can 
be rewritten as
$$r_u(v)=r_{a_1}\circ r_{a_2}\circ\cdots\circ r_{a_n}(b)|{\frak t}^\ell.$$
(Notice that the $r_{a_i}$ commute, so that the composition is independent
of the ordering.)
Clearly, then, $r_u(v)\in\Phi^{\rm res}(\ell)$.

In case $\ov a$  is an exceptional orbit ${\cal O}=\coprod_{i=1}^t{\cal
  O}_i$ with the ${\cal O}_i=\{a_{i,1},a_{i,2}\}$ being exceptional
pairs we have
$$r_{u}(v) =(b- 2\langle b, \sum_{i}
a_{i,1}\spcheck+a_{i,2}\spcheck\rangle a_{i,1})|\frak t^\ell= (b-
\langle b, \sum_ia_{i,1}\spcheck+a_{i,2}\spcheck\rangle 
(a_{i,1}+a_{i,2}))|\frak t^\ell,$$ since $a_{i,1}$ and $a_{i,2}$ have
the same restriction to $\ell$. In this case,  
$c_i=a_{i,1}+a_{i,2}$ is a root and $c_i\spcheck =
a_{i,1}\spcheck+a_{i,2}\spcheck$ since 
$a_{i,1}$, $a_{i,2}$, and $c_i$ have the same length. Thus 
$$r_{u}(v) = r_{c_1}\circ\cdots\circ r_{c_t}( b) |\frak
t^\ell.$$
(Here, the $c_i$ are orthogonal, so that the order of the reflections
is again irrevelant.)
This proves that $r_{u}(v) \in \Phi^{\rm res}(\ell)$ in this
case also. 

The walls in $\frak t^\ell$ defined by the elements of $\Phi^{\rm res}(\ell)$
are the same as the walls of $\Phi^\tau$, viewed as linear hyperplanes. Thus,
the Weyl groups are the same. Finally, if  $\Phi^\tau$ is irreducible, we
cannot divide the set of walls for $\Phi^\tau$ into two nonempty,
mutually orthogonal subsets.  Thus the same is true for $\Phi^{\rm
res}(\ell)$, and hence
$\Phi^{\rm res}(\ell)$ is irreducible.
\end{proof}

For $a\in \Phi$ such that $a|{\frak t}^\ell\neq 0$, the element
$\pi (a\spcheck)$ is given by the formula
$$\pi(a\spcheck) = \frac{1}{n_{\ov a}}\sum _{b\in \ov a}b\spcheck.$$
Define $\Phi^{\rm proj}(\ell)\spcheck\subseteq {\frak t}^\ell$ by 
$$\Phi^{\rm proj}(\ell)\spcheck = \{\, \pi(a\spcheck): a\in \Phi\,\} -\{0\}.$$
Clearly, $\pi(a\spcheck) = 0$ if and
only if $\ov a= 0$ as a linear form. 

We now define a second root system in $\frak t^\ell$ as follows. Its
coroots will be the set $\Phi^{\rm  proj}(\ell)\spcheck$.  Given $a\in
\Phi$ such that $\pi(a\spcheck) \neq 0$, we define the root
inverse to
$\pi(a\spcheck)\in \Phi^{\rm  proj}(\ell)\spcheck$  to be $\epsilon(\ov
a)n_{\ov a}a|{\frak t}^\ell$. As before, this is independent of the
choice of a lift of
$\pi(a\spcheck)$ to $a\in \Phi\spcheck$. Let $\Phi^{\rm
proj}(\ell)\subseteq {\frak t}^\ell$ be the set of all such elements.

Dually to the above results, we have:

\begin{proposition}\label{Wres=Wproj}
$\Phi^{\rm proj}(\ell)$
is a  possibly nonreduced root system in $\frak
t^\ell$. It is irreducible if $\Phi$ is irreducible.
The coroot inverse to $\epsilon(\ov a)n_{\ov
a}a|{\frak t}^\ell$ is $\pi(a\spcheck)$. The coroot lattice of $\Phi^{\rm
proj}(\ell)$ is $\pi(Q\spcheck)$.  The Weyl group of $\Phi^{\rm
proj}(\ell)$ is the same as that of
$\Phi^{\rm  res}(\ell)$ and hence as that of $\Phi^\tau$. 
\end{proposition}

We call $\Phi^{\rm
proj}(\ell)$ the {\sl projection  root system}.

The set of walls defined by  $\Phi^{\rm
proj}(\ell)$ is equal to the set of walls defined by  $\Phi^{\rm
res}(\ell)$, and thus the two systems have the same Weyl groups. In
general, however, there is no one-to-one correspondence  between 
$\Phi^{\rm proj}(\ell)$ and $\Phi^{\rm
res}(\ell)$. Of course, if $\tau ={\rm Id}$, then $\Phi^{\rm
proj}(\ell)=\Phi^{\rm
res}(\ell)=\Phi$. On the other hand, there are examples where one
of the systems is reduced and the other is non-reduced. Note however that
the set
$\widetilde\Delta/\ell$ injects into both 
$\Phi^{\rm res}(\ell)$ and $\Phi^{\rm proj}(\ell)$. Moreover, if $\Phi$ is
simply laced, we can say the following:

\begin{lemma}\label{resdualproj}
 Suppose that $\Phi$ is simply laced. Then, using the Weyl
invariant inner product to identify $\frak t$ and $\frak t^*$,
$\Phi^{\rm proj}(\ell)$ is the inverse system to $\Phi^{\rm res}(\ell)$.
\end{lemma} 
\begin{proof} In case $\Phi$ is simply laced, the inner product identifies
$a$ with $a\spcheck$ and $a|\frak t^\ell$ with $\pi (a\spcheck) \in \frak
t^\ell$. Thus the roots of  $\Phi^{\rm res}(\ell)$ are identified with the
coroots in $\Phi^{\rm proj}(\ell)$.
\end{proof}

\subsection{Generalized Cartan matrices for $\Phi^{\rm res}(\ell)$
and $\Phi^{\rm proj}(\ell)\spcheck$}\label{6.3}

In this section we suppose that $\Phi$ is irreducible and reduced.
We identify $\widetilde \Delta/\ell$  with its image
in ${\frak t}^\ell$ and hence for $a\in\widetilde \Delta$ we write
$\ov a$ both for an orbit in 
$\widetilde \Delta$ and for an element of $\Phi^{\rm res}(\ell)$. 
We denote by $\pi( a\spcheck)$ the element of $\Phi^{\rm
proj}(\ell)\spcheck$ corresponding to $a\spcheck$, but continue to denote
the corresponding orbit of $\widetilde \Delta$ by $\ov a$.

Since $\widetilde \Delta/\ell\subseteq {\frak t}^\ell$ are elements of
the root system $\Phi^{\rm res}(\ell)$, 
and since from Equation~\ref{rooteqn},  $n(\ov a,\ov b)\le 0$ for $\ov
a\not=\ov b$, $\ov a, \ov b\in \widetilde \Delta/\ell$,
the numbers $n(\ov a,\ov
b)$ form a generalized
Cartan matrix.
Two elements $\ov a$ and $\ov b$ are orthogonal if and only if their
orbits span orthogonal subspaces of ${\frak t}$. Since $\Phi$ is
irreducible, 
it follows that this generalized Cartan matrix is indecomposable.

Let $d =\dim \frak t^\ell$. Then the cardinalities of $\widetilde
\Delta/\ell$ and of  $\widetilde \Delta\spcheck/\ell$ are both $d+1$.
Since $\widetilde \Delta/\ell$, resp.\ $\widetilde
\Delta\spcheck/\ell$, spans $({\frak t}^\ell)^*$, resp.\ ${\frak
t}^\ell$, there is a single linear relation among its elements.  
We claim that the one relation has positive integral
coefficients. The root integers $h_a$ only depend on the orbit
$\ov a$. Define $h_{\ov a} =n_{\ov a}h_a$ for any choice of $a\in \ov
a$.  The relation
$\sum _{a\in
\widetilde
\Delta}h_aa=0$ leads to a relation
$$\sum _{\ov a\in\widetilde \Delta/\ell}h_{\ov a}\ov a= 0,$$
Similarly, the relation for the projection coroots is:
\begin{equation}\label{bargaeqn}
\sum_{\ov a}g_{\ov a}\pi(a\spcheck)= 0,
\end{equation}
where $g_{\ov a} = n_{\ov a}g_a$.
Thus, the generalized Cartan matrices determined by $n(\ov a,\ov b)$ and
by $n(\pi(a\spcheck),\pi(b\spcheck))$ are of affine type. It is not in 
general true that
$\widetilde \Delta/\ell$, resp.\ $\widetilde \Delta\spcheck/\ell$
is an extended set of simple roots resp.\ coroots for $\Phi^{\rm
res}(\ell)$ resp.\ $\Phi^{\rm proj}(\ell)$, cf.
Proposition~\ref{setofcoroots}.

 There are then a corresponding affine
diagrams, which we denote by 
$\widetilde D(\widetilde{\Delta}/\ell)$ and $\widetilde D(
\Delta\spcheck/\ell)$.
While these affine Dynkins diagrams will be different in general,
the associated Coxeter graphs will be the same. (Here the Coxeter 
graph of a generalized Dynkin diagram is obtained by keeping
bonds and their multiplicities but forgetting the arrows.) 

Our goal now will be to work out explicitly the Cartan integers for the
coroots $\Phi^{\rm proj}(\ell)\spcheck$ inverse to the projection
root system. We begin with a graph-theoretic lemma:

\begin{lemma} Let $D$ be a finite tree, and let $\ell$ be a 
group of automorphisms of $D$. If $v_1$ and $v_2$ are two vertices of $D$
which are connected by an edge, then either ${\rm Stab}(v_1) \subseteq
{\rm Stab}(v_2)$ or ${\rm Stab}(v_2) \subseteq {\rm Stab}(v_1)$.
\end{lemma}
\begin{proof} First we claim that there is fixed point for the action of
$\ell$ on the topological space $|D|$ associated to $D$. The proof is by
induction on the number of vertices. Clearly $\ell$ has a fixed point if
there are $1$ or $2$ vertices. Otherwise, let $D'\subset D$ be the
subgraph obtained by deleting the leaves. Then $D'$ is a nonempty
contractible proper subgraph on which $\ell$ acts, so by induction
there is a fixed point of the action of $\ell$ on $D'$ and hence on $D$.

Choose a   point $p$ fixed by $\ell$. If $p$ is an interior point of
an edge $e$ whose boundary is $\{v_1, v_2\}$, then it is easy to see
that ${\rm Stab}(v_1) = {\rm Stab}(v_2)$ is the set of $g\in \ell$
such that $g|e ={\rm Id}$. Assume that we are not in
this case. There is a unique path
$\Gamma$ in
$D$ joining $p$ to $v_1$.  Possibly after switching $v_1$ and $v_2$, we
can assume that
$v_2$ does not lie on this path. Hence the unique path $\Gamma '$ from $p$
to $v_2$ is the union of $\Gamma$ with the edge connecting $v_1$ and
$v_2$. If $g\in \ell$ fixes $v_2$, then  $g(\Gamma') = \Gamma'$. Since
$g(p)=p$, $g|\Gamma' ={\rm Id}$. Thus $g(v_1) =v_1$. It follows that  
${\rm Stab}(v_2) \subseteq {\rm Stab}(v_1)$.
\end{proof}

\begin{proposition}\label{Cartanints}
Suppose that $\Phi$ is not of type $A_n$.
Let $\pi(a\spcheck), \pi(b\spcheck) \in \Phi^{\rm proj}(\ell)\spcheck$,
and let $n(\pi(a\spcheck), \pi(b\spcheck))$ be the corresponding Cartan
integer. Then:
\begin{enumerate}
\item If every element of $\ov a$ is orthogonal to every element of $\ov
b$, then $n(\pi(a\spcheck), \pi(b\spcheck)) = 0$.
\item If there exist $a \in \ov a$ and $b \in
\ov b$ such that $a\spcheck$ and $b\spcheck$ are not orthogonal,
then either ${\rm Stab}(a\spcheck)\subseteq {\rm Stab}(b\spcheck)$ or 
${\rm Stab}(b\spcheck)\subseteq {\rm Stab}(a\spcheck)$.  If ${\rm
Stab}(a\spcheck)\subseteq {\rm Stab}(b\spcheck)$, then 
$$n(\pi(a\spcheck),\pi(b\spcheck))=\epsilon(\ov b)n(a\spcheck,b\spcheck).$$
\item If there exist $a \in \ov a$ and $b \in
\ov b$ such that $a\spcheck$ and $b\spcheck$ are not orthogonal
and 
${\rm Stab}(b\spcheck)\subseteq {\rm
Stab}(a\spcheck)$, then 
$$n(\pi(a\spcheck),\pi(b\spcheck))=\epsilon(\ov b)\frac{n_{\ov
    b}}{n_{\ov a}}n(a\spcheck,b\spcheck).$$
\end{enumerate}
\end{proposition}
\begin{proof} 
Let $a\spcheck,b\spcheck\in\widetilde \Delta\spcheck$ have the property
that $n(a\spcheck_i,b\spcheck_j)=0$ for all $a _i\in\ov a$ and
all $b _j\in\ov b$. The two subspaces of
${\frak t}$ spanned by the $a_i\spcheck$ such that $a_i\in \ov a$,
resp.\ the $b_j\spcheck$ such that $b_j\in \ov b$, are
orthogonal, and hence 
$\pi(a\spcheck)$ and
$\pi(b\spcheck)$ in
${\frak t}^{\ell}$ are also orthogonal. Thus
$n(\pi(a\spcheck),\pi(b\spcheck))=0$ and  there is no bond between
$\pi(a\spcheck)$ and $\pi(b\spcheck)$ in  the affine diagram associated
with the generalized Cartan matrix of these elements.

Suppose that there exist $a \in \ov a$ and $b \in
\ov b$ such that $n(a\spcheck,b\spcheck)\not=0$.
By the previous lemma, either ${\rm Stab}(a\spcheck)\subseteq {\rm
Stab}(b\spcheck)$ or 
${\rm Stab}(b\spcheck)\subseteq {\rm Stab}(a\spcheck)$. 
Since the root inverse to $\pi(b\spcheck)$ is $\epsilon(\ov
b)n_{\ov b}b|{\frak t}^\ell$ we have
\begin{eqnarray*}
n(\pi(a\spcheck),\pi(b\spcheck)) & = & \langle
\pi(a\spcheck),\epsilon(\ov b)n_{\ov    b} \ov b\rangle \\
& = & \epsilon(\ov b)\langle \pi(a\spcheck),\sum_{b'\in \ov b}b'\rangle
= \epsilon(\ov b) \langle a\spcheck,\sum_{b'\in \ov b}b'\rangle \\
& = & \epsilon(\ov b) \sum_{b'\in \ov b}\langle a\spcheck,b'\rangle =
\epsilon(\ov b)\sum_{b'\in\ov b}n(a\spcheck,(b')\spcheck).
\end{eqnarray*}

If ${\rm Stab}(a\spcheck)\subseteq {\rm Stab}(b\spcheck)$, we see that
$$\sum_{b'\in\ov
  b}n(a\spcheck,(b')\spcheck)=n(a\spcheck,b\spcheck),$$
and 
$n(\pi(a\spcheck),\pi(b\spcheck))=\epsilon(\ov b)n(a\spcheck,b\spcheck)$
in this case.
On the other hand if ${\rm Stab}(b\spcheck)\subseteq {\rm
Stab}(a\spcheck)$, then  
$$\sum_{b'\in\ov  b}n(a\spcheck,(b')\spcheck)=
\frac{n_{\ov b}}{n_{\ov a}}n(a\spcheck,b\spcheck),$$
and thus
$$n(\pi(a\spcheck),\pi(b\spcheck))=\epsilon(\ov b)\frac{n_{\ov
    b}}{n_{\ov a}}n(a\spcheck,b\spcheck).$$ 
\end{proof}

Similar results handle the case where $\Phi$ is of type $A_n$:

\begin{proposition} Suppose that $\Phi$ is of type $A_n$, and that $\frak
t^\ell \neq 0$. Then the Cartan integers are given by the same formula
as Proposition~\ref{Cartanints} except in the case where $\Phi$ is of
type $A_{2k-1}$, $\ell$ contains a rotation of order $k$ and an
involution fixing two vertices. In this case, the quotient coroot
diagram is of type $\widetilde A_1$.
\end{proposition}
\begin{proof} Since $\ell$ is dihedral, the stabilizer of an
element has either one or two elements. The only case not covered
by Proposition~\ref{Cartanints} is the case where there exist 
two non-orthogonal coroots $a\spcheck$ and $b\spcheck$ such that
${\rm Stab}(a\spcheck)$ and ${\rm Stab}(b\spcheck)$ are both nontrivial.
In this case, the product of the two nontrivial elements is a rotation
which either has order $n+1$, if $n$ is even, or $k=(n+1)/2$, if $n$ is
odd.  In the first case, $\frak t^\ell =\{0\}$, and in the second case
either 
$\frak t^\ell =\{0\}$ or $\Phi$ is of type
$A_{2k-1}$,
the rotation subgroup $\ell'$ of $\ell$ has order exactly $k$, and there
is an involution in $\ell$  fixing two vertices. In this case,
$\frak t^{\ell} =\frak t^{\ell'}$, and the quotient coroot diagram is of
type $\widetilde A_1$.
\end{proof} 

Note that the affine diagrams associated to the Cartan integers we have
calculated here agree with the diagrams given in
Definition~\ref{defofdiagram} in the introduction.

Next we relate the Weyl groups of these roots systems to
$W({\frak t}^\ell,G)$. 

\begin{proposition}\label{prop}
Let $\Phi$ be a reduced root system with $\tau$ and $\ell$ as above.
 Suppose  that, for every component $D$ of $\widetilde D(\Phi)$  which
is  of type $A_n$, the stabilizer $\ell_D$ in $\ell$ of $D$ is either
trivial or is not a cyclic  group of rotations of $D$. Then
the Weyl group 
$W(\frak t^\ell, G)$ is identified with the Weyl group of $\Phi^{\rm
res}(\ell)$ or equivalently with the Weyl group of
$\Phi^{\rm proj}(\ell)$.
\end{proposition}

\begin{proof} 
It suffices to consider the case when $\Phi$ is irreducible.
Clearly, the result holds if $\ell$ is the trivial group.
Thus, we assume that $\ell \neq {\rm Id}$. We have seen that we can
realize the elements of the Weyl group of $\Phi^{\rm res}(\ell)$ or of
$\Phi^{\rm proj}(\ell)$ as elements of the Weyl group $W(\Phi)$ 
normalizing
$\frak t^\ell$. Thus there is a homomorphism from the Weyl group of 
$\Phi^{\rm res}(\ell)$ to $W(\frak t^\ell, G)$, and since the Weyl group
of $\Phi^{\rm res}(\ell)$ acts faithfully on
$\frak t^\ell$, this homomorphism is injective. We must show that its 
image is all of
$W(\frak t^\ell, G)$.

Given $w\in W(\frak t^\ell, G)$, represent $w$ by an element of $W(\Phi)$
which normalizes $\frak t^\ell$. Since $w$ permutes the set $\{\ov a\}$ of
restricted roots and preserves the inner product on $\frak t^\ell$, it 
defines an automorphism of the root system $\Phi^{\rm proj}(\ell)$. We
claim that $\Phi^{\rm proj}(\ell)$ is either non-simply laced or
$A_1$. Assuming this, since every automorphism of a root system which
is either non-simply laced or
$A_1$ is given by a Weyl element, it follows that $w$ is given by an
element of the Weyl group of $\Phi^{\rm proj}(\ell)$, or equivalently
of $\Phi^{\rm res}(\ell)$.

Note that that $\Phi^{\rm proj}(\ell)$ is  non-simply laced if there is an
exceptional orbit. Thus we may assume that there are no exceptional
orbits. First assume that
$\Phi$ is not of $A_n$ type. In particular, the numbers $n_{\ov a}$
cannot all be equal since the diagram $\widetilde D(\Phi)$ is
contractible. In particular there must exist two coroots
$a\spcheck,b\spcheck \in \widetilde D\spcheck(\Phi)$ which are not
orthogonal and such that $n_{\ov a} \neq n_{\ov b}$. If
$\Phi$ is simply laced, it follows from Proposition~\ref{Cartanints} that
$n(\pi(a\spcheck), \pi(b\spcheck))\neq n(\pi(b\spcheck),
\pi(a\spcheck))$, and hence that
$\Phi^{\rm proj}(\ell)$ is not simply laced.  If $\Phi$ is not simply
laced,
$a\spcheck$ is a short coroot, and $b\spcheck$ is a long coroot, it is
easy to see that ${\rm Stab} (a\spcheck )\subseteq {\rm Stab}
(b\spcheck)$, and thus that $n_{\ov a}\geq n_{\ov b}$. Thus
$n(\pi(a\spcheck), \pi(b\spcheck)) = n(a\spcheck, b\spcheck)$, and
$$n(\pi(b\spcheck), \pi(a\spcheck))\geq n(b\spcheck, a\spcheck) >
n(a\spcheck, b\spcheck).$$
Thus $\Phi^{\rm proj}(\ell)$ is  non-simply laced in this case also.

Direct inspection then handles the case where $\Phi$ is of $A_n$ type and
$\ell$ is not cyclic.
\end{proof}

In the next section 
we will prove a related result (Proposition~\ref{prop}) which also covers
the remaining case when  $\Phi$ has a component of type $ A_n$ whose
stabilizer is a group of rotations.

\subsection{The case of an outer automorphism}

For future reference, we want to work out the results of
Subsections~\ref{Linear} and ~\ref{Digression} in the case where 
$\tau =\ell$. In this case, $\tau$ is induced from a group of diagram 
automorphisms of the Dynkin
diagram of $G$. Hence $\tau$ acts on the extended diagram, fixing the
extended root. The results described here are due, for the most part, to
deSiebenthal \cite{deS}.  Notice that if $\tau \neq {\rm Id}$ and
$\Phi$ is irreducible, there are very few possibilities:  $\Phi$
is simply laced and is of type $A_n$, $D_n$, $E_6$ if $\tau$ has
order $2$,  and is of  type
$D_4$ if $\tau$ has order $3$ or $6$.

\begin{lemma}\label{linearcase}
 In the above notation, assuming that $\tau$ is a group of linear
transformations of ${\frak t}$, 
\begin{enumerate}
\item For every $a\in \Phi$, $\ov a = a|\frak t^\tau$ is nonzero.
\item There is a one-to-one correspondence between $\Phi^{\rm res}(\tau)$ 
and the set of orbits $\Phi/\tau$.
\item The subset $\Delta/\tau$ of $\Phi^{\rm res}(\tau)$ is 
a set of simple roots.
\item The highest root for the irreducible factors of $\Phi^{\rm
    res}(\tau)$ corresponding to the above 
set of simple roots are the images  $\ov d$, where $d$ is the highest
root of an irreducible factor of $\Phi$. 
\item The root integer  for $\Phi^{\rm res}(\tau)$ corresponding to  $\ov
a$, for $a\in \widetilde\Delta$, is the integer
$n_{\ov a}h_a$.
\item $\Phi^{\rm res}(\tau)$ is  reduced if
and only if there are no 
exceptional orbits. In this case $\Phi^{\rm res}(\tau)=\Phi^\tau$. 
\item $\Phi^{\rm res}(\tau)$ is  not reduced if and only if $\Phi$ has
  an irreducible factor of type 
$A_{2k}$ and the stabilizer of this component in $\tau$ is non-trivial.
In this case, $\Phi^\tau$ is the
subsystem of $\Phi^{\rm res}(\tau)$ consisting of the roots $a$ such that
$2a$ is not a root. The set of indivisible roots in $\Phi^{\rm res}(\tau)$
is also a root system, and $\Delta/\tau$ is also a set of 
simple roots for this root system.
\end{enumerate}
\end{lemma}

\begin{proof} (1) follows since $\frak t^\tau$ contains a regular element. 

To
see (2), note that, if $a$ and $a'\in \Phi$ are such that $a|\frak t^\tau
= a'|\frak t^\tau$, then $W_{a}\cap \frak t^\tau = W_{a'}\cap \frak
t^\tau$. By Lemma~\ref{samewall}, this can only happen if $a$ and $\pm
a'$ lie in the same orbit or
$a$ lies in an exceptional orbit  and $\pm a' = a + \tau(a)$.
In this case $a'|\frak t^\tau = \pm 2a |\frak t^\tau$. Thus, if $a|\frak t^\tau
= a'|\frak t^\tau$ and $a$ and $a'$ do not lie in the same orbit, then
$ka|\frak t^\tau =0$ for some $k> 0$. This contradicts (1). Moreover, it
follows that $\Phi^{\rm res}(\tau)$ is reduced if and only if there are no
exceptional orbits. 

It suffices to prove (3), (4), (5), (6), and (7)  under the assumption that
$\Phi$ is irreducible.
To see (3), (4), and (5), note that every positive root $b$ can be written as
a positive integral linear combination $\sum_{a\in \Delta}r_aa$. Thus 
$\ov  b = \sum _{\ov a}r_an_{\ov a}\ov a$, and similarly for negative
roots. Since the cardinality of $\Delta/\tau$ is the dimension of $\frak
t^\tau$, it follows that the $\ov a$ are linearly independent and hence a set
of simple roots for $\Phi^{\rm res}(\tau)$. Taking $b=d$, we see that the
coefficients of $\ov b$, in terms of the simple roots, are at most those of
$\ov d$, and hence $\ov d$ is a highest root for $\Phi^{\rm res}(\tau)$.
Hence the root integers are as claimed. 
To see  (6), since $\ov d$ is the highest root for
$\Phi^{\rm res}(\tau)$, it follows that the alcove for $\Phi^\tau$ is the
alcove for $\Phi^{\rm res}(\tau)$. Thus if $\Phi^{\rm res}(\tau)$ is reduced
we must have $\Phi^{\rm res}(\tau)=\Phi^\tau$. If $\Phi^{\rm res}(\tau)$ is not
reduced, it is easy to see that the root system associated to the alcove is
exactly the set of roots $a$ such that $2a$ is not a root. The remaining
statement in (7) follows by a direct inspection.
\end{proof}

Next we determine the Weyl group and the coroot lattice. 

\begin{lemma}\label{Weylcent} If $\tau$ is a group of linear
transformations of
$\frak t$, 
\begin{enumerate}
\item
$W(\frak t^\tau,
G) = W(\Phi^{\rm proj}(\tau)) = W(\Phi^\tau) = Z_W(\tau)$, the subgroup of
elements of
$W$ which commute with $\tau$. 
\item The coroot lattice for $\Phi^\tau=\Phi^{\rm res}(\tau)$ is
$(Q\spcheck)^\tau= Q\spcheck \cap \frak t^\tau$.
\item The coroot lattice for $\Phi^{\rm proj}(\tau)$ is
$\pi(Q\spcheck)$.
\end{enumerate}
\end{lemma}
\begin{proof}
Clearly $Z_W(\tau)$ normalizes $\frak t^\tau$, and since $\frak t^\tau$
contains a regular element, the action of $Z_W(\tau)$ on $\frak t^\tau$
is faithful. Thus $Z_W(\tau) \subseteq W(\frak t^\tau,
G)$. Conversely, if $w\in W(\frak t^\tau,
G)$, choose a regular element $x\in \frak t^\tau$. Then, for all $g\in
\tau$, $w(x) = g(w(x)) = g(w(g^{-1}(x)))$. Thus since $x$ is regular $w =
g\circ w\circ g^{-1}$ for all $g\in
\tau$, so that $w\in Z_W(\tau)$.

To see (2), the coroot lattice for $\Phi^\tau=\Phi^{\rm res}(\tau)$ is
spanned by elements of the form $\epsilon(\ov a)\sum _{b\in \ov
a}b\spcheck$. If $\ov a $ is ordinary, this is just $\sum _{b\in \ov
a}b\spcheck$, and if $\ov a$ is exceptional, corresponding to the pair
$\prod _i\{a_{i,1},a_{i,2}\}$, then 
$a_{i,1}\spcheck+a_{i,2}\spcheck$ is again in the coroot lattice and the
corresponding orbit sum is $\sum _ia_{i,1}\spcheck+a_{i,2}\spcheck$. Thus
the coroot lattice is generated by orbit sums. Since
$\tau$  permutes an integral basis for
$Q\spcheck$, the coroot orbits generate $(Q\spcheck)^\tau$. 

Finally, it follows from the definition that the coroot lattice for 
$\Phi^{\rm proj}(\tau)$ is $\pi(Q\spcheck)$. 
\end{proof}

Using Proposition~\ref{Cartanints} and Lemma~\ref{resdualproj}, we see
that, in case $\Phi$ is irreducible and $\tau$ is nontrivial, $\Phi^{\rm
res}(\tau)$ is given as follows: 
\begin{itemize}
\item If $\Phi$ is of type $A_{2n-2}$, $n> 1$,
then $\Phi^{\rm res}(\tau)$ is of type $BC_n$.
\item If $\Phi$ is of type $A_{2n-1}$, $n> 1$,
then $\Phi^{\rm res}(\tau)$ is of type $C_n$.
\item If $\Phi$ is of type $D_{n+1}$, $n\geq 3$, and $\tau$ has order
$2$ then $\Phi^{\rm res}(\tau)$ is of type $B_n$.
\item If $\Phi$ is of type $D_4$, and $\tau$ has order $3$,
then $\Phi^{\rm res}(\tau)$ is of type $G_2$.
\item If $\Phi$ is of type $E_6$,
then $\Phi^{\rm res}(\tau)$ is of type $F_4$.
\end{itemize}

\section{The torus $\ov S^{w_{\cal C}}$ and the Weyl group
$W(S^{w_{\cal C}},G)$}\label{cmaxtorus}  

We apply the results of the previous section to describe  the 
quotient torus
$\ov S^{w_{\cal C}}$ and the Weyl group $W(S^{w_{\cal C}},G)$ in 
terms of the original root system $\Phi(G)$ on $T$ and the group
${\cal C}$.

\subsection{Further results under an additional hypothesis}

In this subsection, we keep the conventions of the previous section, so
that
$\Phi$ is reduced but not necessarily irreducible and we make one 
further assumption on the group $\tau$: 

\begin{assumption}\label{assumption} The fixed subspace ${\frak
  t}^\ell$ of the linearization $\ell$ of
$\tau$ is written as the intersection of the kernels of a subset of the 
roots of $\Phi$.
\end{assumption}

Let $\Phi^\perp$ be the subset of $\Phi$ consisting of roots vanishing
on ${\frak t}^\ell$. Let ${\frak u}$ be the subspace spanned by the
coroots inverse to the roots in 
$\Phi^\perp$. Then $\Phi^\perp$ is a root system on ${\frak u}$.
The above assumption on ${\frak t}^\ell$ implies that ${\frak u}\oplus
{\frak t}^\ell={\frak t}$.  Let $\pi\colon \frak t \to \frak t^\ell$
denote orthogonal projection.

\begin{lemma}
The intersection ${\frak t}^\tau\cap {\frak u}$ is the barycenter of
an alcove for $\Phi^\perp$.
\end{lemma}

\begin{proof}
Let $\tilde x_0={\frak t}^\tau\cap {\frak u}$.
We claim that no root of $\Phi^\perp$ is integral on $\tilde x_0$ and
hence $\tilde x_0$ is in the interior of an alcove $B\subset \frak u$
of the root system $\Phi^\perp$. In fact,   if $a\in\Phi^\perp$  is
integral  on
$\tilde x_0$, then $a$ is a root of $\Phi$ which vanishes on ${\frak
t}^\ell$ and hence is integral  on
$\{\tilde x_0\}+{\frak t}^\ell={\frak t}^\tau$. But this contradicts our
assumption that ${\frak t}^\tau$ contains an interior point of an alcove
for $\Phi$.
 The group $\tau$
normalizes ${\frak u}$ and the alcove structure for the root system
$\Phi^\perp$. The point $\tilde x_0$ is  the unique fixed point for
$\tau|{\frak u}$ and is in the interior  of $B$. Thus $\tilde x_0$ is the
barycenter of 
$B$.
\end{proof}

\begin{proposition}\label{Weylgroup}
The Weyl group $W({\frak t}^\ell,G)$ is identified with the Weyl group of
$\Phi^\tau$. 
\end{proposition}

\begin{proof}
Let $g\in Z_{W_{\rm aff}(\Phi)}(\tau)$.
Then its differential is an element of the Weyl group of $\Phi$ which
normalizes ${\frak t}^\ell$.
The elements of the Weyl group of $\Phi^\tau$ are exactly the
restrictions to ${\frak t}^\ell$ of the differentials of elements
$g\in Z_{W_{\rm aff}(\Phi)}(\tau)$.
This proves that the Weyl group of $\Phi^\tau$ is identified with a
subgroup of $W({\frak t}^\ell,G)$. 

Let $g\in N_{W(\Phi)}({\frak t}^\ell)$. 
Then $g$ normalizes ${\frak u}$ and the alcove decomposition of ${\frak
u}$ for $\Phi^\perp$. As such, setting $\tilde x_0={\frak
t}^\tau\cap {\frak u}$, the point $g\cdot\tilde x_0$ is the barycenter of
some alcove for this alcove 
decomposition. Thus, there
is an element $h$ of the affine Weyl group $W_{\rm aff}(\Phi^\perp)$ with
$h\cdot\tilde x_0=g\cdot \tilde x_0$. Then $h^{-1}g(\tilde x_0)=\tilde
x_0$. Let $w\in W(\Phi^\perp)$ be the Weyl part of $h$. Since $w|{\frak
t}^\ell={\rm Id}$, the element $w^{-1}g$ normalizes ${\frak t}^\ell$.
Hence,
$h^{-1}g$ normalizes $\{\tilde x_0\}+{\frak t}^\ell={\frak t}^\tau$, and
hence by Lemma~\ref{cent=norm} is an element of $Z_{W_{\rm
aff}(\Phi)}(\tau)$. Thus, the restriction to
${\frak t}^\ell$ of its differential is an element of the Weyl group of
$\Phi^\tau$. The element $h$, and consequently also its differential,
centralize
${\frak t}^\ell$.
Thus, the restriction to ${\frak t}^\ell$ of the differential of
$h^{-1}g$ agrees with that of $g$.
This proves that the Weyl group of $\Phi^\tau$ is all of $W({\frak
t}^\ell,G)$.
\end{proof}

\begin{remark} Unlike the case of a group of linear automorphisms, it is
not in general true that $W({\frak
t}^\ell,G)$ is equal to $Z_W(\ell)$, the set of $w\in W$ which commute
with $\ell$. However, there is a surjection from $Z_W(\ell)$ to $W({\frak
t}^\ell,G)$, and in fact $Z_W(\ell) \cong I\rtimes W({\frak t}^\ell,G)$
for an appropriate subgroup $I$ of $Z_W(\ell)$.
\end{remark}

Proposition~\ref{Weylgroup} gives an extension of
Proposition~\ref{prop} to the case of an arbitrary group $\tau$:

\begin{proposition}\label{extension} Let $\Phi$ be a reduced but not
necessarily irreducible root system on ${\frak t}$. Suppose that $\tau$
is a group of affine isometries normalizing an alcove of $\Phi$. Then
the Weyl group of $W({\frak t}^\ell)$ is
identified with the Weyl group of $\Phi^{\rm res}(\ell)$ or
equivalently with the Weyl group of $\Phi^{\rm proj}(\ell).$
\end{proposition}

\begin{proof}
Clearly, we may assume that $\Phi$ is irreducible.
The only irreducible case not covered by Proposition~\ref{prop} is
when $\Phi$ is of type $A_n$ and $\ell$ is a group of rotations of the
extended Dynkin diagram.
We choose a representation of $\Phi$ on the subspace
$\{(x_1,\ldots,x_{n+1})|\sum x_i=0\}$ in ${\bf R}^{n+1}$ such that the
set of simple roots is $\{e_i-e_{i+1}\}$.
Suppose that $\ell$ is a group of order $m$ and let $k=(n+1)/m$.
Then ${\frak t}^\ell$ is the intersection of the kernels of the roots
$e_i-e_{i+k}$. Hence, this case is covered by
Proposition~\ref{Weylgroup}. 
\end{proof}

Next we describe the coroot lattice of $\Phi^\tau$.

\begin{proposition}\label{lattice}
The coroot lattice of the affine Weyl group associated to $\Phi^\tau$
is identified with $\pi(Q\spcheck(\Phi))$.
\end{proposition}

\begin{proof}
We decompose ${\frak t}= {\frak u}\oplus {\frak t}^\ell$ and let
$\tilde x_0={\frak u}\cap {\frak t}^\tau$.
Suppose that $g\in Z_{W_{\rm aff}(\Phi)}(\tau)$ and that $g|{\frak
t}^\tau$ is translation by an element $\gamma_2\in {\frak t}^\ell$.
We write
$$g(a_1,a_2)=(wa_1+\gamma_1,b+\gamma_2)$$
for $w$ in the Weyl group of $\Phi^\perp$. Since $g\in W_{\rm aff}(\Phi)$,
 $\gamma_1+\gamma_2\in Q\spcheck(\Phi)$, and hence  $\gamma_2\in
\pi(Q\spcheck(\Phi))$.

Conversely, if $\gamma_2\in \pi(Q\spcheck(\Phi)$,  write 
$\gamma=\gamma_1+\gamma_2$
for some $\gamma\in Q\spcheck(\Phi)$, with $\gamma_1\in
{\frak u}$.
Clearly, $\gamma_1$ is contained in the center of the root system
$\Phi^\perp$, and hence $\tilde x_0+\gamma_1$ is the barycenter of an
alcove for $\Phi^\perp$. This means that there is an element
$h\in W_{\rm aff}(\Phi^\perp)$ such that
$h(\tilde x_0+\gamma_1) = \tilde x_0$. The composition of translation by
$\gamma$ followed by $h$ is an element of $Z_{W_{\rm aff}(\Phi)}(\tau)$
whose restriction to ${\frak t}^\tau$ is translation by $\gamma_2$.
\end{proof}

\subsection{The case of a  subgroup of ${\cal C}G$}

We return now to the case where $\Phi=\Phi(T,G)$ is the root system
associated to a simply connected  group $G$.

\begin{lemma}\label{assumptionOK}
Let ${\cal C}\subseteq {\cal C}G$ be a subgroup, and let
$\tau$ be the group of affine automorphisms of ${\frak t}$
normalizing
$A$ determined by the action of ${\cal C}$ on the alcove. Then
$\tau$ satisfies the hypothesis of Assumption~\ref{assumption}.
\end{lemma}
\begin{proof} By Proposition~\ref{torus}, for each $c\in {\cal
C}$, the fixed point subspace
${\frak t}^{w_c}$ is conjugate under the Weyl group to $\frak
t_c$. Since $\frak
t_c$ satisfies the hypothesis of Assumption~\ref{assumption}, so
does ${\frak t}^{w_c}$. Since $w_{\cal C}$ is the linearization of
$\tau$, and ${\frak t}^{w_{\cal C}} = \bigcap _{c\in {\cal
C}}{\frak t}^{w_c}$, the lemma follows.
\end{proof}

\begin{proposition}\label{Cproj=Ctau}
The root system $\Phi^{\tau }$ is identified with the subroot system of 
$\Phi^{\rm proj}(w_{\cal C})$ consisting of all roots whose inverse
coroots are  indivisible elements of $\Phi^{\rm
proj}(w_{\cal C})\spcheck$.
\end{proposition}

\begin{proof}
By Proposition~\ref{lattice} and  Proposition~\ref{Wres=Wproj}, 
$\Phi^{\tau }$ and $\Phi^{\rm proj}(w_{\cal C})$ have the same
coroot lattices. By Proposition~\ref{Wres=Wproj}, their  Weyl
groups are the same. Now there is the following general lemma on
two root systems with the same coroot lattice:

\begin{lemma}\label{2root} Let $\Phi_1$ and $\Phi_2$ be two  root
systems on a vector space $V$, and suppose that $\Phi_1$ is reduced.
Suppose that the coroot lattice
$Q\spcheck(\Phi_1)$ is equal to the coroot lattice $Q\spcheck(\Phi_2)$
and that $W(\Phi_1) = W(\Phi_2)$. Then $\Phi_1$ is the set of all
non-multipliable roots in $\Phi_2$. In particular, if $\Phi_2$ is
reduced, then $\Phi_1 = \Phi_2$.
\end{lemma}
\begin{proof} Let $\Phi_2'$ be the subroot system of $\Phi_2$ consisting
of the non-multipliable roots. Then $(\Phi_2')\spcheck$ is the sub-coroot
system of $\Phi_2\spcheck$ consisting of all indivisible coroots in
$\Phi_2\spcheck$. 
In particular, the coroot lattice of $\Phi_2'$ is equal to that of
$\Phi_2$.
Of course, $\Phi_2$ and $\Phi_2'$ have the same Weyl group. Thus, it
suffices to assume  that $\Phi_2$ is reduced.

Let $Q\spcheck = Q\spcheck(\Phi_1) = Q\spcheck(\Phi_2)$.
For each $a\in \Phi_1$, $a\spcheck$ is an indivisible element of
$Q\spcheck$. Since reflection in the wall defined by $a$ is an element of
$W(\Phi_1) = W(\Phi_2)$, there exists a $b\in \Phi_2$ such that
$a\spcheck = rb\spcheck$ for some real number $r$. But $b\spcheck \in
Q\spcheck$ is also indivisible, so that $r=\pm 1$. Thus $\Phi_1\spcheck
\subseteq \Phi_2\spcheck$. By symmetry $\Phi_2\spcheck
\subseteq \Phi_1\spcheck$, and hence $\Phi_1\spcheck
= \Phi_2\spcheck$.
\end{proof}

Applying the lemma to $\Phi_1 = \Phi^{\tau }$ and $\Phi_2=\Phi^{\rm
proj}(w_{\cal C})$  completes the proof of
Proposition~\ref{Cproj=Ctau}.
\end{proof}

\begin{proposition}\label{setofcoroots} Suppose that $\Phi$ is
irreducible.  Orthogonal projection induces an embedding
$\widetilde
\Delta\spcheck_{\cal C}$ in $\frak t^{w_{\cal C}}$.  Its image is
a set of coroots inverse to an extended set of simple roots either
for
$\Phi^{\rm proj}(w_{\cal C})$ or for the subroot system
consisting of all non-multipliable roots of $\Phi^{\rm
proj}(w_{\cal C})$.
\end{proposition}

\begin{proof}
We have an identification of $\Phi^\tau$ with the non-multipliable
roots of $\Phi^{\rm proj}(w_{\cal C})$. We also know by 
Proposition~\ref{multcoroots} that, up to positive multiples, the
image under orthogonal projection of
$\widetilde \Delta\spcheck_{\cal C}$  forms the set of coroots
inverse to an extended set of simple coroots for $\Phi^\tau$. If
$\Phi^{\rm proj}(w_{\cal C})$ is reduced, then it is equal to
$\Phi^\tau$ and its coroots are indivisible elements in the
coroot lattice. Thus, as before, the multiples are all $+1$ and
the image of $\widetilde
\Delta\spcheck_{\cal C}$  under orthogonal projection forms the
set of coroots inverse to an extended set of simple roots for
$\Phi^{\rm proj}(w_{\cal C})$.  

If $\Phi^{\rm proj}(w_{\cal C})$ is not reduced, then it is of type $BC_n$
for some $n\ge 1$ and $\Phi^\tau$ is the subsystem of type $C_n$.
In the extended set of simple coroots for $\Phi^\tau$ all but two
are neither divisible nor multipliable in $\Phi^{\rm proj}(w_{\cal C})$
and the last two are multipliable by $2$ in $\Phi^{\rm
proj}(w_{\cal C})$. Thus, the image of
$\widetilde \Delta\spcheck_{\cal C}$ contains all the former
coroots and, for each of the latter two, contains either the coroot in
$\Phi^\tau$ or twice it. A priori there are three possibilities: (i)
$\widetilde
\Delta\spcheck_{\cal C}$ is equal to the extended set of coroots
for
$\Phi^\tau$ (i.e. consists of indivisible coroots); (ii) $\widetilde
\Delta\spcheck_{\cal C}$ is equal to the 
extended set of coroots for $\Phi^{\rm proj}(w_{\cal C})$ (i.e., contains
one indivisible but multipliable coroot and one non-multipliable but
divisible coroot); or (iii) contains two non-multipliable but
divisible coroots. In case (iii) the lattice spanned by $\widetilde
\Delta\spcheck_{\cal C}$ is of index two in the coroot lattice of
$\Phi^{\rm proj}(w_{\cal C})$ and hence this case is ruled out by
Proposition~\ref{lattice}. Since by Proposition~\ref{Cproj=Ctau}
$\Phi^\tau$ is the subroot system of $\Phi^{\rm proj}(w_{\cal C})$ consisting
of the non-multipliable roots, cases (i) and (ii) are exactly the two
cases listed in the statement of the proposition.
\end{proof}

\begin{defn}\label{defnof}
We define $\Phi(w_{\cal C})$ to be the subroot system of $\Phi^{\rm
proj}(w_{\cal C})$ such that the image of $\widetilde
\Delta\spcheck_{\cal C}$ in $\frak t^{w_{\cal C}}$ is an extended
set of simple coroots. By the previous proposition, $\Phi(w_{\cal
C})$ is either
$\Phi^{\rm proj}(w_{\cal C})$ or  the subroot system consisting
of all non-multipliable roots of $\Phi^{\rm proj}(w_{\cal C})$. 
\end{defn}

\begin{corollary}\label{bargaremark}
There is a positive integer $n_0$ such that, for all $\ov a\in
\widetilde \Delta_{\cal C}$, 
$g_{\ov a}$ is equal to $n_0$ times the coroot integer $m_{\ov
a}$ for the extended simple root of $\Phi(w_{\cal C})$
corresponding to
$\ov a$.
\end{corollary}
\begin{proof} By Equation~\ref{bargaeqn} in Section~\ref{6.3}, 
$\sum _{\ov a\in \widetilde \Delta_{\cal C}}g_{\ov
a}\pi(a\spcheck) = 0$. On the other hand, by definition $\sum
_{\ov a\in 
\widetilde \Delta_{\cal C}}m_{\ov a}\pi(a\spcheck) = 0$, 
and the $m_{\ov a}$ are relatively prime integers. Thus the
corollary is clear.
\end{proof}

\subsection{Proof of
Theorem~\ref{diagram1}}\label{exceptional} 

Let ${\cal C}$ be a subgroup of ${\cal C}G$ and let $\tau$ be the
corresponding group of affine isometries. 
 By
Lemma~\ref{assumptionOK}, $\tau$ 
satisfies Assumption~\ref{assumption}. Thus, the
results of this section apply to $\tau$. By
Proposition~\ref{setofcoroots},  
$\pi$ embeds
$\widetilde
\Delta\spcheck_{\cal C}$ as an extended set of simple coroots for
either the root system $\Phi^{\rm proj}(w_{\cal C})$ or for the 
subroot system of nonmultipliable roots in $\Phi^{\rm
proj}(w_{\cal C})$. In either case, we let $\Phi(w_{\cal C})$ be the corresponding
root system. The Weyl group of
$\Phi(w_{\cal C})$ is the same as the Weyl group of
$\Phi^{\rm proj}(w_{\cal C})$, and by Proposition~\ref{Weylgroup}, 
this is 
$W(S^{w_{\cal C}}, G)$. The coroot lattice of $\Phi(w_{\cal C})$ is the  lattice
generated by the indivisible coroots in  $\Phi^{\rm proj}(w_{\cal C})\spcheck$,
which is is the coroot lattice of $\Phi^{\rm proj}(w_{\cal C})$. By
Proposition~\ref{Wres=Wproj}, the coroot lattice of $\Phi^{\rm
proj}(w_{\cal C})$ is
$\pi(Q\spcheck)$. By Proposition~\ref{Cartanints}, the Cartan matrix
associated to $\widetilde
\Delta\spcheck_{\cal C}$ agrees with the Cartan matrix
associated to $\widetilde D\spcheck_{\cal C}$ as given in
Definition~\ref{defofdiagram}.  This completes the proof of
Theorem~\ref{diagram1}.

\section{The fixed subgroup of an automorphism}

Our main goal in this section is to describe the centralizer of a
$c$-pair in $G$, and in particular its Lie algebra and group of components.
In order to do so, it will be convenient to study more generally the fixed
subgroups of certain automorphisms of compact, connected groups. In
particular, we shall see  how the fixed subgroup changes as we vary the
given automorphism by composing with an inner automorphism. The fixed
subgroup is completely described by its Lie algebra, its component group,
and the  fundamental group of the identity component. We shall give a
complete description of the Lie algebra and component group of the fixed
subgroup. It is also straightforward to describe the fundamental group of its
identity component, or of the derived subgroup of the identity component, but
we shall not do this here. Many of the results of the first part
of this section, stated in a somewhar different language, can be
found in \cite{deS}.

In this section
$H$ will always denote  a compact, connected  group with maximal torus
$T$ and Lie algebra $\frak h$. For $h\in H$ we denote by $i_h$ the
inner automorphism of $H$ given by conjugating by $h$. Let $\Phi_H$ be the
set of roots of
$H$ with respect to
$T$. We denote by ${\frak t}$ the Lie algebra of $T$ and by ${\frak d}$
the subspace spanned by the coroots inverse to the roots in
$\Phi_H$. Then $\Phi_H$ is a root system on ${\frak d}$.
Furthermore, there is a direct sum decomposition of ${\frak t}$ into 
${\frak d}\oplus {\frak z}$ where ${\frak z}$ is the center of ${\frak
h}$. Let $\Delta_H$ be a set of simple roots for  $\Phi_H$, and let
$A\subseteq {\frak d}$ be the alcove determined by $\Delta_H$.
We denote by $\Delta_H\spcheck\subseteq {\frak t}$ the coroots
inverse to the roots $\Delta_H$.
Let $\Lambda\subseteq {\frak t}$ be the fundamental group of $T$ and let
$Q\spcheck_H\subseteq \Lambda$ be the coroot lattice of $\Phi_H$,
i.e., the lattice with basis $\Delta_H\spcheck$.
The fundamental group of $H$ is $\Lambda/Q\spcheck_H$.

The following lemma is straightforward and its proof is omitted.

\begin{lemma}\label{basicaut}
Suppose that  $\sigma\colon H\to H$ is an automorphism. Then
there exists an inner automorphism
$i_h$ such that, setting $\sigma'=i_h\circ \sigma$, we have:
\begin{enumerate} 
\item The fixed subgroup of $\sigma'$ is conjugate to the fixed subgroup of
$\sigma$.
\item The automorphism $\sigma'$ normalizes $T$ and the
automorphism induced by $\sigma$  on ${\frak d}$ normalizes $A$.
In particular, ${\frak d}^\sigma$ contains an interior point of $A$.
Such a point exponentiates to a regular element of $T$.
\item If
$\sigma''$ is another automorphism normalizing $T$ and $A$ which has the same
image in the outer automorphism group of $H$ as $\sigma$, then   there is an
inner automorphism
$i_t$ with $t\in T$ such that
$\sigma'' = i_t\circ
\sigma'$. 
\item For $t\in T$, the fixed subgroup of $\sigma''= 
i_t\circ
\sigma'$ only depends on the image of $t$ in $T/({\rm Id}-\sigma)(T)$.
\end{enumerate}
\end{lemma}

\begin{defn} An automorphism  $\sigma$ which normalizes $T$ and the alcove
$A$ is said to be {\sl in normal form\/} with respect to $T$ and $A$.
\end{defn} 

By the previous lemma, it suffices to study automorphisms in
normal form.

Let $\sigma$ be an automorphism in normal form  such that the
restriction of $\sigma$  to the center of $H$ has finite order.
We shall always make this assumption in what follows.
Since $\sigma$ normalizes $A$, it normalizes  $\Delta_H$ and
$\Delta_H\spcheck$. 
It follows that the restriction of $\sigma$ to ${\frak t}$ preserves
the decomposition ${\frak t}={\frak d}\oplus {\frak z}$ and that this
restriction is of
finite order. Hence, we can fix an inner product on ${\frak t}$ which
is invariant both under the Weyl group of $H$ and under $\sigma$ and such
that
${\frak d}$ and ${\frak z}$ are orthogonal.

We let $\Delta_H/\langle\sigma\rangle$, resp.\ 
$\Delta_H\spcheck/\langle\sigma\rangle$, denote the set of 
$\sigma$-orbits of $\Delta_H$, resp.\ $\Delta_H\spcheck$.  
The orbit of $a\in \Delta_H$ is denoted $\ov a$. 
All $a_i$ in $\ov a$ have the same restriction to ${\frak t}^\sigma$.
Thus, we can view $\ov a$ as an element of $({\frak t}^\sigma)^*$.
For each $\ov a\in \Delta_H/\langle\sigma\rangle $ we denote by  $n_{\ov
a}$ the cardinality of the orbit $\ov a\subseteq \Delta_H$.

\subsection{Some general facts about group cohomology}

If $J$ is an abelian group, let ${\rm Tor}\, J$ be the torsion 
subgroup of $J$.  Let $\sigma$ be an automorphism of  $J$. Denote by 
$J^\sigma$ the subgroup of invariants
$$J^\sigma = \{\, a\in J: \sigma (a) =a\,\},$$
and by $J_\sigma$ be the group of coinvariants
$$J_\sigma = J/{\rm Im}({\rm Id} -\sigma)J.$$
It is the largest quotient of $J$ on which $\sigma$ acts trivially.
Suppose that $J=\frak t$ is a vector space with a $\sigma$-invariant 
inner product, and that $\Lambda$ is a   sublattice of $\frak t$
such that
$\sigma(\Lambda) = \Lambda$. The  inner product on
$\frak t$ identifies
$t_\sigma$ with
$t^\sigma$ via orthogonal projection $\pi\colon \frak t \to \frak
t^\sigma$. The projection $\pi$ induces a map from $\Lambda$ to
$t^\sigma$ and hence factors through $\Lambda_\sigma$. By
comparing ranks, one checks that this sets up an isomorphism from
$\Lambda_\sigma/{\rm Tor}(\Lambda_\sigma)$ to $\pi(\Lambda)$. In
particular, suppose that $\sigma$ acts as a permutation 
representation on a
${\bf Z}$-basis of $\Lambda$. In this case ${\rm
Tor}(\Lambda_\sigma)  =0$ and $\pi$ induces an isomorphism from
$\Lambda_\sigma$ to
$\pi(\Lambda)$.

Suppose that we have an exact sequence of groups
$$0\to J'\to J \to J''\to 0,$$
such that $\sigma$ acts on all of the groups $J',J,J''$ and the homomorphisms
are $\sigma$-equivariant. Then there is an associated long exact sequence
$$0 \to (J')^\sigma \to J^\sigma \to (J'')^\sigma \to (J')_\sigma \to
J_\sigma \to (J'')_\sigma \to 0.$$
Here the homomorphism $(J'')^\sigma \to (J')_\sigma$ is defined as
follows: given $\xi\in (J'')^\sigma$, lift $\alpha$ to an element
$\tilde \xi \in J$. Since $\sigma(\tilde \xi) -\xi$ projects to $0$ in
$J''$, $\sigma(\tilde \xi) -\xi \in J'$, and its image in $(J')_\sigma$
is independent of the choice of the lift of $\xi$. This defines the
connecting homomorphism. We note that, if $J$ and $J''$ are not assumed
to be abelian but $J'$ is central in $J$, then there is still an exact
sequence
$$0 \to (J')^\sigma \to J^\sigma \to (J'')^\sigma \to (J')_\sigma \to
H^1(J, \langle\sigma \rangle)\to H^1(J'', \langle\sigma \rangle),$$
where $H^1(J, \langle\sigma \rangle)$ and $H^1(J'', \langle\sigma
\rangle)$ are no longer groups in general but only pointed sets, and
the connecting homomorphism is similarly defined.  

As a first application, we have the following:

\begin{lemma}\label{componentgroup} Suppose that $\sigma$ is  an
automorphism of
$H$ in normal form. Then the  component group
$\pi_0(T^\sigma)$ is naturally isomorphic to ${\rm
Tor}(\Lambda_\sigma)$. In particular, if $H$ is simply connected,
then $T^\sigma$ is connected. The fundamental group of the identity
component $(T^\sigma)^0=S^\sigma$ is
$\Lambda^\sigma$. Finally, the coinvariant torus $T_\sigma$ is
connected, and in fact $T_\sigma =\frak t^\sigma/\pi(\Lambda)$, where $\pi\colon \frak
t \to {\frak t}^{\sigma}$ is orthogonal projection.
\end{lemma}
\begin{proof} Consider the exact sequence
$$0 \to \Lambda \to \frak t \to T \to 0.$$
Taking the associated long exact cohomology sequence gives
$$0 \to \Lambda^\sigma \to \frak t^\sigma \to T^\sigma \to
\Lambda_\sigma \to \frak t_\sigma \to T_\sigma \to 0.$$
By the remarks above, $\frak t_\sigma \cong \frak t^\sigma$ and under
this identification the image of $\Lambda_\sigma$ is identified with
$\pi(\Lambda)$. Thus the kernel of the map $\Lambda_\sigma \to
\frak t_\sigma$ is exactly the torsion subgroup. Since $\frak
t^\sigma/\Lambda^\sigma$ is connected, it follows that 
$\pi_0(T^\sigma)\cong{\rm Tor}(\Lambda_\sigma)$. If $H$ is simply
connected, then $\Lambda =Q_H\spcheck$ and $\sigma$ acts as a 
permutation of the set of simple roots $\{a\spcheck\}$. Thus ${\rm
Tor}(\Lambda_\sigma)= 0$. The remaining statements are clear.
\end{proof}

\subsection{The root system $\Phi_H^{\rm proj}(\sigma)$ on ${\frak
t}^{\sigma}$} 

Let  $\sigma$ be an automorphism of $H$ in normal form. Orthogonal
projection
$\pi\colon {\frak t}\to{\frak t}^{\sigma}$ induces a map 
$\Lambda_\sigma\to  {\frak  t}^\sigma$.
Since $\Delta_H\spcheck$ is a $\sigma$-invariant basis of ${\frak d}$
it follows that $\Delta_H\spcheck/\langle\sigma\rangle$ is a basis 
for
${\frak
  d}^\sigma$. Since $\sigma$ permutes a basis for $Q\spcheck_H$, it 
follows that
$(Q\spcheck_H)_\sigma$ is torsion free and that $\pi$ induces an
isomorphism from
$(Q\spcheck_H)_\sigma$  to  a lattice of
maximal rank in ${\frak d}^\sigma$.

The automorphism $\sigma$
generates a cyclic subgroup of linear automorphisms of $\frak d$
which satisfies the hypotheses of Section~\ref{Linear} and 
contains no rotation in case $\Phi_H$ has a factor of type $A_n$. There are
the corresponding resticted and projection root systems, which we shall just
denote in this context by
$\Phi_H^{\rm res}(\sigma)$ and
$\Phi_H^{\rm proj}(\sigma)$. In particular,
$\Delta_H\spcheck/\langle\sigma\rangle$ is a set of simple coroots
for
$\Phi_H^{\rm proj}(\sigma)\spcheck$, and  the coroot lattice for
$\Phi_H^{\rm proj}(\sigma)$ is the lattice
$(Q\spcheck_H)_{\sigma}$ spanned by $\{\pi( a\spcheck)\}_{\ov a\in
\Delta_H/\langle\sigma\rangle}$.
The root inverse to $\pi (a\spcheck)$ is $\epsilon(\ov a)n_{\ov
a}\ov a|{\frak
  t}^{\sigma}$.  Let  
$B\subseteq {\frak d}^{\sigma}$ be the alcove determined by this
set of simple roots.
 We 
denote by ${\cal C}\Phi_H^{\rm proj}(\sigma)$ the quotient of the 
dual of  the root lattice of $\Phi_H^{\rm proj}(\sigma)$ by the
coroot lattice
$(Q\spcheck_H)_{\sigma}$. It is identified with a finite subgroup of
${\frak d}^{\sigma}/(Q\spcheck_H)_{\sigma}$. 

Applying Subsection~\ref{action1}, there is an affine action of
${\cal C}\Phi_H^{\rm
  proj}(\sigma)$ 
on ${\frak d}^{\sigma}$ normalizing $B\subseteq {\frak
  d}^{\sigma}$. We extend this action to one on all of ${\frak
  t}^\sigma$ by letting ${\cal C}\Phi_H^{\rm
  proj}(\sigma)$  act trivially  on ${\frak
  z}^\sigma$.

\subsection{The action of $\pi_0(H^\sigma)$ on ${\frak t}^\sigma$} 

We continue to assume that $\sigma$ is an automorphism of $H$ in
normal form. Let  $S^{\sigma}=(T^\sigma)^0$ be the image under the
exponential mapping of ${\frak t}^{\sigma}$.

\begin{lemma}\label{4321}
The torus $S^{\sigma}$ is a maximal torus for
$H^\sigma$.  
The set of roots of 
$H^\sigma$ with respect to $S^{\sigma}$ is a sub-root system of
$\Phi_H^{\rm res}(\sigma)$.
\end{lemma}

\begin{proof}
First we claim that, if $X\in {\rm Lie}(H^\sigma)\otimes_{\bf R}{\bf
C}$ is fixed by the action of $S^\sigma$, then $X\in \frak
t^\sigma$. Since
${\rm Lie}(H^\sigma) \subseteq \frak h$, it clearly suffices to show
that, if $X\in \frak h\otimes_{\bf R}{\bf C}$ is fixed by $S^\sigma$, then
$X\in
\frak t\otimes_{\bf R}{\bf C}$. The action of
$S^{\sigma}$ on
${\frak h}\otimes_{\bf R}{\bf C}$ normalizes the root spaces of
$\frak h$. The action of
$S^{\sigma}$ on a root space ${\frak h}^a$ is given by the character
$a|S^{\sigma}$. By Lemma~\ref{basicaut}, $S^\sigma$ contains a 
regular element of $H$. Thus, none of these characters is trivial.
Hence, if
$X$ is fixed 
by $S^\sigma$, then $X\in \frak t\otimes_{\bf R}{\bf C}$.
It follows that the fixed subspace of the action of $S^\sigma$ on
${\frak h}$ is ${\frak t}$, and hence that the fixed subspace of the 
action of $S^\sigma$ on
${\frak h}^\sigma$ is ${\frak t}^\sigma$.
Thus $S^{\sigma}$ is a maximal torus of $H^\sigma$. Also,
the roots for the $S^{\sigma}$-action on ${\rm
  Lie}(H^\sigma)\otimes_{\bf 
  R}{\bf C}$ are given by restrictions of roots of $\Phi_H$ to 
$S^{\sigma}$ and so are elements of $\Phi_H^{\rm
res}(\sigma)$. Hence the roots for $H^\sigma$ with respect to
$S^{\sigma}$ form a subset of $\Phi_H^{\rm res}(\sigma)$. To see that
they are a sub-root system, in other words that the corresponding 
coroots are the same, it is enough to show that the inner product on
$\frak d^\sigma$ induced by the Weyl invariant inner product on $\frak t$
is invariant under the Weyl group of $(H^\sigma)^0$. Since $S^\sigma$
contains a regular element of $T$, every element in
the Weyl group of $(H^\sigma)^0$ is the restriction of a Weyl element
of $H$, and thus the last statement is clear.
\end{proof}

There is a natural map from $(N_H(T))^\sigma$ to the Weyl group 
$W(T,H)$. The kernel of this map is $(N_H(T))^\sigma\cap T=
T^\sigma$. We denote  its image by
$\ov N(\sigma)$. This group is contained in the subgroup
$Z_W(\sigma)$ of
$W(T,H)=W$ 
consisting of elements commuting with
$\sigma$. Note further that $N_{(H^\sigma)^0}(S^\sigma)\cap T = S^\sigma$,
since conjugation by $T$ acts trivially on $S^\sigma$ and the Weyl group
$W(S^\sigma, (H^\sigma)^0)$ acts faithfully on $S^\sigma$. Likewise,
$\ov N(\sigma)$ acts faithfully on $S^\sigma$ since $S^\sigma$
contains a regular element.

\begin{lemma}\label{signormal}
 Every component of $H^\sigma$ contains an
element normalizing ${\frak t}^\sigma$ . If $g$ and $g'$ are in the
same component of $H^\sigma$ and normalize $\frak t^\sigma$, then
there exists a $w\in W(S^\sigma, (H^\sigma)^0)$ such that $g'|\frak
t^\sigma = w\circ (g|\frak t^\sigma)$.
\end{lemma}

\begin{proof}
Let $g\in H^\sigma$. Then $gS^{\sigma}g^{-1}$ is a maximal torus of
$H^\sigma$ and 
hence there is $h\in (H^\sigma)^0$ such that
$(gh)S^{\sigma}(gh)^{-1}=S^{\sigma}$.  
Clearly, $gh\in N_H(S^{\sigma})$ and $gh$ and $g$ lie in the same 
component of $H^\sigma$.
If $g$ and $g'$ are in the same component of $H^\sigma$ and both
normalize ${\frak t}^\sigma$, then $g'g^{-1}\in (H^\sigma)^0$
and normalizes ${\frak t}^\sigma$. Thus $g'g^{-1}|\frak t^\sigma =
w$ for some $w\in W(S^\sigma, (H^\sigma)^0)$.
\end{proof}

\begin{lemma}\label{compgrp}  The map $N_H(T)^\sigma\to 
\pi_0(H^\sigma)$ induced   by the inclusion $N_H(T)^\sigma\subseteq
H^\sigma$ factors to give  an  isomorphism 
$$N_H(T)^\sigma/N_{(H^\sigma)^0}(S^\sigma)  \cong \pi_0(H^\sigma).$$
\end{lemma}

\begin{proof}  The induced map  $N_H(T)^\sigma$ to 
$\pi_0(H^\sigma)$ is surjective by Lemma~\ref{signormal}.  Its
kernel is
$N_H(T)^\sigma\cap(H^\sigma)^0$.
Any element in this intersection normalizes $T$ and its action on $T$
commutes with $\sigma$. Thus, it also normalizes
$S^\sigma$. Conversely, since $S^\sigma$ contains a regular element 
of
$T$, any element of $(H^\sigma)^0$ normalizing
$S^\sigma$ normalizes $T$.
\end{proof}

\begin{lemma}\label{compgrp2} There is an exact sequence
$$\{1\} \to \pi_0(T^\sigma) \to
N_H(T)^\sigma/N_{(H^\sigma)^0}(S^\sigma)  \to 
\ov N(\sigma)/W(S^\sigma, (H^\sigma)^0) \to \{1\}.$$
\end{lemma}
\begin{proof} By definition, there is an exact sequence
$$\{1\} \to T^\sigma \to N_H(T)^\sigma \to \ov N(\sigma) \to
\{1\}.$$ Now $N_{(H^\sigma)^0}(S^\sigma)$ is a normal
subgroup of $N_H(T)^\sigma$, and the image of 
$N_{(H^\sigma)^0}(S^\sigma)$ in $W(T,H)$ is $W(S^\sigma,
(H^\sigma)^0)$. Thus $\ov N(\sigma)$ contains
$W(S^\sigma, (H^\sigma)^0)$ as a normal subgroup, and there is a 
surjection 
$$N_H(T)^\sigma/N_{(H^\sigma)^0}(S^\sigma)  
\to \ov N(\sigma)/W(S^\sigma, (H^\sigma)^0).$$
By a diagram chase, the kernel of this map is then
$$T^\sigma/N_{(H^\sigma)^0}(S^\sigma)\cap T =
T^\sigma/S^\sigma=\pi_0(T^\sigma).$$
This establishes the exact sequence of the statement.
\end{proof}

\begin{corollary}\label{isomIused}
There is an exact sequence
$$\{1\} \to \pi_0(T^\sigma) \to
\pi_0(H^\sigma) \to 
\ov N(\sigma)/W(S^\sigma, (H^\sigma)^0) \to \{1\}.$$
\end{corollary}

In the case of an inner automorphism, this result is due to
Steinberg \cite{Stein}.

\subsection{Special automorphisms}

\begin{lemma}\label{sigma1}
Let  $\sigma\colon H\to H$ be an automorphism in normal form.
There is an element $s\in S^{\sigma}$ such that setting
$\sigma_1=i_s\circ \sigma$, the set 
$\Delta_H/\langle\sigma\rangle $ is a set of simple roots for 
$H^{\sigma_1}$ with respect
to the maximal torus $S^{\sigma}=S^{\sigma_1}$.  
\end{lemma}

\begin{proof}
Since $\Delta_H$ is an linearly independent subset of ${\frak t}^*$ 
invariant under $\sigma$, the quotient set $\Delta_H/\langle\sigma\rangle
$ is a linearly independent subset of $({\frak t}^\sigma)^*$.
For each $\ov a\in  \Delta_H\spcheck/\langle\sigma\rangle$, order
$\ov a = \{a_1, \dots, a_{n_{\ov a}}\}$  so that $\sigma\cdot
a_i=a_{i+1}, 1\le i\le n_{\ov a}$ (by convention 
$a_{n_{\ov a}+1}=a_1$).  Then the  action of $\sigma$ on ${\frak
h}\otimes_{\bf R}{\bf C}$ sends the root space ${\frak h}^{a_i}$ to
${\frak  h}^{a_{i+1}}$. Hence, $\sigma^{n_{\ov a}}$ sends ${\frak
h}^{a_i}$ to itself for every $i$. This means that there is
$q_{a_i}\in 
U(1)$ such that   $\sigma^{n_{\ov a}}|{\frak h}^{a_i}$ is
multiplication by 
$q_{a_i}$ on each root space ${\frak h}^{a_i}$.
By equivariance under $\sigma$, we see that
$q_{a_i}=q_{a_j}$ for all $1\le i,j\le n_{\ov a}$. We denote this
common value by $q_{\ov a}$.
Choose an element $\mu_{\ov a}\in U(1)$ such that $\mu_{\ov
a}^{n_{\ov a}}=q_{\ov a}$. 
Since the elements of $ \Delta_H/\langle\sigma\rangle$ are linearly 
independent in $({\frak t}^{\sigma})^*$, it follows that there is an
element $s\in S^{\sigma}$ such that $s^a =\mu_{\ov a}^{-1}$ for
every $a\in I$.

We set $\sigma_1=s\sigma$.
By construction,
if $\ov a = \{a_1, \dots, a_{n_{\ov a}}\}$, then $\sigma_1^{n_{\ov
a}}|\bigoplus_{i}{\frak h}^{a_i}$ is the identity. Choose a non-zero
element
$X_a\in {\frak h}^a$. Set
$$X_{\ov a}=X_a+\sigma_1X_a+\cdots +\sigma_1^{n_{\ov a}-1}X_a.$$
Then $X_{\ov a}$ is  a non-zero element of ${\frak h}\otimes_{\bf
R}{\bf C}$  invariant under $\sigma_1$, and hence $X_{\ov a}
\in {\rm Lie}(H^{\sigma_1})\otimes_{\bf R}{\bf C}$. Moreover,
for all $s\in S^\sigma$, ${\rm Ad}(s)X_{\ov a} = \ov a(s)X_{\ov
a}$. Thus,  each
$\ov a\in \Delta_H/\langle\sigma\rangle$ is a root  of
$H^{\sigma_1}$ with respect to $S^{\sigma}$. Since every root of
$H^{\sigma_1}$ is the restriction of  a root of $H$ to
$S^{\sigma}$ and since $\Delta_H$ is a set of simple roots for
$H$, it follows that  every root of
$H^{\sigma_1}$ can be written uniquely as a linear combination of the
elements in $\Delta_H/\langle\sigma\rangle $ and the coefficients of this 
linear combination are either all positive integers
or all negative integers. Consequently, $\Delta_H/\langle\sigma\rangle $ 
is a set of simple roots of $H^{\sigma_1}$ with respect to
$S^{\sigma}$. 
\end{proof}

\begin{defn}
An automorphism $\sigma_1$ in normal form which satisfies
the conclusions of the previous lemma is called {\sl a special 
automorphism}.
\end{defn}

\begin{proposition}\label{indivis}
Let $\sigma_1$ be a special automorphism. Then 
the roots of
$H^{\sigma_1}$ are the indivisible roots in 
$\Phi_H^{\rm res}(\sigma_1)$. In particular, the Weyl group of
$H^{\sigma_1}$ with respect to $S^{\sigma_1}$ is equal to the Weyl
group of $\Phi_H^{\rm res}(\sigma_1)$, namely $Z_W(\sigma_1)$.
\end{proposition}

\begin{proof}
According to Lemma~\ref{linearcase}, the indivisible roots of
$\Phi_H^{\rm res}(\sigma_1)$ form a reduced root system with
$\Delta_H/\langle\sigma_1\rangle$  a set of simple roots. By 
Lemma~\ref{sigma1}, $\Delta_H/\langle\sigma_1\rangle$ is also a set
of simple roots for $\Phi(S^{\sigma_1},(H^{\sigma_1})^0)$. Since
the latter is also a reduced root system, the two root systems
agree. By Lemma~\ref{Weylcent} and Proposition~\ref{prop},
$W(S^{\sigma_1},(H^{\sigma_1})^0)$ is then $Z_W(\sigma_1)$.
\end{proof}

\begin{corollary} If $\sigma_1$ is a special automorphism,
  then $\pi_0(H^{\sigma_1}) =\pi_0(T^{\sigma_1})$. 
\end{corollary}
\begin{proof} By Lemma~\ref{compgrp2},
$\pi_0(H^{\sigma_1})/\pi_0(T^\sigma) \cong
\ov N(\sigma_1)/W(S^{\sigma_1}, (H^{\sigma_1})^0)$. By the previous
lemma, $W(S^{\sigma_1}, (H^{\sigma_1})^0)= Z_W(\sigma_1)$, and
$\ov N(\sigma_1) \subseteq Z_W(\sigma_1)$. Thus the quotient
$\ov N(\sigma_1)/W(S^{\sigma_1}, (H^{\sigma_1})^0)$ is trivial, so
that
$\pi_0(H^{\sigma_1}) =\pi_0(T^{\sigma_1})$. 
\end{proof}

The following relates the roots of $H^{\sigma_1}$, which are 
restricted roots, to the projection root system.

\begin{corollary}\label{one-to-one}
There is a natural one-to-one correspondence between the roots of
$H^{\sigma_1}$ and the subroot system of the non-multipliable roots 
in
$\Phi_H^{\rm proj}(\sigma_1)$. The corresponding roots in these two
systems are positive multiples of each other.
\end{corollary}

\begin{proof}
According to Proposition~\ref{extension}, these root systems  have
the same Weyl group, which implies that they have the same set of
walls.
Since each of these root systems is reduced, the result follows.
\end{proof}

For example, if $\Phi$ is simply laced, which is always the case if
$\sigma_1$ is nontrivial and
$\Phi$ is irreducible, and if there are no exceptional
orbits, then
$H^{\sigma_1} =
\Phi_H^{\rm res}(\sigma_1)$ and
$\Phi_H^{\rm proj}(\sigma_1)$ is the inverse system, and the
bijection of the corollary simply associates to each element of
$H^{\sigma_1}$ its inverse coroot.

\subsection{The roots of $H^{\sigma}$}

Suppose that 
$\sigma_1$ is a special automorphism, and let $B$ be an
alcove of $\frak d ^{\sigma_1}$ containing the origin for
$\Phi_H^{\rm proj}(\sigma_1)$. Let
$\sigma = i_s\circ
\sigma_1$ for some $s\in S^{\sigma_1}$. Let $\tilde s$ be a lift of
$s$ whose projection to $\frak d^{\sigma_1}$ lies in some alcove
containing the origin. By Part 1 of Lemma~\ref{Weylcent}, the group
$Z_W(\sigma_1)$ acts transitively on the set of all such
alcoves, and hence there is a $w\in Z_W(\sigma)$ such that 
$w(\tilde s) \in B$. Moreover,   by Proposition~\ref{indivis},  we
may lift $w$ to $g\in N_H(T)^{\sigma_1}$. Then
$i_g\circ i_s\circ \sigma_1\circ i_g^{-1} = i_{w(s)}\circ
\sigma_1$, and the fixed subgroups of $i_s\circ \sigma_1$ and $i_g\circ i_s\circ \sigma_1\circ i_g^{-1} = i_{w(s)}\circ
\sigma_1$ are conjugate in $H$. Thus, to analyze the fixed
subgroups, we can always assume that the projection
$\hat s$ of $\tilde s$ to
$\frak d^\sigma$ lies in the alcove $B$. 
Denote by $\ov S^\sigma$ the quotient of ${\frak t}^\sigma$ by the
image under orthogonal projection of $\Lambda$. Note that this
terminology is consistent with that of the introduction. By
Lemma~\ref{basicaut}, the component group of $H^\sigma$ only 
depends on the image of $s$ in $\ov S^\sigma$.

\begin{lemma}\label{rootssy1}
Let $\sigma_1$ be a special  automorphism  and suppose that the
  action of $\sigma_1$ on $\Delta_H$ has no exceptional orbits.
Let $\tilde s\in {\frak  d}^{\sigma_1}$, and set  $s={\rm exp}(\tilde
  s)$. Let  $\sigma=i_s\circ \sigma_1$. Then the roots of
  $H^\sigma$ with respect to the maximal torus $S^{\sigma}$ are
$$\{\ov a\in \Phi_H^{\rm res}(\sigma)\bigl|\ \bigr.
n_{\ov a}\langle \ov a,\tilde s\rangle\in{\bf Z}\rangle\}.$$ 
\end{lemma}

\begin{proof}
Since $\sigma|T=\sigma_1|T$, it follows that $S^\sigma=S^{\sigma_1}$
and  that $\sigma$ and $\sigma_1$ have the same action on $\Phi_H$ so
that
$\Phi_H^{\rm proj}(\sigma)=\Phi_H^{\rm proj}(\sigma_1)$.
Since $\sigma_1$ has no exceptional orbits, it follows from
Lemma~\ref{linearcase}
that 
$\Phi_H^{\rm res}(\sigma)$ is a reduced root system. Thus, by
Lemma~\ref{indivis} each  $\ov a$ in $\Phi_H^{\rm res}(\sigma_1)$ is a
root of of $H^{\sigma_1}$ with respect to ${\frak t}^{\sigma_1}$.

The roots of $H^{\sigma}$ with respect to ${\frak t}^{\sigma}$ are a
subset of $\Phi_H^{\rm res}(\sigma_1)$.
Let $\ov a\in  \Phi_H^{\rm res}(\sigma_1)$ be an orbit
of order $n_{\ov a}$. Let $\{a_1,\ldots,a_{n_{\ov a}}\}$ be this
orbit, ordered so that $\sigma_1\cdot a_i=a_{i+1}$ (by convention
$a_{n_{\ov     a}+1}=\ov a_1$). 
For each $i$ let ${\frak h}^{a_i}\subseteq {\frak h}\otimes {\bf C}$ be
the root space for $a_i$.
Since $\ov a$ is a root of $H^{\sigma_1}$, it follows that the action
of $\sigma_1$ identifies the 
root space ${\frak h}^{a_i}$ with ${\frak h}^{a_{i+1}}$ in such a way
that the  action of $\sigma_1^{n_{\ov a}}$ is the identity on each of
these root spaces. Since $s\in S^{\sigma}$, we have $s^{a_i}
=s^{a_{i+1}}$ for all $i\le n_{\ov a}$. Thus,  
$s^{n_{\ov a}a_1}=1$ if and only if  the  action of
$(s\sigma_1)^{n_{\ov a}}$ is trivial on all these root spaces if and only
if $\ov a$ is a root of $H^{\sigma}$. 
Of course $s^{n_{\ov a}a_1}=1$ if and only if $n_{\ov
  a}\langle \ov a,\tilde s\rangle\in{\bf Z}$.
\end{proof}

\begin{corollary}\label{5678}
With $\sigma_1$, $\tilde s$ and $\sigma$ as in the previous lemma,
 the roots of $H^{\sigma}$ are exactly those of
$\Phi_H^{\rm res}(\sigma_1)$ that correspond under the bijection 
given in Corollary~\ref{one-to-one} to the roots of
$\Phi_H^{\rm proj}(\sigma_1)$ which take integral values on  $\tilde
s$.
\end{corollary}

\begin{proof}
Since $\sigma$ has no exceptional orbits, $\Phi_H^{\rm proj}(\sigma)$
is a reduced root system and the root of this system corresponding to
$\ov a\in \Phi_H^{\rm res}(\sigma)$ is 
$n_{\ov a}\ov a$. Given this the result is immediate from
the previous lemma.
\end{proof}

Suppose that $\sigma_1$ is nontrivial, that $\Phi$ is
irreducible, and that there are no exceptional orbits.
It follows that $\Phi$ is simply laced.
By Lemma~\ref{resdualproj} the systems
$\Phi_H^{\rm res}(\sigma_1)$ and
$\Phi_H^{\rm proj}(\sigma_1)$ are inverse to each other. The
roots of $H^{\sigma}$ are then the elements of $\Phi_H^{\rm
res}(\sigma_1)$ such that the corresponding inverse coroot is
integral on $\tilde s$, or equivalently on $\hat s$. Since $\hat
s$ lies in the alcove $B$, a set of simple roots for
$H^{\sigma}$ is obtained as follows:  take the extended coroot
diagram for
$\Phi_H^{\rm res}(\sigma_1)$, in other words the extended root
diagram for
$\Phi_H^{\rm proj}(\sigma_1)$. Let $\Phi_H^{\rm
proj}(\sigma_1)(\hat s)$ be the set of elements of $\Phi_H^{\rm
proj}(\sigma_1)$ such that the corresponding wall of $B$ contains
$\hat s$. Then the set of elements of $\Phi_H^{\rm
res}(\sigma_1)$ which are inverse to an element of $\Phi_H^{\rm
proj}(\sigma_1)(\hat s)$ is a set of simple roots for $H^\sigma$. 
The possible diagrams obtained in this way will describe all of
the possible root systems for $H^\sigma$. A similar result holds
if there are exceptional orbits, with a slightly more involved
proof. In this case, the root and coroot systems are both of type
$BC_n$, and thus are abstractly isomorphic, and the procedure for 
finding the possible root systems of the $H^\sigma$ is again to 
take the extended coroot diagram, choose a proper subdiagram, and
then pass to the inverse system. Taken together, these results 
generalize a theorem of Kac \cite{Kac1} which dealt with the case
of finite order automorphisms $\sigma$.

\subsection{The component group of $H^\sigma$}

Let $\sigma$ be an automorphism in normal form.
We keep the notation above, and for simplicity assume that there 
are  no exceptional orbits of $\sigma$ (a minor modification
handles the general case). 
We now use the results developed above to give an explicit 
description of the component group of the fixed subgroup 
$H^\sigma$.

\subsubsection{The homomorphism $\delta$}

Let $\widetilde H$ be the
universal covering group of $H$. Then 
there is an exact sequence
$$\{0\}\to \Lambda/Q\spcheck_H\to \widetilde H\to H\to \{1\},$$
where $\Lambda/Q\spcheck_H$ is a central subgroup of $\widetilde H$.
The automorphism $\sigma$ acts on this sequence.
Thus, there is
an exact sequence
$$\{0\}\to\left(\Lambda/Q\spcheck_H\right)^\sigma\to  \widetilde
H^\sigma\to H^\sigma\to \left(\Lambda/Q\spcheck_H\right)_{\sigma}.$$
Since $\widetilde H$ is the product of a vector space and a compact, 
simply connected group, it follows from
\cite{Tohoku} that
$\widetilde H^\sigma$ is connected. (In fact, we shall give
another proof of this fact shortly.) Thus we have:

\begin{lemma}\label{pi0lambda}
Let $g\in H$ and let $\tilde g$ be a lift of $g$ to $\widetilde
 H$.
Then the element $\sigma(\tilde g)\tilde g^{-1}$ is contained in
 $\Lambda/Q\spcheck_H$ and its image in 
$(\Lambda/Q\spcheck_H)_{\sigma}$
 depends only on $g$. This function 
induces an injection $\delta$ from 
$\pi_0(H^\sigma)$ to  ${\rm
Tor}\left(\Lambda/Q\spcheck_H\right)_{\sigma}$.
\end{lemma}

\begin{proof}
   From the long exact cohomology sequence we see that $\delta$ is 
an injective homomorphism from $\pi_0(H^\sigma)$ to
$\left(\Lambda/Q\spcheck_H\right)_{\sigma}$. Since $H^\sigma$ is a
compact group, the image of $\delta$ is finite, and hence
is contained  in 
${\rm Tor}\left(\Lambda/Q\spcheck_H\right)_{\sigma}$.
\end{proof}

Applying the above to the case $H=T$, there is also a homomorphism
from $\pi_0(T^\sigma)$ to ${\rm Tor}(\Lambda_\sigma)$. Of course,
this is just the isomorphism of Lemma~\ref{componentgroup}. By the
functoriality of the connecting homomorphism $\delta$, we have:

\begin{lemma}\label{CD0} There is a commutative diagram with exact
  columns: 
$$\CD
0 @. 0\\
@VVV @VVV\\
\pi_0(T^\sigma) @>{\cong}>> {\rm Tor}(\Lambda_\sigma)\\
@VVV @VVV\\
\pi_0(H^\sigma) @>{\delta}>> {\rm
Tor}\left((\Lambda/Q\spcheck_H)_{\sigma}\right).
\endCD$$
\end{lemma}

\subsubsection{The homomorphism from ${\rm
    Tor}\left((\Lambda/Q\spcheck_H)_{\sigma}\right)$ to ${\cal
    C}\Phi_H^{\rm proj}(\sigma)$}

\begin{lemma} The natural map $\Lambda_{\sigma}\to {\frak t}^{\sigma}$
    induced by 
the orthogonal projection from $\Lambda$ to ${\frak t}^\sigma$ has
kernel equal 
to ${\rm Tor}\left(\Lambda_{\sigma}\right)$ and image contained in the
dual lattice to the root lattice of $\Phi_H^{\rm proj}(\sigma)$. 
\end{lemma}

\begin{proof}
The only thing that needs to be established is that if $r$ is a root
of $\Phi_H^{\rm proj}(\sigma)$ and $\lambda\in \Lambda_{\sigma}$
then $\langle r,\lambda\rangle \in {\bf Z}$. The root lattice 
$\Phi_H$ has as a basis the roots $\{\epsilon(\ov a)n_{\ov a}\ov
a\}_{\ov a\in \Delta_H/\langle\sigma\rangle}$. It suffices to show
that these take integral values on $\lambda$. Let $\{a_1,\ldots,
a_{n_{\ov a}}\}$  be the orbit $\ov a$ and lift $\lambda$ to
$\tilde
\lambda\in \Lambda$. Then
$\langle r,\lambda\rangle=\epsilon(\ov a) \sum_{i=1}^{n_{\ov
a}}\langle a_i,\tilde\lambda\rangle\in{\bf Z}$. 
\end{proof}

Recall that ${\cal C}\Phi_H^{\rm proj}(\sigma)$ is the quotient of 
the  coweight lattice of the root system $\Phi_H^{\rm
proj}(\sigma)$ by its coroot lattice.
The exact sequence
$$0\to Q\spcheck_H \to \Lambda \to \Lambda/Q\spcheck_H\to 0$$
yields an exact sequence
$$(Q\spcheck_H)_\sigma\to\Lambda_\sigma \to
\left(\Lambda/Q\spcheck_H\right)_\sigma\to 0.$$
The image in $\frak t_\sigma$ of any $\mu\in \Lambda_{\sigma}$ which
maps to a torsion element in $(\Lambda/Q\spcheck_H)_{\sigma}$ is
contained in ${\frak
 d}^{\sigma}$. Thus, 
we have constructed a homomorphism from ${\rm
  Tor}\left((\Lambda/Q\spcheck_H)_{\sigma}\right)$ to the quotient
  ${\frak d}^{\sigma}/(Q\spcheck_H)_{\sigma}$. The kernel of this
  homomorphism is the image of ${\rm Tor}(\Lambda_{\sigma})$ in ${\rm
  Tor}\left((\Lambda/Q\spcheck_H)_{\sigma}\right)$ and the
  image is contained in ${\cal C}\Phi_H^{\rm proj}(\sigma)$.
Finally, the map from ${\rm Tor}(\Lambda_{\sigma})$ in ${\rm
  Tor}\left((\Lambda/Q\spcheck_H)_{\sigma}\right)$ is injective by
Lemma~\ref{CD0}. We summarize this picture in the following lemma:

\begin{lemma}\label{usefulexseq} Orthogonal projection from
$\Lambda_\sigma$ to $\frak t_\sigma$ induces an exact sequence
$$0 \to  {\rm Tor}(\Lambda_{\sigma}) \to {\rm
  Tor}\left((\Lambda/Q\spcheck_H)_{\sigma}\right) \to {\cal C}\Phi_H^{\rm
proj}(\sigma).$$
\end{lemma}

Recall from Section~\ref{action1} that is an affine action of  
${\cal C}\Phi_H^{\rm proj}(\sigma)$ on ${\frak  d}^{\sigma}$
normalizing
$B$.
This induces an action of
${\rm Tor}\left((\Lambda/Q\spcheck_H)_{\sigma}\right)$, by affine
isometries,  on ${\frak d}^{\sigma}$  normalizing $B$.
The linearization of the action of  ${\cal C}\Phi_H^{\rm
proj}(\sigma)$ is denoted $\nu_B$.

We have identified $\Lambda_\sigma/{\rm Tor}(\Lambda_{\sigma})$ with
$\pi(\Lambda)$. Moreover $\pi(Q_H\spcheck)\cong (Q_H\spcheck)_\sigma$. Thus
$$\pi(\Lambda)/\pi(Q_H\spcheck)\cong 
(\Lambda/Q\spcheck_H)_{\sigma}/{\rm Tor}(\Lambda_\sigma).$$   The
following lemma, which will be needed in the next section, is then
clear.

\begin{lemma}\label{ratreps}
There is an exact sequence 
$$0\to \pi(Q_H\spcheck) \to (\pi(Q_H\spcheck)\otimes {\bf Q})\cap
\pi(\Lambda) \to {\rm Tor}\left((\Lambda/Q\spcheck_H)_{\sigma}\right)/{\rm
Tor}(\Lambda_\sigma) \to 0.$$
\end{lemma}

\subsubsection{The basic equation and its consequences}

\begin{proposition}\label{basicform} Suppose that $\sigma_1$ is
special, and let $\sigma =i_s\circ \sigma_1$ with  $s =\exp (\tilde
s)$ for some $\tilde s \in \frak t^{\sigma_1}$.  Let $w\in
Z_W(\sigma)$ and let $g$ be a lift of $w$ to
$N_{(H^{\sigma_1})^0}(S^\sigma)$. Let $t\in T$ with $t=\exp (\tilde
t)$ for some $\tilde t\in \frak t$. Then the element $tg$ lies in
$H^\sigma$ if and only if $$({\rm Id} - w)\tilde s + (\sigma - {\rm
Id})\tilde t \in \Lambda,$$ if and only if $({\rm Id} - w)\tilde s
\in \pi (\Lambda)$.  Thus, if $\ov s$ is the image of
$s$ in $T_\sigma=\ov S^\sigma$, $w\in \ov N(\sigma)$ if and only if 
$w(\ov s) = \ov s$ in the induced action of $w$ on $T_\sigma$.
\end{proposition}

\begin{proof} Since $\sigma_1(g) = g$, we see that $i_s\circ \sigma_1(tg) =
tg$ if and only if
$$s\sigma_1(t)\sigma_1(g) s^{-1} = s\sigma_1(t)g s^{-1} =tg,$$
or equivalently
$$1= t^{-1}s\sigma_1(t)g s^{-1}g^{-1}= (t^{-1}\sigma_1(t))(s(w\cdot
s)^{-1}.$$ Equivalently, in additive notation,
$$({\rm Id} - w)\tilde s + (\sigma - {\rm Id})\tilde t \in 
\Lambda.$$ Since every lift of $w$ to an element of $N_H(T)$ is of
the form $tg$ for some $t\in T$, we see that $w$ lifts to an element
of $(N_H(T))^\sigma$ if and only if $s\equiv ws \bmod ({\rm Id} -
\sigma)(T)$, or equivalently 
$w(\ov s) = \ov s$.
\end{proof}

\begin{proposition}\label{fixes}
Suppose that $\sigma =i_s\circ \sigma_1$ with  $s =\exp (\tilde
s)$ for some $\tilde s \in \frak t^{\sigma_1}$ and that the
projection $\hat s$ of $\tilde s$ to $\frak d^\sigma$ lies in the
alcove $B$. Then the group
$\pi_0(H^\sigma)/\pi_0(T^\sigma)$ is isomorphic to the stabilizer of
$\hat s$ in ${\cal C}\Phi^{\rm proj}_H(\sigma)$.
\end{proposition}
\begin{proof} By Proposition~\ref{basicform}, $\ov N(\sigma)={\rm
Stab}_{Z_W(\sigma_1)}(\ov s)$ is the stabilizer in the Weyl group of
$\Phi^{\rm proj}_H(\sigma)$ of the point $\ov s$. By
Corollary~\ref{5678}, the Weyl group $W(S^\sigma, (H^\sigma)^0)$ is
the group $Z_W(\sigma_1)(\ov s)$ generated by  reflections in walls
defined by the roots in
$\Phi_H^{\rm   proj}(\sigma)$ which are integral on $\hat s$.  By 
Lemma~\ref{stabilizer} (and the remarks immediately following it in
case $\Phi^{\rm proj}_H(\sigma)$ is not reduced), the quotient ${\rm
Stab}_{Z_W(\sigma_1)}(\ov s)/Z_W(\sigma_1)(\ov s)$ is isomorphic to
the  stabilizer of $\hat s$ under the action of ${\cal C}\Phi^{\rm
proj}_H(\sigma)$ on $B$, and this is the statement
of the proposition.
\end{proof}

We can now sketch a proof of the theorem of the first author in
\cite{Tohoku}:

\begin{corollary}
If $H$ is simply connected, then $H^\sigma$ is connected.
\end{corollary} 
\begin{proof}  Since $H$ is simply connected, 
$\Lambda/Q_H\spcheck$ is trivial, and hence
$\pi_0(H^\sigma)/\pi_0(T^\sigma)$ is trivial. By 
Lemma~\ref{componentgroup}, $\pi_0(T^\sigma)$ is also trivial. Thus
$\pi_0(H^\sigma)$ is trivial, so that $H^\sigma$ is connected.
\end{proof}

\begin{proposition}\label{CD} With
notation as above,
let $p\colon \pi_0(H^\sigma) \to \ov N(\sigma)/W(S^\sigma,
(H^\sigma)^0)$ be the map defined by the exact sequence of
Lemma~\ref{isomIused}, let $\delta\colon \pi_0(H^\sigma) \to {\rm
Tor}\left((\Lambda/Q\spcheck_H)_{\sigma}\right)$ be the homomorphism
of Lemma~\ref{pi0lambda}, and let $\nu_B$ be the linearization of the
action of  ${\cal C}\Phi_H^{\rm proj}(\sigma)$ on $B$, viewed as a
homomorphism from ${\rm
Tor}\left((\Lambda/Q\spcheck_H)_{\sigma}\right)$ to $Z_W(\sigma_1)$.
Then the image of $\nu_B\circ \delta$ is contained in $\ov
N(\sigma)$, and 
$p=\nu_B\circ \delta$ as homomorphisms from $\pi_0(H^\sigma)$ to 
$\ov N(\sigma)/W(S^\sigma, (H^\sigma)^0)$.
\end{proposition}

\begin{proof} Fix a component of $H^\sigma$, and find, by
Lemma~\ref{signormal}, an element
$g_0\in N_H(T)^{\sigma}$ lying in this component.  Let
$w\in Z_W(\sigma)$ be the image of $g_0$. Then by definition $w\in
\ov N(\sigma)$ and 
$p([g_0]) = w\bmod W(S^\sigma, (H^{\sigma})^0)$. By
Proposition~\ref{indivis}, there is a
$g\in N_H(T)^{\sigma_1}$ whose image in $Z_W(\sigma_1)$ is also $w$.
Thus $g_0 = tg$ for some $t\in T$.   Let $\tilde g$ be a lift of $g$ to 
the  covering $\frak t \rtimes W$ of $N_H(T)$ and let $\tilde
t$ be a lift of $t$ to $\frak t$. Then the element
$\tilde t \tilde g$ is a lift of
$tg$ to  $\frak t \rtimes W$. The
calculation of Proposition~\ref{basicform} shows that
$$\lambda = ({\rm Id} -w)\tilde s + (\sigma -{\rm
Id})\tilde t\in \Lambda$$
and identifies
$\sigma(\tilde t\cdot\tilde g)(\tilde t\cdot\tilde g)^{-1}$ 
(additively) with 
$\lambda$, up to an element in $({\rm Id} - \sigma)(\Lambda)$.
It then follows  from the definition of
$\delta$ that
$\delta([g_0]) = \delta([tg])$ is the image of $\lambda$ in 
$(\Lambda/Q_H\spcheck)_\sigma$. Thus, letting $\pi$ denote the
orthogonal projection onto ${\frak
  t}^\sigma$,  by Lemma~\ref{usefulexseq}, the image $\ov \xi$ of
$\delta([tg])$ in 
${\cal C}\Phi^{\rm proj}_H(\sigma)$ is represented by $\pi(\lambda)$
modulo $\pi(Q_H\spcheck)$.
Of course,  $\pi(\lambda)=\pi({\rm Id}-w)\tilde s=({\rm
  Id}-w)(\hat s)$.
Thus $\xi=\hat s-w(\hat s) \in{\frak t}^\sigma$ lifts $\ov
\xi$. In other words, $w(\hat s) + \xi =\hat s$. Let $B' = w(B) +
\xi$; it is an alcove in $\frak d^\sigma$ for $\Phi_H^{\rm
  proj}(\sigma)$. Now there is a unique
$\gamma\in W_{\rm aff}(\Phi_H^{\rm
  proj}(\sigma))$ such that $\gamma(B') =
B$ and $\gamma(\hat s) = \hat s$. Let $w'$ be the Weyl part of
$\gamma$. It is a product of reflections in walls defined by
roots in $\Phi_H^{\rm
  proj}(\sigma)$ which are integral on $\hat s$. By
Corollary~\ref{5678}, such walls  are walls of roots of
$H^\sigma$. Hence $w'\in W(S^\sigma,
(H^{\sigma})^0)$. Let $\varphi$ be the affine transformation 
$\gamma w +\xi = w'w+ \xi'$, where $\xi' \equiv \xi $ mod
the coroot lattice of $\Phi_H^{\rm
  proj}(\sigma)$. Then by construction $\varphi$ normalizes the
alcove
$B$, and so by definition $w'w = \nu_B(\xi') = \nu_B(\xi) =
\nu_B(\delta([g_0]))$.   Since $w\equiv w'w\bmod W(S^\sigma,
(H^{\sigma})^0)$, we see that $\nu_B(\delta([g_0])) \in \ov
N(\sigma)$ and that
$p([g_0]) =w'w =  \nu_B(\delta([g_0]))$. This completes the proof. 
\end{proof}

\begin{corollary}\label{fullstatement} Let $\sigma = i_s\circ 
\sigma_1$, let
$\hat s$ be the image of $\tilde s$ in $\frak d^\sigma$, and assume 
that
$\hat s$ lies in the alcove $B$. Then the  group
$\pi_0(H^\sigma)$ is isomorphic to the stabilizer of $\hat s$ in the 
group
${\rm Tor}\left((\Lambda/Q\spcheck_H)_{\sigma}\right)$. Moreover, for every
element
$\mu \in {\rm
Tor}\left((\Lambda/Q\spcheck_H)_{\sigma}\right)$, there exists an 
$s$ such that, if $\sigma = i_s\circ \sigma_1$, then $\mu =\delta
([z])$ for some $z\in H^\sigma$.
\end{corollary}
\begin{proof}  
Let $\ov \delta$ be the homomorphism $\pi_0(H^\sigma) \to {\rm
Tor}\left((\Lambda/Q\spcheck_H)_{\sigma}\right) \to {\cal
C}\Phi_H^{\rm proj}(\sigma) $. By Proposition~\ref{CD}, the image of
$\ov \delta$ is exactly the stabilizer of $\hat s$ in ${\cal
C}\Phi_H^{\rm proj}(\sigma)$. Since the image of $\delta$ contains
the kernel of the map from ${\rm
Tor}\left((\Lambda/Q\spcheck_H)_{\sigma}\right)$ to ${\cal
C}\Phi_H^{\rm proj}(\sigma)$, it follows that the image of $\delta$
is the stabilizer of $\hat s$ in ${\rm
Tor}\left((\Lambda/Q\spcheck_H)_{\sigma}\right)$.

The second statement follows immediately from the fact that every
affine automorphism of a simplex has a fixed point. 
\end{proof}

\subsection{The case of $c$-pairs}

We return to the notation of the previous sections, so that $G$ is a
simple and simply connected group with maximal torus $T$, and $A$ is
a fixed alcove in $\frak t$. Let
$c\in {\cal C}G$, and suppose that $(x,y)$ is a $c$-pair. After
conjugation we can assume that
$(x,y)$ is in normal form, so that $x={\rm exp}(\tilde x)$ 
for some $\tilde x\in A^c$.  We will
apply the results above to the compact group $H=Z(x)$ and the
automorphism of $Z(x)$ defined by conjugation by $y$. For
future reference, let us record the following definition:

\begin{defn}\label{specialcpair} A $c$-pair $(x,y_1)$ in normal form
such that conjugation by $y_1$ is a special automorphism of $Z(x)$ is
called a {\sl special $c$-pair}.
\end{defn} 

Note that
$T$ is a maximal torus for $Z(x)$ and that a set of simple roots
$\widetilde I(x)$ for $Z(x)$ is given by the subset of $\widetilde
\Delta$ taking integral values on $\tilde x$.   Let $Q\spcheck(x)$
be the lattice generated by $\widetilde I(x)$. It is the coroot
lattice for the derived subgroup of $Z(x)$, and thus is the lattice
denoted $Q_H\spcheck$ above. Conjugation by
$y$ induces an  automorphism of $Z(x)$ normalizing ${\frak t}$ 
and $A\subseteq {\frak t}$. Its action on ${\frak t}$ is given by the
Weyl element $w_c$. This Weyl element normalizes
$\widetilde \Delta$ and $\widetilde I_c(x)$. Let $\widetilde
\Delta_c$ and $\widetilde I(x)$ be the respective quotients of
these sets by the action of $w_c$. The elements of
$\widetilde
\Delta_c\subseteq {\frak t}^{w_c}$ satisfy one relation  $\sum_{\ov
a\in
\widetilde
\Delta_c}g_{\ov a}\pi( a\spcheck)=0$.

\subsubsection{Description of ${\rm
    Tor}\left((Q\spcheck/Q\spcheck(x))_{w_c}\right)$}

\begin{lemma}\label{torsion1}
${\rm Tor}\left((Q\spcheck/Q\spcheck(x))_{w_c}\right)$ is a
cyclic group of  order $n={\rm gcd}\{\, g_{\ov a}:\ov a\in\widetilde 
\Delta_c-\widetilde I_c(x)\,\}$, generated by 
$$\zeta=\frac{1}{n}\sum_{\ov
a\in \widetilde \Delta_c-\widetilde I_c(x)}g_{\ov a}\pi( a\spcheck)
=-\sum_{\ov a\in \widetilde I_c(x)}\frac{g_{\ov a}}{n}\pi( a\spcheck)\in
Q\spcheck.$$ 
\end{lemma}

\begin{proof}
We have an exact sequence
$$0\to {\bf Z}\left(\sum_{a\in\widetilde \Delta}g_aa\spcheck\right)\to
\bigoplus_{a\in \widetilde \Delta}{\bf Z}(a\spcheck)\to Q\spcheck\to 0.$$
Since $Q\spcheck(x)=\bigoplus_{a\in \widetilde I(x)}{\bf Z}(a\spcheck)$, we
see that there is an exact sequence
$$0\to {\bf Z}\left(\sum_{a\in\widetilde \Delta-\widetilde
    I(x)}g_aa\spcheck\right)\to \bigoplus_{a\in \widetilde\Delta-\widetilde
  I(x)}{\bf Z}(a\spcheck)\to Q\spcheck/Q\spcheck(x)\to 0.$$ 

Conjugation by $y$ induces an action on this sequence.
On the second term it is the action induced by the permutation action
of $\widetilde \Delta-\widetilde I(x)$ with quotient $\widetilde
\Delta_c-\widetilde I_c(x)$.   Hence, taking coinvariants yields an
exact sequence
$$ {\bf Z}\left(\sum_{\ov a\in\widetilde \Delta_c-\widetilde I_c(x)}g_{\ov
a}\pi( a\spcheck)\right)\to \bigoplus_{\ov a\in 
\Delta_c-\widetilde I_c(x)}{\bf Z}(\pi( a\spcheck))\to
\left(Q\spcheck/Q\spcheck(x)\right)_{w_c} \to 0.$$ 
Clearly, then, the torsion subgroup of
$\left(Q\spcheck/Q\spcheck(x)\right)_{w_c}$  is as claimed. Since $n|g_{\ov
a}$ for all $\ov a\in
\widetilde \Delta_c-\widetilde I_c(x)$, we see that
$\zeta\in Q\spcheck$. The two expressions for $\zeta$ are equal since, in
$Q\spcheck$,we have the relation $\sum_{\ov a\in \widetilde
  \Delta}g_{\ov a}\pi( a\spcheck)=0$.
\end{proof}

\begin{corollary}\label{Twc}
The order of $\pi_0(T^{w_c})$ is $n_0 = \gcd\{g_{\ov a}: \ov a \in
\widetilde \Delta_c\}$ and, for every $c$-pair $(x,y)$ in normal
form,  the map from
$\pi_0(T^{w_c})$ to
$\pi_0(Z(x,y))$ is an injection onto the cyclic subgroup of
$\pi_0(Z(x,y))$ of elements whose order divides $n_0$.
\end{corollary}

\begin{corollary}\label{torsion2}
If $(x,y)$ is a $c$-pair in normal form, then $\pi_0(Z(x,y))$ is a 
cyclic group of  order  dividing ${\rm gcd}\{\, g_{\ov a}:\ov
a\in\widetilde 
\Delta_c-\widetilde I_c(x)\,\}$.
\end{corollary}

\subsubsection{The case when all the $g_{\ov a}$ are equal}

There is one case we need to single out for special consideration.

\begin{proposition}\label{equalg}
Suppose that all the  $g_{\ov a}$ are equal,  say to $n$.
Then for every $c$-pair $(x,y)$ in normal form the group $Z(x,y)$ has
$n$ components. Let $T^{w_c}$ be the fixed points of the $w_c$-action
on $T$. The inclusion $T^{w_c}\to Z(x,y)$ induces a  bijection on the group
of components. 
\end{proposition}

\begin{proof}
Fix a $c$-pair $(x,y)$ in normal form.
Suppose that $\tilde x'\in A^c$
exponentiates to a regular element  $x'\in T$. Then $Z(x')=T$,
$Q\spcheck(x') = \{0\}$, and
$Z(x',y)=T^{w_c}$. 
According to Lemma~\ref{torsion1}, 
since $g_{\ov a}=n$ for every $\ov a\in\widetilde \Delta$,
${\rm Tor}(Q\spcheck/Q\spcheck(x))_{w_c}$ is a cyclic group of order
$n$.
Applying this to $x'$ instead of $x$, we see that ${\rm
  Tor}(Q\spcheck)_{w_c}$, which is the group of components of $T^{w_c}$, is a
cyclic group of order
$n$.    
Furthermore, the natural map $Q\spcheck\to Q\spcheck/Q\spcheck(x)$
induces an isomorphism on the torsion subgroups of the
$w_c$-coinvariants. 
It follows from Corollary~\ref{pi0lambda} that 
the inclusion $T^{w_c}=Z(x',y)\to Z(x,y)$ induces a bijection on
components, and in particular that $Z(x,y)$ has $n$ components.
\end{proof}

\subsection{Variation of $\pi_0(Z(x,y))$ as $x$ varies}

\begin{proposition}\label{lambda}
Suppose that $\tilde x, \tilde x'$ are points of $A^c$.
Let $x$, resp.\ $x'$, be the image of $\tilde
x$, resp.\ $\tilde x'$, under the exponential mapping.
Suppose $Z(x')\subseteq Z(x)$. Then the map on fundamental groups  
$Q\spcheck/Q\spcheck(x')\to Q\spcheck/Q\spcheck(x)$ induced by the 
inclusion
$Z(x')\subseteq Z(x)$
descends to a map on coinvariants which gives an 
injection
$${\rm Tor}\left((Q\spcheck/Q\spcheck(x'))_{w_c}\right)\subseteq {\rm
Tor}\left((Q\spcheck/Q\spcheck(x))_{w_c}\right).$$
\end{proposition}

\begin{proof}
The groups $Z(x)$ and $Z(x')$ are connected and both have $T$ as a
maximal torus. Since $Z(x')\subseteq Z(x)$, 
$\widetilde I(x')\subseteq \widetilde I(x)$, and hence
$Q\spcheck(x')\subseteq Q\spcheck(x)$.
The natural surjection $p\colon Q\spcheck/Q\spcheck(x')\to
Q\spcheck/Q\spcheck(x)$ is the map on fundamental groups induced by the
inclusion $Z(x')\subseteq Z(x)$.
It follows immediately from the description of the torsion subgroup
and its generator given in  Lemma~\ref{torsion1} that, if ${\rm
Tor}((Q\spcheck/Q\spcheck(x))_{w_c})$ is cyclic of order $n$ and 
${\rm Tor}((Q\spcheck/Q\spcheck(x'))_{w_c})$ is cyclic of order
$n'$, then $n'|n$ and that the map
$p_{w_c}\colon (Q\spcheck/Q\spcheck(x'))_{w_c}\to
(Q\spcheck/Q\spcheck(x))_{w_c}$ induced by $p$ on coinvariants
sends a generator of ${\rm Tor}((Q\spcheck/Q\spcheck(x'))_{w_c})$
to $n/n'$ times a generator of ${\rm
Tor}((Q\spcheck/Q\spcheck(x))_{w_c})$. Thus it is injective on the
torsion subgroup. 
\end{proof}

\begin{corollary}\label{pi0inj}
Suppose that $(x,y)$ and $(x',y')$ are 
$c$-pairs in normal form.
Suppose
$Z(x')\subseteq Z(x)$ and that $Z(x',y')\subseteq Z(x,y)$.
Then the inclusion $Z(x',y')\to
Z(x,y)$ induces an injective homomorphism $\pi_0(Z(x',y'))\to
\pi_0(Z(x,y))$. 
\end{corollary}

\begin{proof}
This is immediate from Proposition~\ref{lambda} and
Corollary~\ref{pi0lambda}.
\end{proof}

\section{$C$-triples}\label{ctripsect}

Let $C=(c_{ij})_{1\le i,j\le 3}$ be an antisymmetric matrix of
elements of ${\cal C}G$. Let $\langle C\rangle\subseteq {\cal C}G$
be the subgroup generated by all the entries $c_{ij}$ of $C$.
We wish to study the moduli space of conjugacy classes of ordered
$C$-triples in $G$. Quite generally, suppose that ${\bf x} =
(x_1,x_2,x_3)$ is an ordered triple in $G^3$. For an integral
$3\times 3$ matrix $M$, we set ${\bf x}^M =(x_1', x_2', x_3')$,
where $x_i' = x_1^{m_{1i}} x_2^{m_{2i}} x_3^{m_{3i}}$. We have
the following straightforward lemma:

\begin{lemma} Suppose that $M$ is unimodular. The map ${\bf x}
\mapsto {\bf x}^M$ defines an homeomorphism from the moduli space
of
$C$-triples to the moduli space of $C'$-triples, where $C' =
(c_{ij}')_{1\le i,j\le 3}$ and $c_{ij}'=
\prod_{r,s}c_{rs}^{m_{ri}\cdot m_{sj}}$. Moreover, for an
appropriate choice of a unimodular $M$, either $\langle C
\rangle = \langle C'\rangle$ is  cyclic and we can assume
that
$c'_{13}=c'_{23} = 1$ or  $G= Spin (4n)$ and $\langle C
\rangle = \langle C'\rangle$ is not cyclic, and we can assume
that  $c_{12}',
c_{23}', c_{13}'$ are the three nontrivial elements of $\langle C
\rangle$.
\end{lemma}

In case $\langle C
\rangle$ is cyclic, we assume that $c'_{13}=c'_{23} = 1$ and set
$c'_{12} =c$. In other words, $(x,y)$ is a $c$-pair and $z\in
Z(x,y)$. We shall call such a triple a {\sl $c$-triple}. For the
rest of this section, we shall study $c$-triples, postponing the
remaining case until later. We let ${\cal T}_G(c)$ denote the
moduli space of conjugacy classes of $c$-triples in $G$.

\subsection{$c$-triples of rank zero}

\subsubsection{The order of a $c$-triple}

\begin{defn}
We define the {\sl order} of a $c$-triple ${\bf x} = (x,y,z)$ to
be  the order of $[z]$ in $\pi_0(Z(x,y))$. Clearly, the order is a
conjugacy class invariant.
\end{defn}

It is easy to check, using the results of this section, that this
definition coincides with the previous definition given in
Section~\ref{commtripsect} in case
$c=1$. It follows from Corollary~\ref{constant} below that
the order is constant on connected components of ${\cal T}_G(c)$.
Thus we define the {\sl order\/} of a component $X$ of ${\cal
T}_G(c)$ to be the order of ${\bf x}$, where ${\bf x}$ is any
$c$-triple whose conjugacy class lies in $X$.

\subsubsection{The $w_c$-action on $\widetilde \Delta$}

\begin{lemma}\label{c'}
Suppose that   $(x,y,z)$ is a $c$-triple of
order $k$ and rank zero such that the $c$-pair $(x,y)$ is in
normal form.   Then:
\begin{enumerate}
\item $Z^0(x,y)$ is the torus $S^{w_c}$ and the action of $z$ on
this torus has isolated fixed points. 
\item If $k=1$,
then $(x,y)$ is a rank zero $c$-pair in $G$, which  is of type
$A_n$ for some $n$, and $c$ generates the center of $G$.
\item A $c$-triple $(x,y,z')$ is conjugate to $(x,y,z)$ if
  and only if $z$ and $z'$ lie in the same component of $Z(x,y)$. 
\end{enumerate}
\end{lemma}

\begin{proof}
The fact that $(x,y,z)$ is of rank zero means that $Z(x,y)^z$ is
finite.
If $Z^0(x,y)$ is not a torus, then, by \cite{deS} II \S 2,
$\left(Z^0(x,y)\right)^z\subseteq Z(x,y,z)$ is of positive dimension.
Thus $Z^0(x,y)$ is a torus, and since $S^{w_c}$ is a maximal
torus of $Z^0(x,y)$ they are equal.

If $k=1$, then $z\in Z^0(x,y)$ which by Part 1 is a torus.
Thus, the action of $z$ on $Z^0(x,y)$ is trivial. By Part 1
this implies that $Z^0(x,y)$ is a point, which means that $(x,y)$
is a $c$-pair of rank zero.  The statements about $G$ and $c$ now
follow from Proposition~\ref{rank0pair}.

By Part 1 and Corollary~\ref{torsion2},
$Z(x,y)$ is an extension of a cyclic group by a torus.  Thus
$(x,y,z')$ is conjugate to $(x,y,z)$ only if $z'$ is in the same
component as $z$. The converse follows easily from the fact  that
$z$ acts on $S^{w_c}=Z^0(x,y)$ with isolated fixed points.
\end{proof}

Now let us eliminate one exceptional case.

\begin{lemma}\label{excpt}
If there is a $c$-triple of rank zero in $G$, then the action of $w_c$
on $\widetilde \Delta$ has no exceptional orbits as defined in
Section~\ref{exceptional}.
\end{lemma}

\begin{proof}
Suppose that $w_c$ acts on $\widetilde \Delta$ interchanging
two roots 
$a\spcheck_1,a\spcheck_2$ which are not orthogonal.
The extended Dynkin diagram  $\widetilde D\spcheck(G)$ is thus
symmetric about  the bond connecting the nodes corresponding to $a_1$
and $a_2$. This symmetry implies that there are three possible types
of diagrams to consider:
(i) $\widetilde D$ is simply 
laced and has no trivalent vertices; (ii) $\widetilde D$ is simply
laced and has two trivalent vertices; and (iii) $\widetilde D$ has two
double bonds. In the first case $G$ is isomorphic to $SU(2n+1)$ for
some $n\ge 1$. This case is ruled out since the center of 
$SU(2n+1)$ acts freely on the extended Dynkin diagram. In the
second case, the subdiagram of $\widetilde D$ which contains the
chain connecting the trivalent vertices together with all nodes
adjacent to the trivalent vertices is the extended Dynkin diagram
for $D_{2n+1}$, and hence there is a non-trivial linear relation
between the elements of $\widetilde \Delta$ corresponding to the
nodes of this subdiagram.  
Since any proper subset of $\widetilde\Delta$ is linearly independent,
it follows that  this is the entire extended diagram for $G$. Thus,
 $G$ is of type $D_{2n+1}$ for some $n\ge 2$. Direct examination of
the action of the center in this case shows that  $c$ is a
generator of the center and that  the integers $\{g_{\ov a}\}$ are
all equal to $4$. In the third case, $\widetilde D$ has no trivalent
vertices and hence is a chain. Since the highest root is a long
root and since according to Lemma~\ref{root} it must be at one end of
the chain, it  
follows that the $a_1$ and $a_2$ are short roots. The subchain that
contains $a_1$ and $a_2$ and contains one long root on each side of
$a_1$ is the extended diagram for $C_{2n+1}$ for some $n\ge 1$. As
before, since there is a nontrivial linear relation between these
roots, it follows that $G$ is of type $C_{2n+1}$. Direct inspection
shows that the $\{g_{\ov a}\}$ are all equal to $2$ in this case.

What we have seen is that if $w_c$ acts on $\widetilde \Delta$ with an
exceptional orbit, then the integers $\{g_{\ov a}\}$ are all equal and
the number of orbits is at least $2$.

Now suppose that $(x,y,z)$ is a $c$-triple of rank zero in $G$.
Then we know that $Z^0(x,y)=S^{w_c}$.
Since all the $g_{\ov a}$ are equal, by
Proposition~\ref{equalg}  
$Z(x,y)=T^{w_c}$. Consequently, the action of any 
$z\in Z(x,y)$ on $S^{w_c}$ is trivial. That is to say the rank of
$Z(x,y,z)$ is equal to that of $Z(x,y)$.
But the rank of $Z(x,y)$ is one less than the number of orbits
and hence is positive. This contradicts the fact that $(x,y,z)$ is
of rank  zero.
\end{proof}

\subsubsection{The divisibility of the $g_{\ov a}$}

Suppose that  $(x,y,z)$ is a $c$-triple of order $k$ and rank
zero and that the $c$-pair $(x,y)$ is in normal form. Recall that, 
by Lemma~\ref{torsion1}, 
${\rm Tor}\left((Q\spcheck/Q\spcheck(x))_{w_c}\right)$ is a cyclic
group of order $n={\rm gcd}\{\, g_{\ov a}:\ov a\in\widetilde 
\Delta_c-\widetilde I_c(x)\,\}$, and is generated by 
$$\zeta=\frac{1}{n}\sum_{\ov
a\in \widetilde \Delta_c-\widetilde I_c(x)}g_{\ov a}\pi(
a\spcheck) =-\sum_{\ov a\in \widetilde I_c(x)}\frac{g_{\ov
a}}{n}\pi(a\spcheck)\in \pi(Q\spcheck).$$
Let $\Phi(x)\subseteq \Phi$ be the set of roots annihilating $x$
and let $\frak d$ be the subspace of $\frak t$ spanned by the
coroots inverse to $\Phi(x)$. Then $w_c$ acts on $\frak
d$ and on $\Phi(x)$. Let
$\Phi_c(x)$ be the root system
$\Phi^{\rm proj}_{Z(x)}(w_c)$ as defined in the previous section.
Recall that the coroots inverse to the simple roots of
$\Phi_c(x)$ are given by 
$\{\, \pi(a\spcheck): \ov a \in \widetilde I_c(x)\,\}$. Fix an
alcove  $B$ containing the origin in $\frak d^{w_c}$ for the root
system $\Phi_c(x)$.

 We can write
$y=sy_1$ where
$y_1$ is special and
$s\in S^{w_c}$. Let $s =\exp \tilde s$, where $\tilde s \in \frak
t^{w_c}$. After a further conjugation, we can assume that $\tilde
s\in \frak
t^{w_c}$  projects to $\hat s\in \frak d^{w_c}$ lying in the
alcove $B$. The element
$z$ defines a class
$[z] \in \pi_0(Z(x,y))$. Let
$\mu _z$ be the image of $\delta([z])$ in
${\rm Tor}\left((Q\spcheck/Q\spcheck(x))_{w_c}\right)/{\rm
Tor}(Q\spcheck_{w_c})\subseteq {\cal C}\Phi_c(x)$, where $\delta$
is the homomorphism of Lemma~\ref{pi0lambda}. By
Lemma~\ref{ratreps}, we may lift
$\mu_z$ to an element $\tilde \mu_z$ of
$(Q\spcheck(x)_{w_c})\otimes {\bf Q}$, well-defined modulo
$Q\spcheck(x)_{w_c}$. Write 
\begin{equation}\label{mueqn}
\tilde \mu_z = \sum _{\ov a\in
\widetilde I_c(x)}r_{\ov a}\pi( a\spcheck).
\end{equation}

\begin{proposition}\label{c'1} 
Let $(x,y,z)$ be a $c$-triple of order $k$ and rank zero, and
suppose that the $c$-pair $(x,y)$ is in normal form. With
assumptions and  notation as in the previous paragraph,
\begin{enumerate}
\item The group $Z^0(x,y_1)$ is semisimple and hence $\frak
d^{w_c}= \frak t^{w_c}$.
\item The element $\tilde s=\hat s$ is a barycenter of the alcove
$B$.
\item In Equation~\ref{mueqn} no
coefficient $r_{\ov a}$ is integral. 
\item For every set of simple roots of $\Phi_c(x)$, the
coefficients of $\tilde \mu_z$, written as a linear combination
of the inverse coroots, are all non-integral.
\item The element $x$ is the image under the exponential map of a 
vertex of $A^c$ opposite to a face $\{\ov a =0\}$ for some $\ov a\in 
\tilde I_c(x)$.
\item The order $k$ divides $g_{\ov a}$.
\end{enumerate}
\end{proposition}

\begin{proof} 
Since $(x,y,z)$ is of rank zero, it follows from
Lemma~\ref{c'} that $Z^0(x,y) = S^{w_c}$ and that conjugation
by $z$ normalizes $S^{w_c}$ and its action on this torus
has isolated fixed points.
Since $Z^0(x,y)$ is abelian,
its Weyl group is trivial.
Thus, by Lemma~\ref{compgrp2}
there is a well-defined  action of the element
$[z]\in \pi_0(Z(x,y))$   on 
$\frak t^\sigma$, and, by  Proposition~\ref{CD}, it is
given as the Weyl element
$w=\nu_B(\mu_z)$. The element
$w$ preserves $\frak t^{w_c}$ and fixes only the origin there.
By Proposition~\ref{indivis}, $w$ is an element of  the Weyl group
of
$Z^0(x,y_1)$. It follows that $Z^0(x,y_1)$ has  finite center and
thus is semisimple. This proves (1).

Since $w=\nu_B(\mu_z)$, the affine action of
$\mu_z\in{\cal C}\Phi_c(x)$ on the
alcove $B$ has isolated fixed points, and hence fixes only the 
barycenter of $B$. Since $\mu_z$ fixes $\tilde s = \hat s$, by
Proposition~\ref{fixes}, it follows that
$\tilde s$ is the barycenter of $B$, proving Part (2).
Moreover, by
Proposition~\ref{finite}, no coefficient $r_{\ov
a}$ of  $\tilde \mu_z\in
(Q\spcheck(x)_{w_c})\otimes {\bf Q}$ is integral, proving (3).

Any two sets of simple roots of $\Phi_c(x)$ are related by an
element of the Weyl group of $\Phi_c(x)$. Since the Weyl group
acts trivially on the center, and since $\tilde \mu_z$ projects
to an element of the center,  Part (4) 
follows.

Since
$Z^0(x,y_1)$ is semisimple, 
$\widetilde I_c(x)$ has cardinality equal to the dimension of
$\frak t^{w_c}$. Thus $x$ is the image under the exponential map
of a vertex $\tilde x$ of
$A^c$, opposite the face
$\{\ov a =0\}$, say. By Lemma~\ref{torsion1}, ${\rm
Tor}((Q\spcheck/Q\spcheck(x))_{w_c})$ is cyclic of order 
$g_{\ov a}$. Since the class of $z$ is an element of order $k$ in
$\pi_0(Z(x,y))$, and since $\delta$ is injective by
Lemma~\ref{pi0lambda},
$\delta([z])$ is an element of order $k$ in ${\rm
Tor}((Q\spcheck/Q\spcheck(x))_{w_c})$. Consequently, $k|g_{\ov
a}$.
\end{proof}

To see the relationship between $k$ and the remaining $g_{\ov
b}$, we use the following lemma:

\begin{lemma}\label{knotdivide}
Let $x\in A^c$. Suppose that $\mu \in
(Q\spcheck/Q\spcheck(x))_{w_c}$ has order $k$. Let
$\tilde \mu\in
(Q\spcheck(x)_{w_c})\otimes {\bf Q}$ be an element lifting the
image $\ov \mu$ of $\mu$ in ${\rm
Tor}((Q\spcheck/Q\spcheck(x))_{w_c})/{\rm
Tor}(Q\spcheck_{w_c})\subseteq{\cal C}\Phi_c(x)$.
 Then  no coefficient of $\tilde
\mu$ is integral if and only if, for all   $\ov b
\in \widetilde I_c(x)$, $k\not|g_{\ov b}$. 
\end{lemma}

\begin{proof} By Lemma~\ref{torsion1}, a generator for
$(Q\spcheck/Q\spcheck(x))_{w_c}$ is $\zeta =
-\sum_{\ov b\in \widetilde I_c(x)}(g_{\ov b}/n)\pi(
b\spcheck)$, where $n={\rm gcd}\{\, g_{\ov a}:\ov a\in\widetilde 
\Delta_c-\widetilde I_c(x)\,\}$. Thus $k|n$ and we can write 
$\mu = \ell
\zeta$ for some integer $\ell$ of the form $tn/k$, where $t$ and
$k$ are relatively prime. A representative for $\tilde \mu$ is
then given by 
$$\ell \zeta = -\sum_{\ov b\in \widetilde I_c(x)}\frac{tg_{\ov
b}}{k}\pi( b\spcheck).$$
Clearly, then, no coefficient of $\tilde
\mu$ is integral if and only if, for all   $\ov b
\in \widetilde I_c(x)$, $k\not|g_{\ov b}$. 
\end{proof}

\begin{corollary}\label{onlyif} Let $(x,y,z)$ be a rank zero
$c$-triple of order $k$. Then there is a unique $\ov a\in
\widetilde \Delta _c$ such that 
$k|g_{\ov a}$.
\end{corollary} 
\begin{proof} This is immediate from the previous lemma and Parts
(3) and (6) of the previous proposition.
\end{proof}

\subsubsection{Classification of rank zero $c$-triples}

We can now give the classification of rank zero $c$-triples.

\begin{proposition}\label{k} Let $G$ be simple.
\begin{enumerate}
\item There is a $c$-triple of rank zero and order $k$ in $G$ if
and only if
$k$ divides  exactly one of the $g_{\ov a}$.
\item Suppose that $g_{\ov a}$ does not divide $g_{\ov b}$ for
$\ov b \neq \ov a$. Let   $(x,y,z)$ be a
$c$-triple  of rank zero and order $g_{\ov a}$ and let
$(x',y',z')$ be a rank zero $c$-triple whose order $k$ divides
$g_{\ov a}$.  Then
$(x',y',z')$ is conjugate to $(x,y,z^\ell)$ for a unique $\ell$
such that $1\leq \ell < g_{\ov a}$ and moreover $k= g_{\ov
a}/\gcd(\ell, g_{\ov a})$. 
\end{enumerate}
\end{proposition}

\begin{proof}
The only if part of the first statement is Corollary~\ref{onlyif}.

Conversely, suppose that exactly one of the $g_{\ov a}$ is
divisible by
$k$. Choose $\tilde x\in A^c$ so that, for $x =\exp \tilde x$, we
have 
$\widetilde I_c(x) =
\widetilde
\Delta_c-\{\ov a\}$. By Lemma~\ref{torsion1}, the order of ${\rm
  Tor}\left((Q\spcheck/Q\spcheck(x))_{w_c}\right)$ is
divisible by $g_{\ov a}$ and hence is divisible by $k$. Let
$y_1$ be such that
$(x,y_1)$ is a special 
$c$-pair in normal form. Let $\tilde s\in {\frak t}^{w_c}$  be the
barycenter of the alcove $B\subseteq {\frak t}^{w_c}$ for the root
system $\Phi_c(x)$ and let 
$s={\rm exp}(\tilde s)$. We set $y=sy_1$. Fix an element
$\mu\in{\rm
  Tor}\left((Q\spcheck/Q\spcheck(x))_{w_c}\right)$  of order $k$.
Since $\tilde s$ is the barycenter of $B$, $\mu$ fixes $\tilde
s$.  Hence by Corollary~\ref{fullstatement}, there is a
$z\in Z(x,y)$ whose image in ${\rm
  Tor}\left((Q\spcheck/Q\spcheck(x))_{w_c}\right)$ under the map
$\pi_0(Z(x,y)) \to {\rm
  Tor}\left((Q\spcheck/Q\spcheck(x))_{w_c}\right)$ is
$\mu$.  Since $\tilde s$  is a regular element for $\Phi_c(x)$, it
follows that
$Z^0(x,y)=S^{w_c}$. Lifting the image of $\mu$ in ${\rm
  Tor}\left((Q\spcheck/Q\spcheck(x))_{w_c}\right)/{\rm
  Tor}(Q\spcheck)_{w_c}$ to an element $\tilde \mu\in
(Q\spcheck(x))\otimes {\bf Q}$, $\tilde \mu=\sum s_{\ov a}\pi(a\spcheck)$.
It follows by construction that no
coefficient $s_{\ov a}$ of a $\pi( a\spcheck) \in \widetilde
I_c(x)$ is integral. Hence, by 
Proposition~\ref{finite}, 
$\mu$ acts with isolated fixed points on
$S^{w_c}$.   It follows that the conjugation
action of
$z$ on
$S^{w_c}$ also has isolated fixed points, and hence
$(x,y,z)$ is of rank zero. Clearly, it is of order $k$.

Now suppose that $g_{\ov a}$ does not divide $g_{\ov b}$ for
$\ov b \neq \ov a$. Let   $(x,y,z)$ be a
$c$-triple  of rank zero and order $g_{\ov a}$ and let
$(x',y',z')$ be a rank zero $c$-triple whose order $k$ divides
$g_{\ov a}$.
After conjugation 
we can assume that $(x,y)$ and $(x',y')$ are $c$-pairs in normal form.
Then $y'\in S^{w_c}\cdot y$.
Each of $x$ and $x'$ is the image
under the exponential mapping of the
vertex of $A^c$ opposite the face $\{\ov a=0\}$ for the unique $\ov
a\in\widetilde \Delta_c$ for which $k|g_{\ov a}$.
Thus, $x=x'$.
 
By Lemma~\ref{sigma1} there is  $y_1\in S^{w_c}\cdot y$ such that
$(x,y_1)$ is a special $c$-pair in normal form. We write
$y=sy_1$ and $y'=s'y_1$ for elements $s,s'\in S^{w_c}$ which are the
images of $\tilde s$ and $\tilde s'$ in ${\frak t}^{w_c}$.
It follows from Proposition~\ref{c'1} that $\tilde s$ and $\tilde
s'$ are barycenters for the alcove decomposition of ${\frak
t}^{w_c}$ associated with the root system $\Phi_c(x)$.

Thus, there is an element in the Weyl group of $\Phi_c(x)$ on
${\frak  t}^{w_c}$ which conjugates $\tilde s$ to $\tilde s'$.
By Proposition~\ref{indivis}, the Weyl group of $\Phi_c(x)$ is
the Weyl group of
$Z(x, y_1)$. Thus, there exists a $g\in N_T(Z(x,y_1))$
 conjugating $\tilde s$ to $\tilde s'$, and
hence
$(x,sy_1)$ and $(x,s'y_1)$ are conjugate $c$-pairs in $G$.

This allows us to assume further that $y=y'$. By
Proposition~\ref{c'},
$Z^0(x,y)=S^{w_c}$. By Corollary~\ref{torsion2}, the group
$\pi_0(Z(x,y))$ is   cyclic of order
dividing $g_{\ov a}$, and $[z]$ is an element of this group of
order exactly $g_{\ov a}$. Hence  $[z]$ generates
$\pi_0(Z(x,y))$. Thus there is a unique integer
$\ell$ with $1\leq \ell < g_{\ov a}$ such that $z^\ell$ and $z'$
are in the same component of
$Z(x,y)$. Since $z'$ acts on $S^{w_c}$ with isolated fixed
points, this implies that $z^\ell$ and $z'$ are conjugate in
$Z(x,y)$ and hence that $(x,y,z^\ell)$ and $(x,y,z')$ are
conjugate in $G$.
\end{proof}

\subsubsection{Simple groups containing $c$-triples of rank
  zero}\label{listofrankzero}

It follows from the above that a simple group $G$ 
has a
$c$-triple of rank zero and order $k$ if and only if
exactly one of the integers
$g_{\ov a}$  is divisible by $k$. 
Examining the quotient diagrams gives the following list of all
the possibilites for
$c\neq 1$:

\begin{enumerate}
\item $G=A_n$, $c$ a generator of ${\cal C}G$ and $k|n+1$.
\item $G=C_2$, $c$ the non-trivial element of ${\cal C}G$ and 
$k=2$.
\item $G=D_6$, $c$ an exotic element of the center and $k=4$.
\item $G=E_6$, $c$ a generator of ${\cal C}G$ and $k=2$ or $6$.
\item $G=E_7$, $c$ the non-trivial element of ${\cal C}G$ and
$k=3$ or $6$.
\end{enumerate}

\subsubsection{Action of the automorphism group of $G$ on rank 
zero $c$-triples}

\begin{lemma}\label{automorphism}
Let $(x,y)$ be a $c$-pair in normal form in $G$.
Let $g\in G$, let $i_g$ denote conjugation by $g$, and suppose
that $i_g(x) =\zeta x$ for some $\zeta \in {\cal C}G$. Let
$\sigma\colon G\to G$ be an automorphism  normalizing $S^{w_c}$
such that $\sigma(x) =
\zeta_1x$ and
$\sigma (y) =\zeta_2y$, where $\zeta_1, \zeta_2\in {\cal C}G$.
Then both
$i_g$ and
$\sigma$ induce the identity on  ${\rm
Tor}\left((Q\spcheck/Q\spcheck(x))_{w_c}\right)$. Moreover,
$\sigma$ induces the identity on $\pi_0(Z(x,y))$.
\end{lemma}

\begin{proof} 
Since $T$ is a maximal torus for $Z(x)$, there is an inner
automorphism $i_h$,   $h\in Z(x)$, such that $i_h\circ i_g$
normalizes $T$. 
 Since
$Z(x)$ is connected, the inner automorphism $i_h$ of $Z(x)$
induces the identity on its fundamental group.
Thus, without loss of
generality, we can assume that $g$ normalizes $T$.
Let $\tilde x\in {\frak t}$ be the lift of $x$ lying in the alcove
$A$. 
Then $g\tilde x g^{-1}=\tilde x+\tilde \zeta$ for some element $\tilde
\zeta\in{\frak t}$ lifting $\zeta$.
Since $\zeta\in{\cal C}G$, the affine automorphism $t\mapsto
\varphi(t)= gtg^{-1}+\tilde\zeta $ of ${\frak t}$ normalizes the
alcove decomposition associated with the root system of $G$.
Since $\varphi(\tilde x)=\tilde x$, $\varphi(A)=A'$ is an alcove
containing $\tilde x$. Hence  there is an element
$\mu$  of the affine Weyl group of $G$ which fixes $\tilde x$ and
sends $A'$ to $A$. The composition $\mu\circ \varphi$ then 
normalizes
$A$. Since the translational part of this affine transformation 
is congruent to $\tilde \zeta$ modulo the coroot lattice
$Q\spcheck$, we see by Lemma~\ref{char} that this composition is
$\nu_A(\zeta)$, which we denote by $w_\zeta$.
On the other hand, since $\mu$ fixes $\tilde x$, the Weyl part 
$w$ of
$\mu$ is an element of the Weyl group of $Z(x)$. Since $w\circ
i_g|{\frak t}=w_\zeta$, multiplying 
$g$  by an element of $N_{Z(x)}(T)$ allows us to assume, without
loss of generality that $i_g|{\frak t}=w_\zeta$.

As noted in Section~\ref{action1},  we have $g_{w_\zeta(a)}=g_{a}$
for all $a\in \widetilde\Delta$. Since  $w_\zeta$ and
$w_c$ commute, we also have
$g_{w_\zeta(\ov a)}=g_{\ov a}$ for all $\ov a\in
\widetilde\Delta_c$. Now we see directly from the 
expression for a generator of ${\rm
Tor}\left((Q\spcheck/Q\spcheck(x))_{w_c}\right)$ given in
Lemma~\ref{torsion1}, that  $w_\zeta$ acts trivially on ${\rm
Tor}\left((Q\spcheck/Q\spcheck(x))_{w_c}\right)$. 
This completes the proof of the first part of the lemma.

Since $S^{w_c}\cdot x$ contains a generic element of $G$, $\sigma$
normalizes $T$. The equation $\sigma (y) =\zeta_2y$ says that
$\sigma|{\frak t}$ commutes with $i_y=w_c$.
Let $\tilde x\in A^c$ be the lift of
$x$. There is a Weyl element $w\in W(G)$ commuting with $w_c$
such that $w(x) =\zeta_1^{-1}x$ and such that $w\circ\sigma(A)=A$. 
Since $w\circ \sigma(A)=A$, it follows that $g_{w\circ
\sigma(a)}=g_{a}$ for all $a\in \widetilde\Delta$. Since $w\circ 
\sigma$ commutes with $w_c$, it is also the case that
$g_{w\circ\sigma(\ov a)}=g_{\ov a}$ for all $\ov a\in
\widetilde\Delta_c$. Thus, $w\circ\sigma$ acts trivially on ${\rm
Tor}\left((Q\spcheck/Q\spcheck(x))_{w_c}\right)$.  The first part
of the lemma implies that $w$ acts trivially on this group.
Thus, it  follows that $\sigma$ acts trivially on ${\rm
Tor}\left((Q\spcheck/Q\spcheck(x))_{w_c}\right)$.  The final
statement follows since the inclusion $\pi_0(Z(x,y)) \to {\rm
Tor}\left((Q\spcheck/Q\spcheck(x))_{w_c}\right)$ of
Lemma~\ref{pi0lambda} is equivariant with respect to the
automorphism $\sigma$.
\end{proof}

\begin{lemma}\label{outer}
If $(x,y,z)$ is a $c$-triple of rank zero in $G$ and if 
$\sigma\colon G\to G$  is an automorphism of $G$ fixing $c$, then
$\sigma(x,y,z)$ is conjugate  to $(x,y,z)$. 
\end{lemma}

\begin{proof}
We can assume that $(x,y)$ is a $c$-pair
in normal form. Since $\sigma(c)=c$, the triple
$(\sigma(x),\sigma(y),\sigma(z))$ is also a $c$-triple, clearly of
rank zero.
By Proposition~\ref{k}, composing $\sigma$ with an inner
automorphism allows us to assume that $\sigma(x,y,z)=(x,y,z^\ell)$ for
some integer $\ell$.
Since $(x,y,z)$ is of rank zero, by Lemma~\ref{c'},
$Z^0(x,y)=S^{w_c}$. Hence, 
$\sigma(S^{w_c})=S^{w_c}$. 
Applying Lemma~\ref{automorphism} we see that $\sigma$ acts 
trivially on
${\rm Tor}\left((Q\spcheck/Q\spcheck(x))_{w_c}\right)$.  Hence, by
Lemma~\ref{pi0inj}, it follows that $z'=\sigma(z)$ and $z$ are in 
the same component of $Z(x,y)$, and hence that $\sigma(x,y,z)$ is
conjugate to
$(x,y,z)$.   
\end{proof}

\subsubsection{Action of ${\cal C}G$ on the space of rank zero
$c$-triples}

We shall first consider the action on the first
component.

\begin{lemma}\label{x}
Let $\gamma\in {\cal C}G$. Then   $\gamma\cdot
(x,y,z)=(\gamma x,y,z)$ is the trivial action on the set of
conjugacy classes of $c$-triples of rank zero in $G$.
\end{lemma}

\begin{proof}
Note that $Z(\gamma x, y) =Z(x,y)$. It follows that the order of
$(\gamma x, y,z)$ is equal to that of $(x,y,z)$, and it 
clearly has rank zero. We can assume that
$(x,y)$ is in normal form with respect to $A$.  By 
Proposition~\ref{k} there is  $g\in G$ which conjugates
$(\gamma x,y)$ to $(x,y)$. Since $i_g(x) =\gamma^{-1}  x$ and
$i_g(y) = y$, $i_g$ induces an automorphism of $Z(x,y)$. By
Lemma~\ref{automorphism}, the induced action of $i_g$ on
$\pi_0(Z(x,y))$ is trivial. In particular, $i_g(z)$ and $z$ lie
in the same connected component of $Z(x,y)$. Hence they are
conjugate in $Z(x,y)$. It follows that $(\gamma x,y,z)$ and
$(x,y,z)$ are conjugate in $G$.
\end{proof}

\begin{lemma}\label{y}
For any $\gamma \in {\cal C}G$ the action 
$\gamma\cdot(x,y,z)=(x,\gamma y,z)$ is trivial on the space of
conjugacy classes of $c$-triples of rank zero in $G$.
\end{lemma}

\begin{proof}
By symmetry between $x$ and $y$ (at the expense of replacing $c$ by
$c^{-1}$), the result in this case follows from the previous one.
\end{proof}

Lastly, we must consider the action $\gamma \cdot
(x,y,z)=(x,y,\gamma z)$.

\begin{lemma}\label{z}
For $\gamma \in {\cal C}$ the action of ${\cal C}G$ defined by
$\gamma \cdot (x,y,z)=(x,y,\gamma z)$ induces an action of ${\cal
C}G$ on the set of conjugacy classes of rank zero
$c$-triples in $G$. The stabilizer ${\cal K}$ of a conjugacy
class of a
$c$-triple of rank  zero is ${\cal K} ={\cal C}G\cap S^{w_c}$.
The order of each orbit, or equivalently, the index of
${\cal K}$ in
${\cal C}G$
 is ${\rm gcd}(g_{\ov a})=n_0$.
\end{lemma}

\begin{proof}
Let $(x,y,z)$ be a $c$-triple of rank zero such that $(x,y)$ is a
$c$-pair in normal form. Then $Z^0(x,y)=S^{w_c}$. For any $\gamma
\in {\cal C}G$, $\gamma z$ and $z$ are in the same component
of
$Z(x,y)$ if and only if $\gamma \in S^{w_c}$.  This proves that
${\cal K} ={\cal C}G\cap S^{w_c}$.

Lastly, we show that the order of 
${\cal C}G/\left(S^{w_c}\cap {\cal C}G\right)$ is equal to
$n_0$. By Corollary~\ref{Twc},
$n_0$ is the number of components of $T^{w_c}$. Since ${\cal
C}G\subseteq T^{w_c}$ and since
$S^{w_c}$ is the component of the identity of
$T^{w_c}$, this shows that the order of ${\cal 
C}G/\left(S^{w_c}\cap {\cal C}G\right)$ divides $n_0$. 

To complete the proof, we need to show that the inclusion  ${\cal
  C}G\to T^{w_c}$ is onto on the level of components. We state
this as a separate lemma:

\begin{lemma}\label{surjpi0} Let $G$ be a simple group
and let $c\in {\cal C}G$.
The inclusion
${\cal C}G\to T^{w_c}$ induces a surjection ${\cal C}G\to
\pi_0(T^{w_c})$.
\end{lemma}

\proof Applying cohomology to the $w_c$-actions on
the exact sequence
$$0\to Q\spcheck\to P\spcheck\to {\cal C}G\to 0$$
and considering the torsion subgroups,
we have an exact sequence
$${\cal C}G\to {\rm Tor}(Q\spcheck_{w_c})\to {\rm
Tor}(P\spcheck_{w_c})
\to {\cal C}G.$$
Thus, it suffices to show that ${\rm Tor}(P\spcheck_{w_c})\to
{\cal C}G$ is injective. 

Set $h_{\ov a} = n_{\ov a}h_a$. There is the exact sequence
$$0 \to {\bf Z}(\sum _{a\in \widetilde \Delta}h_aa) \to
\bigoplus_{a\in \widetilde \Delta}{\bf Z}(a) \to Q\to 0,$$
where $Q=Q(\Phi)$ is the root lattice of $G$. Dualizing gives 
$$0 \to P\spcheck \to \bigoplus _{a\in \widetilde \Delta}{\bf
Z}a^* \to {\bf Z} \to 0,$$ 
where the second map is obtained by sending $\sum
_ar_aa^*$  to $\sum _ar_ah_a$. (Here the
$a^*$ are the dual basis to the basis $\{a: a\in
\widetilde \Delta\}$ of the free ${\bf Z}$-module
$\bigoplus_{a\in \widetilde \Delta}{\bf Z}(a)$.) Hence, the group
${\rm Tor}(P\spcheck_{w_c})$ is a cyclic group of order   ${\rm
gcd}(h_{\ov a})$.  Clearly, if $\ov a$ contains
the extended root,
$n_{\ov a}$ is the order of $c$ and $h_{\tilde a}=1$ so that
$h_{\ov a}$ is the order of $c$. Thus the order of ${\rm
Tor}(P\spcheck_{w_c})$ divides the order of $c$. On the other
hand, a generator for ${\rm Tor}(P\spcheck_{w_c})$  is
represented by 
$a^*-w_c(a^*)\in {\cal C}G$ for any
$a\in \widetilde \Delta$ for which $h_a=1$. Choose $a$ to be the
root mapped by $w_c$ to $\tilde a$. Then the corresponding
element  of
${\rm Tor}(P\spcheck_{w_c})$ is represented by
$\varpi_a\spcheck\in P\spcheck$, and its image in
${\cal C}G$ is therefore equal to
$c^{-1}$. Since the order of  ${\rm Tor}(P\spcheck_{w_c})$ is
at most that of
$c$, the map ${\rm Tor}(P\spcheck_{w_c})\to
{\cal C}G$ is injective, and in fact an isomorphism. This proves
the lemma.
\end{proof}

\subsection{The maximal torus of a $c$-triple of order $k$}

Fix a non-trivial element $c\in 
{\cal C}G$.
Let $(x,y,z)$ be a $c$-triple in $G$ of order $k$.
We let $S$ be a maximal torus for $Z(x,y,z)$ and let ${\frak s}$ be
its Lie algebra. As usual, $L=DZ(S)$.  By Theorem~\ref{main} 
there is a $c$-triple
$(x_0,y_0,z_0)\in L$ of rank zero and $s_1,s_2,s_3\in S$ such that
$(x,y,z)=(s_1x_0,s_2y_0,s_3z_0)$.

\subsubsection{Determination of ${\frak s}$}

Our  goal is to describe the torus $S$, or equivalently
$\frak s$. We begin with the following definition.

\begin{defn}\label{barg} Let $\ov{\bf g}\colon \widetilde
\Delta_c\to {\bf N}$ be the function defined by $\ov{\bf g}(\ov
a)=g_{\ov a}$. For each integer
$k\ge 1$ dividing at least one of the  integers
$g_{\ov a}$, for $\ov a\in \widetilde \Delta_c$, we define 
$f^c(k)$ to be the maximal face of $A^c$ with the
property that every root $\ov a\in \widetilde
\Delta_c$  for which $k\not|g_{\ov a}$ takes an integral value on 
$f^c(k)$. Let $\frak t^{w_c}(\ov{\bf g},k)$ be the linear subspace
of
${\frak t}^{w_c}$ parallel to $f^c(k)$ and let $S^{w_c}(\ov{\bf
g},k)$ be the torus whose Lie algebra is $\frak t^{w_c}(\ov{\bf
g},k)$.
\end{defn}

We can then state the main result as follows:

\begin{proposition}\label{parallel}
Let $(x,y,z)$ be a $c$-triple of order $k$ and let $S$ be a
maximal torus of $Z(x,y,z)$. Then
$k$ divides at least one of the integers $g_{\ov a}$ and $S$ is
conjugate to $S^{w_c}(\ov{\bf g},k)$. The element
$x$ is conjugate to a point of  
$f^c(k)$. In particular, the dimension of $S$ is equal to one
less than the number of $\ov a\in \widetilde \Delta$ for which
$k$ divides $g_{\ov a}$.
\end{proposition}

First let us show that
$\frak s$ is conjugate to a linear subspace parallel to some face
of $A^c$:

\begin{lemma}
Suppose that $(x,y,z)$ is a $c$-triple and that $(x,y)$ is a 
$c$-pair in normal form. Supoose that $x'\in S\cdot x$ is generic
in the sense that any root of $G$ vanishing on $x'$ vanishes
on $S\cdot x$ and suppose that $x'={\rm exp}(\tilde x')$
for some $\tilde x'$ contained in $A^c$. Let $f$ be the  face
of $A^c$ containing $\tilde x'$  in its interior. Then, up to
conjugation by an  element in $Z(x,y)$, ${\frak s}$ is the linear
subspace of
${\frak t}^{w_c}$ parallel to $f$.
\end{lemma}

\begin{proof}
Since $x'\in S\cdot x$ is generic, it follows that ${\rm dim}
\,{\frak s}
\le {\rm dim}\,f$.  
On the other hand, for every $\tilde x''$ in the interior of $f$,
the triple
$({\rm exp}(\tilde x''),y,z)$ is a $c$-triple with the same 
centralizer as the triple $(x,y,z)$. This means that $ f-\tilde
x'$ is contained in $\frak s$. Comparing this with the dimension
estimate shows that ${\frak s}$ is the linear subspace parallel
to $f$.
\end{proof}

\begin{lemma}\label{subsetfn}
Let $(x,y,z)$ be a $c$-triple of order $k$. Then $k$ divides at
least one of the $g_{\ov a}$ and $x$ is conjugate to the
exponential of an element of
$f^c(k)$. Moreover, the Lie algebra ${\frak s}$ of
$S$ is conjugate to the linear subspace of ${\frak
t}^{w_c}$ parallel to a face of
$A^c$ contained in  $f^c(k)$. 
\end{lemma}

\begin{proof}
After conjugation, we may assume that $(x,y)$ is a $c$-pair in
normal form. Since 
$(x,y,z)$ is of order
$k$, the order of 
$\pi_0(Z(x,y))$ is divisible by $k$, and hence
$\left(Q\spcheck/Q\spcheck(x)\right)_{w_c}$ has order divisible by
$k$. By Lemma~\ref{torsion1}, this means $n={\rm 
  gcd}_{a\in \widetilde  
\Delta-\widetilde I(x)}\{g_{\ov a}\}$ is divisible by $k$, and in
particular $k|g_{\ov a}$ for some $a$. Moreover $\{\ov
a: k\not|g_{\ov a}\} \subseteq\widetilde I(x)$. Thus,
$x$ is the image under the exponential mapping of a point $\tilde
x$ contained in  $f^c(k)$. 

Take  $x'\in S\cdot x$ to be a generic element
with the property that $x={\rm exp}(\tilde x')$ for some
$\tilde x'\in A^c$.
Let $f$ be the face of $A^c$ containing $\tilde x$ in its
interior. By the first part of this lemma, $f$ is a face of
$f^c(k)$. According to the previous lemma
$\frak s$ is parallel to $f$.
\end{proof}

We turn now to the proof of Proposition~\ref{parallel}.
As usual, let $L=DZ(S)$ and let
$(x_0,y_0,z_0)$ be a rank zero $c$-triple in $L$ such that
$(x,y,z) =(s_1x_0, s_2y_0, s_3z_0)$ for elements $s_1,s_2,s_3\in S$.
Let $\frak t_L={\rm Lie}(L)\cap \frak t$; it is the Lie algebra
of a maximal torus for
$L$. Let
$Q\spcheck_L=Q\spcheck\cap  {\frak t}_L$.
Of course, $c\in L$ and the action of $w_c$ on ${\frak t}$ 
normalizes
${\frak t}_L$.
The group $L$ is simply connected and semi-simple, but may not be
simple.
Let $\prod_iL_i$ be its decomposition into simple factors and let
$\Phi_{L_i}$ be the corresponding root system with respect to
${\rm Lie}(L_i)\cap \frak t$. We decompose
$c\in L$ as
$\prod_ic_i$. Write
$x_0 =\prod_i x_{0,i}$, and similarly for $y_0$ and $z_0$. Then
$(x_{0,i},y_{0,i},z_{0,i})$ is a rank zero
$c_i$-triple in $L_i$ for every $i$. Of course, ${\rm
Tor}\left((Q\spcheck_L/Q\spcheck_L(x_0))_{w_c}\right)=
\bigoplus _i{\rm Tor}
\left((Q\spcheck_{L_i}/Q\spcheck_{L_i}(x_{0,i}))_{w_{c_i}}\right)$.
By Lemma~\ref{torsion1} applied to the
$(x_{0,i},y_{0,i})$, the groups
${\rm Tor}
\left((Q\spcheck_{L_i}/Q\spcheck_{L_i}(x_{0,i}))_{w_{c_i}}\right)$
are cyclic.  Let $\mu_i\in {\rm Tor}
\left((Q\spcheck_{L_i}/Q\spcheck_{L_i}(x_{0,i}))_{w_{c_i}}\right)$
be the element
$\delta([z_{0,i}])$, where $[z_{0,i}]\in\pi_0(Z(x_{0,i},y_{0,i}))$
is the class of $z_{0,i}$ and $\delta$ is the homomorphism of
Lemma~\ref{pi0lambda}. By Lemma~\ref{ratreps}, there is 
is an element $\tilde \mu_i\in
(Q\spcheck_{L_i}(x_{0,i})_{w_{c_i}})
\otimes {\bf Q}$ whose image in ${\rm Tor}
\left((Q\spcheck_{L_i}/Q\spcheck_{L_i}(x_{0,i}))_{w_{c_i}}\right)/
{\rm Tor}(( Q_{L_i}\spcheck)_{w_{c_i}})$ is the image of $\mu_i$
under the quotient map.  According to Proposition~\ref{c'1}, the
coefficients of 
$\tilde\mu_i$ are all non-integral when written as a linear
combination of any set of simple coroots for the root system
$\Phi_{Z_{L_i}(x_{0,i})}^{\rm proj}(w_c)$. Clearly
$$\delta([z_0]) = \sum_i\mu_i
\in \bigoplus_i
\left((Q\spcheck_{L_i}/Q\spcheck_{L_i}(x_{0,i}))_{w_c}\right)=
\left((Q\spcheck_L/Q\spcheck_L(x_0))_{w_c}\right).$$ 
Let $\mu = \sum _i\mu _i$ and let $\tilde \mu = \sum _i\tilde
\mu_i \in  (Q\spcheck_L(x_0)_{w_c})\otimes {\bf Q}$. Then $\tilde
\mu$ projects to the image of $\mu$ in ${\rm Tor}
\left((Q\spcheck_{L}/Q\spcheck_L(x_0))_{w_c}\right)/
{\rm Tor}(( Q_L\spcheck)_{w_c})$. For every set of simple roots
for $\Phi_{Z_L(x_0)}^{\rm proj}(w_c)$, if we write $\tilde \mu$ as
a linear combination of the corresponding  coroots, then no
coefficient  is integral.

Fix a generic element $x'\in S\cdot x_0$. Then 
$Z_G(x')=S\cdot Z_L(x_0)\subseteq Z_G(x)$ and so $\Phi(x') =
\Phi_L(x_0)$. Hence $Q\spcheck_L(x_0)=Q\spcheck(x')$. Moreover
$Z(x',y)\subseteq Z(x,y)$. Furthermore,  $Z_G(x',y_0)=Z_G(x',y)
=S\cdot Z_L(x_0,y_0)$.  Clearly,
$z\in Z(x',y)$. Since the order of $[z]\in \pi_0(Z(x,y))$ is $k$,
it follows from Lemma~\ref{pi0inj} that the order of $[z]\in
\pi_0(Z(x',y))$ is
also $k$. Thus the element 
$\mu\in(Q\spcheck/Q\spcheck(x'))_{w_c}$  is of order $k$. 
Since $\Phi(x') =
\Phi_L(x_0)$, $\Phi_c(x') = \Phi_{Z(x')}^{\rm proj}(w_c) =
\Phi_{Z_L(x_0)}^{\rm proj}(w_c)$. Thus, expressing $\tilde\mu$ as
a linear combination of $\pi(a\spcheck)$ for $\ov a\in
\widetilde I_c(x')$,  no 
coefficient  of
$\tilde\mu$ is integral. It
follows by Lemma~\ref{knotdivide} that
$k\not| g_{\ov a}$ for any
$\ov a\in
\widetilde I_c(x')$.
Hence, for a generic element $x'\in S\cdot x_0$ the 
subset
$\widetilde I_c(x')$ is contained in the subset of $\ov a\in 
\widetilde
\Delta_c$ for which $k\not| g_{\ov a}$. This together with
Lemma~\ref{subsetfn} shows 
that the generic $x'\in S\cdot x_0$ is conjugate to a point in the
interior of the face $f^c(k)$ and hence that ${\frak s}$ is
exactly the linear space parallel to $f^c(k)$.
This completes the proof of Proposition~\ref{parallel}.
\relax\ifmmode\expandafter\endproofmath\else
  \unskip\nobreak\hfil\penalty50\hskip.75em\hbox{}\nobreak\hfil\bull
  {\parfillskip=0pt \finalhyphendemerits=0 \bigbreak}\fi

The following is a corollary of the proof:

\begin{corollary}\label{constant}
Let $(x,y,z)$ be a $c$-triple and let $(x_0, y_0, z_0)$ be any
rank zero $c$-triple in $L$ such that $(x,y,z) = (s_1x_0, s_2y_0,
s_3z_0)$, where the $s_i\in S$. Then the order of $(x,y,z)$ is
the order of $[z_0] \in \pi_0(Z_L(x_0, y_0))/\pi_0(S\cap Z_L(x_0,
y_0))$, and hence the order is constant on connected components
of ${\cal T}_G(c)$. Finally, the order of the $c$-triple
$(x,y,z)$ of
$G$ divides the order of the $c$-triple $(x_0, y_0, z_0)$ in $L$. 
\end{corollary}

The order of $(x,y,z)$ in $G$ is not always the
order of $(x_0, y_0, z_0)$ in $L$.

By Proposition~\ref{parallel}, if $G$ has a $c$-triple of order
$k$, then $k$ divides at least one of the $g_{\ov a}$. Just as in
the case of commuting triples, there is a converse to this
statement:

\begin{proposition} Let $k$ be a positive integer. There is a
$c$-triple of order
$k$ in $G$ if  and only if $k$
divides at least one of the
$g_{\ov a}$. 
\end{proposition}

\begin{proof} The ``only if" direction follows from
Proposition~\ref{parallel}. Conversely, suppose that  $k$
divides at least one of the
$g_{\ov a}$. Choose $x$ to
be any element in the face 
$f^c(k)$. Then $k$ divides the order of ${\rm
Tor}\left((Q\spcheck/Q\spcheck(x))_{w_c}\right)$. It then follows
by Corollary~\ref{fullstatement} that there is a choice of $y$
such that $(x,y)$ is a $c$-pair and such that $k$ divides the
order of $\pi_0(Z(x,y))$. Choose a $z \in Z(x,y)$ mapping to
an element of order $k$ in $\pi_0(Z(x,y))$. Then $(x,y,z)$ is a
$c$-triple in $G$ of order $k$.
\end{proof}

\subsubsection{A normal form for $c$-triples of order $k$}

Suppose that ${\bf x}$ is a $c$-triple and that $S$, $L$ are as
above. Write ${\bf x} =(s_1x_0,s_2y_0, s_3z_0)$ where the $s_i\in
S$ and $(x_0,y_0,z_0)\in L$.  Let
$L=\prod _iL_i$ be the decomposition of
$L$ as a product of simple groups and let $c= c_1\cdots c_r$ be
the corresponding decomposition of $c$. We can write $x_0$ as a
product of elements $x_{0,i}$, and similarly for $y_0$, $z_0$. An
alcove for
$L$ is a product of alcoves for the $L_i$. Likewise, the root
system
$\Phi_c(x_0)$ on $(\frak t_L)^{w_c}$ is a product of the root 
systems
$\Phi_{c_i}(x_{0,i})$ on $(\frak t_{L_i})^{w_{c_i}}$. If
$(x_0, y_0)$ is a
$c$-pair in $L$, we say that $(x_0, y_0)$ is in normal form for
the product of the $B_i$ if each $c_i$-pair $(x_{0,i}, y_{0,i})$
is in normal form for $B_i$.

With this notation, we have the following:

\begin{proposition}\label{ctripnormform} 
Suppose that $k$ divides at least one of the $g_{\ov a}$
for $\ov a\in\widetilde \Delta_c$. Let $S\subseteq T^{w_c}$
be the torus whose Lie algebra ${\frak s}$
is parallel to the face
$f^c(k)$ of $A^c$. Let $L=DZ(S)$, and let ${\frak t}_L={\frak
t}\cap {\rm Lie}(L)=\frak s^\perp$. Let $\tilde x_0$ be the
unique point such that $({\frak s}+f^c(k))\cap{\frak
t}_L=\{\tilde x_0\}$, and let $x_0={\rm exp}(\tilde x_0)\in L$.
Finally let $y_0=sy_1$, where $y_1$ is such that $(x_0, y_1)$ is
a special $c$-pair in normal form and $s$ is the exponential of
the products of the barycenters of the alcoves $B_i$ for the
root systems $\Phi_{c_i}(x_{0,i})$ in the simple factors of $L$.
If 
${\bf x}$ is a
$c$-triple of order $k$ in $G$, then there are elements
$s_1,s_2,s_3$ in
$S$ and $z_0\in L$ such that
$(s_1x_0,s_2y_0,s_3z_0)$ is a $c$-triple conjugate to ${\bf x}$.  
\end{proposition}

\begin{proof} Write ${\bf x}= (x,y,z)$. Let $S$ be a maximal torus
for
$Z({\bf x})$. By Proposition~\ref{parallel}, possibly after
conjugating
${\bf x}$, we can assume that 
$S$ is the torus whose Lie algebra is parallel to $f^c(k)$ and
that $x$ is the exponential of an element of $f^c(k)$. Thus $x =
s_1x_0$ for some $s_1\in S$. There are $s_2,s_3\in S$ such
that $y=s_2y_0$ and $z=s_3z_0$, where $y_0, z_0\in L$. After a
further conjugation in
$L$, we can assume that $(x_0, y_0)$ is a $c$-pair in normal
form for $L$. Note that $(x_0, y_0, z_0)$ has rank zero in $L$. It
then follows from Proposition~\ref{c'1}, applied to the simple
factors of $L$, that
$y_0$ is as described in the statement of the proposition.
\end{proof}

\subsubsection{More on the group $L$}

Let ${\bf x}$ be a $c$-triple in $G$, let $S$ be a maximal
torus for $Z({\bf x})$, and let $L=DZ(S)$.

\begin{proposition}\label{atmost}
 Let $n_0= \gcd\{g_{\ov a}\}$.
The following are equivalent:
\begin{enumerate}
\item The order $k$ of the $c$-triple $(x,y,z)$ divides
$n_0$.
\item $S$ is conjugate to $S^{w_c}$.
\item $L$ is conjugate to $L_c = DZ(S^{w_c})$.
\item Every simple factor of $L$ is of type $A_n$ for some $n$.
\end{enumerate}
Finally,  there is at most one simple factor of $L$
which is not of type $A_n$.
\end{proposition}

\begin{proof} As in the case of commuting triples, since the
Dynkin diagram of $L$ is a proper subdiagram of the Dynkin
diagram of $G$, $L$ can have at most one simple factor which is
not of $A_n$-type. This proves the last statement.

To prove the equivalences of the proposition, first note
that, by
Proposition~\ref{parallel}, $S=S^{w_c}$ if and only if
$k$  divides $g_{\ov a}$ for every ${\ov a}$.  Thus (1) is
equivalent to (2). If $S=S^{w_c}$, then by definition $L= 
DZ(S^{w_c}) = L_c$. Thus (2) implies (3). Conversely, if
$L=L_c$, then $S$ is the torus associated to the intersections of
the kernels of the roots of $L_c$,  and thus $S=S^{w_c}$.

Since, by Proposition~\ref{finite}, all simple factors of
$L_c$ are of type
$A_n$, (3) implies (4). 
Finally, we show that (4) implies (3). Let $(x_0, y_0, z_0)$ be a
$c$-triple of rank zero in $L$ such that $(x,y,z) = (s_1x_0,
s_2y_0, s_3z_0)$, where the $s_i\in S$. Write $x_0 =
\prod_ix_{0,i}$, where the $x_{0,i}$ lie in the simple factors
$L_i$ of $L$, and similarly for
$y_0, z_0$. Since $L_i$ is of type $A_n$ for some $n$,
$(x_{0,i},y_{0,i})$ is a rank zero $c_i$-pair in $L_i$. It
follows that $L$ contains the rank zero $c$-pair $(x_0, y_0)$.
Since $c\in L$, it follows by Lemma~\ref{c} that $L_c\subseteq
L$. Since $L$ contains a rank zero $c$-pair, it follows that
$L\subseteq L_c$. Thus $L=L_c$.
\end{proof}

We can say   more about $L$ in case there is a simple
factor not of type $A_n$: 

\begin{proposition} Suppose that $L$ contains a simple factor,
not of type $A_n$, and write $L=L_0\times \prod_{i=1}^sL_i$,
where $L_0$ is simple and not of type $A_n$, and the $L_i$ are
simple and of type $A_n$ for $i\geq 1$. Let $c =\prod _{i=0}^sc_i$ be
the corresponding decomposition of $c$. Finally let $L_{0,c_0}$
be the subgroup of $L_0$ associated to the element $c_0\in {\cal
C}L_0$ as in Section~\ref{center}. Then $L_c = L_{0,c_0} \times
\prod_{i=1}^sL_i$.
\end{proposition}

\begin{proof} Since $c_0 \in L_{0,c_0}$, we have
$c\in L_{0,c_0} \times
\prod_{i=1}^sL_i$. Moreover, $L_{0,c_0} \times
\prod_{i=1}^sL_i$ is a product of groups of $A_n$-type, and $c$
projects to a generator of the center of each factor. The result 
now follows from  Lemma~\ref{LisLc}.
\end{proof}

\subsection{The number of components}

Our goal here is  to  describe the number of components of ${\cal
T}_G(c)$ of order $k$ and to identify each such component
explicitly. We will postpone the explicit determination of the
Weyl group
$W(S,G)$ to Section~\ref{sd}, however.

Suppose that $k$ divides at least one of the $g_{\ov a}$
for $\ov a\in\widetilde \Delta_c$.
We keep the notation of Proposition~\ref{ctripnormform}.
In particular, let $\tilde x_0$ be the
unique point such that $({\frak s}+f^c(k))\cap{\frak
t}_L=\{\tilde x_0\}$, and let $x_0={\rm exp}(\tilde x_0)\in L$,
and  let $y_0=sy_1$, where $y_1$ is such that $(x_0, y_1)$ is a
special $c$-pair in normal form and $s$ is the exponential of the
products of the barycenters of the alcoves $B_i$ for the root
systems $\Phi_{c_i}(x_{0,i})$ in the simple factors of $L$. Then
$(x_0, y_0)$ is a $c$-pair in $L$.

\begin{proposition}\label{whoknows2}
Let ${\bf x}$ be a $c$-triple of order $k$ in $G$. Then
${\bf x}$ is conjugate to a
$c$-triple of the form $(s_1x_0,s_2y_0,s_3z_0)$ where $x_0, y_0$
are as above,
$z_0\in L$ and $s_1,s_2,s_3\in S$. Denote by
$\psi({\bf x})$ the class
$[z_0]\in \pi_0(Z_G(x_0,y_0))$.
\begin{enumerate}
\item $\psi({\bf x})$ is well-defined and depends only on the
conjugacy class of ${\bf x}$.
\item $\psi$ is constant on the components of ${\cal T}_G(c)$ of
order $k$ in $G$.
\item $\psi({\bf x})$ is of order $k$.
\item $\psi$ induces a bijection from the set of components of
${\cal T}_G(c)$ of order $k$ in $G$ to the set of
elements of order $k$ in $\pi_0(Z_G(x_0,y_0))$.
\end{enumerate}
\end{proposition}
\begin{proof} By Proposition~\ref{ctripnormform}, ${\bf x}$ is 
conjugate to  $(s_1x_0,s_2y_0,s_3z_0)$ as claimed. To prove Part
(1), it is clearly sufficient to show the following: suppose that
$(s_1x_0,s_2y_0,s_3z_0)$ and
$(t_1x_0,t_2y_0,t_3z_0')$ are two $c$-triples which are conjugate
in $G$. Then $[z_0] = [z_0']$ in $\pi_0(Z_G(x_0,y_0))$. By
Theorem~\ref{main}, there is a $g\in N_G(S)$ such that
$$i_g(x_0,y_0, z_0)  =(u_1x_0,u_1y_0, u_3z_0'),$$ where
$u_i\in S\cap L$. Moreover we can assume that
$g\in N_G(T)$. Clearly $[z_0'] = [i_g(z_0)]$, and we must show
that $[z_0] = [i_g(z_0)]$.  

First note that, since $[g, y_0]\in {\cal C}L\subseteq T$, $i_g$
commutes with
$w_c$, and thus $i_g$ acts on $T^{w_c}$ and hence on
$\pi_0(T^{w_c})$. In particular $i_g$ acts on $L_c$. By
Lemma~\ref{surjpi0}, the center
${\cal C}G$ surjects onto
$\pi_0(T^{w_c})$, and so the induced action of $i_g$ on
$\pi_0(T^{w_c})$ is trivial.  

Since $S\cap
L\subseteq {\cal C}L$, the inner automorphism $i_g$ defines an
automorphism of
$Z_L(x_0, y_0)$. Moreover, $[z_0]$ is the image of the
corresponding element of $\pi_0(Z_L(x_0,y_0))$ under the natural
homomorphism $\rho\colon\pi_0(Z_L(x_0,y_0)) \to
\pi_0(Z_G(x_0,y_0))$. We must show that $\rho\circ i_g = \rho$.
We may write $L\cong L_0\times \prod _{i\geq 1}L_i$, where the
$L_i$ are simple groups of $A_n$ type and $L_0$ is either trivial
or a simple group not of type $A_n$. Clearly $i_g$ preserves the
factors $L_0$ and
$\prod _{i\geq 1}L_i$. There is a corresponding direct sum
decomposition
$$\pi_0(Z_L(x_0,y_0))  = \pi_0(Z_{L_0}(x_{0,0},y_{0,0}))  \oplus
\bigoplus _{i\geq 1} \pi_0(Z_{L_i}(x_{0,i},y_{0,i})) ,$$
where the $x_{0,i}, y_{0,i}$ are the projections  of $x_0$, $y_0$
to the factor $L_i$. We analyze the action on each of these
factors separately.

The group $\prod _{i\geq 1}L_i$ is a subgroup of $L_c$. Thus the
map $\bigoplus _{i\geq 1} \pi_0(Z_{L_i}(x_{0,i},y_{0,i})) \to
\pi_0(Z_G(x_0,y_0))$ factors through $\pi_0(Z_{L_c}(x_0,y_0))$.
On the other hand, $\pi_0(Z_{L_c}(x_0,y_0)) = Z_{L_c}(x_0,y_0) =
{\cal C}L_c$. Since ${\cal C}L_c \subseteq T^{w_c}$, the map
$\pi_0(Z_{L_c}(x_0,y_0)) \to \pi_0(Z_G(x_0,y_0))$ factors through
$\pi_0(T^{w_c})$. But as we have seen, $i_g$ acts trivially on
$\pi_0(T^{w_c})$. Thus, for any element $\xi'\in \bigoplus _{i\geq
1} \pi_0(Z_{L_i}(x_{0,i},y_{0,i}))$, we have $\rho(i_g(\xi')) =
\rho(\xi')$.

Now consider the action of $i_g$
on $\pi_0(Z_{L_0}(x_{0,0},y_{0,0}))$. Since $i_g$ defines an
automorphism of $Z_{L_0}$ which sends $(x_{0,0},y_{0,0})$ to an
element of the form $(v_1x_{0,0},v_2y_{0,0})$, where the $v_i\in
{\cal C}L_0$, it follows from Lemma~\ref{automorphism} that the
action of $i_g$
on $\pi_0(Z_{L_0}(x_{0,0},y_{0,0}))$ is trivial. Thus, for
$\xi_0\in \pi_0(Z_{L_0}(x_{0,0},y_{0,0}))$, $\rho(i_g(\xi_0)) =
\rho(\xi_0)$. It follows that $\rho\circ i_g = \rho$, and finally
that $[z_0] = [i_g(z_0)]$ as claimed.

To see (2), if the conjugacy class of ${\bf x}'=(x', y', z')$ is
in the same component of
${\cal T}_G(c)$ as that of ${\bf x}$, then after conjugation we
can assume that
$(x', y', z') = (t_1x_0, t_2y_0, t_3z_0)$ for   $t_i\in S$. Then 
by the definition of $\psi$ and the fact that it is well-defined,
we see that $\psi({\bf x}') =[z_0] =
\psi({\bf x})$.

To see (3), set ${\bf x}_0= (x_0, y_0, z_0)$. Then $\psi({\bf
x}) = \psi({\bf x}_0) =[z_0]$. On the other hand, ${\bf x}_0$ is
a $c$-triple in the same component of ${\cal T}_G(c)$ as ${\bf
x}$, so by Corollary~\ref{constant}, the order of ${\bf x}_0$ is
$k$. By definition, this is the order of $[z_0] = \psi({\bf x}_0)
=
\psi({\bf x})$ in $\pi_0(Z_G(x_0, y_0))$.

Finally, we prove (4). Suppose that $\psi({\bf x}_1) = \psi({\bf
x}_2)$. Then the conjugacy class of ${\bf x}_1$ is in the same
connected component of ${\cal T}_G(c)$ as that of $(x_0, y_0,
z_1)$, say, and the conjugacy class of
${\bf x}_2$ is in  the same connected component of ${\cal
T}_G(c)$ as that of $(x_0, y_0, z_2)$, where $z_1$ and $z_2$ lie
in
$L$ and differ by an element of
$Z^0_G(x_0, y_0)$. Hence the classes of $(x_0, y_0, z_1)$ and
$(x_0, y_0, z_2)$ are in the same connected component of ${\cal
T}_G(c)$, as are those
${\bf x}_1$ and ${\bf x}_2$. Thus $\psi$ is injective on the set
of components, and its image is contained in the set of elements
in
$\pi_0(Z_G(x_0, y_0))$ of order $k$. To see that it is surjective,
choose any element of order $k$ in
$\pi_0(Z_G(x_0, y_0))$, represented by
$[z_0]$, say. Then $(x_0, y_0, z_0)$ is a $c$-triple of order
$k$, and by construction $\psi(x_0, y_0, z_0) = [z_0]$.
\end{proof}

\begin{corollary}\label{noofcomps}
If $k$ divides at least one of the $g_{\ov a}$ for $\ov
a\in\widetilde \Delta_c$, then there are exactly $\varphi(k)$
components of ${\cal T}_G(c)$ of order $k$ in $G$.
\end{corollary}

\subsection{Proof of Parts 1,2,3   of
Theorem~\protect{\ref{ctrip}} for
$\langle C\rangle$ cyclic}

We assume that $\langle C\rangle$ is cyclic.
Part 1 of Theorem~\protect{\ref{ctrip}} is contained in
Proposition~\ref{parallel}. Part 2 of
Theorem~\protect{\ref{ctrip}} is contained in
Proposition~\ref{constant}. The first statement of Part 3 of
Theorem~\protect{\ref{ctrip}} is Corollary~\ref{noofcomps}. Let
$X$ be a component of ${\cal T}_G(c)$ of order $k$. Let $S =
S^{w_c}(\ov {\bf g}, k)$. Since
$X$ is the quotient of $S^3$ by a finite group,  Part (1) implies
that
$d_X=\frac13\dim X + 1$ is equal to the number of $\ov a$ such
that
$k|g_{\ov a}$. By the first statement of Part (3) and
Lemma~\ref{standardfact},
$$\sum _Xd_X= \sum _{\ov a\in \widetilde \Delta_c}g_{\ov a} =
\sum _{\ov a\in \widetilde \Delta_c}n_{\ov a}g_a = g.$$

\subsection{Proof of Part 4 of Theorem~\protect{\ref{ctrip}} for
$\langle C\rangle$ cyclic}

\begin{theorem}\label{ccomp}
Let $X$ be a component of ${\cal T}_G(c)$.
 Associated to $X$ is the torus $S$ and  the group
$L=DZ(S)$.
Decompose $L=L_0\times L'$ where $L_0$ is either
trivial or is a simple group not of $A_n$-type for any $n\ge 1$.
If $L_0$ is trivial, the map $\ov p$ of Theorem~\ref{main}
induces a homeomorphism from 
$$(\ov S\times \ov S \times S)/W(S,G)$$
to $X$.
If $L_0$ is not trivial, then the map $\ov p$ of
Theorem~\ref{main} induces a homeomorphism from 
$$(\ov S\times \ov S \times (S/{\cal K}))/W(S,G)$$
to $X$ where ${\cal K}$ is a subgroup of order at most $2$ in
${\cal C}L_0$.
\end{theorem}

\begin{proof}
Let $F = S\cap L$. There is an action of $F\times F\times F$ on the
set of all conjugacy classes of rank zero $c$-triples in $L$. By
Lemma~\ref{x} and Lemma~\ref{y}, the action of $F$ on the first
two factors is trivial.
We analyze the action of $F$ on the
last factor. Write $c=c_0c'$ with $c_0\in {\cal C}L_0$ and
$c'\in{\cal C}L'$. 
There is an inclusion $S\cap L \subseteq {\cal C}L$.
Since $L'$ decomposes as a product of groups of $A_n$-type.
It is easy to see that ${\cal C}L'$ acts freely on the moduli
space of $c'$-triples in $L'$. 
According to Lemma~\ref{z} the stabilizer of a conjugacy class of
$c_0$-triples in $L_0$ is the intersection  ${\cal C}L_0\cap
S^{w_{c_0}}$, which is a  subgroup  ${\cal K}_0\subseteq {\cal
C}L_0$ of index given by the gcd of the quotient coroot integers
for $L_0$. The two possibilities for $L_0$
when ${\rm dim}(S^{w_c})
>0$ are $L_0$ of type $C_2$ or $D_6$.
In both these cases, ${\cal K}_0$ is of order $2$ and may be
described explicitly. We let ${\cal K}$ be the intersection
$S\cap {\cal K}_0 \subseteq S\cap L$. Clearly, ${\cal K}$ has
order at most $2$, and is trivial if $L_0$ is trivial.

Thus, the stabilizer in $F^3$ of any $c$-triple of rank zero in 
$L$ is equal to $F\times F\times {\cal K}$ where ${\cal K}$ is
trivial if $L_0$ is trivial, and ${\cal K}$ has order at most $2$
if $L_0$ is non-trivial. Finally, it follows from 
Lemma~\ref{outer} that the Weyl group of $S$ acts trivially on
the set of conjugacy classes of rank zero $c$-triples in $L$. The
theorem now follows from Theorem~\ref{main}. 
\end{proof}

The possibilities for ${\cal K}$ and the component $X$ in case
$L_0$ is not trivial are as follows:
\begin{itemize}
\item If $G$ is of type $B_n$, and hence $L_0=L$ is of type $C_2$,
then ${\cal K}={\cal C}L = S\cap L$ has order $2$, and $X= (\ov
S\times \ov S \times \ov S)/W(S,G)$.
\item If $G$ is of type $C_{2n}$, and hence $L_0$ is of type
$C_2$, then ${\cal K}$ is trivial and $X= (\ov
S\times \ov S \times S)/W(S,G)$.
\item If $G$ is of type $D_{2n}$, and hence $L_0$ is of type
$D_6$, then ${\cal K}$ is trivial and $X= (\ov
S\times \ov S \times S)/W(S,G)$.
\item If $G$ is of type $E_7$, and hence $L_0=L$ is of type $D_6$,
then ${\cal K}=S\cap L$ has order $2$ and $X= (\ov
S\times \ov S \times \ov S)/W(S,G)$.
\end{itemize}

This explicit list completes the proof of
Theorem~\ref{ctrip} in case $\langle C\rangle$ is cyclic.

\section{The tori $\ov S(k)$ and $\ov S^{w_C}(\ov {\bf g}, k)$ and
their Weyl groups}\label{sd}  

\setcounter{theorem}{0}

Let $\Phi$ be an irreducible but possibly non-reduced 
root system on a 
vector space ${\frak t}$ with set of simple roots $\Delta$.
Let $A\subset {\frak t}$ be the alcove associated with $\Delta$,
and let $d$ be the highest root.

\begin{defn}\label{10}
Fix a positive integer $n_0$ and  let ${\bf
n}\colon
\widetilde \Delta\to {\bf N}$ be the function $n_0{\bf g}$.
Fix an integer $k> 1$ dividing at least one of the
integers $\{{\bf n}(a)\}_{a\in\widetilde\Delta}$.
Let $\widetilde I({\bf n},k)=\{a\in \widetilde \Delta: k\not|
{\bf n}(a)\}$.  
Note that $\widetilde I({\bf g},k)=\widetilde I(k)$ as
previously defined.
More generally, if we factor $k$ as $dk'$ where $d={\rm gcd}(k,n_0)$,
then $\widetilde I({\bf n},k)=\widetilde I( k')$. 
Let $f({\bf n},k)\subset A$ be defined by 
$$f({\bf n},k) = \{\,\tilde x \in A: a(\tilde x ) \in {\bf Z} {\rm \,\,
for\,\,all\,\,} a\in \widetilde I({\bf n},k)\,\}.$$
Let
$\hat f({\bf n},k)$ be the affine space generated by $f({\bf n},k)$ and let
${\frak t}({\bf n},k)$ be the  linear space parallel to $\hat f({\bf n},k)$. 
Notice that if $n_0=1$, then $f({\bf n},k)=f(k)$ and  ${\frak
t}({\bf
  n},k)={\frak t}(k)$ as defined in Theorem~\ref{commuttrip}.
More generally, if we factor $k$ as $dk'$ where $d={\rm gcd}(k,n_0)$,
then $f({\bf n},k)=f(k')$ and ${\frak t}({\bf n},k)={\frak t}(k')$.
\end{defn} 

\subsection{A root system on $\hat f(k)$}

Let $k\ge 1$ be an integer dividing at least one of the $g_a$.
We denote by $\Phi(k)$ the subroot system of $\Phi$ consisting of all
roots which annihilate ${\frak t}(k)$.
Suppose that $\Phi$ is reduced and let $\Phi^+$ be the set of
positive roots for $\Phi$.
Recall from Proposition~\ref{kroots} and
Claim~\ref{kroots2} that $\Phi^+(k)=\Phi^+\cap  
\Phi(k)$ is a set of positive roots for $\Phi(k)$ and that
 $\Delta(k)$, 
the set of simple roots  determining this set  
of positive roots, is given by $\Delta(k)=
I(k)\cup \{b\}$ for some root $b\in \Phi(k)$,
where $I(k)=\widetilde I(k)\cap \Delta$. Furthermore,
the root system
$\Phi(k)$ is irreducible and $d$
is its highest root.
 Writing $d\spcheck=\sum_{a\in
I(k)}g'_aa\spcheck +g'_bb\spcheck$, we have that $g'_b=k$.
Furthermore, for every $a\in I(k)$ the 
coroot integer $g'_a$ is not  divisible by $k$.

We let ${\frak u}$ denote the subspace of ${\frak t}$ spanned by the
coroots inverse to the roots in $\Phi(k)$. Then
${\frak u}={\frak t}(k)^\perp$.

\begin{lemma}
 Let  $\tilde x_0=\hat f(k)\cap {\frak u}$ and let $A(k)$
be the alcove in ${\frak u}$ determined by the set of
simple roots $\Delta(k)$ for $\Phi(k)$.
Then $\tilde x_0$ is the vertex of the alcove $A(k)$ for $\Phi(k)$
 opposite the
face defined by $\{b=0\}$.
\end{lemma}

\begin{proof}
Each  $a\in I(k)$ vanishes on $\hat f(k)$
and hence on $\tilde x_0$. Similarly, $d(\tilde x_0)=1$.
This proves that $\tilde x_0$ is the vertex of $A(k)$ opposite
the face $\{b=0\}$.
\end{proof}

\begin{lemma}
Suppose that $g\in
W_{\rm aff}(\Phi)$ and the differential $w$ of $g$ 
normalizes ${\frak t}(k)$. Then there 
is an element  $w'\in W(\Phi(k))$ such that $w'w$
normalizes $\hat f(k)$.
\end{lemma}

\begin{proof}
Suppose first that $\Phi$ is reduced.
Since $w$ normalizes ${\frak t}(k)$, it normalizes ${\frak u}$ and
hence it also normalizes the root system $\Phi(k)$ on ${\frak
u}$.  Thus, $w(\tilde x_0)$ is the vertex of an alcove $A'(k)$
for $\Phi(k)$ containing the origin. 
There is an element $w'\in
W(\Phi(k))$ such that $w'(A'(k))=A(k)$. Then $w' w$ is an
automorphism of $A(k)$. Since $b$ is the unique element of $\widetilde
\Delta(k)$ for which the coroot integer $g'_b$ is divisible by $k$, the
root $b$ is fixed by any
automorphism of $\widetilde \Delta(k)$. Thus, $\tilde x_0$ is fixed by
$w' w$. Since $w' w$ fixes $\tilde x_0$ and normalizes
${\frak t}(k)$ it normalizes $\hat f(k)$.

If $\Phi$ is not reduced, then $\Phi$ is of type $BC_n$ for some $n$.
The fundamental relation among the coroots in this case is of the form
$$2\tilde a\spcheck+2a_1\spcheck+\cdots+2a_{n-1}\spcheck+a_n\spcheck$$
where $\tilde a=-d$. Thus, $I(2)=\widetilde I(2)=\{a_n\}$.
(The difference here is that the coefficient of $\tilde a\spcheck$ is
not one, as it is for reduced root systems.)  This
implies that $f(2)$ is the face of $A$ given by the equation
$\{a_n=0\}$, and so $\hat f(2)$ is a linear
subspace.
Hence $\hat f(2)={\frak t}(2)$, and the statement of the
lemma is obvious.
\end{proof}

\begin{corollary}\label{repren}
Every element  $\ov w\in W({\frak t}(k),\Phi)$ has a 
representative
$w\in N_{W(\Phi)}(\hat f(k))$.  
\end{corollary}

\begin{lemma}
Suppose that $\Phi$ is reduced.
Any $\tilde x'\in \hat f(k)$ is equivalent under the action of  $W_{\rm
aff}(\Phi)$  to a point in $f(k)$.
\end{lemma}

\begin{proof}
We may assume that $\Phi =\Phi(T, G)$ for some simple group
$G$ with maximal torus $T$. We set
$S(k)={\rm exp}({\frak t}(k))$ and let
$L(k)=DZ(S(k))$. Let $x_0={\rm exp}(\tilde x_0)\in L(k)$ and let
$(x_0,y_0,z_0)$ be a commuting triple of order $k$ and rank zero
 in $L(k)$. Then $Z_G(x_0,y_0,z_0)$
has $S(k)$ as a maximal torus.
Let $\tilde x'\in\hat f(k)$ and let $x'={\rm exp}(\tilde x')$. Then 
$(x',y_0,z_0)$ is a commuting triple of order $k$ in $G$. Then 
according to Part 3 of Proposition~\ref{component}, the point $x'$
is conjugate in
$G$ and hence in $W(\Phi)$ to a point which is the
exponential of a point in $f(k)$. Hence, $\tilde x'$ is
conjugate under $W_{\rm aff}(\Phi)$ to a point of $f(k)$.
\end{proof}

\begin{corollary}
Suppose that $\Phi$ is irreducible but not necessarily reduced.
The alcove decomposition of ${\frak t}$ induces a decomposition of
$\hat f(k)$ into compact convex regions with disjoint interiors. Let
${\cal A}(k)$ be this decomposition. 
Then $f(k)\in  {\cal A}(k)$.
Moreover, the action of $N_{W_{\rm
aff}(\Phi)}(\hat f(k))$ on $\hat f(k)$ preserves the
decomposition ${\cal A}(k)$ and is transitive on
${\cal A}(k)$.
\end{corollary}

\begin{proof}
First let us consider the case when $\Phi$ is reduced.
Clearly, $N_{W_{\rm aff}(\Phi)}(\hat f(k))$ acts on $\hat f(k)$
normalizing ${\cal A}(k)$ and $f(k)$ is one of the elements of 
this decomposition.
 Suppose that $A_1$ is
an element of this decomposition. Let $\tilde x'$ be an interior point
of $A_1$. According to the previous lemma, there is an
element $g\in W_{\rm aff}(\Phi)$ which conjugates $\tilde x'$ to a
point of $f(k)$. Since there are only finitely many elements of
$W_{\rm aff}(\Phi)$ which conjugate $A_1$ so as to meet
$f(k)$, it follows that there is some element $g\in W_{\rm aff}(\Phi)$
which conjugates an open subset of $A_1$ into $f(k)$. This
element normalizes $\hat f(k)$ and hence sends $A_1$ onto
$f(k)$. This proves the transitivity statement in the reduced case.

If $\Phi$ is not reduced, then it is of type $BC_n$ for some $n$,
$k=2$  and $\hat f(2)$ is the linear subspace spanned by the coroots of the
subsystem $BC_{n-1}$. The induced alcove decomposition of 
$\hat f(2)$ is exactly the alcove decomposition for $BC_{n-1}$. Since
$N_{W_{\rm aff}(\Phi)}(\hat f(2))$ contains
the affine Weyl group of $BC_{n-1}$, the lemma is
clear in this case.
\end{proof}

\begin{lemma}
If $g\in W_{\rm aff}(\Phi)$ has the property that there is an element
of the decomposition ${\cal A}(k)$ normalized by $g$, then
$g|\hat f(k)$ is the identity.
\end{lemma}

\begin{proof}
Let $A_1$ be an element of ${\cal A}(k)$ such that
$g(A_1)= A_1$. Let $A$ be an alcove for $\Phi$ containing
$A_1$ in its closure, and let $B=g(A)$. Then $B$ also contains
$ A_1$, and hence $A\cap B$ contains $A_1$. But the unique
element of $W_{\rm aff}(\Phi)$ taking $A$ to $B$ fixes pointwise the
intersection $A\cap B$. Thus, $g|A_1$ is the identity, and
consequently $g|\hat f(k)$ is the identity.
\end{proof}

\begin{corollary}
Define $W_{\rm aff}(\hat f(k))$ be $N_{W_{\rm aff}(\Phi)}(\hat f(k))/Z_{W_{\rm
aff}(\Phi)}(\hat f(k))$. Then $ W_{\rm aff}(\hat f(k))$ acts on $\hat
f(k)$ as a group of affine isometries. It acts simply transitively on the
elements of the decomposition ${\cal A}(k)$.
\end{corollary}

\begin{proposition}\label{generateW}
The reflections of $\hat f(k)$ in the walls of $f(k)$ are realized
by elements  in
$W_{\rm aff}(\hat f(k))$. These reflections generate
$W_{\rm aff}(\hat f(k))$ which is thus a Coxeter group
with fundamental domain  $f(k)$.
\end{proposition}

\begin{proof}
Consider an element $A_1$ for the decomposition ${\cal
A}(k)$ which shares a codimension-one wall with $f(k)$. Let $g$
be  an element of $N_{W_{\rm aff}(\Phi)}(\hat f(k))$ which sends
$f(k)$ to
$A_1$. This affine isometry of $\hat f(k)$ fixes the
intersection of $f(k)$ and $A_1$, which is a codimension-one
affine subspace. Thus, it is a reflection in this subspace.
This shows that the reflections in the walls of $f(k)$ are 
elements of
$W_{\rm aff}(\hat f(k))$. Since $W_{\rm aff}(\hat f(k))$ acts 
simply transitively on the 
elements in the decomposition ${\cal A}(k)$, it follows that
these reflections generate $W_{\rm aff}(\hat f(k))$.
\end{proof}

\begin{corollary}
There is a reduced root system $\Phi({\frak t}(k))$ on ${\frak t}(k)$
and  a vertex $v$ of $f(k)$ such that, using $v$ to identify $\hat
f(k)$ with ${\frak t}(k)$, the affine Weyl group of $\Phi({\frak t}(k))$
is identified with $W_{\rm aff}(\hat f(k))$.
\end{corollary}

\begin{proof}
The follows immediately from the fact, established in the previous
proposition, that $W_{\rm aff}(\Phi(\hat f(k))$ is a Coxeter
group.
\end{proof}

\begin{corollary}\label{reflectionresult}
We define an embedding of $\widetilde \Delta\spcheck-\widetilde
I\spcheck(k)$ into 
${\frak t}(k)$ by sending $a\spcheck$ to $\pi_k(a\spcheck)$, where
$\pi_k$ is the orthogonal projection.
Up to positive multiples $\widetilde
\Delta\spcheck-I\spcheck(k)\subset {\frak t}(k)$ is the set of
extended coroots for the root system $\Phi({\frak t}(k))$.
In particular, the Weyl group of ${\frak t}(k)$ in $\Phi$ is the group
generated by reflections 
in the $\widetilde \Delta\spcheck-\widetilde
I\spcheck(k)$.
\end{corollary}

\begin{proof}
Since the walls of the alcove $f(k)\subset \hat f(k)$ are the
subspaces of $\hat f(k)$ orthogonal to
the $\pi_k(a\spcheck)$ for $a\spcheck\in\widetilde
\Delta\spcheck-\widetilde I(k)$, the first statement is clear. 
Clearly, the image under the differential of $W_{\rm aff}(\hat
f(k))=N_{W_{\rm 
aff}(\Phi)}(\hat f(k))/Z_{W_{\rm aff}(\Phi)}(\hat f(k))$ is contained
in the Weyl group of ${\frak t}(k)$. By Corollary~\ref{repren} this
map is onto. By Proposition~\ref{generateW}, $W_{\rm aff}(\hat
f(k))$ is generated 
by the reflections in the walls of $f(k)$. The image under the
differential of these reflections is the set of reflections of ${\frak t}(k)$
in the $\pi_k(a\spcheck)$ for $a\in \widetilde \Delta-\widetilde
I(k)$.
\end{proof}

\begin{lemma}\label{latticeresult}
The lattice generated by $\pi_k(a\spcheck)$ for $a\in\widetilde
\Delta-\widetilde I(k)$ is the image under the orthogonal projection
of the coroot lattice $Q\spcheck$.
\end{lemma}

\begin{proof}
Since $\widetilde \Delta$ spans $Q\spcheck$ and since
$\pi_k(a\spcheck)=0$ if $a\in\widetilde I(k)$, this is clear.
\end{proof}

\begin{proposition}\label{XI.1.13}
The coroot lattice of $\Phi({\frak t}(k))$ is $\pi_k(Q\spcheck)$.
\end{proposition}

\begin{proof}
The coroot lattice of $\Phi({\frak t}(k))$ is identified with the
group of translations of $\hat f(k)$ which occur as restrictions of
elements of $W_{\rm aff}(\Phi)$ normalizing $\hat f(k)$.

Let $\lambda\in Q\spcheck$. Translation by $\lambda$ carries $\hat
f(k)$ to an affine subspace $\hat f_1$ of ${\frak t}$ which meets
${\frak u}$ in a point of the form $\tilde x_0+\zeta$ where $\zeta$ is
an element in the dual to the root lattice of $\Phi(k)$. In
particular, $[\zeta]\in {\frak u}/Q\spcheck(\Phi(k))$ is an element of
${\cal C}\Phi(k)=P\spcheck(\Phi(k))/Q\spcheck(\Phi(k))$. But
since $\tilde x_0$ is the vertex of
$A(k)$ opposite the wall $\{b=0\}$, and $b$ is the unique element
in
$\widetilde \Delta(\Phi(k))$ whose coroot integer is divisible by $k$,
it follows that the action of ${\cal C}\Phi(k)$ on $A(k)$ fixes
$\tilde x_0$. Thus there is an element $w\in W_{\rm aff}(\Phi(k))$
such that $w(\tilde x_0+\zeta)=\tilde x_0$.
The composition of translation by $\lambda$ followed by $w$ is an
element of $W_{\rm aff}(\Phi)$ normalizing $\hat f(k)$ and acting on
it by translation by $\pi_k(w\lambda)= \pi_k(\lambda)$. This shows
that
$\pi_k(Q\spcheck)$ is contained in the coroot lattice of
$\Phi({\frak t}(k))$. 

Conversely, suppose that $w\in W_{\rm aff}(\Phi)$ normalizes $\hat
f(k)$ and acts on it by a pure translation. We write
$w(x)=w_0(x)+\lambda$ where $\lambda\in Q\spcheck$ and $w_0$ is in the
Weyl group of $\Phi$. The element $w_0$ normalizes ${\frak t}(k)$ and
hence ${\frak u}$ so that $w_0(\tilde x_0)=\tilde x_0-\zeta$ for some
$\zeta\in {\frak u}$. Thus, $w(\tilde x_0)=\tilde
x_0+(\lambda-\zeta)$.
Since this element restricts to $\hat f(k)$ to give  a translation,
its restriction is  translation by $\pi_k(\lambda-\zeta)$. Since 
$\zeta\in {\frak u}={\rm Ker}(\pi_k)$, we have
$\pi_k(\lambda-\zeta)=
\pi_k(\lambda)\in\pi_k(Q\spcheck)$. 
\end{proof}

\begin{theorem}\label{kthm}
The reduced root system $\Phi({\frak t}(k))$ on ${\frak t}(k)$ has
Weyl group equal to 
the Weyl group of ${\frak t}(k)$ in $\Phi$. It has coroot lattice
equal to $\pi_k(Q\spcheck)$.  If
$\dim {\frak
  t}(k)\geq 1$, then
$\{\pi_k(a\spcheck)\}_{a\in\widetilde \Delta-\widetilde I(k)}$ as
an extended set of coroots for $\Phi({\frak t}(k))$.
\end{theorem}

\begin{proof}
Everything except the last statement is contained in
Corollary~\ref{reflectionresult} and Proposition~\ref{XI.1.13}.
Since $\pi_k(a\spcheck)\in \pi_k(Q\spcheck)$ and since these form,
up to positive multiples, a set of extended coroots, to complete
the proof we need only see that the $\pi_k(a\spcheck)$ are
indivisible elements of
$\pi_k(Q\spcheck)$. 
This will follow from the next lemma.

\begin{lemma}\label{numerical}
Writing $b\spcheck=\sum_{a\in\Delta-I(k)}m_aa\spcheck+\sum_{a\in
  I(k)}n_aa\spcheck$, one of the 
following holds:
\begin{enumerate}
\item The cardinality of $\Delta-I(k)$ is one and the unique $m_a$ is
  one;
\item There are at least two $a\in \Delta-I(k)$ for which $m_a=1$.
\end{enumerate}
\end{lemma}

\begin{proof}
We have
$$d=\sum_{a\in I(k)}g_a'a\spcheck+
k\left(\sum_{a\in \Delta-I(k)}m_aa\spcheck+
\sum_{a\in I(k)}n_aa\spcheck\right)=\sum_{a\in
\Delta} g_aa\spcheck.$$ Thus, for $a\in \Delta-I(k)$, $g_a
=km_a$, and the remaining $g_a$ are not divisble by $k$. The
result then follows from Lemma~\ref{atleast2}. 
\end{proof}

Now let us return to the proof of the theorem.
Since $I(k)\cup \{b\}$ is a set of simple roots for $\Phi(k)$, 
we have an exact sequence
$$0\to \sum_{a\in I(k)}{\bf Z}(a\spcheck)\oplus {\bf
  Z}(b\spcheck)\to 
\sum_{a\in \Delta}{\bf Z}(a\spcheck)\to \pi_k(Q\spcheck)\to 0.$$
We can rewrite this sequence as
$$0\to {\bf Z}(b\spcheck)\to \sum_{a\in \Delta-I(k)}{\bf
Z}(a\spcheck)\to
\pi_k(Q\spcheck)\to 0.$$
According to Lemma~\ref{numerical}, the image of
$b\spcheck$ in $\pi_k(Q\spcheck$ is $\sum_{a\in\Delta-
  I(k)}m_aa\spcheck$    where either  the
cardinality of $\Delta-I(k)$ is one or there are at least two
$a\in
\Delta-I(k)$ for which $m_a=1$.   In the first case ${\frak
t}(k)$ is a point. In the second
case, it is easy to see  that the image of each
$a\spcheck\in
\Delta-I(k)$ is indivisible in
$\pi_k(Q\spcheck)$.
\end{proof}

Let us return to the case of a general function
${\bf n}=n_0{\bf g}$.
Then $\widetilde I({\bf n},k)=\widetilde I(k')$
and  ${\frak t}({\bf n},k)={\frak t}(k')$
where $k'=k/{\rm gcd}(n_0,k)$, and thus $\Phi({\frak t}({\bf
  n},k))=\Phi({\frak t}(k'))$.
The following is then an immediate corollary of
Theorem~\ref{kthm}.

\begin{corollary}\label{kcor}
The Weyl group of ${\frak t}({\bf n},k)$ in $\Phi$ is the Weyl
group of the reduced root system $\Phi({\frak t}({\bf n},k))$.
The coroot lattice of $\Phi({\frak t}({\bf n},k))$ is
$\pi_k(Q\spcheck)$.  
If $\dim {\frak  t}({\bf n},k)\geq 1$, then
$\{\pi_k(a\spcheck)\}_{a\in\widetilde \Delta-\widetilde I({\bf
n},k)}$ is an extended set of simple coroots of $\Phi({\frak
t}({\bf n},k))$.
\end{corollary}

\subsection{Completion of the proof of
Theorem~\protect{\ref{commuttrip}}}

To complete the proof of Theorem~\ref{commuttrip} we must show that
the torus $\ov S(k)$ and the Weyl group $W(S(k),G)$ are as described
in Part 5 of the statement of that theorem.
This is immediate by applying Theorem~\ref{kthm}
to the root system $\Phi=\Phi(G)$.

\subsection{Completion of the proof of Theorem~\protect{\ref{ctrip}}
in case $\langle C\rangle$ is cyclic}

Let us assume  that $\langle C\rangle$ is cyclic and generated by
$c$. It remains to establish that $\ov S^{w_c}(\ov {\bf g},k)$ and
$W(S^{w_c}(\ov{\bf g},k),G)$ 
are as given in Part 5 of the statement of that theorem.

\begin{lemma}
Every element in $W(S^{w_c}(\ov{\bf g},k),G)$ has a
representative in the Weyl group of $G$ which normalizes
$S^{w_c}$. 
\end{lemma}

\begin{proof}
Set $S=S^{w_c}(\ov{\bf g},k)$ and set $L=DZ(S)$.
According to Theorem~\ref{main}, there is a  $c$-triple 
$(x_0,y_0,z_0)$ of rank zero in $L$, and $S$ is a maximal
torus of $Z_G(x_0,y_0,z_0)$.
By Corollary~\ref{normalform} we can assume that
 $(x_0,y_0)$ is a $c$-pair in $L$ in weak normal form
with respect to $T\cap L$.
 
Let $g\in N_G(T)$ normalize $S$. Then $g$ also normalizes $L$ and
$g(x_0,y_0,z_0)g^{-1}$ is another $c$-triple of rank zero and
order
$k$ in $L$. Thus, by Proposition~\ref{k} this triple is
conjugate by an element $h\in L$ to a 
triple of the form $(x_0,y_0,z_0^\ell)$ for some $\ell\in {\bf Z}$. 
In particular, $hg(x_0,y_0)g^{-1}h^{-1}=(x_0,y_0)$. Since
$(x_0,y_0,z_0)$ is a triple of rank zero in $L$, the
centralizer $Z^0_L(x_0,y_0)$ is  a torus. Suppose that a  root of
$G$ with respect to $T$ annihilates $S\cdot x_0$ and
$S\cdot y_0$. This root then 
annihilates $S$, and hence is a root of $L$ with respect to
$T\cap L$ annihilating
$x_0$ and $y_0$. But we have just seen that there are no such roots.
Thus, $Z_G^0(S\cdot x_0,S\cdot y_0)$ is also a torus. Since
$(x_0,y_0)$ is a $c$-pair in weak normal form,
Lemma~\ref{4321} implies that
$S^{w_c}\subseteq Z^0_G(x_0,y_0)$ is a maximal torus.
Because $hg$ fixes $(x_0,y_0)$ and normalizes $S$, it
normalizes $Z^0_G(S\cdot x_0,S\cdot y_0)=S^{w_c}$. Clearly, the
image in $W(S,G)$ of $hg$ is equal to that of $g$.
\end{proof}

Let us consider the root system $\Phi^{\rm proj}(w_c)$ on ${\frak
t}^{w_c}$. 
According to Proposition~\ref{Wres=Wproj} it is irreducible,
but possibly not reduced. Let $\ov {\bf g}\colon \widetilde
\Delta_c \to {\bf N}$ be the function defined by $\ov {\bf
g}(\ov a) = g_{\ov a}$. Applying Definition~\ref{10} to 
$\frak t^{w_c}$ and the function $\ov {\bf g}$ produces the
subspace ${\frak t}^{w_c}(\ov{\bf g}, k)$ as given in
Definition~\ref{barg}.

We know by Proposition~\ref{setofcoroots}
that the elements $\widetilde \Delta\spcheck_c$ are
the extended set of coroots for 
the root system   $\Phi(w_c)$ of
Definition~\ref{defnof}. Of course, $\Phi(w_c)$ has 
the same Weyl group and coroot lattice as $\Phi^{\rm proj}(w_c)$.
In particular, its coroot lattice is the orthogonal projection
$Q\spcheck$ into ${\frak t}^{w_c}$.
Corollary~\ref{kcor} applied to $\Phi(w_c)$
and $\widetilde \Delta_c$ then implies that the lattice
generated by the images under orthogonal projection to ${\frak
t}^{w_c}(\ov {\bf g},k)$ of $\widetilde \Delta_c-\widetilde I_c(k)$ is
exactly the image under orthogonal projection of $Q\spcheck$, and
that the Weyl group
of ${\frak t}^{w_c}(\ov {\bf g},k)$ with respect to the root system $\Phi^{\rm
proj}(w_c)$, or equivalently with respect to the root system
$\Phi(w_c)$,  is the group generated by reflections in the images
under orthogonal projection of the $a\spcheck\in \widetilde
\Delta_c-\widetilde I_c(k)$.
By Proposition~\ref{extension},
 the Weyl group of
$\Phi^{\rm proj}(w_c)$ is equal to the Weyl group of ${\frak
t}^{w_c}$ in $G$. Thus, we see that the group generated by the
reflections in the images under orthogonal projection of the
$a\spcheck\in
\widetilde \Delta_c-\widetilde I_c(k)$  is equal to the subgroup
of the Weyl group of ${\frak t}^{w_c}(\ov {\bf g},k)$ realized by elements
normalizing both ${\frak t}^{w_c}(\ov {\bf g},k)$ and ${\frak t}^{w_c}$.
By the previous lemma, this is the entire Weyl group of ${\frak
t}^{w_c}(\ov {\bf g},k)$ in $G$.

\subsection{The generalized Cartan matrix associated to $\widetilde
\Delta\spcheck-\widetilde I\spcheck({\bf n},k)\subset {\frak t}({\bf
  n},k)$}

Let $\Phi$ be an irreducible, but possibly
non-reduced, root system on ${\frak t}$ with extended set of simple
coroots $\widetilde \Delta$. Let $\widetilde D$ be the extended
coroot diagram of $\Phi$. Fix a  function ${\bf n}\colon \widetilde
\Delta\to {\bf N}$ of the form $n_0{\bf g}$ for some positive integer
$n_0$ and fix $k\ge 1$ dividing at least one of the integers $n_0g_a$. 
Let $\ell\colon\widetilde \Delta\spcheck\to {\bf
  R}^+$ be the length function determined by the inner product on
${\frak t}$. 
Fix $k\ge 1$ dividing one of the integers $n_0g_a$. Let $\pi_k$
denote orthogonal projection from $\frak t$ to $\frak t({\bf
n},k)$.  Consider the image under orthogonal projection of
$\widetilde
\Delta\spcheck-\widetilde I\spcheck({\bf n},k)$. 
According to Theorem~\ref{kthm}, $\pi_k$ embeds $\widetilde
\Delta\spcheck -\widetilde I\spcheck({\bf n},k)\subset {\frak
t}({\bf n},k)$ as a set of extended coroots for a reduced root
system
$\Phi({\frak t}({\bf n},k))$. In Definition~\ref{newlengths} we
defined a diagram  $\widetilde D({\bf n},k)$ with nodes 
$\widetilde \Delta\spcheck-\widetilde I\spcheck({\bf n},k)$.
On the other hand, by Corollary~\ref{kcor}, the set 
$\{\pi_k(a\spcheck): a\spcheck\in
\widetilde
\Delta\spcheck-\widetilde I\spcheck({\bf n},k)\}\subset \Phi\spcheck({\frak
  t}({\bf n},k))$ is an extended set of simple coroots for the
root system $\Phi(\frak t({\bf
n},k))$, whose Cartan integers are given by
$$n(\pi_k(a\spcheck), \pi_k(b\spcheck)) = 2\frac{\langle
\pi_k(a\spcheck), \pi_k(b\spcheck)\rangle}{\langle
\pi_k(b\spcheck), \pi_k(b\spcheck)\rangle}.$$
Let $\widetilde D_0({\bf n}, k)$ be the corresponding extended
coroot diagram. Orthogonal projection identifies the
nodes of $\widetilde
D({\bf n},k)$ with those of $\widetilde
D_0({\bf n},k)$. 

\begin{theorem}\label{samediags}
Under the above identification of the nodes, the  diagrams
$\widetilde D({\bf n},k)$ and $\widetilde
D_0({\bf n},k)$ coincide.
\end{theorem}

Proposition~\ref{newcorootdiag} is an immediate
corollary of Theorem~\ref{samediags}. We will prove
Proposition~\ref{typesofconnections} in the course of proving 
Theorem~\ref{samediags}.

Since $\widetilde I({\bf n},k)=\widetilde I(k')$ where $k'=k/{\rm
  gcd}(k,n_0)$, without loss of generality, for the rest of this
subsection,  we assume that ${\bf  n}={\bf g}$ and drop it from
the notation. In particular, we denote $\widetilde D({\bf g},k)$
and $\widetilde D_0({\bf g},k)$ by $\widetilde D(k)$ and
$\widetilde D_0(k)$.

If $k=1$, then
$\widetilde I(k)=\emptyset$, $\pi_k$ is the identity and there is
nothing to prove. Thus, 
from now on we assume that $k>1$.  This implies that $\Phi$ is
not of type $A_n$ and hence $\widetilde D$ is contractible. 
Let
$\widetilde D'(k)$ be the sub-diagram of $\widetilde
D(\Phi) = \widetilde D$ spanned by the nodes of
$\widetilde I(k)$.

\begin{lemma} Let $a\spcheck,b\spcheck\in \widetilde
  \Delta\spcheck-\widetilde I\spcheck(k)$. 
\begin{enumerate}
\item The node of $\widetilde D$
  corresponding to $a\spcheck$ is   not connected in $\widetilde D$ to any
  node of $\widetilde D'(k)$ if and only if  $a\spcheck\in {\frak
    t}(k)$. 
\item If $a\spcheck\in{\frak t}(k)$, then 
$\langle a\spcheck, b\spcheck\rangle =\langle
\pi_k(a\spcheck),\pi_k(b\spcheck)\rangle$.
\item Suppose that $a\spcheck,b\spcheck$ are adjacent nodes of
  $\widetilde D$. Then 
$\langle \pi_k(a\spcheck),\pi_k(b\spcheck)\rangle =\langle
a\spcheck,b\spcheck\rangle$.
\item If $a\spcheck$ and $b\spcheck$ are not adjacent nodes of
  $\widetilde D$ and if 
  they are not connected to a common component of $\widetilde
D'(k)$,
  then $\langle \pi_k(a\spcheck),\pi_k(b\spcheck)\rangle=0$. 
\item  If $a\spcheck$ and $b\spcheck$ are not adjacent in $\widetilde
  D$ , but $a\spcheck$ and 
  $b\spcheck$ are  connected
  to a common component of $\widetilde D'(k)$, then $\langle
  \pi_k(a\spcheck),\pi_k(b\spcheck)\rangle<0 $.  
\end{enumerate}
\end{lemma}

\begin{proof}
The subspace ${\frak u}$ is the span of the coroots inverse to
the roots represented by nodes of $\widetilde D'(k)$. As such it
decomposes as an orthogonal sum of the subspaces ${\frak u}_i$
indexed by the connected components of $\widetilde D'(k)$. The
factor ${\frak u}_i$ corresponding to a component is the subspace
of ${\frak u}$ spanned by the coroots inverse to the roots
represented by the nodes of that component.

The coroot $a\spcheck$ is contained in ${\frak t}(k)$ if and only if
it is orthogonal to all the coroots inverse to the roots
corresponding to nodes of  $\widetilde D'(k)$.  This is
equivalent to  the node of
$\widetilde D$ corresponding to $a\spcheck$  not being connected in
$\widetilde D$ to any node of $\widetilde D'(k)$. 

The second item is clear. If $a\spcheck\in {\frak t}(k)$, then
$\pi_k(a\spcheck)=a\spcheck$ and $\langle
\pi_k(a\spcheck),\pi_k(b\spcheck)\rangle=\langle
a\spcheck,b\spcheck\rangle$. 

If $a\spcheck$ and $b\spcheck$ are adjacent nodes of $\widetilde D$,
then since 
$\widetilde D$ is contractible, $a\spcheck$ and $b\spcheck$ are not
connected to a 
common component of $\widetilde D'(k)$. This means that the
projections $\pi_{\frak u}(a\spcheck)$ and $\pi_{\frak
u}(b\spcheck)$ are orthogonal, where $\pi_{\frak u}$ denotes
orthogonal projection to $\frak u$. The third item follows. 

Suppose that $a\spcheck$ and $b\spcheck$ are not adjacent and are not
connected in 
$\widetilde D$ to a common component of $\widetilde D'(k)$. 
Then there is no factor ${\frak u}_i$ of ${\frak u}$ with the property
that the orthogonal projections of both $a\spcheck$ and 
$b\spcheck$ into ${\frak u}_i$ are both non-trivial.  Thus,
$\langle
\pi_k(a\spcheck),\pi_k(b\spcheck)\rangle = \langle
a\spcheck,b\spcheck\rangle=0$. 

Lastly, if $a\spcheck$ and $b\spcheck$ are not adjacent in $\widetilde
D$ then 
$\langle a\spcheck,b\spcheck\rangle=0$ and hence 
$\langle \pi_k(a\spcheck),\pi_k(b\spcheck)\rangle =-\langle
\pi_{\frak
  u}(a\spcheck),\pi_{\frak u}(b\spcheck)\rangle$.  
Thus, we complete the proof by showing that,
under the hypothesis of Part 5, the inner product
$\langle \pi_{\frak u}(a\spcheck),\pi_{\frak
u}(b\spcheck)\rangle>0$. Let ${\frak u}_i$ be the subspace
spanned by the coroots inverse to the roots corresponding to the
nodes of the component of
$\widetilde D'(k)$ connected to both $a\spcheck$ and $b\spcheck$.
By Theorem~\ref{Antype}, the component of $\widetilde D'(k)$
connected to both
$a\spcheck$ and $b\spcheck$ is of type $A_n$ for some $n\ge
1$. Furthermore, if 
$\{a\spcheck_1,\ldots,a\spcheck_n\}$ is the set of simple coroots for
this component 
given by the nodes of $\widetilde D'(k)$,
then $a\spcheck$ is connected to a unique $a\spcheck_i$ and
$b\spcheck$ is connected to a unique
$a\spcheck_j$. Thus, $\pi_{{\frak u}_i}(a\spcheck)$ is a
negative multiple  of the fundamental coweight
$\varpi\spcheck_{a_i}$ for the root system corresponding to
$\widetilde D'(k)$, and likewise
$\pi_{{\frak u}_i}(b\spcheck)$ is a negative multiple of
$\varpi\spcheck_{a_j}$. The following  computation in
$A_n$ shows then that these vectors have positive inner product.

\begin{lemma}\label{someAncomps}  Let $a_1, \dots, a_n$ be the
simple roots in
$A_n$, ordered so that $\langle a_i, a_{i+1}\rangle =-1$ for
$1\leq i\leq n-1$, where $\langle \cdot, \cdot \rangle $ is the
standard Weyl invariant inner product. Let
$\varpi_{a_i}$ be the fundamental weight corresponding to $a_i$.
Then for $i\leq j$,
$$\langle \varpi_{a_i}, \varpi_{a_j}\rangle
=\frac{i(n+1-j)}{n+1}.$$
\end{lemma}
\proof This is a straightforward computation.
\end{proof}

\begin{corollary}\label{connection}
Let $a\spcheck,b\spcheck\in\widetilde \Delta\spcheck-\widetilde
I\spcheck(k)$. 
Then  $\langle
\pi_k(a\spcheck), \pi_k(b\spcheck)\rangle < 0$  if and only if either
$a\spcheck$ and
$b\spcheck$ correspond to nodes of $\widetilde D$ which are adjacent in
$\widetilde D$, or correspond to nodes of $\widetilde D$ which are
connected in $\widetilde D$ to a common component of $\widetilde
D'(k)$. 
\end{corollary}

In particular, the corollary tells us which nodes in $\widetilde
D_0(k)$ are connected by a bond in the diagram. This agrees with
the recipe given in Part 1 of Definition~\ref{newlengths} for $\widetilde
D(k)$. To determine the 
multiplicities and directions of the bonds in $\widetilde D_0(k)$, we
compute the lengths of the
$\pi_k(a\spcheck)$. Suppose that $a\spcheck\in \widetilde
\Delta\spcheck-\widetilde 
I\spcheck(k)$. Let
$\{a\spcheck_1,\ldots,a\spcheck_r\}\subset \widetilde 
I\spcheck(k)$ be the nodes of $\widetilde
D'(k)$ that $a$ meets. 
Since by hypothesis $\widetilde D$ is contractible, the $a\spcheck_i$ 
are simple coroots of distinct irreducible factors of $\Phi(k)$, and hence 
are mutually orthogonal.
We have $\sum_bg_bb\spcheck=0$ and $g_b\equiv 0
\pmod k$ for all $b\spcheck\notin \widetilde I\spcheck(k)$.
Since $n(a\spcheck,b\spcheck)=0$ for all $b\spcheck\in \widetilde
I\spcheck(k) -\{a\spcheck_1,\ldots,a\spcheck_r\}$, we have
\begin{equation}\label{congruence}
\sum_{i=1}^rn(a_i\spcheck, a\spcheck) g_{a_i}\equiv 0
\pmod k.
\end{equation}

Of course, in the above congruence, each $g_{a_i}$ is not divisible by
$k$. 

If $\widetilde
\Delta\spcheck-\widetilde I\spcheck(k)$ is a single node, then clearly
this node is of Type $\infty$ in the terminology of
Proposition~\ref{typesofconnections}. Thus, we can 
assume that the cardinality of $\widetilde
\Delta\spcheck-\widetilde I\spcheck(k)$ is at least $2$, and hence
that $\pi_k(a\spcheck)\not=0$ for every $a\spcheck\in\widetilde
\Delta\spcheck-\widetilde I\spcheck(k)$.
This rules out the 
root systems  $G_2$
and $BC_1$.
Thus, from now on we assume that $\widetilde D$ has
only single and double bonds.

We shall now complete the proof of Proposition~\ref{typesofconnections} and the
proof that $\widetilde D_0(k)=\widetilde D(k)$ by examining the
various possibilities for $k$. The projection of $a\spcheck
\in\widetilde
\Delta\spcheck-\widetilde I\spcheck(k)$ is a sum of coweights of
the form $\varpi_{a_i}\spcheck$, where $a_i$ is a simple root in
a root system of type $A)n$, and we shall tacitly use
Lemma~\ref{someAncomps} in the calculations below in the case
$i=j$.

\noindent
{\bf The case $k=2$.}

By Theorem~\ref{Antype}, or direct inspection in the case of
$BC_n$, all the   components of
$\widetilde D'(2)$ are of $A_1$-type. 
Number the nodes $a\spcheck_1, \dots, a\spcheck_r$ of $\widetilde
D'(2)$ connected to  $a\spcheck$  in such a way that the nodes
$a\spcheck_1,\ldots,a\spcheck_s$ are of different length from $a\spcheck$ and  
$a\spcheck_{s+1},\ldots,a\spcheck_{r}$ are of same length as $a\spcheck$.
We set $t=r-s$.
The $a\spcheck_i$ are mutually orthogonal.
Thus, the length squared of the orthogonal projection
of $a\spcheck$ into ${\frak u}$ is
$$\left(\frac{t}{4}+\frac{s}{2}\right)|a\spcheck|^2.$$ 
Since the projection of $a\spcheck$ into ${\frak u}$ has length
less than that of $a\spcheck$, it follows
that
\begin{equation}\label{s+t}
t+2s<4.
\end{equation}

If $s=0$, then  $\ell(a\spcheck) =\ell(a_i\spcheck)$ for all $i$, and
hence $n(a_i\spcheck,a\spcheck)=-1$ 
for all $i$.
  From Equation~\ref{congruence} we see that $t$ is  
even, and,  by Inequality~\ref{s+t}, that $t$ is either $0$ or
$2$. If $s=t=0$, then $a\spcheck$ is of Type 1 in the terminology
of Proposition~\ref{typesofconnections} and
$a\spcheck=\pi_k(a\spcheck)$ so that
$|\pi_k(a\spcheck)| =|a\spcheck|=\ell_2(a\spcheck)$.
If $s=0$ and $t=2$, then  $a\spcheck$ is of Type 2(i) 
and
$|\pi_k(a\spcheck|^2=|a\spcheck|^2/2=\ell^2_2(a\spcheck)$.

If $s=1$ and $a\spcheck$ is long, then 
$n(a_i\spcheck,a\spcheck)=-1$ for all $i$ and by 
Equation~\ref{congruence} we see that $s+t$ is even. Thus by
Inequality~\ref{s+t}, $s=t=1$.  Thus, $a\spcheck$ is of Type 4(i)
In this case the length squared of
the orthgonal 
projection of
$a\spcheck$ into ${\frak u}$ is $(3/4)|a\spcheck|^2$ and thus 
$|\pi_k(a\spcheck)|^2=|a\spcheck|^2/4=\ell^2_4(a\spcheck)$. 

If $s=1$ and $a\spcheck$ is short then
$n(a_1\spcheck,a\spcheck)=-2$ and $n(a_i\spcheck,a\spcheck)=-1$
for $i>1$ and hence by Equation~\ref{congruence} it follows that
$t$ is even. Inequality~\ref{s+t} implies that 
$t=0$. Thus, $a\spcheck$ is of Type 2(ii) and the length 
squared of the projection of $a\spcheck$ into ${\frak u}$ is
$(1/2)|a\spcheck|^2$ and thus 
$\pi_k(a\spcheck)|^2=|a\spcheck|^2/2=\ell^2_2(a\spcheck)$.

This shows that all $a\spcheck\in\widetilde
\Delta\spcheck-\widetilde I\spcheck(2)$ are of Type 1,
Type 2, Type 4 or Type $\infty$, and that $\widetilde
D_0(2)=\widetilde D(2)$.

\noindent
{\bf The case $k=3$.}

In this case, again by Theorem~\ref{Antype}, all the
  components of 
$\widetilde D'(3)$  are of $A_2$-type. Suppose that
  $a\spcheck\in\widetilde \Delta\spcheck-\widetilde I\spcheck(3)$ and
  that it is 
  connected to  
nodes
$\{a\spcheck_1,\ldots,a\spcheck_r\}$ of $\widetilde D'(3)$. 
Number the nodes $a\spcheck_1, \dots, a\spcheck_r$ of $\widetilde
D'(3)$ connected to  $a\spcheck$  in such a way that the nodes
$a\spcheck_1,\ldots,a\spcheck_s$ are of different length from $a\spcheck$ and  
$a\spcheck_{s+1},\ldots,a\spcheck_{r}$ are of same length as $a\spcheck$.
We set $t=r-s$.
Since the $a\spcheck_i$ are roots of  distinct irreducible components of
$\widetilde D'(3)$ and hence mutually orthogonal.
The length squared of the orthogonal projection of
$a\spcheck$ into  ${\frak u}$ is
$$\left(\frac{2s+t}{3}\right)|a\spcheck|^2.$$
Since this length squared is less than $|a\spcheck|^2$, the
possibilities are $s=1, t=0$ or $s=0, t\le 2$. By
Equation~\ref{congruence}, we see that  $s=0$ and that   $t$ is
either $0$ or $2$. When $s=t=0$ $a\spcheck$ is of Type 1 and
$\pi_k(a\spcheck)=a\spcheck$ so that
$\ell_3(a\spcheck)=|\pi_k(a\spcheck)|$.  If $s=0$ and $t=2$, then
$a\spcheck$ is of Type 3  and
$|\pi_k(a\spcheck)|^2=|a\spcheck|^2/3=\ell^2_3(a\spcheck)$.

This shows that all $a\spcheck\in\widetilde
\Delta\spcheck-\widetilde I\spcheck(3)$ are of Type 1,
Type 3, or  Type $\infty$, and that $\widetilde
D_0(3)=\widetilde D(3)$.

\noindent
{\bf The case $k=4$.}
Of course, the same kind of general arguments as above can be
made, using the fact that the components of $\widetilde D'(4)$
are all of type $A_1$ or $A_3$. However, the only case where
there is a node not of Type $\infty$   is
when
$\Phi$ is of  type
$E_8$. It can be checked directly in this case that all nodes are
of Type 1 or  Type 4(ii) or (iii) and that the lengths are as
stated. 

\noindent
{\bf The case $k>4$.}
In this case only nodes of Type $\infty$ arise.

This completes the proof that $\widetilde D_0(k)=\widetilde D(k)$, 
and hence of Theorem~\ref{samediags}. 
In the course of the proof we showed that every node of 
$\widetilde
\Delta\spcheck-\widetilde I\spcheck(k)$ is of one of the types 
listed in Proposition~\ref{typesofconnections}, thus proving that
result.

\subsection{Proof of Theorem~\ref{quotdiag}}

According to Theorem~\ref{diagram1}
there is a root system $\Phi(w_C)$ on ${\frak t}^{w_C}$ such that
the image under orthogonal projection of 
$\widetilde \Delta\spcheck_C$ is an  extended set of simple
coroots, and such that the extended coroot diagram is $\widetilde
D\spcheck/\langle C\rangle$.  By
Corollary~\ref{bargaremark}, $\ov {\bf g }$ is a positive
integral multiple of the coroot integer function on $\widetilde
\Delta_C$. Let $\Phi(w_C, k) = \Phi(\frak t^{w_C}(\ov {\bf g},
k))$.    Theorem~\ref{quotdiag} now follows by applying 
Corollary~\ref{kcor} and Theorem~\ref{samediags} to the root
system $\Phi(w_C)$, the function $\ov {\bf g }$, and the integer
$k$.

\section{The Chern-Simons invariant}

\subsection{An algebraic invariant of $c$-triples}

We introduce an invariant $CS_G({\bf x})$ of a $c$-triple
${\bf x}$, which is a refinement of the order, and which we will
 relate to the Chern-Simons invariant of the
corresponding flat bundle over the three-torus later in this
section. 

First let us record the following useful lemma concerning simple,
non-simply laced groups.

\begin{lemma}\label{nonsimply} Let $G$ be simple and non-simply
laced. Then the order of ${\cal C}G$ is at most $2$. If the
order of ${\cal C}G$ is $2$, then $\widetilde D(G)$ is
either a chain with   two multiple bonds at the ends or has one
multiple bond meeting one leaf and one trivalent vertex which
meets the remaining two leaves.
\end{lemma}
\begin{proof} The Dynkin diagram $D(G)$ of a simple, non-simply
laced group is a chain with a single multiple bond. Therefore,
$\widetilde D(G)$ is
either a chain with at most two multiple bonds or has one
multiple bond and one trivalent vertex.
The proof is then an elementary argument
involving the possible diagram automorphisms of
$\widetilde D(G)$. 
\end{proof}

Let $G$ be   simple and let $I_0^G$ be
the  unique Weyl
invariant positive definite inner product on ${\frak t}$  with
the property that
$I_0^G(a\spcheck,a\spcheck)=2$ for every short coroot
 $a$. It is easy to check that, for all roots $a$
and all $t\in
\frak t$,
$$I_0^G(a\spcheck, t) =\frac{h_a}{g_a}a(t).$$ In particular, if
$a$ is a long root, then $I_0^G(a\spcheck, t) = a(t)$.
  There is
an induced inner product on the Lie algebra of any maximal torus
of
$G$, which we also denote by $I_0^G$. Now suppose that
$G=\prod_{i=1}^rG_i$, where the
$G_i$ are simple. Set ${\frak t}_i={\frak
t}\cap {\rm Lie} (G_i)$ and define
$I_0^G=\sum_iI_0^{G_i}$.

\begin{lemma}\label{IGandIH}
Suppose that $G$ is simple and that $\widetilde \Delta$ is 
an extended set of simple roots for $G$.
Suppose that $\widetilde I\subset \widetilde\Delta$ is a proper
subset. Let $H\subseteq G$ be the semi-simple subgroup
whose complexified Lie algebra is generated by the root spaces
${\frak g}^a$ for $\pm a\in\widetilde I$, and let $\widetilde H$
be the universal covering of $H$.  Suppose that $\widetilde
H=\prod_iH_i$ is the decomposition into simple factors, and let
$\widetilde I=\coprod_iI_i$ be the corresponding decomposition.
Let ${\frak t}_H={\frak t}\cap {\rm Lie}(H)$ and
${\frak t}_i={\frak t}\cap {\rm Lie}(H_i)$. Then
${\frak t}_H=\bigoplus_i{\frak t}_i$ and  
$I_0^G|{\frak t}_H=\sum_i\epsilon_iI_0^{H_i}$,
where $d_i$ is the highest
root of $H_i$ with respect to $I_i$
and $\epsilon_i=I_0^G(d\spcheck_i,d\spcheck_i)/2$. In
particular, if $G$ is simply laced or if
$H_i$ contains a root which is a long root of $G$, then
$\epsilon_i=1$.
\end{lemma}  

\begin{proof}
The inner product $I_0^G|{\frak t}_i$ is invariant under the Weyl 
group of
$H_i$, and thus is a multiple of $I_0^{H_i}$. Clearly, this
multiple is
$I_0^G(d\spcheck_i,d\spcheck_i)/2$.  
\end{proof}

Suppose that ${\bf x}=(x,y,z)$ is a $c$-triple and that the
$c$-pair $(x,y)$ is in normal form. Let $\widetilde Z(z)$ be the
universal covering of
$Z(z)$. Choose a maximal torus $T(z)$ for $Z(z)$. Of course,
$T(z)$ is a maximal torus of $G$. Let $\frak t(z)$ be its Lie
algebra. Let
$\tilde x$, resp. $\tilde y$,  be a lift of $x$, resp. $y$ to
$\widetilde Z(z)$, and let $\tilde c= [\tilde x,\tilde y]\in 
\widetilde Z(z)$. In fact, $\tilde c\in {\cal C}\widetilde Z(z)$,
and hence there is a 
$\zeta\in {\frak t}(z)$ which projects under the
exponential mapping to
$\tilde c\in\widetilde Z(z)$.  Since $z\in {\cal C}Z(z)$, there
is an element $\hat z\in{\frak t}(z)$  whose image under the
exponential mapping is $z$.

\begin{defn} We define $CS_G({\bf x})=[I_0^G(\zeta,\hat
z)]\in{\bf R}/{\bf Z}$.  The {\sl order\/} of $CS_G({\bf x})$
is its order as an element of ${\bf
R}/{\bf Z}$.
\end{defn}

\begin{lemma}
The value of $CS_G({\bf x})$ depends only on the conjugacy
class of ${\bf x}$.
 \end{lemma}

\begin{proof}
We begin by showing that $CS_G({\bf x})$ only depends on ${\bf
x}$ and not any of the choices made above.  First  fix the maximal
torus
$T(z)$. Then the choice of
$\hat z\in \frak t(z)$ is determined up to an element
$\lambda\spcheck$
 in the coroot lattice of $G$.
Since $\zeta$ projects to an element of the center of $G$, we see 
that
$I_0^G(\lambda\spcheck,\zeta)\in{\bf Z}$. Thus $[I_0^G(\zeta,\hat
z)]$ is independent of the choice of $\hat
z$. Now fix $\hat z$ and vary $\zeta$. 
If  $\zeta'$ is another lift of $\tilde c$, then  $\zeta'-\zeta$ is an
element of the  coroot lattice of $\widetilde Z(z)$. Since $\hat
z$ projects under the exponential mapping to a central element of
$Z(x)$, we see that
$I_0^G(\zeta'-\zeta,\hat z)\in{\bf Z}$. 

Finally consider another maximal torus $T'(z)$ in $Z(z)$.
There is an element $g\in \widetilde Z(z)$ conjugating $T'(z)$ to
$T(z)$. Since $c$ and $z$ are central in $Z(z)$,
conjugation by $g$ fixes
$\tilde c$ and $z$. This establishes that $CS_G({\bf x})$ is
well-defined. Clearly, then, it is a
conjugacy class invariant.
\end{proof}

\begin{lemma}\label{CScomp}
Suppose that ${\bf x} =(x,y,z)$ and ${\bf x}'= (x',y',z')$ are
$c$-triples in
$G$ which lie in the same component of the moduli space. Then
$CS_G({\bf x})=CS_G({\bf x}')$.
\end{lemma}

\begin{proof}
Choose a maximal torus $S\subseteq T$ of
$Z^0(x,y,z)$. Choose a maximal torus $T(z)$ for $Z(z)$ with
$S\subseteq T(z)$. Let
$L=DZ(S)$. Then
$c\in L$ and there is a rank  zero
$c$-triple $(x_0,y_0,z_0)={\bf x}_0$ in $L$ and elements
$(s,t,u)\in S\times S\times S$ such that
$(x,y,z)=(sx_0,ty_0,uz_0)$. It suffices by Theorem~\ref{main} to
show that
$CS_G({\bf x})=CS_G({\bf x}_0)$.  We fix a maximal torus $T_L$
for $L$ such that
$z_0\in T_L$. Clearly
$Z_L(z_0)\subseteq Z_G(z)$.
Let $\tilde x_0,\tilde y_0\in \widetilde Z_L(z_0)$ be lifts of 
$x_0,y_0$, let $\tilde c_0=[\tilde x_0,\tilde y_0]$, and let
$\zeta_0\in {\frak t}_L(z_0)$ project under the exponential
mapping to $\tilde c_0$. If $\hat z_0\in \frak t(z_0)$
exponentiates to $z_0$, then $CS_G(x_0, y_0, z_0) =
[I_0^G(\zeta_0, \hat z_0)]$.  We can lift
$x$ to 
$\tilde x\in
\tilde S\cdot\tilde x_0$, where $\tilde S$ is the identity
component of the inverse image of $S$ in $\widetilde Z(z)$, and
similarly for
$\tilde y$. Thus $[\tilde x, \tilde y] = [\tilde x_0, \tilde y_0]
=\tilde c_0$ and hence
$CS_G(x,y,z_0)=[I_0^G(\zeta_0, \tilde z_0)] =CS_G(x_0,y_0,z_0)$.

Lastly, replace $z_0$ by $uz_0\in Sz_0$. Clearly, a lift of $uz_0$
to
${\frak t}(z)$ is given by $\hat u+\hat z_0$ where $\hat
u\in{\frak s}={\rm Lie}(S)$ exponentiates to $u$. Thus,
$CS_G({\bf x})=[I_0^G(\zeta_0,\hat u+\hat z_0)]$. But ${\frak
s}={\frak t}_L^\perp$ under the pairing $I_0^G$, so that 
$I_0^G(\zeta_0,\hat u)=0$. Hence $CS_G({\bf x})=CS_G(x,y,z_0) = 
CS_G(x_0,y_0,z_0)$. This completes the proof.
\end{proof}

Given a component $X$ of the moduli space ${\cal T}_G(c)$ of
$c$-triples in
$G$, we define $CS_G(X)$ to be $CS_G({\bf x})$, where ${\bf x}$
is any
$c$-triple whose conjugacy class lies in $X$.

\begin{proposition}\label{orderone} 
Let $(x,y)$ be a $c$-pair. Then the function
$Z(x,y)\to {\bf R}/{\bf Z}$ defined by $z\mapsto CS_G(x,y,z)$
induces a homomorphism $\pi_0(Z(x,y))\to {\bf R}/{\bf Z}$.
Hence the order of $CS_G({\bf x})$
divides the order of ${\bf x}$.
\end{proposition}

\begin{proof} The fact that $z\mapsto CS_G(x,y,z)$ descends to a
function
$\pi_0(Z(x,y))\to {\bf R}/{\bf Z}$ is immediate from
Lemma~\ref{CScomp}. To see that it is a homomorphism, since 
$\pi_0(Z(x,y))$ is cyclic  by Corollary~\ref{torsion2}, it
suffices to show that
$CS_G(x,y,z^\ell)=\ell CS_G(x,y,z)$ for all $z\in Z(x,y)$.
Let $\tilde x,\tilde y\in \widetilde Z(z)$ be lifts of $x,y$, let
$\tilde c=[\tilde x,\tilde y]$, and let $\zeta\in{\frak t}(z)$
project to $\tilde c$.  Finally, let $\tilde z$ be a lift of
$z$ to ${\frak t}_z$. Then $CS_G(x,y,z)=I_0^G(\zeta,\tilde
z)\bmod{\bf Z}$. Since $Z(z)\subseteq Z(z^\ell)$, the element
$\ell \hat z \in \frak t(z) =\frak t(z^\ell)$ is
a lift of $z^\ell$. Clearly, then, $CS_G(x,y,z^\ell)=
[I_0^G(\zeta,\ell\tilde z)]=\ell CS_G(x,y,z)$. 
\end{proof}

We determine the order of $CS_G$ in case $G=SU(n+1)$:

\begin{lemma}\label{SUn}
Suppose that $(x,y,z)$ is a $c$-triple in $SU(n+1)$, where $c$
generates the center of $SU(n+1)$. Then the order of $CS_G(x,y,z)$
is the order of
$z\in {\cal C}SU(n+1)$.
\end{lemma}
\begin{proof} This follows easily from  Lemma~\ref{someAncomps}.  
\end{proof}

\subsubsection{More on the structure of $L_c$}

The existence of the invariant $CS_G$ leads to more detailed,
classification-free information on the structure of $L_c$:

\begin{proposition}\label{moreonLc}
Suppose that  $G$ is  simple and let $c\in{\cal C}G$ be of
order
$o(c)>1$. Let $n_0={\rm gcd}\{g_{\ov a}\}$. Then:
\begin{enumerate}
\item $L_c=\prod_{i=1}^rL_i$
where each $L_i$ is a simply connected, simple group of
type $A_{n_i}$ for some $n_i\ge 1$;
\item
${\rm lcm}\,\{n_i+1: i=1,\ldots,r\}=o(c)$, and there is an
$i$ for which $n_i+1=o(c)$.
\item If $G$ is simply laced, then $o(c)=n_0$. 
\item If $G$ is not simply laced,
then either $n_0=1$, in which case each $L_i$ is of type $A_1$ and
the roots of $L_i$  are short roots of $G$,
or $n_0=2$, in which case each $L_i$ is of type
$A_1$, and exactly one of the $L_i$ has a simple root
which is a long root of $G$.
\item If $G$ is not simply laced and $n_0=1$, then $\Delta(c)$
contains the unique short simple root which is not perpendicular
to at least one long simple root.
\item If $G$ is not simply laced and $n_0=2$, then $\Delta(c)$
contains a long simple root $a$, ${\rm exp}(\varpi_a\spcheck)=c$
and all other simple roots of
$G$ are short. 
\item There is   a
$c$-triple $(x,y,t)$ in $G$ of order $n_0$, where $(x,y)$ is a
$c$-pair in normal form and $t\in T^{w_c}$, such that the order of
$CS_G(x,y,t)$ is also $n_0$.
\end{enumerate}
\end{proposition}

\begin{proof}
By Theorem~\ref{Antype}, we know that $L_c=\prod_i L_i$ is
isomorphic to a product of simply connected groups of type
$A_{n_i}$ and that $c$ is a product of elements $c_i\in{\cal
C}L_i$ generating the center of
$L_i$. It follows immediately that $o(c)={\rm
lcm}\{n_i+1:i=1,\ldots,r\}$. If $G$ is of type $A_N$ for some
$N\ge 1$, then $L_c$ is a product of simple groups of type
$A_{n-1}$ where
$n=o(c)$, and so Parts 1 and 2 hold in this case. Assume that $G$
is not of type
$A_N$ for any $N$. Then
$\widetilde D(G)$ is contractible and has at most two vertices of
order $>2$. Furthermore, if it has a vertex of order $>3$, then
it is
$\widetilde D_4$. Any diagram automorphism of such a diagram has
order $1,2,3$, or $4$.   Since the center acts faithfully on
$\widetilde D(G)$, 
$o(c)$ is divisible by at most one prime.  It
follows that  $n_i+1=o(c)$ for some $i$. This proves Parts 1 and
2. 

We number the
$L_i$ so that $n_1+1=o(c)$, and we let ${\frak t}_i={\frak t}\cap
{\rm Lie}(L_i)$. 
Suppose $G$ is simply laced. Let $t\in{\cal C}H_1$ be an element of
order $o(c)$ and let
$(x,y)$ be a $c$-pair in $L_c$. Clearly, ${\cal C}L_c\subseteq
T^{w_c}$ and hence
$t\in T^{w_c}$. Let $\tilde t\in{\frak t}_1$ lift $t$ and let
$\zeta_1\in \frak t_1$ lift $c_1$. Since
$G$ is simply laced, it follows from Lemma~\ref{IGandIH} that
$CS_G(x,y,t)=CS_{L_c}(x,y,t)=[I_0^{H_1}(\tilde t,\zeta_1)]$.
Since $t$ and $c_1$ each generate ${\cal C}H_1$, Lemma~\ref{SUn} 
shows that
$[I_0^{H_1}(\tilde t,\zeta_1)]$ is of order $n_1+1 =o(c)$. By
Lemma~\ref{orderone}, $o(c)$ divides the order of $(x,y, t)$ in
$G$. Since $t\in T^{w_c}$, the order of $(x,y, t)$ in $G$ divides
$n_0$ by Proposition~\ref{atmost}. If $\ov a$ is the orbit
containing the extended root, then $n_0$ divides $g_{\ov a} =
o(c)$. It follows that all of the above divisibilities are in fact
equalities. In particular $o(c) = n_0$ and the order of
$CS_G(x,y,t)$ is $n_0$. Thus we have proved Part 3, as well as 
Part 7 in the simply laced case.

To treat the non-simply laced case, we need the
following.

\begin{claim}
Suppose that $G$ is non-simply laced and that $c\in {\cal C}G,
c\neq 1$. Let $\ov I^G_0(c,c)\in{\bf R}/{\bf Z}$ be defined as
follows: Choose  
$\zeta,\mu\in{\frak t}$ with ${\rm exp}(\zeta)={\rm exp}(\mu)=c$,
and set
$\ov I^G_0(c,c) = [I_0^G(\zeta,\mu)]$. Then $\ov I^G_0(c,c)$ is
well-defined.  Finally, $CS_G({\bf x})=0$ for every $c$-triple
${\bf x}$ in
$L_c$ if and only if 
$\ov I^G_0(c,c)=0$ if and only if $c\in S^{w_c}$ if and only if
$T^{w_c}$ is connected.
\end{claim}

\begin{proof}
Since $\zeta,\mu\in P\spcheck$, it follows immediately that
varying  $\zeta$ and $\mu$ by elements in $Q\spcheck$ changes
$I^G_0(\zeta,\mu)$ by an integer, showing that $\ov I^G_0(c,c)$
is well-defined.

By Lemma~\ref{nonsimply} $o(c)=2$. Hence, by Part 2, all the
$L_i$ are  of type
$A_1$. Let $\{a_1,\ldots,a_t\}$ be the simple roots of $L_c$.
Then a representative for ${\rm log}(c)$ is 
$\zeta_0=\sum_i(1/2)a_i\spcheck$. Another representative for $c$
is $\mu_0=\varpi_a\spcheck$ for some simple root
$a$ of $G$. Since $\varpi_a\spcheck$ represents an element
of ${\cal C}G$, $h_a=1$, and thus
$g_a=h_a$ and  $a$ is a long root of $G$.

Let ${\bf x}_0=(x,y,c)$ be a $c$-triple in $L_c$. By definition,
for any lifts
$\zeta$ and $\mu$ of $c$ to $\frak t$, we have
$CS_G({\bf x}_0)=[I^G_0(\zeta,\mu)]$.
Thus $\ov I^G_0(c,c)=CS_G({\bf x}_0)$. Hence, if $CS_G({\bf
x})=0$ for every $c$-triple ${\bf x}$ in $L_c$, then $\ov
I^G_0(c,c)=0$. By Lemma~\ref{IGandIH},
$$I^G_0(\zeta_0,\mu_0)=
I_0^G(\sum_i(1/2)a_i\spcheck,\varpi_a\spcheck)=
(1/2)\sum_i\epsilon_i\delta_{a,a_i}.$$
Thus $I^G_0(\zeta_0,\mu_0)$
is   zero if $a$ is
distinct from all the $a_i$. Using the fact that $a$ is a  long
root and hence
$\epsilon_i=1$ if $a=a_i$, it follows that $I^G_0(\zeta_0,\mu_0)$
is equal to $1/2$ if
$a$ is equal to one
of the $a_i$.
In particular,   if $\ov I^G_0(c,c)=0$,  then
$\varpi_a\spcheck$ is orthogonal to all
the $a_i\spcheck$ and thus is orthogonal to ${\frak t}_{L_c}$.
In this case, $\varpi_a\spcheck\in{\frak t}^{w_c}$ and hence
$c\in S^{w_c}$. Thus, $\ov I_0^G(c,c)=0$ implies that
$c\in S^{w_c}$. 

 By Lemma~\ref{nonsimply},  ${\cal
C}G =\langle c\rangle$. On the other hand, by Lemma~\ref{surjpi0},
the map ${\cal C}G \to \pi_0(T^{w_c})$ is surjective. If
$c\in S^{w_c}$, then $\pi_0(T^{w_c}) = 0$ and hence $T^{w_c}$ is
connected. 

Finally suppose that $T^{w_c}$ is connected. Then, by
Corollary~\ref{Twc}, $n_0 =1$. Let ${\bf x}$ be a $c$-triple in
$L_c$. Then by Proposition~\ref{atmost}, the order of ${\bf x}$
divides $n_0$ and hence is $1$. Since the order of $CS_G({\bf
x})$ divides the order of ${\bf x}$, by
Proposition~\ref{orderone}, it follows that  $CS_G({\bf
x}) =0$.
\end{proof}

Returning to the proof of the proposition, let us suppose that $G$
is non-simply laced and that $o(c)=2$, so that each $L_i$ is of
type $A_1$.  Let $a_i\in
\Delta$ be the simple root of
$L_i$. Let $a\in \Delta$ be the unique simple root such
that $c=\exp \varpi_a\spcheck$. Since $h_a=1$, the root $a$ is
long. Let
$\zeta\in \frak t_{L_c}$ be any lift of $c\in L_c$ and let 
$z_i={\rm exp}((1/2)a_i\spcheck)$. Then
$$CS_G(x,y,z_i)=[I_0^G((1/2)a_i\spcheck,\zeta)]=
[\epsilon_i a_i(\zeta)/2].$$
Since $\zeta$ projects to a
generator of the center of $L_i$, we see that
$$CS_G(x,y,z_i)= [\epsilon_i/2].$$
On the other hand, using
the lift $\varpi_a\spcheck$ for $c$ shows that
$$CS_G(x,y,z_i)=[I_0^G((1/2)a_i\spcheck,\varpi_a\spcheck)]=
[\epsilon_i\delta_{a_i,a}/2].$$ 
Since $a$ is long, $\epsilon_i\delta_{a_i,a} = \delta_{a_i,a}$. 
It follows that $\epsilon_i\equiv \delta_{a_i,a}\bmod 2$. In
other words,  every $a_i$ distinct
from $a$ is short. Hence either no root of $L_c$ is a long root of
$G$, in which case $\ov I_0^G(c,c)=0$, or exactly one simple root
of $L_c$ is a long root of $G$ and $\ov I_0^G(c,c)\neq 0$. In the
first case, 
$T^{w_c}$ is connected, and hence $n_0=1$,
and in the second case $T^{w_c}$ has two components and so
$n_0=2$. The argument also shows that, if $a_i = a$, then the
order of
$CS_G(x,y, z_i)$ is exactly $2=n_0$. Of course, if $n_0=1$, then
the order of $CS_G(x,y, z_i)$ is also $1$. This proves Part 4, as
well as Part 7 in the non-simply laced case.

The proofs of Parts 5 and 6 are very similar, and we shall just
prove Part 6. Suppose that $G$ is not simply laced and that
one of the $a_i$, say $a_1$, is a long root of $G$. Since the
Dynkin diagram for
$G$ has a unique multiple bond, the long roots form a connected
chain. Thus,  if there is another long simple root $b\neq a_1$ of
$G$, there is a long simple root $b$ of $G$ with
$I_0^G(a_1\spcheck,b\spcheck)=-1$.
Since the $a_i$ are short roots for $i>1$, and hence $\epsilon
_i = 2$, we have that
$I_0^G(a_i\spcheck,b\spcheck)=\epsilon _ia_i(b\spcheck)\equiv
0\pmod 2$ for all
$i>1$. It follows that 
$ b(\sum_i(1/2)a_i\spcheck)=
I_0^G(b\spcheck,\sum_i(1/2)a_i\spcheck) \equiv(1/2) \bmod {\bf
Z}$. This is impossible since
${\rm exp}(\sum_i(1/2)a_i\spcheck)=c\in {\cal C}G$.
Thus $a$ is the unique long simple root  in $G$.
\end{proof}

The proof of Parts 5 and 6 actually shows the following: The
simple roots in $\Delta(c)$ are given as follows. The Dynkin
diagram of $G$ is a chain, and the node at one end is a short
root. Begin with this node, and then take every other node in
the diagram until you reach the node of the double bond. Thus,
for $G$ of type $B_n$, $\Delta (c)$ is a single node
corresponding to the short simple root, and for $G$ of type $C_n$,
$\Delta(c)$ consists of $n/2$ short simple roots if $n$ is even
and consists of $(n-1)/2$ short simple roots plus the long simple
root if $n$ is odd.

\subsubsection{Order of $CS_G({\bf x})$ in the rank zero case}

\begin{proposition}\label{CSorderzero} Let $G$ be simple and
suppose that
${\bf x}$ is a  
$c$-triple of rank zero in  $G$.
Then the order of
$CS_G({\bf x})$ equals that of ${\bf x}$. For every $k$
dividing exactly  one of the
$\{g_{\ov a}:\ov a\in\widetilde \Delta_c\}$, the function
 $CS_G$
defines a bijection between the set of conjugacy classes
of $c$-triples of rank zero and order $k$ in $G$ and the points
of order $k$ in ${\bf R}/{\bf Z}$.
\end{proposition}

\begin{proof}
Fix $k$ dividing exactly one of the integers $g_{\ov a}$. By
Proposition~\ref{k} and Proposition~\ref{orderone} it suffices to
exhibit a single $c$-triple ${\bf x}$ such
that $CS_G({\bf x})$  has order $g_{\ov a}$.

Let  $(x,y,z)$ be a 
$c$-triple of rank zero and order $g_{\ov a}$.  
By Lemma~\ref{vertex}, $z$ is conjugate to the exponential of a
vertex of the alcove
$A$. Moreover, $Z(z)$ is semi-simple and contains the rank zero
$c$-pair $(x,y)$. Thus by Proposition~\ref{rank0pair}, the
universal cover
$\widetilde Z(z)$ is a product of groups $H_i$ of type
$A_{n_i}$ and there is  a lift $\tilde c$ of $c$ to
$\widetilde Z(z)$ such that the image of
$\tilde c$ generates  the center of every simple factor.
Let ${\frak t}(z)=\bigoplus_i{\frak t}_i$ be the orthogonal
direct sum decomposition induced by the decomposition of
$\widetilde Z(z)=\prod_{i=1}^rH_i$.
We write $\zeta=\sum_i\zeta_i$ where $\zeta_i\in \frak t_i$
exponentiates to a generator of ${\cal C}H_i$.
 
 Let $\tilde z$ be a
lift of $z$ to $\widetilde Z(z)$. Then $\tilde z$ lies in the
center of $\widetilde Z(z)$. Hence, if $m$ is the least common
multiple of the integers $n_i+1$, then $\tilde z^m = 1$. Since
the image of $z$ in the group $\pi_0(Z(x,y))$ has order $g_{\ov
a}$, it follows that $g_{\ov a}|m$.
For every $z'\in {\cal C}Z(z)$ the triple $(x,y,z')$ is a 
$c$-triple. We shall find a $z'\in {\cal C}Z(z)$ such that
$CS_G(x,y,z')$ has order $m$. Supposing this, by
Corollary~\ref{orderone}, it follows that $m$ divides the order
of $(x,y,z')$. Thus, $g_{\ov a}$ divides the order of $(x,y,z')$,
which in turn divides at least one of the $g_{\ov b}$. It follows
that $g_{\ov b} = g_{\ov a}=m$ and that the order of $(x,y,z')$
is $g_{\ov a}$. Thus, the order of $CS_G(x,y,z')$ is equal to
$g_{\ov  a}$.

It remains to construct the required element $z'\in{\cal C}Z(z)$.
For any $z'$ in ${\cal C}Z(z)$,
let $\hat z'\in {\frak t}(z)$ be a lift of $z'$.
We write $\hat z'=\sum_i\hat z'_i$ with $\hat z'_i\in{\frak
t}_i$.  Clearly, 
$$CS_G(x,y,z')=[I_0^G(\zeta,\hat z')]=
\sum_i\epsilon_i[I_0^{H_i}(\zeta_i,\hat z'_i)],$$
Since $\zeta_i$ projects to a generator of the center of $H_i$ and 
$H_i$ is isomorphic to $SU(n_i+1)$, it follows from
Lemma~\ref{SUn} that, for every 
$r_i\in {\bf R}/{\bf Z}$ of order dividing $n_i+1$,
there is an element $\hat z'_i\in{\frak t}_i$,
exponentiating to an element contained in  the center
of $H_i$, such that 
$[I^{H_i}_0(\zeta_i,\hat z'_i)[=r_i\bmod{\bf Z}$.
For appropriate choices of elements  $r_i$ of order
dividing $n_i+1$, the element $\sum_ir_i$ has order $m={\rm
lcm}\,\{n_i+1: i=1,\ldots,r\}$. Consequently, 
there is an element $\hat z'=\sum_i\hat z_i'$ such that
$\sum_iI_0^{H_i}(\zeta_i\hat z'_i)$
is of order $m$ modulo ${\bf Z}$.

If $G$ is simply laced, then all the roots of $H_i$ are long
roots  of
$G$, and hence the factors $\epsilon_i$ are all one. In this 
case, the element $\hat z'\in{\frak t}(z)$ constructed
 in the last paragraph exponentiates to  an element $z'$
in the center of $Z(z)$ and $I_0^{G}(\zeta,\hat z')$
is of order $m$ modulo ${\bf Z}$.

Now suppose that $G$ is non-simply laced and $c=1$, so that
$g_{\ov a} =g_a$. We write $Z(z)=\prod H_i/\langle \zeta\rangle$
where $H_1$ is the factor  containing the highest
root of $G$.
Then, by Corollary~\ref{values},
the factor $H_1$ is of type
$A_{g_a-1}$ and the image of $\pi_1(Z(z))$ in this factor
is the center of $A_{g_a-1}$. Since one of the roots of this
factor is  the highest root of $G$, the roots of this factor are
long roots of $G$. Choose a generator $c$ of $\pi_1(Z(z))$, and
let $c_i\in H_i$ be the image of $c$ under the projection
$\widetilde Z(z)\to H_i$.
Since $H_1$ is of type $A_{g_a-1}$ and $c_1$ generates the center
of
$H_1$, there is a $c_1$-triple ${\bf x}_1$ in $H_1$ of
order $g_a$.
By Lemma~\ref{SUn},  $CS_{H_1}({\bf x}_1)$ is also of 
order $g_{\ov a}$.
 For each $i>1$ there is   a $c_i$-triple ${\bf x_i}$
in $H_i$  of order one. Clearly, the product $\prod_i{\bf x}_i$
is a $\prod c_i$-triple in $\widetilde Z(z)$, automatically
of rank zero. It  projects
to a commuting triple ${\bf x}$ of rank zero in $Z(z)$.
Clearly, $CS_G({\bf x})=\sum_i\epsilon_iCS_{H_i}({\bf x}_i)$.
Since the order of ${\bf x_i}$ is one for all $i>1$, and since
the roots of $H_1$ are long roots of $G$  we see that
$CS_G({\bf x})=CS_{H_1}({\bf x}_1)$ and hence has order $g_a$.
This proves the proposition  in   case   $c=1$ or   $G$ is simply
laced. 

There remains the possibility that $G$ is not  simply laced
and $c\not= 1$. However, as the next lemma shows, there is just
one possible $G$ in this case:

\begin{lemma}\label{C2}
 Suppose that $G$ is not  simply laced
and $c\not= 1$. If there is a rank zero $c$-triple $(x,y,z)$ in
$G$, then
$G$ is of type $C_2$.
\end{lemma}
\begin{proof} 
 By Lemma~\ref{nonsimply}, $o(c) =2$.  Since
$\widetilde Z(z)$ is a product of groups of type $A_n$, $z\not=
1$, and hence $z$ is the exponential of a vertex of the alcove
$A$ contained in a wall of $A$ corresponding to the highest root.
Let $b$ be the simple
root of $G$ so that the face
$\{b=0\}$ of
$A$ is opposite to the vertex exponentiating to $z$. By
Lemma~\ref{nonsimply}, $\widetilde D(G)$ is
either a chain with  two multiple bonds at the ends or has one
multiple bond meeting one leaf and one trivalent vertex, two of
whose ears are the remaining two leaves. Moreover, the complement
of the node corresponding to $b$ is a diagram which is a union of
diagrams of type $A_n$. It is easy to see that the only
possibilities for such extended diagrams $\widetilde D(G)$ are
$\widetilde C_2$,
$\widetilde B_3$, or $\widetilde B_4$. The three possibilities
for the quotient coroot integers are $1,2$ in the case  of 
$\widetilde C_2$, $1,2,2$ in the case  of 
$\widetilde B_3$, and $1,2,2,2$ in the case  of 
$\widetilde B_4$. Thus, by Proposition~\ref{k}, only in case
$G=C_2$ does
$G$ contain  a rank zero $c$-triple.
\end{proof}

Returning to the proof of Proposition~\ref{CSorderzero} in the
case where $G$ is of type $C_2$, let ${\bf x} =(x,y,z)$ be any
rank zero $c$-triple. Up to conjugation, it follows that
$z$ is the exponential of the vertex opposite the wall defined by
$\{b=0\}$, where $b$ is the unique short simple root of $C_2$.
Hence $Z(z) =
\widetilde Z(z)\cong H_1\times H_2$, where each $H_i$ is of type
$A_1$ and the extended root $-d$ is a simple root for one of the
$H_i$, say $H_1$.   Moreover
$c=c_1c_2$, where the $c_i$ is the nontrivial central element of
$H_i$, and 
$(x,y)$ is a product of  rank zero $c_i$-pairs $(x_i, y_i)$ in  
$H_i$. Furthermore $z$ is the exponential of $(1/2)d\spcheck$.
By Lemma~\ref{SUn}, $CS_{H_1}(x_1, y_1, z)  =1/2$. 
Since $z\in H_1$, $CS_G(x,y,z) = CS_{H_1}(x_1, y_1, z)  =1/2$. 
\end{proof}

\subsubsection{The order of $CS_G({\bf x})$}

\begin{proposition}\label{thereexists}
Let $k$ be a positive integer dividing at least one of the
$g_{\ov a}$ and suppose that
$k\not| n_0$. Let $S=S^{w_c}(k)$  and let $L=DZ(S)$.  
Then there exists a rank zero
$c$-triple ${\bf x}$ in $L$ of order $k$ in $G$ such that
$CS_G({\bf x})$ has order $k$.
\end{proposition}

\begin{proof}  We begin with the following lemma on the
structure of $L$:

\begin{lemma}\label{unmarked} Suppose that $c\neq 1$. With $L$ as
above,
$L_c$ is properly contained in $L$ and
$L$ has a unique component $L_0$ which is not of $A_n$-type for
some $n$. Write $L=L_0\times L'$. 
If $G$ is not simply laced,
then $L_0$ is of type $C_2$ and all
simple factors of
$L'$ are of type $A_1$ whose roots are short roots of $G$.
\end{lemma}

\begin{proof}
 By Proposition~\ref{atmost}, since $k$ does not divide $n_0$,
$L_c$ is properly contained in $L$ and $L$ has a unique simple
factor
$L_0$ not of
$A_n$-type. Clearly if $G$ is
not simply laced, then $L_0$ is also not simply laced. In
particular the two nodes of the double bond for $G$ are simple
roots for $L_0$. According to Proposition~\ref{moreonLc}, one
of these nodes lies in $\Delta(c)$. It follows that the
projection $c_0$ of $c$ to $L_0$ is nontrivial. Since $L_0$
contains a rank zero $c_0$-pair, it follows from Lemma~\ref{C2}
that $L_0$ is of type $C_2$.   By Proposition~\ref{moreonLc},
$L_c$ is a product of groups of type $A_1$ and has at most one
simple root which is a long root of
$G$. If there is   a long simple root in $L_c$, the
corresponding node is a node of the double bond of $D(G)$. Hence
it is a root of $L_0$. 
It follows that all of the simple factors of $L'$
are of type $A_1$ whose roots are short roots of $G$.
\end{proof}

Returning to the proof of Proposition~\ref{thereexists}, first
assume that
$G$ is simply laced. Let
$L=\prod_{i=0}^rL_i$ be the product decomposition of $L$ in
simple factors, where $L_0$ is the factor which is not of type
$A_n$. Let
$c =\prod_ic_i$ be the corresponding decomposition of $c$. For
each
$i$, suppose that we are given a $c_i$-triple of rank zero ${\bf
x}_i$ in $L_i$ of order $k_i$. Then $\prod_i{\bf
x}_i$ is a $c$-triple in $L$. Moreover, since $G$ is simply
laced, $CS_G({\bf x}) = \sum _iCS_{L_i}({\bf
x}_i)$. 

Let ${\bf x}$ be a $c$-triple in $L$ whose   order in $G$ is
$k$. Let ${\bf x}_i$ be the image of ${\bf x}$ in $L_i$ and
let $k_i$ be its order as a $c_i$-triple in $L_i$. Note that
${\bf x}_i$  has rank zero in $L_i$. Thus by 
Proposition~\ref{CSorderzero}, the order of 
$CS_{L_i}({\bf x}_i)$  is the order of ${\bf x}_i$ as a
$c_i$-triple in $L_i$, and hence
$CS_{L_i}({\bf x}_i) = [r_i/k_i]$ for some integer $r_i$
relatively prime to $k_i$.  The order of ${\bf
x}$, as a $c$-triple in $L$, is the least common multiple $\ell$
of the
$k_i$.
By Corollary~\ref{constant}, $k|\ell$.  It is an elementary
number-theoretic argument that there exist $a_i$
such that 
$[r_0/k_0 +\sum _{i\geq 1}a_i/k_i]\in {\bf R}/{\bf Z}$ has order
$\ell$. 
 For $i\geq 1$, $L_i$ is of type $A_{n_i}$, and thus, given the
integer $a_i
\bmod k_i$, there exists a rank zero $c_i$-triple ${\bf x}_i'$ in
$L_i$ such that 
$CS_{L_i}({\bf x}_i') = [a_i/k_i]$. Note that, by
Proposition~\ref{CSorderzero}, the order of
${\bf x}_i'$ in $L_i$ divides $k_i$. Replace
${\bf x}$ by the 
$c$-triple ${\bf x}' = {\bf x}_0\cdot \prod_{i\geq 1}{\bf
x}_i'$. Then 
$CS_G({\bf x}') = \sum_{i\geq 0}CS_{L_i}({\bf x}_i')$ is of
order $\ell$. On the other hand, the order of ${\bf x}'$ in $L$,
which is the least common multiple of the orders of the ${\bf
x}_i'$, divides
$\ell$. By Proposition~\ref{orderone}, the order of $CS_G({\bf
x}')$ divides the order of ${\bf x}'$ as a $c$-triple in $G$,
which in turn by Corollary~\ref{constant} divides the order of
${\bf x}'$ as a $c$-triple in $L$ which divides $\ell$ which is
the order of
$CS_G({\bf x}')$.
Therefore, the order of ${\bf x}'$ in $G$ is the order of
$CS_G({\bf x}')$, namely $\ell$.
We write ${\bf x}'=(x',y',z')$. Then ${\bf
x}''=(x',y',(z')^{\ell/k})$ is a
$c$-triple of order $k$ in $G$ such that  the order of $CS_G({\bf
x}'')$ is also $k$.

Next suppose that $c=1$. In this case $L=L_0$ is simple. By
Claim~\ref{whoknows}, there is a simple root for $L$ which is a
long root of $G$. Choose a rank zero commuting triple ${\bf x}$ in
$L$ of order $k$ in $G$. By 
Lemma~\ref{IGandIH} and Proposition~\ref{CSorderzero}, the order
of
$CS_G({\bf x}) = CS_L({\bf x})$ is $k$. 

We may thus assume that $G$ is not simply laced and that $c\neq
1$. By Lemma~\ref{unmarked}, $L_0$ is of type $C_2$ and all
simple factors of
$L'$ are of type $A_1$ whose roots are short roots of $G$. Thus
$k=2$. Let
${\bf x}= {\bf x}_0\cdot {\bf x}'$ be a $c$-triple in $L$. By
Lemma~\ref{IGandIH},
$CS_G({\bf x}) = CS_{L_0}({\bf x}_0)$. Choose a rank zero
$c_0$-triple
${\bf x}_0$ in $L_0$. Then its order is $2$ and by 
Proposition~\ref{CSorderzero}, the order of $CS_{L_0}({\bf x}_0)$
is also $2$. Choose any rank zero $c'$-triple ${\bf x}'$ in $L'$.
Its order divides $2$. Thus the order of ${\bf x}= {\bf
x}_0\cdot {\bf x}'$ in $L$ is $2$. Since the order of $CS_G({\bf
x})$ is $2$, it follows that the order of ${\bf x}$ in $G$ is
also $2$.  This concludes the proof.
\end{proof}

\begin{theorem}\label{bijection}
Let $G$ be  simple, let $c\in{\cal C}G$,
 and let ${\bf x}$ be a $c$-triple. Then the order of ${\bf x}$
is equal to the order of $CS_G({\bf x})$.
For any $k\ge 1$ dividing at least one of the $g_{\ov a}$,
the function $CS_G$ induces a bijection between the components
of ${\cal T}_G(c)$ of order $k$
and the points in ${\bf R}/{\bf Z}$ of order $k$.
\end{theorem}
\begin{proof}
Fix a $k$ dividing at least one of the $g_{\ov a}$. There is the
corresponding group $L$ containing a rank zero $c$-triple of
order $k$ in $G$. It follows from Proposition~\ref{whoknows2} that
there exists a rank zero  $c$-triple $(x_0, y_0, z_0)$ in $L$,
such that every component $X$ of ${\cal T}_G(c)$
of order $k$ in $G$ contains the conjugacy class of $(x_0, y_0,
z_0^\ell)$ for exactly one $\ell$ where $1\leq \ell \leq k$ and
$\ell$ is relatively prime to
$k$. By  Proposition~\ref{orderone}, $CS_G(x_0, y_0,
z_0^\ell) = \ell CS_G(x_0, y_0, z_0)$. Since $CS_G$ is
constant on connected components,   it
clearly suffices to find, for every $k$, a $c$-triple ${\bf x}$
of order $k$ in $G$ such that the order of $CS_G({\bf x})$ is
also $k$. This follows from
Part 7 of Proposition~\ref{moreonLc} and
Proposition~\ref{thereexists}.
\end{proof}

\subsection{Flat connections and the Chern-Simons Invariant}
\def\u{\underline}

\subsubsection{Relation of $C$-triples and flat
$G/\langle C\rangle$-connections}

Let $\Gamma$ be a flat connection on a principal $G$-bundle
$\xi$ over the three-torus $T^3 = S^1\times S^1\times S^1$.
Since $G$ is connected and simply connected, 
$\xi$ is trivial.  Choosing a trivialization of $\xi$,
the holonomy of $\Gamma$ around
the three coordinate circles is a
commuting triple $(x,y,z)$ in $G$. 
Varying the trivialization conjugates $(x,y,z)$. Thus, the
isomorphism class of the
$G$-bundle and flat connection
determines the conjugacy class of the commuting
triple. This sets up an isomorphism between ${\cal T}_G$ and
the moduli space of isomorphism classes of flat connections on
principal
$G$-bundles over the three-torus.

If $C$ is not the identity, then a $C$-triple $(x,y,z)$
in $G$ does not determine
a flat connection on a principal $G$-bundle over the three-torus.
It does determine a flat connection on a principal
$K$-bundle over the three-torus, where $K=G/\langle C\rangle$, but
the isomorphism  class of this flat connection determines and is
determined by the conjugacy class of the image
commutative triple $(\ov x,\ov y,\ov z)$ in $K$, which is not
the same as the conjugacy class in $G$ of $(x,y,z)$.
To deal with this incompatibility, we shall consider an enhanced
notion of flat connections up to isomorphism.

Fix a compact connected  group $K$. Denote its universal cover by
$\widetilde K$. Let  $M$ be a  manifold, let $\xi$ be a principal
$K$-bundle over $M$  and let $X\subset M$ be a subspace with
the homotopy type of a connected one-complex
and which carries the fundamental group. A {\sl lifting of $\xi$
over $X$\/} is a pair
$(\tilde \xi, f)$, where
\begin{enumerate}
\item $\tilde \xi$ is a principal $\widetilde K$-bundle over $X$;
\item $f\colon \tilde \xi \times _{\widetilde K}K \to  \xi|X$ is an
isomorphism of principal $K$-bundles.
\end{enumerate}
An {\sl enhanced $K$-bundle\/} $(\xi, \tilde \xi, f)$ over
$(M,X)$ consists of an underlying
$K$-bundle
$\xi$ together with a lifting $(\tilde \xi, f)$ over $X$.  An {\sl
isomorphism\/} between two enhanced $K$-bundles
$(\xi, \tilde
\xi, f)$ and $ (\xi',\tilde \xi', f')$ consists of an underlying 
$K$-bundle isomorphism
$\sigma\colon \xi \to \xi'$
and  an isomorphism
$\tilde
\sigma
\colon \tilde \xi \to \tilde \xi'$  of
$\widetilde K$-bundles such that $(\sigma|X)\circ f = f'\circ
(\tilde
\sigma \times_{\widetilde K}{\rm Id}_K)$.

\begin{lemma}\label{enhantotriv}
Given a principal $K$-bundle $\xi$ over $M$ and
a trivialization $\tau\colon X\times K \to \xi|K$,  there is an 
enhanced $K$-bundle  $\Xi_\tau=(\xi,  X\times \widetilde K,
f_\tau)$, where
$f_\tau$ is the composition $(X\times \widetilde K)\times
_{\widetilde K}K = X\times K
\buildrel{\tau}\over{\longrightarrow}\xi|K$. 
\begin{enumerate}
\item  Given an enhanced $K$-bundle $\Xi=(\xi,\tilde\xi, f)$
over $(M,X)$ there is a trivialization $\tau$ of $\xi|X$ such
that $\Xi_\tau$ is isomorphic to $\Xi$ by an isomorphism
whose underlying $K$-bundle isomorphism is the identity.
\item Given two trivializations $\tau$ and $\tau'$ of $\xi|X$,
the  enhanced $K$-bundles $\Xi_\tau$ and
$\Xi_{\tau'}$ are isomorphic if and only if the function $\kappa
\colon X \to K$ corresponding to the automorphism $\tau^{-1}\circ
\tau '$ of the trivial bundle lifts to a function from $X$ to
$\widetilde K$. 
\end{enumerate}
\end{lemma}
\begin{proof} Since $\widetilde K$ is connected, every principal
$\widetilde K$-bundle over $X$ is isomorphic to the trivial
bundle. Part 1 follows. A lifting of the $K$-bundle isomorphism 
$\tau^{-1}\circ
\tau '$ to an automorphism of the trivial $\widetilde K$-bundle
$X\times \widetilde K$ is the same as a lifting of $\kappa$ to
$\widetilde K$. This proves Part 2.
\end{proof}

Suppose that $\xi\to M$ is  $K$-bundle, that
$\tau\colon X\times K\to \xi|K$ is a trivialization, and that
$\Gamma$ is a flat conection on $\xi$. Then the flat connection
$\tau^*(\Gamma|X)$ lifts uniquely to a flat connection
$\widetilde \Gamma$ on $X\times \widetilde K$. The holonomy of
$\widetilde \Gamma$ is a homomorphism $\rho_{(\Gamma,\tau)}
\colon\pi_1(X)\to \widetilde K$, called the
{\sl $\widetilde K$-holonomy\/} of $(\Gamma,\tau)$.

\begin{lemma}\label{holconj}
Let $\xi$, resp.\ $\xi'$ be a $K$-bundle over $M$ with a flat connection
$\Gamma$, resp.\ $\Gamma'$, and let $\tau$, resp.\ $\tau'$, be a 
trivialization of $\xi|X$, resp.\ $\xi'|X$. Suppose that there
is an isomorphism $(\sigma,\tilde\sigma)$ from 
$\Xi_\tau$ to $\Xi'_{\tau'}$ with
$\sigma^*\Gamma'=\Gamma$. Then the
$\widetilde K$-holomonies
$\rho_{(\Gamma,\tau)}$ and $\rho_{(\Gamma',\tau')}$ are
conjugate by an element of $\widetilde K$.
\end{lemma}

\begin{proof}
Let $\widehat\sigma\colon X\times K\to X\times K$ be the $K$-bundle
map obtained from $\tilde \sigma$ by dividing out by
$\pi_1(K)\subseteq {\cal C}
\widetilde K$.
The connection $\tilde\sigma^*\widetilde \Gamma'$ on $X\times \widetilde K$
lifts the connection $\widehat\sigma^*f_{\tau'}^*\Gamma$ on $X\times K$. 
Since $f_{\tau'}\circ \tilde \sigma=\sigma\circ f_\tau$ and since
$\sigma^*\Gamma'=\Gamma$, 
it follows that $\widehat\sigma^*f_{\tau'}^*\Gamma'=f_\tau^*\Gamma$.
Thus, $\tilde\sigma^*\widetilde \Gamma'$ lifts $f_\tau^*\Gamma$.
But $\widetilde \Gamma$ is the unique lifting of $f_\tau^*\Gamma$.
It follows that $\widetilde \Gamma=\tilde\sigma^*\widetilde\Gamma'$,
 and hence the holonomies of $\widetilde\Gamma$ and $\widetilde 
\Gamma'$ are conjugate in $\widetilde K$.
\end{proof}

Let $(\Xi, \Gamma)$ be a pair consisting of an enhanced
$K$-bundle over $(M,X)$ and a flat connection on the
underlying $K$-bundle. 
We define the {\sl $\widetilde K$-holonomy\/} of $(\Xi,\Gamma)$ to
be the  conjugacy class of the homomorphism $\rho_{(\Gamma,\tau)}
\colon\pi_1(X)\to \widetilde K$ where $\tau\colon X\times K\to \xi|K$
is any trivialization with the property that 
there is an isomorphism from $\Xi_\tau$ to $\Xi$ which is the
identity on the underlying $K$-bundles.
According to Lemmas~\ref{enhantotriv} and~\ref{holconj}
 the $\widetilde K$-holonomy of
$(\Xi,\Gamma)$ is well-defined.
Two  pairs $(\Xi, \Gamma)$ and $(\Xi', \Gamma')$ are {\sl
isomorphic\/} if there exists   an isomorphism $(\sigma,\tilde \sigma)$
from $\Xi$ to $\Xi'$ such that
$\sigma^*\Gamma' = \Gamma$.
We have established the following.

\begin{proposition}\label{tildeKhol}
Let $M$ be a manifold and $X\subset M$ a subset with the homotopy
type of a one-complex carrying the fundamental group of $M$. The
$\widetilde K$-holonomy of a pair 
$(\Xi,\Gamma)$ consisting of an enhanced $K$-bundle over $(M,X)$
and a flat connection $\Gamma$ on the underlying $K$-bundle is a
conjugacy class of homomorphisms from $\pi_1(X)$ to $\widetilde
K$. This $\widetilde K$-holomony  depends only on the isomorphism
class of the pair.
\end{proposition}

 In case when $M$ is the $N$-torus $T^N =
(S^1)^N$ and  $X\subset T^N$ deformation retracts onto the union
$\bigvee_NS^1$ of the
coordinate circles, the $\widetilde K$-holonomy of a pair $(\Xi,
\Gamma)$ is a conjugacy class of homomorphisms from
$\pi_1(X)$ to
$\widetilde K$. We identify this conjugacy class with a conjugacy
class of $N$-tuples $(x_1,
\dots, x_N)
\in (\widetilde K)^N$, where $x_i$ is the image of the element in
$\pi_1(X)$ represented by the $i^{\rm th}$ coordinate circle. 

\begin{proposition}\label{enhantotrip}
Suppose that $G$ is the universal cover of $K$.
Let $T^N = (S^1)^N$ be the $N$-torus and let $X\subset T^N$
deformation retract onto 
the union $\bigvee_NS^1$ of the coordinate circles.
Suppose that $\Xi$ is an
enhanced $K$-bundle on $(T^N,X)$ with underlying bundle $\xi$,
and suppose
that $\Gamma$ is a flat connection on $\xi$. Then the  
$G$-holonomy of $(\Xi, \Gamma)$ is a conjugacy class of  almost
commuting $N$-tuples $(x_1, \dots, x_N) \in G^N$. If $c_{ij}=
[x_i, x_j]$, then $c_{ij}\in \pi_1(K) \subseteq {\cal C}G$ and is
equal to the $2$-dimensional characteristic class of $\xi$
evaluated on the
$(ij)^{\rm th}$-coordinate two-torus.
The map which associates to  the pair $(\Xi, \Gamma)$,
consisting of an enhanced $K$-bundle over $(T^N,X)$ and a flat 
connection on the underlying $K$-bundle, its
$G$-holonomy is a bijection from the set of isomorphism classes
of such pairs $(\Xi, \Gamma)$ to the space of conjugacy classes of
almost commuting  $N$-tuples in $G$. 
\end{proposition}

\begin{proof}
Given a flat $K$-connection on the two-torus with holonomy  $\ov
x,\ov y$ around the coordinate circles, let $x,y\in G$ be lifts
of these elements. Then the commutator $[x,y]\in G$ lies in
$\pi_1(K)\subseteq{\cal C}G$ and is equal to the  characteristic
class
$w(\ov\xi)\in H^2(T^2;\pi_1(K))$. Applying this to the various
coordinate two-tori inside the $N$-torus, shows that the
$N$-tuple associated to $(\Xi, \Gamma)$ proves the first
statement.  The rest of the proposition is a straightforward
exercise.
\end{proof}

Next we show that two enhanced $K$-bundles over
$(T^3, X)$ whose underlying $K$-bundles are isomorphic are
isomorphic as enhanced $K$-bundles.

\begin{lemma}\label{autoenhance}
 Let $\xi$ be a $K$-bundle over
$T^3$. Suppose that $\tau_1$ and $\tau_2$ are two trivializations
of
$\xi |X$. Then there is a bundle automorphism  $\sigma$ of $\xi$
such that $\sigma^*\tau _1 = \tau _2$.
\end{lemma}
\begin{proof} Since the question only involves the homotopy type
of the pair $(M, X)$, we may assume that
$X=\bigvee_3S^1$. Over
$X$, we can take $\sigma = \tau _1\circ \tau_2^{-1}$. Let
$\kappa\colon X \to K$ be the corresponding map. Since $\xi|X$ is
trivial, so is ${\rm Aut} (\xi)|X$, and hence the fundamental
group of $T^3$ acts trivially on the fundamental group of the
fiber of
${\rm Aut}(\xi)$, which is isomorphic to $K$.  Thus the 
obstructions to extending
$\sigma$ to $T^3$ lie in
$H^i(T^3, X; \pi_{i-1}(K))$. The only nonzero such group is
$H^2(T^3, X; \pi_{1}(K))$. The value of the obstruction on a
relative $2$-cell
$e$ is $\kappa_*[\partial e]\in \pi_1(K)$. For  the standard
relative cell decomposition of $(T^3, X)$ for which the
$2$-skeleton is the union of the $T_{ij}$ and the corresponding
cells are  $e_{ij}$, the elements
$[\partial e_{ij}]$ are commutators in $\pi_1(X)$. Since
$\pi_1(K)$ is abelian, $\kappa_*[\partial e_{ij}]$ is trivial.
Hence the automorphism extends.
\end{proof}

\begin{corollary}\label{Xibiject}
Let $C$ be an anti-symmetric $3\times 3$ matrix with values in
$\pi_1(K) \subseteq {\cal C}G$. Let $\Xi$ be an enhanced
$K$-bundle over $(T^3, X)$ with $C(\Xi) = C$. Then there is a
bijection from the set of flat connections modulo
automorphism of $\Xi$ to ${\cal T}_G(C)$.
\end{corollary}

\subsubsection{The Chern-Simons invariant}

Let $G$ be simple and let $\hat I_0(t) =
I_0^G(t,t)$. 

\begin{lemma}\label{form} For all $t\in {\frak
t}$,
$$\frac12\hat I_0(t) =
\frac{1}{2}I_0^G(t,t)=-\frac{1}{16\pi^2g}{\rm Tr}({\rm
ad}(t)^2),$$ where  
$g$ is the dual Coxeter number of $G$.
Equivalently,
$$\hat I_0(s,t) =-\frac{1}{16\pi^2g}{\rm Tr}(2{\rm ad}(s)\cdot
{\rm ad}(t)).$$
\end{lemma}
\begin{proof} By a result of Looijenga \cite{Loo} Lemma (1.2),
$$\sum _{a\in \Phi} a(t)^2 = 2gI_0(t,t).$$
On the other hand, the action of ${\rm ad}\,t$ on $\frak g$ has
eigenvalue $2\pi\sqrt{-1}a(t)$ on the root space $\frak g^a$. The
result is then a direct computation.
\end{proof}

The  Weyl invariant
 quadratic form $\hat I_0$ on ${\frak t}$ has a unique extension
to an
 ${\rm ad}\,G$-invariant quadratic form on ${\frak g}$, also
denoted by $\hat I_0$. It follows from Chern-Weil theory that the
quadratic form $\frac12\hat I_0$ defines a cohomology class
$\tilde c_2
\in H^4(BG;{\bf R})$. By \cite{MT}, VI \S 6 Theorem 6.23 (see also
\cite{Bott}),  
$H^4(BG;{\bf Z})\cong {\bf Z}$ and $\tilde c_2$ is an integral
generator for $H^4(BG;{\bf Z})$. Let  $M$ be a manifold and let
$\xi$ be  a principal $G$-bundle over $M$. There is a
corresponding classifying map $p\: M \to BG$, and we set
$c_2(\xi) =p^*\tilde c_2$. Given a connection $A$  on $\xi$,
the image of $c_2(\xi)$ in $H^4(M; {\bf R})$ is represented by
the  closed $4$-form $\frac12\hat I_0(F_A)$.

There is a secondary characteristic   class 
associated to $c_2$, the {\sl Chern-Simons invariant}. Let $M$ be
a three-manifold and
$\xi\to M$ a principal $G$-bundle. Let $\Gamma$ and $A$ be
connections on
$\xi$. Then  
\begin{equation}\label{CSdefn}
{\rm CS}_\Gamma(A)=-\frac{1}{16\pi^2g}\int _M{\rm
  Tr}\left(2a\wedge F_\Gamma+  a\wedge 
  d_\Gamma(a)+
\frac{2}{3}a\wedge (a\wedge a)\right),
\end{equation}
where $a=A-\Gamma\in \Omega^1(M;{\rm ad}\,\xi)$.
Defined in this manner, ${\rm CS}_\Gamma$ is a real-valued
function on the space of connections on a given bundle. The flat
connections are the critical points of ${\rm CS}_\Gamma$, and
hence ${\rm CS}_\Gamma$ is constant on continuous path of
flat connections.

Since ${\rm CS}_\Gamma$ is a secondary class associated
to the primitive integral class $c_2$, if $\Gamma$ and $\Gamma'$
are gauge equivalent connections then
${\rm CS}_\Gamma(A)-{\rm CS}_{\Gamma'}(A)$ is an integer.
As we vary over all connections $\Gamma'$ gauge equivalent to
$\Gamma$, we can vary the Chern-Simons invariant by an arbitrary
integer. Similarly, if $A$ and $A_1$ are gauge equivalent
connections then
${\rm CS}_\Gamma(A)-{\rm
  CS}_\Gamma(A_1)$ is an integer, and as we vary $A_1$ over
all connections gauge equivalent to $A$ we can vary 
${\rm CS}_\Gamma$ by an arbitrary integer.
In other words,  ${\rm CS}_\Gamma$ is a  well-defined
function from the space of gauge equivalence classes of
connections into
${\bf R}/{\bf Z}$. 

More generally, suppose that $K$ is a compact group whose
universal cover is
$G$ and that $\xi\to M$ is a principal $K$-bundle with connection
$A$. The form
$\frac12\hat I_0(F_A)$ is a closed $4$-form representing a
characteristic class of $\xi$, which we will also denote by
$c_2(\xi)$. The only difference is that this class need not be
integral. In fact, it lies in $(1/n){\bf Z}$, where the image
of $H^4(BK; {\bf Z})$ in $H^4(BG; {\bf Z})$ is generated by
$n\tilde c_2$.  Given
$K$-connections $\Gamma$ and $\Gamma'$ on $ \xi$,  we
can still define ${\rm CS}_{ \Gamma}( \Gamma')\in{\bf R}$ 
by Equation~\ref{CSdefn}.  As before, this is a secondary
characteristic class associated to the four-dimension class $c_2$
of $K$-bundles. The difference is that,  in general, if
$\Gamma''$ is a $K$-connection gauge equivalent to $\Gamma'$, then
$${\rm CS}_{ \Gamma}( \Gamma')-{\rm CS}_{
  \Gamma}( \Gamma'') \in\frac{1}{n}{\bf Z},$$
and hence ${\rm CS}_{\Gamma}$ is only  well-defined
modulo $\frac1n {\bf Z}$ for isomorphism classes of
$K$-connections. As the next lemma shows,
${\rm CS}_{\Gamma}$ is well-defined modulo
${\bf Z}$ on enhanced isomorphism classes of flat $K$-connections. 

\begin{lemma}\label{integerdiff}
Let $\Xi$  be an enhanced $K$-bundle over $(M, X)$, where $M$ is a
closed, oriented three-manifold, and suppose that there is an 
automorphism of
$\Xi$ such that the induced automorphism  of the underlying
$K$-bundle $\xi$ is
$\sigma$. Then 
$${\rm CS}_{\Gamma}(\Gamma')-{\rm CS}_{
\Gamma}(\sigma^*\Gamma')\in {\bf Z}.$$
\end{lemma}

\begin{proof}
Let $\hat\xi\to M\times S^1$ be the $K$-bundle obtained from
$\xi\times I\to M\times I$ by gluing the ends together by
$\sigma$. 
The difference ${\rm CS}_{\Gamma}(\Gamma')-{\rm CS}_{
  \Gamma}(\sigma^*\Gamma')$ is equal to $\int_{M\times
S^1}c_2(\hat\xi)$.
We choose a trivialization $\tau$ of $\xi|X$ such that there
is an isomorphism $\Xi_\tau$ to $\Xi$ whose underlying $K$-bundle
isomorphism is the identity. Then in this trivialization 
$\sigma|X$ is given by a continuous map $\kappa\colon X\to K$.
The fact that 
$\sigma$ comes from an automorphism of the enhanced bundle means
that 
$\kappa$ lifts to a map $X\to G$. This means that 
$\hat\xi|X\times S^1$ lifts to a $G$-bundle, and hence, since $G$
is connected and simply connected,
$\hat \xi|X\times S^1$ is trivial.
A standard obstruction theory argument then shows that
$\hat\xi$ is isomorphic as a $K$-bundle to the connected sum of the
product bundle $\xi\times S^1\to M\times S^1$ and a $K$-bundle
$\eta\to S^4$.  Thus, $\int_{M\times S^1}c_2(\hat
\xi)=\int_{S^4}c_2(\eta)$. But since $S^4$ is simply connected,
$\eta\to S^4$ lifts to a $G$-bundle $\tilde \eta$ and thus
$c_2(\eta)=c_2(\tilde\eta)$ takes an integral value on $S^4$.
\end{proof}

\begin{corollary}
Fix an enhanced $K$-bundle $\Xi$ over $(M, X)$, where $M$ is a
closed, oriented three-manifold, and a flat
connection
$\Gamma_1$ on  $\Xi$. Then the
function
${\rm CS}_{\Gamma_1}$ induces a well-defined   function from
the set of isomorphism classes of pairs $(\Xi, \Gamma)$, where 
$\Gamma$ is a flat connection on $\Xi$,  to ${\bf R}/{\bf Z}$.
\end{corollary}

One important property of the Chern-Simons invariant is the following
additivity property:

\begin{lemma}\label{CSadd}
Let
$\Gamma_0$, $\Gamma_1$ and $A$ be connections on a $K$-bundle
$\xi$ over  a three-manifold.  Then
$${\rm CS}_{\Gamma_0}(A)= {\rm CS}_{\Gamma_1}(A)+{\rm CS}_{\Gamma_0}(\Gamma_1).$$ 
\end{lemma}

\begin{corollary}\label{CSconstant}
 Let $\{\Gamma_t\}_{t\in
[0,1]}$ and
$\{A_t\}_{t\in [0,1]}$ be two continuous paths of flat
connections. Then ${\rm CS}_{\Gamma_0}(A_0) = {\rm
CS}_{\Gamma_1}(A_1)$. 
\end{corollary}
\begin{proof}  Since the critical points of the
functional $A\mapsto {\rm CS}_{\Gamma_0}(A)$ are the flat
connections, ${\rm CS}_{\Gamma_0}(A_0) = {\rm
CS}_{\Gamma_0}(A_1)$. In particular, ${\rm
CS}_{\Gamma_0}(\Gamma_1) =0$. The corollary now follows from the
additivity formula of the previous lemma.
\end{proof}

\subsection{The basic computation}

Let $(x,y,z)$ be a $c$-triple in a simply connected group $G$.
Then $Z(z)$ is a connected group containing $x,y$ and $c$.
We denote by $\ov
Z(z)$  the quotient $Z(z)/\langle c\rangle$ and we denote by
$\widetilde Z(z)$ the universal covering of $Z(z)$.
Let $\tilde x,\tilde y$ be lifts of $x,y$ to the universal covering
$\widetilde Z(z)$.  Let $\zeta =[\tilde x, \tilde y]$ and let
$\tilde \zeta\in \frak z(z)$ be an element whose exponential is
$\zeta$.   Let
$T^2$ be the torus
${\bf R}^2/{\bf Z}^2$. Let $D\subseteq T^2$ be a closed disk and
let
$T_0$ be the closure  of $T^2-D$. Let $\partial$ be the boundary
of $T_0$, with coordinate $\theta$.

\begin{lemma} With notation as above, there is a flat connection
$\tilde A_0$ on the trivial $\widetilde Z(z)$-bundle over $T_0$
such that the holonomy of $\tilde A_0$ along the two coordinate
circles in
$T_0$ is given by $(\tilde x, \tilde y)$, and $\tilde A_0|\partial
=
\tilde
\zeta \, d\theta$.
\end{lemma}
\begin{proof} Clearly, there exists a flat connection $A_0'$ with
the required holonomy on the  bundle $T_0\times \widetilde
Z(z)$.  The connections $A_0'|\partial$ and $\tilde
\zeta \, d\theta$ have the same holonomy. It is then easy to see
that there is an automorphism $\sigma$ of the trivial bundle,
supported near $\partial$, such that $\tilde A_0 =\sigma^*A_0'$
satisfies the conclusions of the lemma.
\end{proof} 

There is an induced flat $\ov Z(z)$-connection $\ov A_0$ on
$\xi_0 =T_0
\times
\ov Z(z)$, whose holonomy along $\partial$ is the identity. Thus
there is a trivialization $\tau_\partial$ of $\xi_0|\partial$ for 
which 
$\ov A_0|\partial$ is the product connection, and hence we can extend
$\xi_0$ and the flat connection $\ov A_0$ to a $\ov Z(z)$-bundle
$\xi'$ with a flat connection $A_0$ over $T^2$. By construction,
there is a trivialization $\tau_D$ of  $\xi'|D$
extending $\tau_\partial$ and with respect to  which
$\ov A_0|D$ is a trivial connection.
We denote by $\tau_0$ the given trivialization  of $\xi'|T_0
=\xi_0$.

\begin{lemma}
 The composition $(\tau_0|\partial)\circ (\tau_D^{-1}|\partial)$
  is the map 
$S^1\times \ov Z(z)\to S^1\times \ov Z(z)$ given by $(\theta,g)\mapsto
(\theta,\zeta^\theta\cdot g)=(\theta,{\rm exp}(\theta
\tilde\zeta)\cdot g)$.  
\end{lemma}

\begin{proof}
The composition $(\tau_0|\partial)\circ (\tau_D^{-1}|\partial)$ 
pulls back the trivial connection to $\tilde \zeta d\theta$.  
The composition is given by left multiplication by $\gamma\colon
S^1\to
\ov Z(z)$, and moreover
$\gamma^{-1}d\gamma=\tilde \zeta d\theta$. 
By changing the trivialization $\tau_D$ by a constant
change of gauge, we can arrange that $\gamma(0)=1$. It follows
that
$\gamma(\theta)={\rm exp}(\theta \tilde \zeta)=\zeta^\theta$.   
\end{proof}

Now consider $T^3= T^2\times S^1$ where  $S^1= {\bf R}/{\bf
Z}$ has coordinate $u$.  The trivialization 
$\tau_0\times S^1$ restricts to a trivialization
$\tau$ of $\xi\times S^1$ over
$X=T_0\bigvee S^1\subset T_0\times S^1$.

There is a maximal torus
$T$ of
$Z(z)$, and hence  of
$G$, whose Lie algebra contains $\tilde \zeta$.  
Of course, $T$ also contains the central element $z$ of $Z(z)$.
Let
$\hat z$ be an element of ${\frak t}={\rm Lie}(T)$ such that
$\exp \hat z =z$. For $u\in{\bf R}$, define
$z^u={\rm exp}(u\hat z)$.  
For
all $u,\theta \in{\bf R}$, the elements $z^u$ and $\zeta^\theta$
lie in
$T$ and hence commute.  

  Define a
connection $A$ on $\xi'\times S^1$ over
$T^2\times S^1$ as follows. 
Over $D\times S^1$, $A$ has connection $1$-form $\hat z\,du$ with
respect to the trivialization
$\tau_D\times {\rm
  Id}_{S^1}$. 
Over $T_0\times S^1$, $A$ has connection $1$-form
$z^{-u}\tau_0^*A_0z^u+\hat zdu$ with respect to the trivialization
$\tau_0\times S^1$. Direct computation shows that the
connections given 
on $D\times S^1$ and on $T_0\times S^1$ are flat. Since $z^u$ and
$\zeta^\theta$ commute for all $u$ and $\theta$, it follows that 
these two partial connections glue together to define a
connection $A$ on $\xi'\times S^1$. Clearly, $A$ is flat, the
restriction of
$A$  to $\xi'\times \{0\}$ is isomorphic to $A_0$, and the
$\widetilde Z(z)$ holonomy of $(A,\tau)$ 
 around the three circle factors gives
the commuting triple $(\tilde x,\tilde y,\exp \hat z)$ in
$\widetilde Z(z)$. 

Let us define the $K$-bundle $\xi = \xi'
\times _{\ov Z(z)}K$ over $T^2$. The connection $A$ induces a flat
connection on $\xi\times S^1$, which we continue to denote by
$A$. The trivialization $\tau$ of $\xi'\times S^1|X$
induces a trivialization 
of $\xi\times  S^1|X$ which we will also denote by $\tau$.
 The $G$-holonomy of $(A,\tau)$
 around the three coordinate circles is $(x,y,z)$.

\begin{lemma}\label{calculate} 
Fix a connection $\Gamma$  on $\xi'\times S^1$ with the
following properties:
\begin{enumerate}
\item If $p\colon \xi'\times S^1\to \xi'$ is the natural
  projection, then $\Gamma=p^*\Gamma_0$ for some connection
$\Gamma_0$ on
  $\xi'$. 
\item The restriction of $\Gamma_0$ to $\xi_0=\xi'|T_0$ is trivial
in the
  trivialization $\tau_0$. 
\item The restriction of $\Gamma_0$ to $\xi'|D$ is a $T$-connection
in
  the trivialization $\tau_D$. 
\end{enumerate}
Also denote by $\Gamma$ the resulting $K$-connection on
$\xi\times S^1$. Then
$${\rm CS}_\Gamma(A)=I_0^G(\hat z, \tilde
\zeta).$$    
\end{lemma}

\begin{proof}
Let $A'$ be the connection on $\xi\times S^1$ which is trivial 
over
$D\times S^1$ in the trivialization $\tau_D\times S^1$ and which 
is given by the one-form $z^{-u}A_0z^u$ over
$T_0\times S^1$  in the
trivialization $\tau_0\times S^1$. As before, one sees easily that
these two descriptions match over $\partial D\times S^1$.

Let us compute ${\rm CS}_{A'}(A)$. First notice that $a=A-A'=\hat
zdu$  over both $D\times S^1$ and $T_0\times S^1$. Clearly, then
$a\wedge a=0$. Over $D\times S^1$ the connection $A'$ is trivial.
Thus, on this patch $d_{A'}(a)=F_{A'}=0$, and consequently the
Chern-Simons integrand vanishes on this patch. Over $T_0\times
S^1$ we have
$d_{A'}(a)=[z^{-u}A_0z^u,\hat zdu]$ so that $a\wedge
d_{A'}(a)=0$. Lastly, on this patch 
$F_{A'}=du\wedge [z^{-u}A_0z^u ,\hat z]$ so that $a\wedge
F_{A'}=0$. This shows that ${\rm CS}_{A'}(A)=0$ and hence by
Lemma~\ref{CSadd}
${\rm CS}_\Gamma(A)={\rm CS}_\Gamma(A').$

Now we compute ${\rm CS}_\Gamma(A')$ First let us show that the
Chern-Simons integrand for this invariant is zero over $D\times 
S^1$. Over $D\times 
S^1$, and with the  trivialization $\tau_D\times S^1$, $\Gamma$ is
the pullback of a connection $B$ on
$\xi|D$ and $A'$ is trivial. 
Thus $a=A'-\Gamma=-B$, 
$F_B$ and $d_B(a)$ are all pulled back from forms on $D$.
Thus, we see that over $D\times S^1$ the
Chern-Simons integrand vanishes identically.

Now we compute the Chern-Simons integral over $T_0\times S^1$ using
the trivialization $\tau_0\times S^1$. The connection $\Gamma$ is
trivial on this patch and the 
one-form $a$ is $z^{-u}\tau_0^*A_0z^u$.
It follows that $d_{\Gamma}(a)=du\wedge
[z^{-u}\tau_0^*A_0z^u,\hat z]+z^{-u}\tau_0^*dA_0z^u$. 
Thus, 
\begin{eqnarray*}
{\rm CS}_\Gamma(A') & = &-\frac{1}{16\pi^2g}\int_{T_0\times S^1}{\rm
  Tr}\left(z^{-u}\tau_0^*A_0z^u\wedge (du\wedge 
  [z^{-u}\tau_0^*A_0z^u,\hat z]+z^{-u}\tau_0^*dA_0z^u)\right)\\
 & = & -\frac{1}{16\pi^2g}\int_{T_0\times S^1}{\rm
  Tr}\left(z^{-u}\tau_0^*A_0z^u\wedge
  (du\wedge[z^{-u}\tau_0^*A_0z^u,\hat z])\right) \\
& = &  -\frac{1}{16\pi^2g}\int_{T_0\times S^1}{\rm
  Tr}\left(-2(z^{-u}\tau_0^*A_0z^u)\wedge (z^{-u}\tau_0^*A_0z^u)\wedge \hat zdu\right)
\end{eqnarray*}
Since $z^{-u}\hat zz^u=\hat z$, we can rewrite this as
$${\rm CS}_\Gamma(A')= -\frac{1}{16\pi^2g}\int_{T_0\times S^1}{\rm
  Tr}\left(-2\tau_0^*A_0\wedge \tau_0^*A_0\wedge \hat zdu\right).$$
Doing the $u$-integration and using the fact that $dA_0 +
A_0\wedge A_0 =0$ yields
\begin{eqnarray*}
{\rm CS}_\Gamma(A') & = & -\frac{1}{16\pi^2g}\int_{T_0}{\rm
  Tr}\left(-2(\tau_0^*A_0\wedge \tau_0^*A_0)\cdot \hat z\right) \\
& = & -\frac{1}{16\pi^2g}\int_{T_0}{\rm 
Tr}\left(2\tau_0^*dA_0\cdot \hat z\right) \\ & = &
-\frac{1}{16\pi^2g}\int_{\partial T_0}{\rm
  Tr}\left(2\tau_0^*A_0 \cdot \hat z\right).
\end{eqnarray*}  
Since $\tau_0^*A_0|\partial T_0=\tilde \zeta d\theta$, and
using Lemma~\ref{form},  we have
$${\rm CS}_\Gamma(A)={\rm CS}_\Gamma(A')=-\frac{1}{16\pi^2g}{\rm
Tr}(2\tilde \zeta\cdot \hat z)=[I_0(\hat
z,\tilde\zeta)]=CS_G(x,y,z).$$ 
\end{proof}

\begin{proposition}\label{enhCS}
Let $\Theta$ be a connection on $\xi\times
S^1$ which is pulled back from a connection on $\xi$
via the natural projection mapping. Then
$${\rm CS}_\Theta(A)=I_0^G(\hat z,\tilde\zeta).$$ 
\end{proposition}

\begin{proof}
By Lemma~\ref{CSadd} it suffices to show that if $\Theta$ and
$\Theta'$ are connections on $(\xi\times _{\ov Z(z)}K)\times S^1$
pulled back from connections on $(\xi\times _{\ov Z(z)}K)$ then
${\rm CS}_\Theta(\Theta')=0$. But this is clear -- under this
hypothesis the Chern-Simons integrand (which is a three-form)
is pulled back from a form on  the two-torus $T^2$ and hence
vanishes identically.
\end{proof}

\begin{corollary}\label{CScorol}
Let $T^3=T^2\times S^1$  and let $X\subset T^3$ 
be $T_0\vee S^1$.
Given $c\in {\cal C}G$ let $K=G/\langle c\rangle$.
Fix a $K$-bundle $\xi_c\to T^2$ with
$w(\xi_c)=c\in H^2(T^2;\pi_1(K))\subseteq {\cal C}G$.
Fix a trivialization $\tau_0$ of $\xi_c|T_0$ and
let $\tau_c$ be the trivialization of
$(\xi_c\times S^1)|X$
obtained from $\tau_0$ and the given product structure
on the last $S^1$-direction. Let $\Xi_c$ be the corresponding
enhanced $K$-bundle over $(T^3, X)$.
\begin{enumerate}
\item 
For every $c$-triple ${\bf x}$, there is a flat connection
$A_0({\bf x})$ on $\xi_c\times S^1$ such that the $G$-holonomy of
$(A_0({\bf x}), \tau_c)$ is the conjugacy class of  ${\bf x}$.
\item 
For every flat connection $A$ on $\xi_c\times S^1$ such
that the 
$G$-holonomy of $(A, \tau_c)$ is the conjugacy class of  ${\bf
x}$, and for  every connection $\Theta$ on $\xi_c\times S^1$
which is pulled back by the natural projection from a connection
$\Theta_0$  on
$\xi_c$,
$${\rm CS}_{\Theta}(A)\equiv  CS_G(x,y,z)\bmod {\bf Z}.$$
\end{enumerate} 
\end{corollary}

\begin{proof}
The above construction shows that, given a $c$-triple ${\bf x}$,
 there is a $K$-bundle $\xi\to T^2$, a trivialization
$\tau'_0$ of $\xi|T_0$  and a flat connection $A'({\bf x})$ on
$\xi\times S^1$ such that the $G$-holonomy of $A'({\bf x})$
measured using the  trivialization $\tau'$, which is the union of
the trivialization $\tau'_0$ on $T_0$ with the product
trivialization around the third
coordinate circle, is ${\bf x}$.
Since $(x,y)$ is a $c$-pair, there is an isomorphism
$\psi\colon\xi_c\to\xi$.  By Lemma~\ref{autoenhance}, there is a
bundle isomorphism from $\xi_c\times S^1$ to $\xi \times S^1$
which carries the trivialization $\tau'$ to $\tau$.   

Let $A_0({\bf x})$ be the pullback of the connection $A'({\bf
x})$.  Then $A_0({\bf x})$ is a flat connection on $\xi_c\times
S^1$ whose
$G$-holonomy  is the conjugacy class of  ${\bf x}$, and hence the
$G$-holonomy of $(\Xi_c,A_0({\bf x}))$ is the conjugacy class of 
${\bf x}$.

By Proposition~\ref{enhCS}
 ${\rm CS}_\Theta(A_0(x,y,z))\equiv CS_G(x,y,z)\bmod {\bf Z}$.
 More generally, 
suppose that $A$
is a flat connection on $\xi_c\times S^1$ such that the
$G$-holonomy of
$(\Xi_c,A(x,y,z))$ is the conjugacy class of ${\bf x}$.
Then by Lemma~\ref{enhantotrip}
 there is an enhanced automorphism of $\Xi_c$
carrying $A$ to $A_0(x,y,z)$. By Lemma~\ref{integerdiff}
 this implies that
$${\rm CS}_\Gamma(A)-{\rm CS}_\Gamma(A_0(x,y,z))\in {\bf Z},$$
and the result follows.
\end{proof}

Notice that the enhancement over $T_0\subset T^2$ is irrelevant --
it is only the enhancement in the last $S^1$-direction that is important.
This is the connection analogue of the fact that if ${\bf
x} = (x,y,z)$ is a $c$-triple, then multiplying $x$ and
$y$ by powers of $c$ does not change the $G$-conjugacy class of 
${\bf x}$ and hence does not change its $CS_G$-invariant,
but   multiplying $z$ an a power of $c$
will change the conjugacy class of ${\bf x}$ in general, and
will even change the connected component of ${\cal T}_G(c)$
containing the conjugacy class of ${\bf x}$, or equivalently will
change 
$CS_G({\bf x})$.

The case of commuting triples is worth stating separately.

\begin{proposition}
Suppose that $(x,y,z)$ is  a commuting triple.
Let $A(x,y,z)$ be a flat connection on a principal 
 $G$-bundle $\xi$ over $T^3$ with holonomy around the
three coordinate circles equal to the conjugacy class of
$(x,y,z)$. Then $\xi$ is a trivial bundle.
Let $\Theta$ be a connection on this bundle
isomorphic to the
 trivial connection. Then
$${\rm CS}_{\Theta}(A)\cong CS_G(x,y,z)\pmod{\bf Z}$$.
\end{proposition}

\begin{proof}
Since $c=1$, the $G$-bundle $\xi_c$ given in 
the statement of Corollary~\ref{CScorol} is trivial.
Furthermore, an enhancement of a $G$-bundle is no extra
information.
The result is now immediate from this corollary and
Lemma~\ref{integerdiff}.
\end{proof}

Let $c\in {\cal C}G$ be given and let $K=G/\langle c\rangle$.
We are now ready to define the Chern-Simons invariant of an
isomorphism class of a pair $(\Xi,A)$ consisting of an 
enhanced $K$-bundle over $(T^3,X)$ whose underlying
$K$-bundle  is isomorphic
to $\xi_c\times S^1$ and a flat connection.
The $G$-holonomy of such a pair is the conjugacy class
of a $c$-triple, and thus by Corollary~\ref{CScorol} there
is an enhanced isomorphism from $\Xi_c$ to $\Xi$. Let $\sigma$ be
the underlying
$K$-bundle isomorphism.
Then we define
$${\rm CS}(\Xi,A)=[{\rm CS}_\Theta(\sigma^*A)]\pmod{\bf Z},$$
where $\Theta$ is any flat connection on $\xi_c\times S^1$
pulled back from
a flat connection on $\xi_c$. Note that, if $\Theta$ is such a
connection, then its $G$-holonomy is the conjugacy class of
$(x,y,1)$, where $(x,y)$ is a $c$-pair, and, in particular, the
$G$-holonomy lies in the component $X_1$ of $c$-triples of order
$1$. If $c=1$, then we can take $\Theta$ to be the trivial
connection and are computing the usual Chern-Simons invariant. 

\subsection{Proof of Theorem~\ref{CS} and Theorem~\ref{clock} in
the case where $\langle C\rangle $ is cyclic}

By Lemma~\ref{CSconstant}, the value of
${\rm CS}(\Xi,A)$ only depends on the component $X$ of ${\cal
T}_G(c)$ containing the 
$G$-holonomy of $(\Xi , A)$.
Given this, 
the proof
 of Theorem~\ref{CS} for $c$-triples is immediate from
Corollary~\ref{CScorol} and Theorem~\ref{bijection}. The proof of
Theorem~\ref{clock} in the cyclic case  follows from the
dimension statement in Theorem~\ref{ctrip}, the computation of
the Chern-Simons invariant contained in  Theorem~\ref{bijection}
and Corollary~\ref{CScorol}, and the numerology of 
Theorem~\ref{clocked}, via the natural identification of ${\bf
Z}/2g{\bf Z}$ with $\displaystyle\frac{1}{2g}{\bf Z}/{\bf Z}$.

\section{The case when $\langle C\rangle$ is not cyclic}

In this section we assume that $G=Spin(4n)$ for some $n\ge 2$.
Then ${\cal C}G\cong {\bf Z}/2{\bf Z}\times {\bf Z}/2{\bf Z}$.
Choose an identification of a maximal torus $T$ for $G$ with the
quotient of ${\bf R}^{2n}$ with basis $\{e_i\}$ by the even
integral lattice so that the roots for $G$ are $\{\pm
e_i\pm e_j\}$. A set of simple roots for $G$ is then
$\Delta= \{a_1, \dots, a_{2n-2}, a_{-},a_{+}\}$, 
where $a_i=e_i-e_{i+1}$ for $1\le i\le 2n-2$, $a_{-}=e_{2n-1}-e_{2n}$,
and $a_{+}=e_{2n-1}+e_{2n}$.
Label the
non-trivial elements of ${\cal C}G$ as $c_0,c_1,c_2$ so
that $c_0$ is represented by $e_{2n}$, $c_1$ is represented by 
$(\sum_ie_i)/2$,  and $c_2$ is represented by $(e_1+\cdots
+e_{2n-1}-e_{2n})/2$.  In terms of the simple roots, $c_0
\equiv \frac12a_++\frac12a_- \bmod Q\spcheck$, $c_1 \equiv
\sum _{i=1}^{n-1}\frac12a_{2i-1} + \frac12a_+\bmod Q\spcheck$,
and $c_2 \equiv
\sum _{i=1}^{n-1}\frac12a_{2i-1} + \frac12a_-\bmod Q\spcheck$.
Thus
$\Delta(c_0) =
\{a_{+}, a_{-}\}$ and $L_{c_0}=L_+\times L_-$, where $L_{\pm}$
is  a group of type $A_1$ whose simple root is
$a_{\pm}$.  Let
$\Delta({\cal C})$ be the subset of $\Delta$ consisting of all
roots
$a\in\Delta$ such that $\varpi_a$ does not annihilate ${\cal C}G$.
Thus
$$\Delta({\cal C}) = \Delta(c_0) \cup \Delta(c_1)\cup
\Delta(c_2) =\{a_1, a_3, \dots, a_{2n-3}, a_-,a_+\}.$$ 
Then $L_{\cal C} = 
L_{\Delta({\cal C})} =\prod_{i=1}^{n-1}L_i\times L_+\times L_-$
is  a product of
$n+1$ groups of type $A_1$, where $a_{2i-1}$ is the simple root
of $L_i$. We also let
$S_{\cal C}= S_{\Delta({\cal C})}$.

We have the following analogue of Corollary~\ref{LisLc}, whose
proof is left to the reader:

\begin{lemma}\label{L=LC} Let  $I\subseteq \Delta$ and let
$L=L_I$. Then $L_{\cal C}\subseteq L$ if and only if ${\cal C}
\subset L$, and $L=L_{\cal C}$ if and only if 
\begin{enumerate}
\item ${\cal C} \subset L$;
\item  $L$ is a product of simple factors $L_i\cong SU(n_i)$ for 
some $n_i$;
\item The projection of ${\cal C}$ to $L_i$ generates 
${\cal C}L_i$.
\end{enumerate}
\end{lemma}

With our choice of $T$ and $\Delta$, the action of $w_{c_0}$ on
$\frak t$ is given by:
$$w_{c_0}(t_1, t_2, \dots, t_{2n-1}, t_{2n}) = (-t_1, t_2, \dots,
t_{2n-1}, -t_{2n}),$$
so that $\frak t^{w_{c_0}} = \{(t_1, t_2, \dots, t_{2n-1},
t_{2n})\in \frak t: t_1 = t_{2n} = 0\}$. Of course, $\frak
t^{w_{c_0}}$ is conjugate to 
$$\frak t_{c_0} ={\rm Ker}(a_+) \cap
{\rm Ker}(a_-) = \{(t_1, t_2, \dots, t_{2n-1},
t_{2n})\in \frak t: t_{2n-1} = t_{2n} = 0\}.$$
Note that $c_0 \in S^{w_{c_0}}$ but that $c_1, c_2\notin
S^{w_{c_0}}$. 

We shall consider 
triples $(x,y,z)$ in $G$ with 
$[x,y]=c_0$, $[x,z]=c_1$, and $[y,z]=c_2$, and let $C$ be the
corresponding antisymmetric matrix. 

\subsection{Rank zero $C$-triples}

\begin{lemma}\label{Crankzero}
There is a rank zero $C$-triple in $G$ if and only if 
$G=Spin(8)$. For $G=Spin(8)$ there are exactly two conjugacy
classes of rank zero
$C$-triples. If $(x,y,z)$ is a rank zero
$C$-triple in
$Spin(8)$, then
$(x,y,z^{-1})$ is a rank zero $C$-triple not conjugate to
$(x,y,z)$.  
\end{lemma}

\begin{proof}
Suppose that $(x,y,z)$ is a rank zero $C$-triple in $G$.
After conjugation, we may
assume that
$(x,y)$ is a
$c_0$-pair in normal form. Thus $S^{w_{c_0}}$ is a maximal
torus of $Z^0(x,y)$.  The element
$z$ commutes with $x$ and
$y$ up to an element of the center.   Since 
the fixed subgroup of conjugation by $z$ on $Z(x,y)$ has rank
zero, it follows  from \cite{deS} II \S 2 that 
$Z^0(x,y)$ is a torus and hence $Z^0(x,y)=S^{w_{c_0}}$. 
 By inspection,
all of the integers $g_{\ov a}$, for the action of $w_{c_0}$ on
the coroot diagram of $G$  are $2$. Hence by
Proposition~\ref{equalg}, $Z(x,y) = T^{w_{c_0}}$. The triple
$(x,y,z^2)$ is a $c_0$-triple in $G$. By Proposition~\ref{atmost},
$S^{w_{c_0}}$ is a maximal torus of 
$Z(x,y,z^2)$.  This implies
that $z^2$ acts trivially on
$S^{w_{c_0}}$.
Since $z$ acts on this torus without fixed points, it follows
that $z$ acts by $-1$ on ${\frak t}^{w_{c_0}}$, and hence
$zsz^{-1} = s^{-1}$ for all $s\in S^{w_{c_0}}$. 

By Proposition~\ref{torus}, the torus $S^{w_{c_0}}$ is conjugate
to $S_{c_0}$. It will be convenient to conjugate $x,y,z$ so that
$Z^0(x,y) = S_{c_0}$. Of course, in this case $zsz^{-1} = s^{-1}$
for all $s\in S_{c_0}$. Write
$x=s\cdot x_0$ and $y=s'\cdot y_0$ for elements
$s,s'\in S_{c_0}$ and $x_0,y_0\in L_{c_0}$.  
Since
$c_1sx_0={}^z(sx_0)=(s^{-1}){}^zx_0$, it follows  that
$s^2=c_1({}^zx_0x_0^{-1})=c_1\cdot u$ for some $u\in S_{c_0}\cap
L_{c_0}$. This implies that $s^2$ is the exponential of an
element of the form
$(\frac{1}{2},\frac{1}{2},\ldots,\frac{1}{2},0,0)$ 
mod ${\bf Z}^n$.
The same computation  for $y$   shows that $(s')^2$ is also the
exponential of such an element. Thus, both $s$ and $s'$ are the
images under the exponential mapping of elements of the form
$(\pm\frac{1}{4},\pm\frac{1}{4},\ldots,\pm\frac{1}{4},0,0)$ 
mod ${\bf Z}^n$.
It  follows easily that, for $4n>8$,  there is a root of
$Spin(4n)$ annihilating $L_{c_0}$ and annihilating $s$ and $s'$.
This root 
then annnihilates both $x$
and $y$, contradicting the fact that $Z^0(x,y)$ is a torus.
Hence $G=Spin(8)$.

Now let us describe all conjugacy classes of rank zero
$C$-triples  in $Spin(8)$. Let $(x,y,z)$ be a rank zero $C$-triple
in
$Spin(8)$. As above,  we arrange that $(x,y)=(s_1x_0,s_2y_0)$
with $s_i\in S_{c_0}$ and
$(x_0,y_0)\in L_{c_0}$. Then
$Z^0(x,y)=S_{c_0}$.  According to Corollary~\ref{maxtorus},  the
conjugacy class of
$(x,y)$ depends only on the the pair $(\ov s_1,\ov s_2)\in \ov
S_{c_0}=S_{c_0}/S_{c_0}\cap
L_{c_0}$. It is easy to see that $\ov S_{c_0}$ is the quotient
of
${\bf R}^2\times\{0\}\times\{0\}$ by ${\bf Z}^2$. The image in
the Weyl group of $S_{c_0}$ of elements in $N_G(S_{c_0})$ which
act trivially on $L_{c_0}$ is 
${\bf Z}/2{\bf Z}\times {\bf Z}/2{\bf Z}$, where one of the
factors acts by
$-1$ on $\frak t_{c_0}$ and the other acts by switching the
coordinates.  The argument above shows that
$\ov s_i$ is the  exponential of an element of the form
$(\pm\frac{1}{4},\pm\frac{1}{4},0,0)$ mod ${\bf Z}^2$. Let 
$\tilde s_i$ be a lift of $\ov s_i$ to
$\frak t_{c_0}$. 

The element $z$ normalizes $S_{c_0}$ and acts by
$-1$ on $\frak t_{c_0}$.Thus $z$  sends $S_{c_0}\cdot x$ 
to $S_{c_0}\cdot c_1x$.  Since $S_{c_0}\cdot x$ contains a
regular element of $T$, the element $z$ normalizes $T$ and hence
normalizes the unique maximal torus of $L_{c_0}$ containing
$x_0$.    Since $x_0$ is the  
exponential  of
$\frac12(a_++a_-)= (0,0,\frac{1}{2},0)$, it follows that
$zx_0z^{-1}$ is the  exponential of one of
$\pm\frac12a_++\pm\frac12a_-$ and hence of one of
$(0,0,\pm\frac{1}{2},0)$ or $(0,0,0,\pm\frac{1}{2})$. The fact
that
$zxz^{-1}=xc_1$ implies that $zx_0z^{-1}$ is  the
exponential   of $(0,0,0,\pm\frac{1}{2})$. By conjugation
by an element commuting with $x$ we can arrange that $zx_0z^{-1}$
is the   exponential  of
$(0,0,0,-\frac{1}{2})$. A computation shows that $\tilde s_1$ must
be congruent to $\pm (\frac{1}{4},\frac{1}{4},0,0)$ mod ${\bf
Z}^2$.  After conjugating by an element of the normalizer of
$S_{c_0}$, we may assume that $\tilde s_1$ is congruent to
$(\frac{1}{4},\frac{1}{4},0,0)$ mod ${\bf Z}^2$. Since there is
no root of
$Spin(8)$ annihilating both
$s_1$ and
$s_2$, this means that  $\tilde s_2$
is either $(\frac{1}{4},-\frac{1}{4},0,0)$ or
$(-\frac{1}{4},\frac{1}{4},0,0)$ mod ${\bf Z}^2$. After
conjugation by an element of the normalizer of $S_{c_0}$ which 
interchanges the two coordinates, and hence fixes $\tilde s_1$,
we may  assume that
$\tilde s_2$ is congruent to
$(\frac{1}{4},-\frac{1}{4},0,0)$  mod ${\bf Z}^2$.
This proves that, given a rank zero $C$-triple $(x,y,z)$ in
$Spin(8)$, the $c_0$-pair $(x,y)$ is determined up to conjugation
in $Spin(8)$.
The representatives that we have chosen are:
$x={\rm exp}(\frac{1}{4},\frac{1}{4},\frac{1}{2},0)$ and 
$y={\rm exp}(\frac{1}{4},-\frac{1}{4},0,0)\cdot w$
where $w$ is an element of $L_{c_0}$ normalizing
the unique maximal torus of  $L_{c_0}$ containing $x_0={\rm
  exp}(0,0,\frac{1}{2},0)$, and the image of $w$ in the Weyl group
of $L_{c_0}$ is the product of the non-trivial  elements  in each
factor of $L_{c_0}$. 

Next we show that there is $z\in G$ such that $(x,y,z)$ is a
rank zero $C$-triple. Write
$x_0 = x_0^+x_0^-$ as a product, where $x_0^\pm\in L_\pm$.
Similarly write
$y_0=y_0^+y_0^-$. Let $z=\epsilon\cdot x^+_0y^+_0$ where
$\epsilon\in  N_{Spin(8)}(S _{c_0})$ commutes with $L_{c_0}$ and
represents the Weyl element which is multiplication by $-1$ on
${\frak t}_{c_0}$. 
Direct computation shows that $(x,y,z)$ is a rank zero
$C$-triple.

Suppose that $(x,y,z)$ is a rank zero $C$-triple in $Spin(8)$. Any
other rank zero
$C$-triple is conjugate to $(x,y,z')$, where $z'=zg$ for some
$g\in Z(x,y)$. Since $\pi_0(Z(x,y))$ is cyclic of order two and
since $i_z =-{\rm Id}$ on $Z^0(x,y)$,  there are
exactly two conjugacy classes of rank zero $C$-triples in
$Spin(8)$, one of which is represented by $(x,y,z)$ and the other
by $(x,y,zt)$, where $t\in Z(x,y)$ but $t\notin Z^0(x,y)$. 
Direct computation shows that
$z$ is of order $4$ and that
$z^2$ is in the non-trivial component of $Z^0(x,y)$. Thus,
$(x,y,z^{\pm 1})$ represent the   two conjugacy classes of
rank zero $C$-triples. 
\end{proof}

\subsection{Action of the center and of the outer automorphism
group}

The following lemma is easy and its proof is left to the reader.

\begin{lemma} An automorphism of $Spin(8)$ which acts trivially on
${\cal
  C}Spin(8)$ is   an inner automorphism and hence acts
  trivially on the space of of conjugacy classes of rank zero
  $C$-triples in $Spin(8)$.
\end{lemma}

Let us consider the action of $({\cal C}Spin(8))^3$ on the space
of rank zero
$C$-triples.

\begin{lemma}\label{noname}
The action of ${\cal C}Spin(8)$ on the moduli space of conjugacy
classes of rank zero $C$-triples in $Spin(8)$ defined by 
$\mu\cdot (x,y,z)=(x,y,\mu\cdot z)$ is transitive.
The isotropy subgroup of this action is $\langle c_0\rangle$.
Likewise, the action  defined by 
$\mu\cdot (x,y,z)=(\mu\cdot x,y,  z)$ is transitive,
and its isotropy is
$\langle c_2\rangle$, and the action defined by 
$\mu\cdot (x,y,z)=(x,\mu\cdot y, z)$ is
transitive, with isotropy $\langle c_1\rangle$.
\end{lemma}
\begin{proof}
This is immediate from the explicit description given above.
\end{proof}

\subsection{The general case}

Following the general pattern, we begin by determining the
possibilities for the  maximal torus of
$Z(x,y,z)$.
Recall from Theorem~\ref{main} that, for every ${\cal C}G$-triple
$(x,y,z)= {\bf x}$, there is a subset 
$I({\bf x})\subset \Delta$ such that $S_I$ is conjugate in $G$ to
a maximal torus of $Z(x,y,z)$.

\begin{lemma}\label{firstlemma}
Let $I'=\Delta({\cal C})\cup \{a_{2n-2}\}$.
If $(x,y,z)={\bf x}$ is a
${\cal C}G$ triple, then either  $I({\bf x})=\Delta({\cal C})$ or
$I({\bf x})=I'$.
\end{lemma}
\begin{proof}
Let $L = L_{I({\bf x})}$. Then by Lemma~\ref{L=LC}, since ${\cal
C}G\subset L$, $\Delta({\cal C}) \subseteq I({\bf x})$. We may
write $L=L_0\times \prod_{i=1}^rL_i$, where each $L_i$ is of type
$A_n$ and $L_0$ is either trivial or   of type $D_{2k}$ for some
$k\geq 2$. For $0\leq i\leq r$, let
$\pi_i$ denote projection to the $i^{\rm th}$ factor. The
projection of
${\bf x}$ to
$L_i$ is a rank zero $\pi_i(C)$-triple. Thus, if $i\geq 1$, then
$L_i$ is of type $A_1$. If
$L_0$ is not trivial, then 
$\pi_0({\cal C})$ is the full center of $L_0$. By
Lemma~\ref{Crankzero},
$k=2$. Since $\Delta({\cal C}) =  \{a_1, a_3, \dots, a_{2n-1},
a_-,a_+\}$, it is clear that the only possibilities for $I({\bf
x})$ which satisfy the above conditions are either $I({\bf x})=
I_{\cal C}$ or $I({\bf x})=I'$.
\end{proof}

The tori $S_{\cal C}$ and $S_{I'}$ corresponding to the sets
$\Delta({\cal C})$ and $I'$ have the following Lie algebras:
\begin{eqnarray*}
\frak t_{\cal C} &=& \bigcap _{i=1}^{n-1}{\rm Ker}(a_{2i-1}) \cap
{\rm Ker}(a_+)\cap {\rm Ker}(a_-) =\{(t_1, t_1, t_3, t_3, \dots,
t_{2n-3}, t_{2n-3}, 0,0)\};\\
\frak t_{I'} &=& \bigcap _{i=1}^{n-1}{\rm Ker}(a_{2i-1}) \cap
{\rm Ker}(a_+)\cap {\rm Ker}(a_-)\cap {\rm Ker}(a_{2n-2})
\\
&=& \{(t_1, t_1, t_3, t_3, \dots, t_{2n-5}, t_{2n-5},0,0, 0,0)\}.
\end{eqnarray*}

In the quotient diagram $\widetilde D\spcheck(G)/{\cal C}G$, one
of the quotient coroot integers $g_{\ov a}$   is $2$ and all the
others are $4$.
Thus, in terms of the quotient diagram and the quotient coroot
integers, we have:
\begin{eqnarray}\label{firsttorus}
\frak t^{w_C} &=& \{(0, s_2, \dots, s_{n-1}, s_n, -s_n, -s_{n-1},
\dots, -s_2, 0)\};\\ \label{secondtorus}
\frak t^{w_C}(\ov {\bf g},4) &=& \{(0, s_2, \dots, s_{n-1}, 0, 0,
-s_{n-1},
\dots, -s_2, 0)\}.
\end{eqnarray}

    From this, it is easy to establish the following:

\begin{lemma}\label{toriOK}
The torus $S_{\cal C}$ is conjugate to $S^{w_C} = S^{w_C}(\ov
{\bf g},1) = S^{w_C}(\ov {\bf g}, 2)$. The torus $S_{I'}$ is
conjugate to 
$S^{w_C}(\ov {\bf g}, 4)$.
\end{lemma}

Having described the possible maximal tori, we need to determine  
how many components of the moduli space correspond to a given
maximal torus. We begin with
the following elementary computation.

\begin{lemma}\label{A_1}
Let $\gamma$ be the non-trivial central element in $SU(2)$.
Consider triples $(p,q,r)$ in $SU(2)$ with $[p,q]=1$ and
$[p,r]=[q,r]=\gamma$. 
Then, up to conjugation, there are exactly two such triples.
We have $p^2=q^2=r^2=\gamma$.
The element $pq^{-1}$ lies in ${\cal C}SU(2)$ and
is a complete invariant of the conjugacy class of $(p,q,r)$.
 The
action of $({\cal C}SU(2))^3$ on the conjugacy classes of such
triples given by $(a,b,c)\cdot (p,q,r)=(ap,bq,cr)$ is transitive
and the stabilizer of any conjugacy class is $\{(a,b,c): a=b\}$.
\end{lemma}

\begin{proof}
The pair $(p,r)$ is a $\gamma$-pair in $SU(2)$. Also, $pq^{-1}\in
Z(p,r)$ and hence lies in ${\cal C}SU(2)$. The result follows
from Proposition~\ref{rank0pair} and Corollary~\ref{orders}.
\end{proof}

Suppose that ${\bf x}$ is a $C$-triple in $L_{\cal C}=L_{\cal
C}=\prod_{i=1}^{n-1}L_i\times L_+\times L_-$.
The triple ${\bf x}$
is  a product of triples $\prod _i{\bf x}_i\cdot {\bf x}_+\cdot
{\bf x}_-$ in the factors.
We write ${\bf x}_i=(x_i,y_i,z_i)$ and ${\bf x}_\pm=(x_\pm,
y_\pm,z_\pm)$.  Let $\epsilon_i\colon {\cal C}L_i\to \{\pm 1\}$,
$\epsilon_\pm\colon {\cal C}L_\pm\to \{\pm 1\}$ be the unique
isomorphisms. Then define $\epsilon({\bf x})=\prod_{i=1}^{n-1}
\epsilon_i(x_iy_i^{-1})
\epsilon_+(x_+z_+^{-1})\epsilon_-(y_-z_-^{-1})$.

\begin{lemma}\label{Ccompbig}
There are exactly two components of ${\cal T}_G(C)$ with maximal
torus conjugate to $S_{\cal C}$. Two $C$-triples ${\bf x}$ and
${\bf x}'$ in $L_{\cal C}$ lie in the same component of ${\cal
T}_G(C)$ if and only if $\epsilon({\bf x})=\epsilon({\bf x}')$.
\end{lemma}

\begin{proof}
By Theorem~\ref{main}, the number of
components of ${\cal T}_G(C)$ with  maximal torus conjugate to 
$S_{\cal C}$ is given by the number of conjugacy classes of rank
zero 
$C$-triples  in $L_{\cal C}$ modulo the action of $F_{\cal C}^3$
and of
$W(S_{\cal C}, G)$, where 
$F_{\cal C}=S_{\cal C}\cap L_{\cal C}I\subseteq {\cal C}L_{\cal
C}\cong ({\bf Z}/2{\bf Z})^{n+1}$. Let us first consider the
action of $F_{\cal C}^3$. It is easy to see that 
\begin{equation}\label{firstcrazy}
F_{\cal C}=\{(\mu_1,\ldots,\mu_{n-1},\mu_+,\mu_-):
\mu_1\cdots \mu_{n-1}=\mu_+=\mu_-\}.
\end{equation}
Straightforward computation shows that the action of $F_{\cal
C}^3$ on the space of conjugacy classes of rank zero $C$-triples
in $L_{\cal C}$ preserves the invariant
$\epsilon$ and acts transitively on the set of conjugacy classes
with a given $\epsilon$. 

The image of $W(S_{\cal C}, G)$ in the outer automorphism group
of $L_{\cal C} = \prod_iL_i\times L_+\times L_-$ is easily
checked to be the subgroup of all permutations of the $L_i$
factors. Thus the action of $W(S_{\cal C}, G)$ fixes $\epsilon$. 
This completes the proof of the lemma.
\end{proof}

\begin{lemma}
Under the action of $F_{\cal C}^3$ on the space of conjugacy
classes of rank zero $C$-triples in $L_{\cal C}$ the stabilizer
of any class is the set of $({\mbox{\boldmath${\mu}$}} =
(\mbox{\boldmath${\mu}$}_1,\ldots,
\mbox{\boldmath${\mu}$}_{n-1},\mbox{\boldmath${\mu}$}_+,
\mbox{\boldmath${\mu}$}_-)\in F_{\cal C}^3$ such that
\begin{equation}\label{crazy}
\mu_j^{(1)}=\mu_j^{(2)}, 1\le  j\le n-1,
\mu_-^{(2)}=\mu_-^{(3)},  \mu_+^{(1)}=\mu_+^{(3)},
\end{equation}
where
$\mbox{\boldmath${\mu}$}_j=(\mu_j^{(1)},\mu_j^{(2)},
\mu_j^{(3)})$ and similarly for
$\mbox{\boldmath${\mu}$}_\pm$.
\end{lemma}

\begin{proof}
This is immediate from  Lemma~\ref{A_1}.
\end{proof}

\begin{corollary}\label{order1,2}
There are two components of ${\cal T}_G(C)$
associated  with the torus $S_{\cal C}$, and each is homeomorphic
to 
$$\left((S_{\cal C}\times S_{\cal C}\times
S_{\cal C})/F\right)/W(S_{\cal C},G),$$
where $F\subset
 ({\cal C}L_{\cal C})^3$ is the set of {\boldmath${\mu}$}
satisfying  Equations~\ref{firstcrazy} and~\ref{crazy}. 
\end{corollary}

Now let us consider the case when $I=I'$.
Recall that
$L_{I'}=\prod_{i=0}^{n-2}L_i$, where $L_i$ is
of type $A_1$ for $i> 0$ and $L_0\cong Spin (8)$.  If ${\bf x}$ is
a rank zero $C$-triple in $L_{I'}$, then we write ${\bf x}
=\prod_{i=0}^{n-2}{\bf x}_i $. The group
$L_i$ has two rank zero
$\pi_i(C)$-triples, up to conjugation, where $\pi_i$ is the
projection from $L$ to $L_i$. It follows that $L_{I'}$ has
$2^{n-1}$ rank zero
$C$-triples up to conjugation. If ${\bf x}_i$ and  ${\bf
x}_i'$ are two rank zero
$\pi_i(C)$-triples in $L_i$, define $\delta_i({\bf x}_i, {\bf
x}_i') =1$ if ${\bf x}_i$ is conjugate in $L_i$ to ${\bf x}_i'$
and $-1$ otherwise. If $i>0$, then it is easy to see
that $\delta_i({\bf x}_i, {\bf x}_i') =
\epsilon_i(x_iy_i^{-1})\epsilon_i(x_i'(y_i')^{-1})$, in the
notation introduced prior to Lemma~\ref{Ccompbig}. Define
$$\delta({\bf x}, {\bf x}') =\prod_{i=0}^{n-2}\delta_i({\bf x}_i,
{\bf x}_i').$$

\begin{lemma}\label{Ccompsmall} Two rank zero $C$-triples ${\bf
x}$ and
${\bf x}'$ in $L_{I'}$ lie in the same component of ${\cal
T}_G(C)$ if and only if $\delta({\bf x}, {\bf x}') =1$. Hence
there are exactly two components of ${\cal T}_G(C)$ with  maximal
torus conjugate to  $S_{I'}$.
\end{lemma}
\begin{proof} By Theorem~\ref{main}, the number of
components of ${\cal T}_G(C)$  with  maximal torus conjugate to 
$S_{I'}$ is given by the number of conjugacy classes of rank
zero 
$C$-triples  in $L_{I'}$ modulo the action of $F_{I'}^3$
and of
$W(S_{I'}, G)$, where 
$F_{I'}=S_{I'}\cap L_{I'}\subseteq {\cal C}L_{I'}\cong ({\bf
Z}/2{\bf Z})^{n-2}\times {\cal C}L_0$. Let us first consider the
action of
$F_{I'}^3$. Let $\nu\colon\{\pm 1\} \to {\cal C}L_0$ be the
emdedding which sends $-1$ to $\pi_0(c_0)$. Then it is easy to
check that
\begin{equation}\label{secondcrazy}
F_{I'} = \left\{(\mu_0, \dots,
\mu_{n-2}):
\mu_0 =\nu \left(\prod_{i=1}^{n-2}\epsilon
_i(\mu_i)\right)\right\}\subseteq
\prod_{i=0}^{n-2}{\cal C}L_i.
\end{equation} 
Note that, by Lemma~\ref{noname}, if {\boldmath${\mu}$} is an
element of
$F_{I'}^3$ and
${\bf x}$ is a rank zero $C$-triple in $L_{I'}$, then 
$\delta ({\bf x}, \mbox{\boldmath${\mu}$}\cdot {\bf x}) = 1$.
Conversely, it is clear that, if $\delta({\bf x}, {\bf x}') =1$,
then ${\bf x}$ and ${\bf x}'$ are in the same  $F_{I'}^3$-orbit.
Finally, the image of $W(S_{I'}, G)$ in the outer automorphism
group of $L_{I'}$ is the permutation group of the factors
$L_i$ for $i\geq 1$. Hence $\delta ({\bf x}, w\cdot {\bf x}) =1$
for all $w\in W(S_{I'}, G)$. The lemma follows.
\end{proof}

\begin{lemma}
Under the action of $F_{I'}^3$ on the space of conjugacy
classes of rank zero $C$-triples in $L_{I'}$ the stabilizer
of any class is the set of $\mbox{\boldmath${\mu}$}=
(\mbox{\boldmath${\mu}$}_0,\ldots,
\mbox{\boldmath${\mu}$}_{n-0})\in F_{I'}^3$ such that
\begin{equation}\label{crazier}
\mu_j^{(1)}=\mu_j^{(2)}, 1\le  j\le n-2 ,
\end{equation}
where
$\mbox{\boldmath${\mu}$}_j=(\mu_j^{(1)},\mu_j^{(2)},
\mu_j^{(3)})$.
\end{lemma}
\begin{proof}
This is immediate from  Lemma~\ref{A_1} and Lemma~\ref{noname}.
\end{proof}

\begin{corollary}\label{order4}
There are two components of ${\cal T}_G(C)$
associated  with the torus $S_{I'}$, and each is homeomorphic
to 
$$\left((S_{I'}\times S_{I'}\times
S_{I'})/F'\right)/W(S_{I'},G),$$
where $F'\subset
({\cal C}L_{I'})^3$ is the subgroup of elements satisfying 
Equations~\ref{secondcrazy} and~\ref{crazier}. 
\end{corollary}

To each component $X$ of ${\cal T}_G(C)$, we assign a positive
integer which we shall call the {\sl order\/} of $X$, in the
following way. Suppose that
$X$ is a component containing the conjugacy class of ${\bf x}$,
where ${\bf x}$ is a rank zero $C$-triple in $L_{\cal C}$. Then
the order of $X$ is $1$ if $\epsilon ({\bf x}) =1$ and the order
of $X$ is $2$ otherwise.   Each of the remaining components is
associated to
$S_{I'}$, and we define each of them to have order $4$.
Thus, for every positive integer $k$ dividing $4$, the number of
components of ${\cal T}_G(C)$ is $\varphi(k)$.

Parts 1,2,3,4 of Theorem~\ref{ctrip}, in the case where $\langle
C\rangle$ is not cyclic, are now contained in the statements of
Lemma~\ref{firstlemma}, Corollaries~\ref{order1,2}
and~\ref{order4}, and the definition of the order.

Lastly we prove Part 5 of
Theorem~\ref{ctrip} in the non-cyclic case. The explicit
coordinates on the Lie algebras
$\frak t^{w_C}$ and $\frak t^{w_C}(\ov {\bf g}, 4)$ given in
Equations~\ref{firsttorus} and~\ref{secondtorus} identify these
vector spaces with ${\bf R}^{n-1}$ and ${\bf R}^{n-2}$,
respectively. Via these identifications, the projections of
$Q\spcheck$ are the full integral lattices.   In each case, the
Weyl group is the full isometry group of the lattice. The images
under projection of the extended simple coroots of $G$ are a set
of extended simple coroots for a root system of type
$BC_{n-1}$, resp.\ $BC_{n-2}$. Part 5 of
Theorem~\ref{ctrip} is then clear.

\subsection{Chern-Simons invariants}

Let $G=Spin(4n)$ and let $\ov G=G/{\cal C}G$.
Fix a $\ov G$-bundle $\xi$ over the three-torus whose  second
Stiefel-Whitney class $w(\xi)$ evaluates over the three
coordinate two-tori $T_{12}, T_{13}, T_{23}$ to give
$c_0,c_1,c_2\in
\pi_1(\ov G)={\cal C}G$ and let $\Xi$ be an enhanced $\ov
G$-bundle whose underlying bundle is $\xi$.  If
${\bf x} =(x,y,z)$ is a
$C$-triple, then, by Corollary~\ref{Xibiject}, there is a  flat
 connection on
$\Xi$ whose $G$-holonomy around the coordinate circles is the
conjugacy class of $(x, y,z)$.  
Let $\Gamma_1$  be a flat connection on $\Xi$ 
whose $G$-holonomy is the
conjugacy class of a $C$-triple ${\bf x}=(x,y,z)$ in 
$L_{\cal C}$ with $\epsilon({\bf x}) =1$, i.e. a $C$-triple whose
conjugacy class lies in the component $X_1$ of ${\cal T}_G(C)$  of
order one. Write ${\bf x} =\prod_i{\bf
x}_i\cdot {\bf x}_+\cdot {\bf x}_-$, where ${\bf x}_i = (x_i,
y_i, z_i)$ lies in
$L_i$, and similarly for ${\bf x}_\pm$. Given a flat connection
$A$ on
$\Xi$,
${\rm CS}_{\Gamma_1}(A)$ is independent mod ${\bf Z}$ of the
choice of
$\Gamma_1$, by Lemma~\ref{integerdiff} and
Corollary~\ref{CSconstant}.  We will denote its class in
${\bf R}/{\bf Z}$ by 
${\rm CS}(A)$. 

\begin{lemma}\label{oppositebig}
Let $\Gamma$ be a flat connection on $\Xi$. 
If the $G$-holonomy of $\Gamma$  is contained in $X_1$, then  
${\rm CS}(\Gamma)=0
\bmod {\bf Z}$. If the $G$-holonomy of $\Gamma$  is contained in
the component of order $2$, then ${\rm CS}(\Gamma)=1/2\bmod {\bf
Z}$. 
\end{lemma}

\begin{proof}
Since ${\rm CS}$ is constant on components of ${\cal T}_G(C)$, 
the  first statement is clear. To prove the second it suffices to
compute
${\rm CS}( \Gamma)$ for one flat connection $\Gamma$
whose $G$-holonomy is of order $2$.  By Lemma~\ref{Ccompbig},   a
$C$-triple ${\bf x}'$ of order $2$ is given by replacing 
$x_1$ in ${\bf x}$ by $\gamma_1x_1$ where $\gamma_1$ is  
the non-trivial   element in ${\cal C}L_1$. Since the
inclusion of each of the $A_1$-factors of $L_{\cal C}$ into $G$
induces an isomorphism on $\pi_3$, to compute the
Chern-Simons invariant ${\rm CS}(\Gamma)$, it
suffices to compute the 
relative Chern-Simons of the $A_1$-connections obtained taking the
images of $\Gamma_1$ and $\Gamma$ under projection the first
$A_1$-factor of $L_{\cal C}$. Lemma~\ref{SUn} and
Corollary~\ref{CScorol} applied to $SU(2)$ show that this
relative Chern-Simons invariant is $1/2\bmod {\bf Z}$. This shows
that
${\rm CS}( \Gamma)=1/2\bmod {\bf Z}$.
\end{proof}

\begin{lemma}
Let $\Gamma$ be a flat connection on $\Xi$ whose $G$-holonomy has
order $4$. Then ${\rm CS}(\Gamma)=\pm 1/4$ modulo ${\bf Z}$. 
If $\Gamma'$ represents a point in the other component of ${\cal
T}_G(C)$ of order $4$, then ${\rm CS}(\Gamma') = -{\rm
CS}(\Gamma)$.
\end{lemma}

\begin{proof}
Let ${\bf u} =(u,v,w)\in L_{I'}$ be a
rank zero $C$-triple. Write ${\bf u}=\prod_{i=0}^{n-2}{\bf u}_i$,
where ${\bf u}_i=(u_i, v_i, w_i)$ is a rank zero
$\pi_i(C)$-triple in $L_i$.  Let $\Gamma$ be a
flat connection on $\Xi$ whose $G$-holonomy is the conjugacy
class of ${\bf u}$. The
$C$-triple
${\bf u}'$ obtained by replacing ${\bf u}_i$ by ${\bf u}_i' =
(u_i, v_i, w_i^{-1})$ is also of
  rank zero. Clearly, for $i>0$, $\delta_i({\bf u}_i, {\bf
u}_i') =1$. By Lemma~\ref{Crankzero}, $\delta_0({\bf u}_0, {\bf
u}_0') =-1$, and hence, by  Lemma~\ref{Ccompsmall},
the conjugacy classes of ${\bf u}'$
and
  ${\bf u}$ 
lie in different components of order $4$
of ${\cal T}_G(C)$.
On the other hand, given the $C$-triple ${\bf x} $  defined
above,   let
${\bf x}'=\prod_i{\bf
x}_i'\cdot {\bf x}_+'\cdot {\bf x}_-'$, where ${\bf x}_i' = (x_i,
y_i, z_i^{-1})$, and similarly for ${\bf x}_\pm'$. Clearly
$$\epsilon ({\bf x}) =
\epsilon_+(\gamma_+)\epsilon_-(\gamma_-)\epsilon ({\bf
x}')=\epsilon ({\bf x}').$$ (Here $\gamma_\pm$ are the nontrivial
elements of ${\cal C}L_\pm$.) Thus by Lemma~\ref{Ccompbig}, ${\bf
x}$ and ${\bf x}'$ represent points  in the component $X_1$.

Let $r$ be the diffeomorphism of the three-torus which is
the identity on the first two factors and is inversion on
the third. Note that $r(X) = X$ so that $r^*\Xi$ is an enhanced
bundle over $(T^3, X)$. Then
${\bf x}'$ is the
$G$-holonomy of
$r^*\Gamma_1$ on $r^*\Xi$, and hence by the above computations 
$r^*\Gamma_1$ represents a point in $X_1$. Likewise
$r^*(\Gamma)$ represents a point in a different component from
$\Gamma$. 
Since $r$  is an
  orientation-reversing diffeomorphism, we have
$${\rm CS}(\Gamma) \equiv {\rm CS}_{\Gamma_1}(\Gamma ) =
-{\rm CS}_{r^*\Gamma_1}(r^*(\Gamma))
\equiv -{\rm CS}(r^*(\Gamma)) \bmod {\bf Z}.$$  Thus, 
${\rm CS}$ takes opposite values modulo ${\bf Z}$ on the two
components of order $4$ of ${\cal T}_G(C)$. 

With ${\bf u}$ as above, let $\hat{\bf u}$ be the triple  $(u,
v, w^2)$. Then $\hat{\bf u}$ is a $c_0$-triple.
Likewise, let
$\hat{\bf x}$ be the $c_0$-triple  $(x, y, z^2)$.

\begin{claim}
The $c_0$-triple $\hat{\bf u}$ is in the
non-trivial component of the modulo space of $c_0$-triples in $G$.
The $c_0$-triple $\hat {\bf x}$ is in the
trivial component of the modulo space of $c_0$-triples in $G$.
\end{claim}

\begin{proof}
It follows from the proof of Lemma~\ref{Crankzero} that the square 
of $w_0$ in $L_0$ is $\gamma_+$ which is in the
  non-trivial component of $T_{L_0}^{w_{c_0}}$.
Direct computation shows that   $S^{w_{c_0}}\cap L_0$
 is connected, and hence is equal to $(T_{L_0}^{w_{c_0}})^0$. 
Thus  $w^2$ lies in the non-trivial
  component of $T^{w_{c_0}}$. Since $\pi_0(T^{w_{c_0}}) \cong
\pi_0(Z(u,v))$, by Proposition~\ref{equalg},  the
$c_0$-triple
$\hat{\bf u}$ represents a point in the non-trivial component of
  the moduli space of $c_0$-triples in $G$.

By Lemma~\ref{A_1},
$z^2=\gamma_1\cdots \gamma_{n-1}\cdot \gamma_+\cdot \gamma_-$.
Since the $\gamma_i, 1\le i\le n-1$ lie in $S^{w_{c_0}}$, it
follows that $z^2$ is congruent modulo $S^{w_{c_0}}$ to
$\gamma_+\cdot \gamma_- =c_0$ which lies in $S^{w_{c_0}}$. Hence 
$z^2\in S^{w_{c_0}}$, and so the $c_0$-triple
$\hat{\bf x}$ lies in the trivial component of
the moduli space of
$c_0$-triples in $G$.
\end{proof}

Let $t$ be the map on the three-torus which double covers the last
coordinate. Note that $t$ induces a map from the pair $(T^3,
X)$ to itself. Thus $t^*\Xi$ is an enhanced $\ov G$-bundle.  The
$c_0$-triple
$\hat{\bf x}$ is the $G$-holonomy of 
$t^*\Gamma_1$ while
$\hat{\bf u}$ is the $G$-holonomy of $t^*\Gamma$. Since the
Chern-Simons invariant is given by integrating a local expression
involving connections and  curvature,  
$${\rm
CS}_{t^*(\Gamma_1)}(t^*\Gamma)=2{\rm CS}_{\Gamma_1}(\Gamma) \equiv
2{\rm CS}(\Gamma)\bmod {\bf Z}.$$  Since the conjugacy classes of
$\hat{\bf u}$ and
$\hat{\bf x}$  lie in opposite components of the moduli space of
$c_0$-triples, by
Theorem~\ref{CS} applied to the case of $c_0$-triples, 
$${\rm CS}_{t^*(\Gamma_1)}(t^*\Gamma)\equiv  1/2 \bmod {\bf Z}.$$
It follows that 
$${\rm CS}(\Gamma)=\pm 1/4 \bmod{\bf Z}.$$

    From the fact that the value of this invariant on one
component of order $4$ is the negative of its value on
the other such component, it now follows
that on one of these components the value is $1/4$ modulo ${\bf
Z}$ and on the other it is $-1/4$ modulo ${\bf Z}$.
\end{proof}

\begin{corollary}\label{CCS}
For every positive integer $k$ dividing $4$, the function  ${\rm
CS}$ defines a bijection between the components of order $k$ of
${\cal T}_G(C)$ and
the points of order $k$ in
${\bf R}/{\bf Z}$.
\end{corollary}

\subsection{Proof of Theorem~\ref{CS} and Theorem~\ref{clock}
when $\langle C\rangle$ is not cyclic}

Theorem~\ref{CS} in the case where $\langle C\rangle$ is not
cyclic is contained in Corollary~\ref{CCS}.  The number of
components of
${\cal T}_G(C)$ and their dimensions are given by
Lemma~\ref{firstlemma} and Corollaries~\ref{order1,2}
and~\ref{order4}. These together with Corollary~\ref{CCS} now
give all the necessary  computations in order to apply
Theorem~\ref{clocked} to establish  Witten's clockwise symmetry,
Theorem~\ref{clock}, in the case when
$\langle C\rangle$ is non-cyclic.

\pagebreak

\section*{Diagrams and tables}

\subsection*{Extended coroot diagrams and extended coroot
integers}

In the diagrams, $\circ$ represents the extended coroot.

\vskip.4in

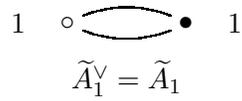
$$\xymatrix@C=.25cm{
1&{\circ}\ar@/^/ @{-}[rrr]\ar@/_/ @{-}[rrr]&&&{\bullet}&1
}$$
\centerline{$\widetilde A\spcheck_1=\widetilde A_1$}

\vfill

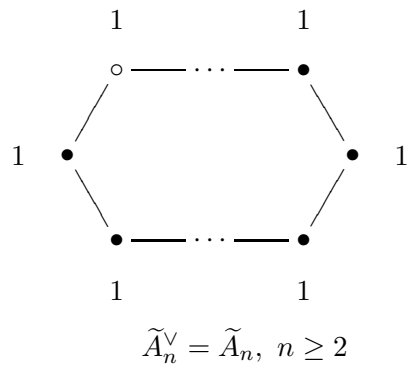
$$\xymatrix@R=.25cm@C=.25cm{
&&1&&&&1&&\\
&&{\circ}\ar@{-}[rr]&&\cdots\ar@{-}[rr]& &{\bullet}&&\\
&&&&&&&&&&\\
1&{\bullet}\ar@{-}[uur]\ar@{-}[ddr]&&&&&&{\bullet}
\ar@{-}[uul]\ar@{-}[ddl]&1\\ &&&&&&&&&&\\
&&{\bullet}\ar@{-}[rr]&&\cdots\ar@{-}[rr]& &{\bullet}&&\\
&&1&&&&1&&
}$$
\centerline{$\widetilde A\spcheck_n=\widetilde A_n,\ n\ge 2$}

\vfill

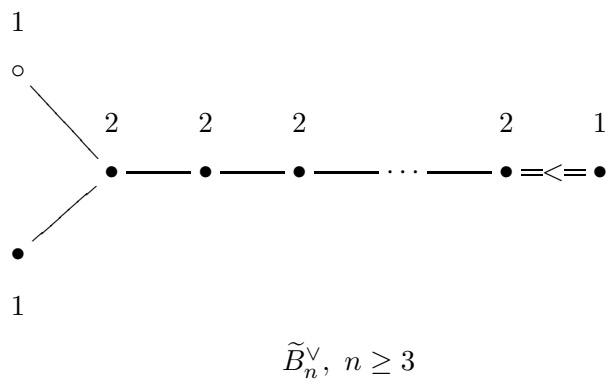
$$\xymatrix@R=.25cm{
1&&&&&&\\
{\circ}&&&&&&\\
&2&2&2&&2&1\\
&{\bullet}\ar@{-}[uul]\ar@{-}[ddl]\ar@{-}[r]&{\bullet}\ar@{-}[r]&{\bullet}
\ar@{-}[r]&\cdots\ar@{-}[r]& {\bullet}\ar@2{-}[r] |{\textstyle <}&{\bullet}\\
&&&&&&&\\
{\bullet}&&&&&&\\
1&&&&&&
}$$
\centerline{$\widetilde B\spcheck_n,\ n\ge 3$}

\pagebreak

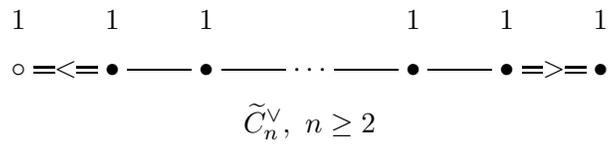
$$\xymatrix@R=.25cm{
1&1&1&&1&1&1\\
{\circ}\ar@2{-}[r]
|{\textstyle <}&{\bullet}\ar@{-}[r]&{\bullet}\ar@{-}[r]&\cdots
\ar@{-}[r]& {\bullet}\ar@{-}[r]&{\bullet}\ar@2{-}[r] |{\textstyle >}&{\bullet}
}$$
\centerline{$\widetilde C\spcheck_n,\ n\ge 2$}

\vfill

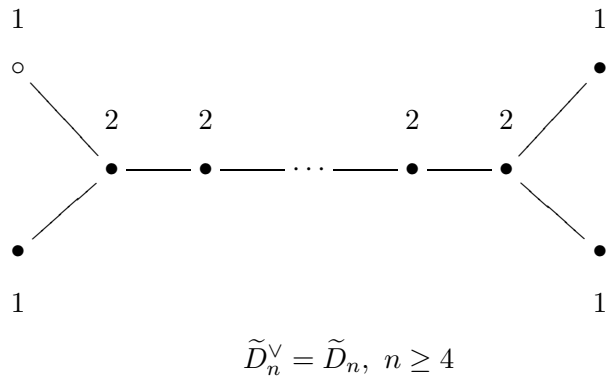
$$\xymatrix@R=.25cm{
1&&&&&&1\\
{\circ}&&&&&&{\bullet}\\
&2&2&&2&2&\\
&{\bullet}\ar@{-}[uul]\ar@{-}[ddl]\ar@{-}[r]&{\bullet}\ar@{-}[r]&\cdots
\ar@{-}[r]&{\bullet}\ar@{-}[r]& {\bullet}\ar@{-}[uur]\ar@{-}[ddr]&\\
&&&&&&&\\
{\bullet}&&&&&&{\bullet}\\
1&&&&&&1
}$$
\centerline{$\widetilde D\spcheck_n=\widetilde D_n,\ n\ge 4$}

\vfill

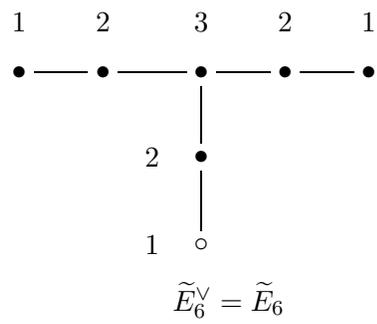
$$\xymatrix@R=.25cm@C=.25cm{
1&&2&&3&&2&&1\\
{\bullet}\ar@{-}[rr]&&{\bullet}\ar@{-}[rr]&&{\bullet}\ar@{-}[dd]\ar@{-}[rr]&&
{\bullet}\ar@{-}[rr] &&{\bullet}\\
&&&&&&&&&&\\
&&&2&{\bullet}\ar@{-}[dd]&&&&\\
&&&&&&&&&&\\
&&&1&{\circ}&&&&
}$$
\centerline{$\widetilde E\spcheck_6=\widetilde E_6$}

\vfill

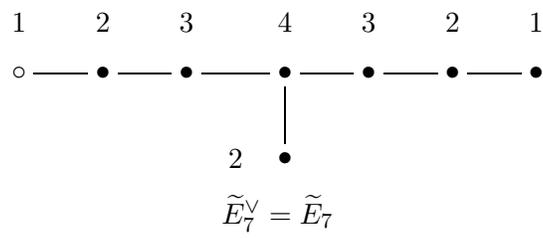
$$\xymatrix@R=.25cm@C=.25cm{
1&&2&&3&&4&&3&&2&&1\\
{\circ}\ar@{-}[rr]&&{\bullet}\ar@{-}[rr]&&{\bullet}\ar@{-}[rr]&&{\bullet}
\ar@{-}[dd]\ar@{-}[rr]&& {\bullet}\ar@{-}[rr] &&{\bullet}\ar@{-}[rr]
&&{\bullet}\\ &&&&&&&&&&&&\\
&&&&&2&{\bullet}&&&&&&}$$
\centerline{$\widetilde E\spcheck_7=\widetilde E_7$}

\pagebreak

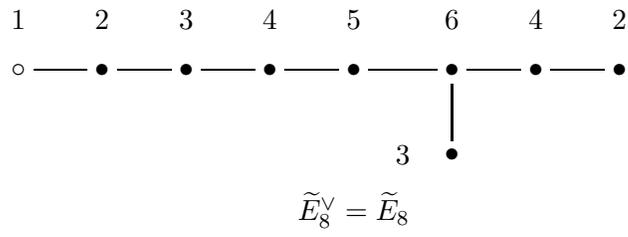
$$\xymatrix@R=.25cm@C=.25cm{
1&&2&&3&&4&&5&&6&&4&&2\\
{\circ}\ar@{-}[rr]&&{\bullet}\ar@{-}[rr]&&{\bullet}\ar@{-}[rr]&&{\bullet}
\ar@{-}[rr]&&{\bullet}
\ar@{-}[rr]&&{\bullet}\ar@{-}[dd]\ar@{-}[rr]&& {\bullet}\ar@{-}[rr]
&&{\bullet}\\ &&&&&&&&&&&&&&&\\
&&&&&&&&&3&{\bullet}&&&&&&
}$$
\centerline{$\widetilde E\spcheck_8=\widetilde E_8$}

\vfill

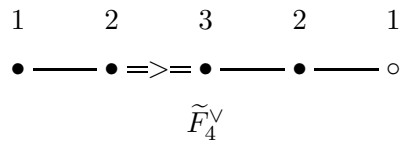
$$\xymatrix@R=.25cm{
1&2&3&2&1\\
{\bullet}\ar@{-}[r]&{\bullet}\ar@2{-}[r]
|{ \textstyle >}&{\bullet}\ar@{-}[r]&{\bullet}\ar@{-}[r]&{\circ} }$$
\centerline{$\widetilde F\spcheck_4$}

\vfill

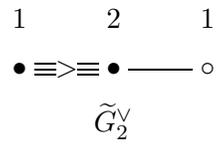
$$\xymatrix@R=.25cm{
1&2&1\\
{\bullet}\ar@3{-}[r] |{\textstyle >}&{\bullet}\ar@{-}[r]&{\circ}
}$$
\centerline{$\widetilde G\spcheck_2$}

\vfill

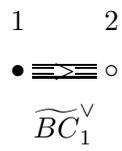
$$\xymatrix@R=.25cm{
1&2\\
{\bullet}\ar@{-}[r]<2.4pt>\ar@{=}[r]|{ \textstyle >}\ar@{-}[r]<-2.4pt>
&{\circ} }$$
\centerline{$\widetilde{BC}\spcheck_1$}

\vfill

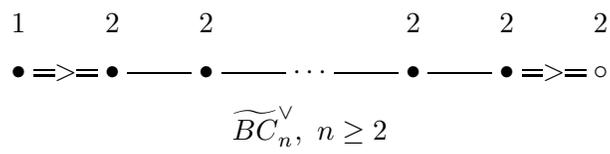
$$\xymatrix@R=.25cm{
1&2&2&&2&2&2\\
{\bullet}\ar@2{-}[r]
|{\textstyle >}&{\bullet}\ar@{-}[r]&{\bullet}\ar@{-}[r]&\cdots
\ar@{-}[r]& {\bullet}\ar@{-}[r]&{\bullet}\ar@2{-}[r] |{\textstyle >}&{\circ}
}$$
\centerline{$\widetilde{BC}\spcheck_n,\ n\ge 2$}

\pagebreak

\subsection*{Quotient extended coroot diagrams and quotient
coroot integers}

 In the diagrams below, $c$ denotes an element of the center and
$o(c)$ is its order. Except in the case of $\widetilde A_n$, $c$
will always denote a generator of the center. 
$c_{\rm SO}$ is the non-trivial element
    in $\pi_1(SO(2n))\subset {\cal C}Spin(2n)$ and  $c_{\rm
exotic}$ is
    any other non-trivial element in ${\cal C}Spin(4n)$.

\vfill

$$\xymatrix@R=.25cm{
o(c)\\{\bullet}}$$
\centerline{$\widetilde A\spcheck_n/c$ when $o(c)=n+1$}

\vfill

$$\xymatrix@C=.25cm{
o(c)&{\bullet}\ar@/^/ @{-}[rrr]\ar@/_/ @{-}[rrr]&&&{\bullet}&o(c)
}$$
\centerline{$\widetilde A\spcheck_{n}/c$ when $o(c)=(n+1)/2$}

\vfill

$$\xymatrix@R=.25cm@C=.25cm{
&&o(c)&&&&o(c)&&\\
&&{\bullet}\ar@{-}[rr]&&\cdots\ar@{-}[rr]& &{\bullet}&&\\
&&&&&&&&&&\\
o(c)&{\bullet}\ar@{-}[uur]\ar@{-}[ddr]&&&&&&{\bullet}
\ar@{-}[uul]\ar@{-}[ddl]&o(c)\\ &&&&&&&&&&\\
&&{\bullet}\ar@{-}[rr]&&\cdots\ar@{-}[rr]& &{\bullet}&&\\
&&o(c)&&&&o(c)&&}$$
\centerline{$\widetilde A\spcheck_n/c$ when $o(c)<(n+1)/2$}

\vfill

$$\xymatrix@R=.25cm{
2&2&2&&2&2&1\\
{\bullet}\ar@2{-}[r]
|{\textstyle <}&{\bullet}\ar@{-}[r]&{\bullet}\ar@{-}[r]&\cdots
\ar@{-}[r]& {\bullet}\ar@{-}[r]&{\bullet}\ar@2{-}[r] |{\textstyle <}&{\bullet}
}$$
\centerline{$\widetilde B\spcheck_n/c$}

\pagebreak

$$\xymatrix@R=.25cm{
1&2\\
{\bullet}\ar@{-}[r]<2.4pt>\ar@{=}[r]|{ \textstyle >}\ar@{-}[r]<-2.4pt>
&{\circ} }$$
\centerline{$C_2\spcheck/c$}

\vfill

$$\xymatrix@R=.25cm{
1&2&2&&2&2&2\\
{\bullet}\ar@2{-}[r]
|{\textstyle >}&{\bullet}\ar@{-}[r]&{\bullet}\ar@{-}[r]&\cdots
\ar@{-}[r]& {\bullet}\ar@{-}[r]&{\bullet}\ar@2{-}[r] |{\textstyle >}&{\bullet}
}$$
\centerline{$\widetilde C\spcheck_{2n}/c,\ n\ge 2$}

\vfill

$$\xymatrix@C=.25cm{
2&{\bullet}\ar@/^/ @{-}[rrr]\ar@/_/ @{-}[rrr]&&&{\bullet}&2
}$$
\centerline{$\widetilde C\spcheck_3/c$}

\vfill

$$\xymatrix@R=.25cm{
2&2&2&&2&2&2\\
{\bullet}\ar@2{-}[r]
|{\textstyle <}&{\bullet}\ar@{-}[r]&{\bullet}\ar@{-}[r]&\cdots
\ar@{-}[r]& {\bullet}\ar@{-}[r]&{\bullet}\ar@2{-}[r] |{\textstyle >}&{\bullet}
}$$
\centerline{$\widetilde C\spcheck_{2n+1}/c,\ n\ge 2$}

\vfill

$$\xymatrix@R=.25cm{
2&2&2&&2&2&2\\
{\bullet}\ar@2{-}[r]
|{\textstyle <}&{\bullet}\ar@{-}[r]&{\bullet}\ar@{-}[r]&\cdots
\ar@{-}[r]& {\bullet}\ar@{-}[r]&{\bullet}\ar@2{-}[r] |{\textstyle >}&{\bullet}
}$$
\centerline{$\widetilde D\spcheck_n/c_{\rm SO}$}

\pagebreak

$$\xymatrix@R=.25cm{
2&2&2 \\
{\bullet}\ar@2{-}[r]|{\textstyle <}&{\bullet}\ar@2{-}[r] |{\textstyle
  >}&{\bullet} 
}$$
\centerline{$D\spcheck_4/c_{\rm exotic}$}

\vfill

$$\xymatrix@R=.25cm{
2&&&&&&\\
{\bullet}&&&&&&\\
&4&4&4&&4&2\\
&{\bullet}\ar@{-}[uul]\ar@{-}[ddl]\ar@{-}[r]&{\bullet}\ar@{-}[r]&{\bullet}
\ar@{-}[r]&\cdots\ar@{-}[r]& {\bullet}\ar@2{-}[r] |{\textstyle <}&{\bullet}\\
&&&&&&&\\
{\bullet}&&&&&&\\
2&&&&&&
}$$
\centerline{$\widetilde D\spcheck_{2n}/c_{\rm exotic},\ n\ge 3$}

\vfill

$$\xymatrix@R=.25cm{
4&4&4&&4&4&4\\
{\bullet}\ar@2{-}[r]
|{\textstyle <}&{\bullet}\ar@{-}[r]&{\bullet}\ar@{-}[r]&\cdots
\ar@{-}[r]& {\bullet}\ar@{-}[r]&{\bullet}\ar@2{-}[r] |{\textstyle
>}&{\bullet}}$$
\centerline{$\widetilde D\spcheck_{2n+1}/c$}

\vfill

$$\xymatrix@R=.25cm{
2&4\\
{\bullet}\ar@{-}[r]<2.4pt>\ar@{=}[r]|{ \textstyle
  >}\ar@{-}[r]<-2.4pt>&{\bullet} 
}$$
\centerline{$D_4/{\cal C}D_4$}

\pagebreak

$$\xymatrix@R=.25cm{
2&4&4&&4&4&4\\
{\bullet}\ar@2{-}[r]
|{\textstyle >}&{\bullet}\ar@{-}[r]&{\bullet}\ar@{-}[r]&\cdots
\ar@{-}[r]& {\bullet}\ar@{-}[r]&{\bullet}\ar@2{-}[r] |{\textstyle >}&{\bullet}
}$$
\centerline{$\widetilde D\spcheck_{2n}/{\cal C}D_{2n},\ n\ge 3$}

\vfill

$$\xymatrix@R=.25cm{
3&6&3\\
{\bullet}\ar@3{-}[r] |{\textstyle >}&{\bullet}\ar@{-}[r]&{\bullet}
}$$
\centerline{$\widetilde E\spcheck_6/c$}

\vfill

$$\xymatrix@R=.25cm{
2&4&6&4&2\\
{\bullet}\ar@{-}[r]&{\bullet}\ar@2{-}[r]
|{ \textstyle >}&{\bullet}\ar@{-}[r]&{\bullet}\ar@{-}[r]&{\bullet} }$$
\centerline{$\widetilde E\spcheck_7/c$}

\vfill

\subsection*{Root systems on $\frak t^{w_{\cal C}}$}

\vskip.2in

\hskip-1in\vbox{\begin{tabular}{||c|c|c|c|c|c|c|c||}
\hline
$G$ & ${\cal C}$ & $L_{\cal C}$ & $\Phi^{w_{\cal C}}$ &
$\Phi^{\rm res}$ &
$\Phi^{\rm   proj}$ & $\Phi(w_{\cal C})$ & $g_{\overline a}$ \\ 
\hline\hline
$B_n$ & ${\cal C}B_n$  & $A_1$ (short root) & $C_{n-1}$  & $B_{n-1}$ &
$BC_{n-1}$ & $BC_{n-1}$ & $1,2,\ldots,2$ \\ \hline  
$C_{2n+1}$ & ${\cal C}C_{2n+1}$ & $\prod_{i=1}^{n+1}A_1$ &  $C_n$ &
$BC_n$  & $BC_n$ &  $C_n$ & $2,2,\ldots,2$ \\ \hline
$C_{2n}$ & ${\cal C}C_{2n}$ & $\prod_{i=1}^{n}A_1$ & $C_n$  & $C_n$ &
$BC_n$ & $BC_n$ & $1,2,\ldots,2$ \\ \hline 
$D_n$ & $\pi_1(SO(2n))$ & $A_1\times A_1$ & $C_{n-2}$  & $B_{n-2}$ &
$C_{n-2}$  & $C_{n-2}$  & $2,2,\ldots,2$ \\ \hline
$D_{2n}$ & $\langle c_{\rm exotic}\rangle$ & $\prod_{i=1}^{n}A_1$ &
$B_n$  & $C_n$ & $B_n$ & $B_n$ & $2,4,\ldots,4,2,2$ \\ \hline 
$D_{2n+1}$ & ${\cal C}D_{2n+1}$ & $\prod_{i=1}^{n-1}A_1\times A_3$ &
$C_{n-1}$ & 
$BC_{n-1}$ & $BC_{n-1}$ & $C_{n-1}$ & $4,4,\ldots,4$ \\ \hline  
$D_{2n}$ & ${\cal C}D_{2n}$ & $\prod_{i=1}^{n+1}A_1$ & $C_{n-1}$ & $BC_{n-1}$ &
$BC_{n-1}$ & $BC_{n-1}$ & 
$2,4,\ldots, 4$  \\ \hline
$E_6$ & ${\cal C}E_6$ & $A_2\times A_2$ & $G_2$  & $G_2$ & $G_2$ & $G_2$ &
$3,6,3$ \\ \hline 
$E_7$ & ${\cal C}E_7$ & $A_1\times A_1\times A_1$ & $F_4$ & $F_4$  & $F_4$
& $F_4$ & $2,4,6,4,2$ \\ \hline\hline
\end{tabular}}

\pagebreak

\subsection*{Root systems on ${\mathfrak t}(k)$ for $k>1$}

\vskip.2in

\hskip.6in\vbox{\begin{tabular}{||c|c|c|c|c||}
\hline
$G$ &  $k$ & $L$ & $\Phi(\frak t(k))$ & $g_a$ divisible by $k$ 
\\ 
\hline\hline
$B_n$ &  $2$  & $B_3$ & $C_{n-3}$ & $2,2,\ldots,2$      \\ \hline  
$D_{n}$ & $2$ & $D_4$ & $C_{n-4}$ & $2,2,\ldots,2$ \\ \hline 
$E_6$ &  $2$ & $D_4$  & $A_2$ & $2,2,2$ \\ \hline 
$E_6$ &  $3$ & $E_6$  & trivial & $3$ \\ \hline 
$E_7$ &  $2$ & $D_4$ & $B_3$ & $2,4,2,2$   \\ \hline
$E_7$ &  $3$ & $E_6$ & $A_1$ & $3,3$    \\ \hline
$E_7$ & $4$ &  $E_7$ & trivial & $4$    \\ \hline
$E_8$ & $2$ & $D_4$ & $F_4$ & $2,4,6,4,2$ \\ \hline
$E_8$ & $3$ & $E_6$ & $G_2$ & $3,6,3$ \\ \hline
$E_8$ & $4$ & $E_7$ & $A_1$ & $4,4$ \\ \hline
$E_8$ & $5$ & $E_8$ & trivial & $5$ \\ \hline
$E_8$ & $6$ & $E_8$ & trivial & $6$  \\ \hline
$F_4$ & $2$ & $B_3$ & $A_1$ & $2,2$ \\ \hline
$F_4$ & $3$ & $F_4$ & trivial & $3$  \\ \hline
$G_2$ & $2$ & $G_2$ & trivial & $2$  \\ 
\hline\hline
\end{tabular}}

\vskip.7in

\subsection*{Root systems on $\frak t^{w_C}(\overline {\bf
g},k)$ for
$\langle C\rangle \not=1$
 and $k\not| n_0$}

\vskip.2in

\begin{tabular}{||c|c|c|c|c|c||}
\hline
$G$ & $\langle C\rangle $ & $k$ & $L$ &
$\Phi(\frak t^{w_C}(\overline {\bf g},k))$ & $g_{\overline a}$
divisible by $k$  \\  
\hline\hline
$B_n$ & ${\cal C}B_n$  & $2$  & $C_2$ & $C_{n-2}$ & $2,2,\ldots,2$
\\ \hline   
$C_{2n}$ & ${\cal C}C_{2n}$ & $2$  &  $\prod_{i=1}^{n-1}A_1\times C_2$ &
$C_{n-1}$ & $2,2,\ldots,2$ 
 \\ \hline 
$D_{2n}, n\ge 3$ & $\langle c_{\rm exotic}\rangle $ & $4$ &
$\prod_{i=1}^{n-3}A_1\times D_6$ & $C_{n-3}$ & $4,4,\ldots, 4$  \\ \hline 
$D_{2n}$ & ${\cal C}D_{2n}$ & $4$ & $\prod_{i=1}^{n-2}A_1\times
D_4$ & $C_{n-2}$ & $4,4\ldots,4$ \\ \hline
$E_6$ & ${\cal C}E_6$ & $2$ & $E_6$  & trivial& $6$  \\ \hline 
$E_6$ & ${\cal C}E_6$ & $6$ & $E_6$  & trivial & $6$  \\ \hline 
$E_7$ & ${\cal C}E_7$ & $4$ & $D_6$ & $A_1$ & $4,4$    \\ \hline
$E_7$ & ${\cal C}E_7$ & $3$ & $E_7$ & trivial & $6$    \\ \hline
$E_7$ & ${\cal C}E_7$ & $6$ & $E_7$ & trivial & $6$   \\ 
\hline\hline
\end{tabular}

\end{document}